%% file: morseftlch.tex
\documentclass[12pt]{amsart}
\usepackage{amscd}
\usepackage[dvips]{graphicx}

\newtheorem{thm}{Theorem}[section]

\newtheorem{dfn}[thm]{Definition}

\newtheorem{prp}[thm]{Proposition}

\newtheorem{cor}[thm]{Corollary}

\newtheorem{lma}[thm]{Lemma}

\newtheorem{rem}[thm]{Remark}
\newenvironment{rmk}{\begin{rem}\rm}{\end{rem}}

\newtheorem*{clm}{Claim}

\newenvironment{pf}{\begin{proof}}{\end{proof}}

\numberwithin{equation}{section}

\newcommand{\R}{{\mathbb{R}}}
\newcommand{\C}{{\mathbb{C}}}

\newcommand{\Z}{{\mathbb{Z}}}
\newcommand{\conf}{{\mathcal{C}}}
\newcommand{\B}{{\mathcal{B}}}

\newcommand{\Ordo}{{\mathbf{O}}}
\newcommand{\ordo}{{\mathbf{o}}}
\newcommand{\W}{{\mathcal{W}}}
\newcommand{\V}{{\mathcal{V}}}
\newcommand{\M}{{\mathcal{M}}}

\newcommand{\gdim}{\operatorname{gdim}}

\newcommand{\Exp}{\operatorname{Exp}}
\newcommand{\ev}{\operatorname{ev}}
\newcommand{\la}{\langle}
\newcommand{\ra}{\rangle}
\newcommand{\pa}{\partial}
\newcommand{\id}{\operatorname{id}}

\newcommand{\ind}{\operatorname{Index}}
\newcommand{\Spa}{\operatorname{Span}}
\newcommand{\area}{\operatorname{Area}}
\newcommand{\ix}{\operatorname{Index}}
\newcommand{\krn}{\operatorname{Ker}}

\newcommand{\lk}{\operatorname{lk}}

\newcommand{\sblv}{{\mathcal{H}}}
\newcommand{\lr}{\lrcorner}
\newcommand{\inn}{\mathbf{in}}
\newcommand{\out}{\mathbf{out}}
\newcommand{\0}{\mathbf{0}}

\title[Morse flow trees and Legendrian contact homology\dots]
{Morse flow trees and Legendrian contact homology in $1$-jet spaces}
\author{Tobias Ekholm}

\begin{document}

\begin{abstract}
Let $L\subset J^1(M)$ be a Legendrian submanifold of the $1$-jet space of a Riemannian $n$-manifold $M$. A correspondence is established between rigid flow trees in $M$ determined by $L$ and  boundary punctured rigid pseudo-holomorphic disks in $T^\ast M$, with boundary on the projection of $L$ and asymptotic to the double points of this projection at punctures, provided $n\le 2$, or provided $n>2$ and the front of $L$ has only cusp edge singularities. This result, in particular, shows how to compute the Legendrian contact homology of $L$ in terms of Morse theory.
\end{abstract}

\subjclass[2000]{57R17; 53D40}
\thanks{Part of the work on this paper was done while the author was a Research Fellow of the Royal Swedish Academy of Sciences sponsored by the Knut and Alice Wallenberg foundation. The author is an Alfred P. Sloan Research Fellow and acknowledges support from NSF-grant DMS-0505076.}
\maketitle

\input{Sec/1intr}

\input{Sec/2defs}

\input{Sec/3dim}

\input{Sec/4metric}

\input{Sec/5disktotree}

\input{Sec/6treetodisk}

\input{Sec/7expl}

\input{Sec/ref}
\end{document}

%% file: Sec/1intr.tex
\section{Introduction}\label{S:intr}
Let $M$ be a smooth $n$-manifold, let $J^0(M)=M\times\R$ be its $0$-jet space, and let $J^1(M)=T^\ast M\times\R$ be its $1$-jet space endowed with the standard contact structure $\xi$ which is the kernel of the $1$-form $dz-p\,dq$, where $z$ is a coordinate in the $\R$-direction and where $p\,dq$ is the canonical $1$-form on $T^\ast M$. An $n$-dimensional submanifold $L\subset J^1(M)$ is {\em Legendrian} if it is everywhere tangent to $\xi$. This paper concerns Legendrian submanifolds of $1$-jet spaces, and in particular their contact homology. Legendrian contact homology is a part of Symplectic Field Theory \cite{EGH}. It is a framework for finding isotopy invariants of Legendrian submanifolds of contact manifolds by "counting" rigid (pseudo-)holomorphic disks. The analytical foundations of Legendrian contact homology in $1$-jet spaces were established in \cite{EES4}, see also \cite{EES2, EES3}.

Contact homology has proved very useful, see e.g. \cite{Ch, El, EES1}, especially for Legendrian submanifolds of dimension $1$, where the Riemann mapping theorem can be used to give a combinatorial and computable description of the theory, see \cite{Ch, ENS}. In the higher dimensional case, finding holomorphic disks involves solving a non-linear first order partial differential equation. This is an infinite dimensional problem and it is therefore often difficult to compute the contact homology of a given Legendrian submanifold. The main result of this paper (Theorem \ref{t:main}) reduces the computation of Legendrian contact homology in $1$-jet spaces to a finite dimensional problem in Morse theory. One of its applications is as a step in showing that Ng's powerful invariants of knots in $3$-space \cite{Ng1, Ng2, Ng3}, constructed in a combinatorial fashion by considering the contact homology of co-normal lifts of knots (which are Legendrian tori in the $1$-jet space of the $2$-sphere) indeed agree with the invariants arising from holomorphic disk counts, see \cite{EESN}. Other applications (e.g. to contact surgery and to Legendrian embeddings) will be described in forthcoming papers.

Similar results, relating Morse theory and holomorphic disks, were obtained in the (slightly different) setting of Floer homology of Lagrangian intersections. In \cite{Fl}, Floer showed that holomorphic strips ($2$-gons) in a cotangent bundle $T^\ast M$ of a Riemannian manifold $M$, with one boundary component mapping to the $0$-section and the other to the graph $\Gamma_{df}$ of the differential of a function $f\colon M\to\R$, correspond to gradient lines of $f$ connecting critical points. In \cite{FO} Fukaya and Oh generalized Floer's result showing that holomorphic $n$-gons with boundary components mapping to distinct graphical Lagrangian submanifolds $\Gamma_{df_1},\dots,\Gamma_{df_n}$ correspond to gradient flow trees. The situation under study in this paper is more involved than that in \cite{FO} because of the appearance of non-graphical Lagrangian submanifolds and this is reflected in the corresponding Morse theory.

In order to state Theorem \ref{t:main}, we give brief descriptions of the holomorphic- and Morse theoretic objects involved. Let $L\subset J^1(M)$ be a closed Legendrian submanifold. The image of $L$ under the {\em Lagrangian projection} $\Pi_\C\colon J^1(M)\to T^\ast M$ is an immersed Lagrangian submanifold. After small perturbation, the only self intersections of $\Pi_\C(L)$ are transverse double points. Let $J\colon T(T^\ast M)\to T(T^\ast M)$ be an almost complex structure tamed by the standard symplectic form $\omega=d(p\,dq)$ on $T^\ast M$. (An almost complex structure $J$ is tamed by $\omega$ if $\omega(X, J X)>0$ for all nonzero $X\in T(T^\ast M)$. A Riemannian metric on $M$ induces an almost complex structure tamed by $\omega$, see Subsection \ref{s:acs}.)

If $S$ is a Riemann surface with complex structure $j\colon TS\to TS$ then a map $u\colon S\to T^\ast M$ is {\em $J$-holomorphic} if $du+J\circ du\circ j=0$. We study boundary punctured $J$-holomorphic disks with boundary mapping to $\Pi_\C(L)$, which are asymptotic to double points at the punctures, and which have restrictions to the boundary which admit a continuous lift to $L$. We call such disks {\em $J$-holomorphic disks with boundary on $L$}. Two of their key properties are as follows. First, the punctures come equipped with signs, see Definition \ref{d:phol}. Second, associated to a $J$-holomorphic disk is its formal dimension, see Proposition \ref{p:tfdim=dfdim}, which measures the expected dimension of the space of nearby $J$-holomorphic disks. We say that a disk is {\em rigid} if it has formal dimension $0$ and if it is transversely cut out by its defining differential equation. Legendrian contact homology is defined using disks with only one positive puncture and we say that an almost complex structure $J$ is {\em regular for $L$}, if there are no $J$-holomorphic disks with one positive puncture of negative formal dimension and if all such disks of formal dimension $0$ are rigid.

Consider the {\em front projection} $\Pi_F\colon J^1(M)\to J^0(M)$ and the {\em base projection} $\Pi\colon J^1(M)\to M$. After small perturbation of $L$ the map $\Pi\colon L\to M$ is an immersion outside a codimension $1$ submanifold $\Sigma\subset L$ and at points in $\Sigma$ outside a codimension $1$ subset $\Sigma'\subset\Sigma$, $\Pi$ has a standard fold singularity. The corresponding singularity of $\Pi_F$ is a {\em cusp edge}. For any point $m\in M$, the smooth sheets of $\Pi_F(L)$ over $m$ are graphs of a finite number of locally defined functions.

Equip $M$ with a Riemannian metric. A {\em flow tree} determined by $L$ is a continuous map of a tree $\Gamma$ to $M$ such that the restriction of the map to any edge of $\Gamma$ parameterizes a part of a gradient flow line of some local function difference, see Definition \ref{d:ftree}. Here we point out three  properties. First, the vertices of a flow tree which map to critical points of a local function difference are called punctures and come equipped with a sign, see Definition \ref{d:tpunctsign}. Second, associated to a flow tree is its formal dimension, see Definition \ref{d:tfdim}, which measures the expected dimension of the space of nearby flow trees. We say that a tree is {\em rigid} if its formal dimension is $0$ and if it is transversely cut out, see Proposition \ref{p:ttv}.
Third, a vertex of a tree is either simple or a multiple cover and a tree without vertices which are multiple covers is called {\em simple}, see Definition \ref{d:simple}. All trees with only one positive puncture are simple, see Lemma \ref{l:simple}, and all rigid trees of Legendrian submanifolds of dimension $\le 2$ are simple, see \S\ref{ss:many+}. We say that a metric is {\em $P$-regular for $L$}, where $P>0$ is an integer, if the following holds for flow trees with $\le P$ positive punctures: there exist no simple flow trees of formal dimension $<0$, all simple trees of formal dimension $0$ are rigid, and the set of simple rigid flow trees is finite.

For $0<\lambda\le 1$, consider the fiber scaling $s_\lambda\colon J^1(M)\to J^1(M)$, $s_\lambda(m,v,z)=(m,\lambda v,\lambda z)$, where $m\in M$, $v\in T^\ast_m M$, and $z\in\R$. If $L\subset J^1(M)$ is a Legendrian submanifold, let $L_\lambda=s_\lambda(L)$. Since $s_\lambda^\ast(dz-p\,dq)=\lambda(dz-p\,dq)$, $L_\lambda$ is a Legendrian submanifold. In \S\ref{ss:LegisoCE} we Legendrian isotope $L_\lambda$ by a $C^0$-small Legendrian isotopy supported near the cusp edge of $L_\lambda$ to $\tilde L_\lambda$.

\begin{thm}\label{t:main}
Let $L\subset J^1(M)$ be an $n$-dimensional Legendrian submanifold and fix $P>0$. (If $n>2$, assume that the front of $L$ has only cusp-edge singularities.)
\begin{itemize}
\item[{\rm (a)}]
After a small perturbation of $L$ the set $\Omega$ of $P$-regular metrics for $L$ is open and dense.
\item[{\rm (b)}] For any $g\in\Omega$ there exists $\lambda_0>0$ and almost complex structures $J_\lambda$, regular for $\tilde L_\lambda$, $0<\lambda<\lambda_0$, with the following property. The rigid $J_\lambda$-holomorphic disks with boundary on $\tilde L_\lambda$ and with at most one positive puncture are in 1-1 correspondence with the rigid flow trees with one positive puncture determined by $L$.
\end{itemize}
In particular, the contact homology of $L$ can be computed counting rigid flow trees (defined using any $1$-regular metric) instead of holomorphic disks.
\end{thm}

Theorem \ref{t:main} (a) is proved in Subsection \ref{s:ttv} and (b) in Subsection \ref{s:unifinv}. In fact, Theorem \ref{t:main} (b) follows from (a) in combination with two other results which we state in Subsection \ref{s:outline} and which give more information about the nature of the 1-1 correspondence and about non-rigid disks. In many cases the extra assumptions on the front singularities for Legendrian submanifolds of dimension $>2$ in Theorem \ref{t:main} can be arranged to hold: there exists h-principles which guarantee that unless there is a homotopy reason for front singularities of codimension larger than $1$ to exist, they can be removed by Legendrian isotopy, see \cite{En}.

Theorem \ref{t:main} (b) fails for disks with many positive punctures in dimensions higher than $2$. This is related to the phenomenon that $J$-holomorphic disks may converge to a multiply covered tree. To understand rigid holomorphic disks with several positive punctures in the higher dimensional setting it is in general not enough to study rigid trees, also higher dimensional spaces of trees must be considered. We give a brief discussion of limits of disks with several positive punctures in \S\ref{ss:many+}.

\subsection{Outline and statements of two theorems}\label{s:outline}
The paper basically constitutes a proof of Theorem \ref{t:main}. Since the proof is rather long and involved we give, in this subsection, a brief sketch of its two main parts. This sketch simultaneously serves as an outline of the paper. Let $L\subset J^1(M)$ be a Legendrian submanifold. We use notation as above.

In Section \ref{S:defs}, we define the basic objects: $J$-holomorphic disks with boundary on $L$ and flow trees determined by $L$. In particular, we associate to each flow tree its {\em $1$-jet lift}. This is a curve in $L$ with a finite-to-one map to the tree.

In Section \ref{S:dimcount} we discuss dimension concepts for trees and for disks. We associate two dimension concepts to a flow tree: its {\em formal-} and its {\em geometric dimension}, see Definitions \ref{d:tfdim} and \ref{d:tgdim}, respectively. The formal dimension of a given tree is intended to measure the dimension of the space of nearby trees with $1$-jet lifts which are homotopic to the $1$-jet lift of the given tree, whereas the geometric dimension is intended to measure the dimension of the space of nearby trees with the exact same geometric properties as the given one. We show that the formal dimension is at least as large as the geometric dimension, see Lemma \ref{l:fdim>gdim}, and prove, under mild genericity conditions on the metric or on $L$, that the total number of edges of a flow tree can be bounded in terms of its number of positive punctures and its formal dimension, see Lemma \ref{l:tfinite}. This latter result allows us to prove a transversality result, see Proposition \ref{p:ttv}, which implies Theorem \ref{t:main} (a).

The formal dimension of a holomorphic disk is a function of its punctures and of the homotopy class (in $L$) of its restriction to the boundary. We present the dimension formula for disks in Proposition \ref{p:tfdim=dfdim} and note that the formal dimension associated to the boundary condition which arises as the $1$-jet lift of a flow tree $\Gamma$ agrees with the formal dimension of $\Gamma$.

In Section \ref{S:metrpert}, we describe how to alter a metric $g$ on $M$, regular for $L$, and the family $L_\lambda$ itself in order to gain better control of the behavior of holomorphic disks. The alteration consists of a Legendrian isotopy taking $L_\lambda$ to $\tilde L_\lambda$ (mentioned above) and a further (ambient) isotopy of $\Pi_\C(\tilde L_\lambda)$ which results in an immersed totally real submanifold $\hat L_\lambda\subset T^\ast M$, close to $\Pi_\C(L_\lambda)$, and converging to it as certain deformation parameters approach $0$. In particular, for sufficiently small deformation parameters and sufficiently small $\lambda>0$ there is a 1-1 correspondence between rigid flow trees determined by $\hat L_\lambda$ and rigid flow trees determined by $L$, see Lemma \ref{l:ftrel}. (The flow trees of $\hat L_\lambda$ are not completely independent of $\lambda$ but converges as $\lambda\to 0$. We call the limiting trees {\em flow trees of $\hat L_0$}.) Moreover, there exists fiber preserving diffeomorphisms, arbitrarily close to the identity,  $\Phi_\lambda\colon T^\ast M\to T^\ast M$ covering the identity on $M$ such that $\Phi_\lambda(\Pi_\C(\tilde L_\lambda))=\hat L_\lambda$. If $J$ is the almost complex structure induced by the metric $g$ then $J_\lambda$ in Theorem \ref{t:main} (b) is given by $J_\lambda=d\Phi_\lambda^{-1}\circ J\circ d\Phi_\lambda$. Thus, $J_\lambda$-holomorphic disks with boundary on $\Pi_\C(\tilde L_\lambda)$ correspond in a 1-1 fashion (composing with $\Phi_\lambda$) to $J$-holomorphic disks with boundary on $\hat L_\lambda$. We will work in the latter setting, keeping the almost complex structure standard and deforming the Lagrangian projection.

In Section \ref{S:disktotree} we show that any sequence of $J$-holomorphic disks with boundary on $\hat L_\lambda$ has a subsequence which converges to a flow tree of $\hat L_0$. More precisely, we show the following.
\begin{thm}\label{t:disktotree}
Let $L\subset J^1(M)$ be a Legendrian submanifold of dimension $n$. If $n>2$, assume that the front of $L$ has only cusp-edge singularities. If $u_\lambda$, $\lambda\to 0$ is a sequence of $J$-holomorphic disks with one positive puncture, of formal dimension $d$, and with boundary on $\hat L_{\lambda}$, then there exists a subsequence $\lambda\to 0$ such that $u_{\lambda}$ converges to a (possibly broken) flow tree $\Gamma$ of $\hat L_0$ of formal dimension $d$ as $\lambda\to 0$.
\end{thm}
Theorem \ref{t:disktotree} is proved in Subsection \ref{s:totavlkbu} (a similar result holds for disks with many positive punctures, see \S\ref{ss:many+}). We employ a specific conformal model for the source spaces of $J$-holomorphic disks, see Lemma \ref{l:confmod}, which in particular allows us to identify the space of conformal structures on a disk with $m$ punctures on the boundary with $\R^{m-3}$ in a natural way. The convergence statement in Theorem \ref{t:disktotree} involves two things: on the one hand, the coordinates of the conformal structure on the domain of $u_{\lambda}$ has, when multiplied by $\lambda$, a limit which is determined by the flow tree (see Remark \ref{r:domaincontrolI}), on the other hand, the image of the holomorphic map lies close to the flow tree. The main step in the proof is to show that after addition of a uniformly finite (as $\lambda\to 0$) number of punctures in the domains of $u_{\lambda}$, the derivatives $|Du_\lambda|$ satisfy a $K\lambda$-bound for some $K>0$. To demonstrate the uniform finiteness we employ blow-up arguments in combination with a total average linking number, see \S\ref{ss:deftotavlk}.

In Section \ref{S:treetodisk} we construct $J$-holomorphic disks near rigid flow trees. More precisely we show the following.

\begin{thm}\label{t:treetodisk}
For every $P>0$, there exists $\lambda_0>0$ and $\delta>0$ such that the following hold for all $0<\lambda<\lambda_0$. If $\Gamma$ is a simple rigid flow tree of $\hat L_\lambda$ with $p\le P$ positive punctures then there exists exactly one (up to holomorphic automorphisms of the source) $J$-holomorphic disk with boundary on $\hat L_\lambda$, of formal dimension $0$, and with $p$ positive punctures in a $\delta$-neighborhood of $\Gamma$ and this disk is rigid. Moreover, if a sequence of holomorphic disks $u_{\lambda}$, $\lambda\to 0$, converges to a simple rigid tree of $\hat L_0$ then there exists a simple rigid tree $\Gamma$ of $\hat L_\lambda$ such that $u_\lambda$ lies in a $\delta$-neighborhood of $\Gamma$.
\end{thm}

The proof of Theorem \ref{t:treetodisk} is given in Subsection \ref{s:unifinv}. The notion of a $\delta$-neighborhood of a rigid tree involves on the one hand a neighborhood of conformal structures on domains naturally associated to a rigid tree and on the other hand a norm in a certain function space on this domain. A rough sketch of the main steps of the proof is as follows. We first construct approximately $J$-holomorphic disks $u_\lambda$ with boundary on $\hat L_\lambda$ near $\Gamma$. To infer the existence of actual holomorphic disks we use Floer's Picard lemma, see Lemma \ref{l:FloerPicard}, which in particular requires the existence of a uniformly bounded inverse of the linearization of the $\bar\pa_J$-operator at $u_\lambda$. To establish the existence of such an inverse is somewhat complicated for the following two reasons. First, the boundary conditions for the linearized $\bar\pa_J$-operator degenerates as $\lambda\to 0$. Second, the conformal structures on the domains of $u_\lambda$ converge to the deepest stratum of the boundary of the moduli space of conformal structures. (Expressed in another way: if the domain has $m$ punctures, all components of the corresponding point in $\R^{m-3}$ go to $\infty$.) To deal with these problems we introduce a Sobolev space with a weight function, which grows exponentially as punctures- and as certain marked points on the boundary of the domain are approached. To compensate the drop in Fredholm index that these weights cause near the limit we augment the Sobolev space by cut-off local solutions of the $\bar\pa_J$-equation. In these augmented function spaces, the kernel of the linearized $\bar\pa_J$-operator is closely related to the tangent space of the space of flow trees. This allows us to find the inverse of the linearization and thereby to prove existence of a solution. To prove the last statement we use Theorem \ref{t:disktotree} together with knowledge of the local form of solutions to the $\bar\pa_J$-equation.

In Section \ref{S:expl}, we consider a simple example of a $1$-parameter family of flow trees which illustrates many of the properties of rigid flow trees discussed in earlier sections of the paper.

\subsection*{Acknowledgements} The author would like to thank J. Etnyre, L. Ng, and M. Sullivan for many discussions and for reading earlier versions of the paper. He would also like to thank K. Cielebak, Y. Eliashberg, and K. Zhu for discussions. 

%% file: Sec/2defs.tex
\section{Definitions of disks and trees}\label{S:defs}
In Subsection \ref{s:phol} we define $J$-holomorphic disks with boundary on $L\subset J^1(M)$ and  give specific (flat) models for the source spaces of such disks. In Subsection \ref{s:ftree} we define flow trees and derive some of their elementary properties.

\subsection{Holomorphic disks}\label{s:phol}
Let $L\subset J^1(M)=T^\ast M\times\R$ be a closed Legendrian submanifold and let $z$ be a coordinate in the $\R$-direction. We assume that $L$ is sufficiently generic so that the Lagrangian projection $\Pi_\C\colon J^1(M)\to T^\ast M$ restricted to $L$ has only transverse double points. If $c$ is a double point of $\Pi_\C(L)$ then we write $\Pi_\C^{-1}(c)\cap L=\{c^+,c^-\}$, where $z(c^+)>z(c^-)$. (Note that the Reeb field of the contact form $dz-p\,dq$ is $\pa_z$, therefore there exists a 1-1 correspondence between double points of $\Pi_\C(L)$ and Reeb chords on $L$ and we will use these two notions interchangeably.)

Let $D_{m+k}$ be the unit disk $D$ in the complex plane $\C$ with punctures  $x_1,\dots,x_m,y_1,\dots,y_k$ on the boundary and let $J$ be an almost complex structure on $T^\ast M$ which is tamed by the standard symplectic form $\omega$ on $T^\ast M$.
\begin{dfn}\label{d:phol}
A $J$-holomorphic disk with positive punctures $p_1,\dots,p_m$ and negative punctures $q_1,\dots,q_k$ and with boundary on $L$, is a map $u\colon D_{m+k}\to T^\ast M$ with the following properties.
\begin{itemize}
\item $du + J\circ du\circ i=0$ (where $i$ is the complex structure on the complex plane).
\item The restriction $u|\pa D_{m+k}$ has a continuous lift $\tilde u\colon\pa D_{m+k}\to L\subset J^1(M)$.
\item $lim_{\zeta\to x_j}u(\zeta)=p_j$ and $\lim_{\zeta\to x_j\pm}{\tilde u}(\zeta)=p_j^{\pm}$, where $\lim_{\zeta\to x_j+}$ means that $\zeta$ approaches $x_j$ from the region in $\pa D_{m+k}$ in the positive direction as seen from $x_j$ and $\lim_{\zeta\to x_j-}$ means it approaches $x_j$ from the region in the negative direction.
\item $lim_{\zeta\to y_j}u(\zeta)=q_j$ and $\lim_{\zeta\to x_j\pm}{\tilde u}(\zeta)=q_j^{\mp}$.
\end{itemize}
\end{dfn}

We note that if $\Phi\colon \Delta\to D_{m+k}$ is a conformal map and if $u$ is $J$-holomorphic then so is $u\circ\Phi$.

\subsubsection{Conformal models}\label{ss:confmod}
We will consider disks with $m\ge 3$ punctures (if the original disk has less than $3$ punctures we add marked points on the boundary). For technical reasons we shall distinguish one of the positive punctures of the disks considered. (In the case of disks with only one positive puncture this is automatic). Let $\conf_m$ denote the space of conformal structures on $D_m$ with one distinguished puncture. Then $\conf_m$ is a simplex of dimension $m-3$: fixing the position of the distinguished puncture at $1\in \pa D$ and the two immediately following punctures at $i$ and $-1$, respectively, standard coordinates on $\conf_m$ are given by the positions of the $m-3$ punctures in the lower hemisphere of the circle. From another perspective the complex structure on a disk with punctures is an endomorphism $j\colon TD_m\to TD_m$ with $j^2=-1$. We will view the $m-3$ dimensional tangent space of $T_j\conf_m$ at a complex structure structure $j$ in the following way. Consider $m-3$ vector fields $v_1,\dots,v_{m-3}$ on $D$ where $v_r$ has support in a small neighborhood of the $r^{\rm th}$ moving puncture, agrees with the holomorphic vector field $z\mapsto i\cdot z$ generating a rotation moving this puncture along the boundary, and is tangent to the boundary along the boundary. Then it follows from standard properties of the Riemann-Hilbert problem that if $\gamma$ is any linearized variation of the complex structure $j$. (That is, $\gamma$ is an endomorphism of $TD_m$ such that $j\gamma+\gamma j=0$.) Then there exists a (unique) vector field $u$ along $D$ tangent to the boundary along the boundary and vanishing at all the punctures of $D_m$ and unique constants $\lambda_1,\dots,\lambda_{m-3}$ such that
$$
j\cdot\gamma=\bar\pa u + \sum_j \lambda_j(\bar\pa v_j).
$$
We thus think of $(\bar\pa v_1,\dots,\bar\pa v_{m-3})$ as a basis of $T_j\conf_m$.

\begin{figure}[htbp]
\begin{center}
\includegraphics[angle=0, width=8cm]{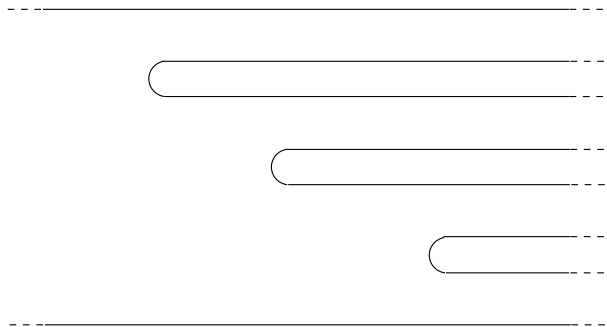}
\end{center}
\caption{A standard domain}
\label{f:stdom}
\end{figure}

We construct different coordinates on $\conf_m$ as follows. Consider $\R^{m-2}$ with coordinates $a=(a_1,\dots,a_{m-2})$. Let $t\in\R$ act on $\R^{m-2}$ by $t(a)=a+t(1,\dots,1)$. The orbit space of this action is $\R^{m-3}$. Define $\Delta_m(a)$ as the subset of $\R\times[0,m-1]$ obtained by removing $m-2$ horizontal slits of width $\epsilon$ , $0<\epsilon\ll 1$, starting at $(a_j,j)$, $j=1,\dots,m-1$ and going to $\infty$. All slits have the same shape, ending in a half-circle, see Figure \ref{f:stdom}.

Endowing $\Delta_m(a)$ with the flat metric we get a conformal structure $\kappa(a)$ on the $m$-punctured disk. Moreover, using the fact that translations are biholomorphic we find that $\kappa(a)=\kappa(t(a))$ for all $a$. If $I$ is a boundary component of $\Delta_m$ which has both of its ends at $+\infty$, we will call the point with smallest real part along $I$ a {\em boundary minimum} of $\Delta_m$.

\begin{lma}\label{l:confmod}
The map $a\mapsto \kappa(a)$ induces a diffeomorphism $\R^{m-3}\to\conf_m$.
\end{lma}

\begin{pf}
We prove the lemma by induction. For $m=3$ it holds trivially. Assume it holds for $m=n-1$. It then follows that the map induced from $\R^{n-3}\to\conf_n$ on the boundary of the space of conformal structures has degree $1$. Thus it suffices to show that the map is an immersion for $m=n$. To this end we argue as follows. Fix the coordinate of the last boundary minimum. Let $\epsilon=(\epsilon_1,\dots,\epsilon_{m-3})$ and let $\phi_{t\epsilon}\colon \Delta_m(a)\to\Delta_m(a+(t\epsilon,0))$ be a diffeomorphism which shifts the $j^{\rm th}$ boundary minimum a distance $t\epsilon_j$. Then the linearized variation of the conformal structures induced by this diffeomorphism is $-i\bar\pa w_\epsilon$, where the vector field $w_\epsilon$ is horizontal along the $j^{\rm th}$ bend if $\epsilon_j\ne 0$.

Consider now the representation of the tangent space of conformal variations on the unit disk with punctures generated by $\bar\pa v_j$ as discussed above.
The map taking a non-compact end $[0,-\infty)\times[0,1]$ to a compact one is $z\mapsto e^{\pi z}$. Thus the pull backs of the vector fields $v_j$ to $\Delta_m(a)$ have support around the punctures and have a Fourier expansion of the form $z\mapsto c_1 e^{-\pi z} + \sum_{n\ge0} c_ne^{n\pi z}$ near their puncture of support. The differential of the map $\R^{m-3}\to\conf_m$ can thus be described as follows. Given $\epsilon$ there exists $u\colon \Delta_m(a)\to\C$ which lies in a Sobolev space with small negative exponential weights at all punctures and which is tangential to the boundary along the boundary, and constants $\lambda_1,\dots,\lambda_{m-3}$ such that
$$
\bar\pa w_\epsilon = \bar\pa u + \sum_j\lambda_j\bar\pa v_j.
$$
The image of the tangent vector $\epsilon$ is then $(\lambda_1,\dots,\lambda_{m-3})$. To show that this map is injective we must thus show that
$$
\bar\pa w_\epsilon =\bar\pa u
$$
implies $\epsilon=0$. To see this we note that the Maslov index of the boundary condition satisfied by $w_\epsilon-u$ is at most $-(m-2)+(m-3)$. (Each boundary minimum contributes $-1$ to the Maslov index of the boundary condition of $u$ and the addition of a locally constant section along the bend raises the Maslov index by at most $1$.) Thus the index of the $\bar\pa$-problem with boundary condition of $w_\epsilon-u$ and with small negative exponential weights at all punctures equals $1-(m-2)+(m-3)=0$. Since the $1$-dimensional Riemann-Hilbert problem has only kernel or only cokernel we find that $w_\epsilon-u=0$ which implies that $\epsilon=0$. The lemma follows.
\end{pf}

We will often write $\Delta_m$ for a representative $\Delta_m(a)$ of the $m$ punctured disk with some conformal structure as described above, suppressing the conformal structure from the notation.
\begin{rmk}\label{r:metricconf}
Occasionally we will discuss distances between different conformal structures on $D_m$. When doing that we will use the norm on $\conf_m$ induced by the isomorphism in Lemma \ref{l:confmod} from the standard norm on $\R^{m-3}$.
\end{rmk}

\begin{rmk}
We will call domains of the form $\Delta_m$ discussed above {\em standard domains}. Except for boundary minima of a standard domain explained above we will also talk about vertical line segments. A {\em vertical line segment} of a standard domain $\Delta_m$ is a line segment of the form $\{\tau\}\times[a,b]$ contained in $\Delta_m$ and with $(\tau,a)$ and $(\tau,b)$ in $\pa\Delta_m$.
\end{rmk}

\subsection{Flow trees}\label{s:ftree}
Before we define flow trees we must describe the local functions induced by a Legendrian submanifold $L\subset J^1(M)$.

\subsubsection{Legendrian submanifolds as multi-valued functions}\label{ss:legmvfunc}
For $m\in M$ we write $J^0_m(M)$ and $J^1_m(M)$ to denote the fibers over $m$ of the fibrations $J^0(M)\to M$ and $J^1(M)\to M$, respectively. Let $L\subset J^1(M)$ be a closed Legendrian submanifold. We say that $L$ is {\em front generic} if it satisfies the following conditions.
\begin{itemize}
\item The map $\Pi\colon L\to M$ is an immersion outside a codimension $1$ submanifold $\Sigma\subset M$. (Recall $\Pi\colon J^1(M)\to M$ is the base projection.)
\item For points $s\in\Sigma\setminus\Sigma'$, where  $\Sigma'$ is a codimension $1$ subset $\Pi_F\colon L\to J^0(M)$ has a standard cusp edge singularity. (Recall $\Pi_F\colon J^1(M)\to J^0(M)$ is the front projection.) That is, there are coordinates $u=(u_1,\dots,u_n)$ around $s$ on $L$ and coordinates  $x=(x_1,\dots,x_n)$ around $\Pi(s)\in M$ such that if $z$ is a fiber coordinate in $J^0(M)$, $\Pi_F(u)=(x(u),z(u))$, where
\begin{align}\notag
x_1(u) &=\tfrac12 u_1^2,\\\notag
x_j(u) &= u_j,\quad j=2,\dots,n,\\\label{e:cusp}
z(u) &= \tfrac13u_1^3 + \beta\tfrac12 u_1^2 +\alpha_2 u_2 +\dots+ \alpha_n u_n,
\end{align}
where $\alpha_j$, $j=2,\dots,n$, and $\beta$ are constants.
\item The map $\Pi\colon(\Sigma\setminus\Sigma')\to M$ is an immersion. We require that this immersion is self-transverse.
\end{itemize}
Note that if $L\subset J^1(M)$ is any Legendrian submanifold it can be made front generic after an arbitrarily small Legendrian isotopy. If $L\subset J^1(M)$ is a Legendrian submanifold such that $\Sigma'\subset\Sigma$ is empty we say $L$ has {\em simple front singularities}. If the dimension of $L$ is $1$ then $L$ can be deformed into a Legendrian submanifold with simple front singularities by an arbitrarily small deformation. When $\dim(L)>2$ we will restrict attention to Legendrian submanifolds with simple front singularities. Henceforth, we will thus assume that all Legendrian submanifolds of dimension $>2$ discussed are front generic and have simple front singularities unless otherwise is explicitly stated. We will treat the case of Legendrian submanifolds with simple front singularities as the main case, and point out the differences in the $2$-dimensional case as needed.

\begin{rmk}\label{r:2Dswt}
When $\dim(L)=2$, a generic Legendrian submanifold has, except for cusp-edge singularities, also swallow-tail singularities. The normal form of such singularities is the following.
\begin{align*}
x_1(u,q)&=u\\
x_2(u,q)&=q\left(\frac13 q^2 +\frac12 u\right),\\
y_1(u,q)&=\alpha+\frac12 q^2+\beta q,\\
y_2(u,q)&=q+\beta,\\
z(u,q)&=\frac14 q^4 + \frac12  u q^2 +\beta\left(q\left(\frac13 q^2+u\right)\right) +\alpha u.
\end{align*}
where $(x_1,x_2)$ are coordinates on $M$, $(y_1,y_2)$ corresponding coordinates in the cotangent fibers, $z$ is a coordinate in the $\R$-direction, $(u,q)$ are coordinates on $L$, and where $\alpha$ and $\beta$ are constants.
\end{rmk}

If $L\subset J^1(M)$ is a Legendrian submanifold with simple front singularities then the projection of the singular set is naturally stratified
$$
\Pi(\Sigma)=\Sigma_1\supset\Sigma_2\supset\dots\supset\Sigma_k.
$$
Here $\Sigma_j$ is the set of self-intersection points of $\Pi\colon\Sigma\to M$ of multiplicity at least $j$. Note that $k\le n$ and that $\Sigma_k$ is a closed submanifold of $M$ of codimension $n-k$. Moreover, $\Sigma_j^\circ=\Sigma_j-\Sigma_{j+1}$ is a submanifold of $M$ which is non-closed if $j<k$.

\begin{rmk}
In the case that $\dim(L)=2$ we consider a similar stratification
$$
\Pi(\Sigma)=\Sigma_1\supset (\Sigma_2\cup\Sigma_2^{\rm sw}),
$$
where $\Sigma_2$ consists of the transverse double points of $\Pi(\Sigma)$ and where $\Sigma_2^{\rm sw}$ consists of all swallow-tail points.
\end{rmk}

Let $m\in M-\Sigma_1$ and assume that $\Pi^{-1}(m)\cap L$ consists of $r$ points. Then there exists some neighborhood $U\subset M-\Sigma_1$ such that $\Pi^{-1}(U)\cap L$ consists of $r$ disjoint open subsets $U_1,\dots, U_r$ of $L$ such that $\Pi\colon U_j\to U$, $1\le j\le r$, are embeddings. In particular, $\Pi_F(U_j)\subset J^0(U)$ is the graph of some function $f_j\colon U\to\R$, and $\Pi_\C(U_j)\subset T^\ast U$ is the graph of $df_j$, $1\le j\le r$. Thus, associated to each point in $\Pi^{-1}(m)\cap L$ is a local function defined in some neighborhood of $m$ in $M$.

If $m\in\Sigma_s$ and if there are $r$ points in $\Pi^{-1}(m)\cap (L-\Sigma)$ then exactly as above there are $r$ local functions defined in some neighborhood of $m$ in $M$. Furthermore, there exists a neighborhood $V$ of $m$ in $\Sigma_1$ and $s$ disjoint open subsets $V_1,\dots, V_s$ of $\Sigma$ such that $\Pi^{-1}(V)\cap \Sigma=V_1\cup\dots \cup V_s$ and such that $\Pi\colon V_j\to M$ are embeddings. In particular, if $W$ is a sufficiently small neighborhood of $m$ in $M$ then $\Pi(V_j)$ subdivides $W$ into two components $W_+^j$ and $W_-^j$, where $W_-^j$ is the component into which the image of any vector in $\krn(d\Pi)\subset T_v V_j$, $v\in V_j$ under the second derivative of $\Pi$ points. Associated to a point in $\Pi^{-1}(m)\cap V_j$ are two local functions $f_j^1$ and $f_j^2$ defined on $W_+^j$. Note that these functions have natural extensions to the closure of $W_+^j$, that they agree on $\pa W_+^j$, and that the limit of $df_j^1$ and $df_j^2$ as $\pa W_+^j$ is approached from points in $W_+^j$ agree. These properties are all straightforward consequences of the local form in \eqref{e:cusp}. For example, for the standard cusp edge in \eqref{e:cusp} the local functions are
\begin{align*}
f^1(x)&=\tfrac13(2x_1)^{\tfrac32} + \beta x_1 +\alpha_2 x_2 +\dots+ \alpha_n x_n,\\
f^1(x)&=-\tfrac13(2x_1)^{\tfrac32} + \beta x_1 +\alpha_2 x_2 +\dots+ \alpha_n x_n,
\end{align*}
defined on $x_1>0$ and the common limit of the differentials along $\{x_1=0\}$ is
$$
\beta dx_1 + \alpha_2 dx_2+\dots+\alpha_n dx_n.
$$

If $m\in M$ then points in $\Pi^{-1}(m)\cap(L-\Sigma)$ will be called {\em smooth points over} $m$ and points in $\Pi^{-1}(m)\cap\Sigma$ will be called {\em cusp points over} $m$.

\begin{figure}[htbp]
\begin{center}
\includegraphics[angle=0, width=3cm]{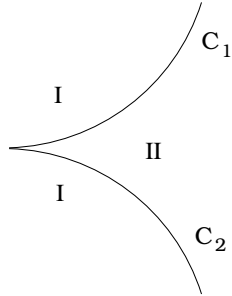}
\end{center}
\caption{The base projection near a swallow tail point.}
\label{f:swall}
\end{figure}

\begin{rmk}
In the case that $\dim(L)=2$ the above description holds outside a neighborhood of $\Sigma_2^{\rm sw}$. A neighborhood of a swallow-tail point is depicted in Figure \ref{f:swall}. The local sheet containing the swallow tail point determines $3$ local functions $f_1$, $f_2$, $f_3$ in the region marked ${\bf II}$ and $1$ function $f$ in the region marked ${\bf I}$. The functions $f_1$ and $f_2$ ($f_2$ and $f_3$) are newborn functions corresponding to the cusp edge $C_1$ ($C_2$).
\end{rmk}

\subsubsection{Local gradients and lifts}\label{ss:lgranddlift}
In order to discuss gradient flows we need a metric. Thus, let $g$ be a Riemannian metric on $M$. For $m\in M$, we consider the vector fields
$$
\nabla(f_1-f_2),
$$
where $f_1$ and $f_2$ are local functions determined by $L$ at $m$ and where $\nabla$ denotes the $g$-gradient. These vector fields are defined on the subsets of $M$ around $m$ where both local functions $f_1$ and $f_2$ are defined. We call them the {\em local gradients at $m$}.

Consider a path $\gamma\colon [0,1]\to M$. A {\em $1$-jet lift} of $\gamma$ is an (unordered) pair $\{\tilde\gamma_1,\tilde\gamma_2\}$ of continuous lifts $\tilde\gamma_j\colon [0,1]\to L\subset J^1(M)$, $j=1,2$ of $\gamma$ such that either $\tilde\gamma_1(t)\ne \tilde\gamma_2(t)$, or $\tilde\gamma_1(t)=\tilde\gamma_2(t)$ and their common value is a cusp point over $\gamma(t)$. We write $\bar\gamma_j=\Pi_\C\circ\tilde\gamma$ and call $\bar\gamma_j$ a {\em cotangent lift} of $\gamma$.

Let $\gamma\colon (-\epsilon,\epsilon)\to M$ be a path with $1$-jet lift $(\tilde\gamma_1,\tilde\gamma_2)$, let $f_1$ and $f_2$ be the local functions specified by the $1$-jet lift, and assume that $\gamma$ satisfies the gradient flow equation
$$
\dot\gamma(t)=-\nabla(f_1-f_2)(\gamma(t)).
$$
Then we say such a path is a {\em flow line} of $L$.

If $\gamma$ is a flow line of $L$ defined on $(-\epsilon,\epsilon)$ and if $\lim_{t\to\pm\epsilon}\tilde\gamma_j(t)\notin\Sigma$, $j=1,2$, then $\gamma$ can be continued as a flow line of $L$ defined on a larger interval. It follows that any flow line of $L$ has a maximal interval of definition. (Note that this maximal interval may consist of only one point.)
\begin{lma}\label{l:maxintvl}
If the maximal interval of $\gamma$ has a non-compact end $(a,\infty)$ or $(-\infty,b)$ then
$$
\lim_{t\to\pm\infty}\gamma(t)=m,
$$
where $m$ is a critical point of some local function difference. If it has a compact end $(a,l]$ or $[l,b)$ then
$$
\lim_{t\to l}\tilde\gamma_j(t)\in\Sigma,\quad j=1\text{ or }j=2.
$$
\end{lma}
\begin{pf}
The proof of this lemma is standard. The only point to check is that the flow time for the flow of two new-born functions near $\Pi(\Sigma)$ to reach $\Pi(\Sigma)$ is finite. To this end let $f_1$ and $f_2$ be the new-born functions at the cusp edge. Pick coordinates as in \eqref{e:cusp} and let the metric be $g_{ij}$ with inverse $g^{ij}$. Then
$$
\nabla (f_1-f_2)=g^{ij}(x)\frac{\pa (f_1-f_2)}{\pa x_j}\pa_{x_i}= \frac32g^{i1}(x)x_1^{\frac12}\pa_{x_i}.
$$
In particular, the solution curve $x(t)$ to the gradient equation satisfies
$$
\frac{d}{dt} \sqrt{x_1}=-\frac34 g^{11}(x).
$$
Since $g^{11}(x)>\delta>0$ for some fixed $\delta$ it follows that the flow time to $\Pi(\Sigma)$ is finite.
\end{pf}

Let $\gamma\colon (a,b)\to M$ be a flow line of $L$ and let $\{\tilde\gamma_1,\tilde\gamma_2\}$ be its $1$-jet lift. We orient the components of the $1$-jet lift as follows.
\begin{dfn}\label{d:flowori}
The {\em flow orientation} of $\tilde\gamma_1$ at $p$ is given by the unique lift of the vector
$$
-\nabla(f_1-f_2)(\Pi(p))\in T_{\Pi(p)} M
$$
to $T_p L$, where $f_1$ is the local function determined by $p$.
\end{dfn}

\subsubsection{Definition and elementary properties of flow trees}\label{ss:ftree}
We first define the source spaces of flow trees. A {\em source-tree} is a tree $\Gamma$ with finitely many edges. Its extra structure is the following. At each $k$-valent vertex $v$, $k\ge 2$ the $k$ edges adjacent to $v$ are cyclically ordered.

\begin{dfn}\label{d:ftree}
A {\em flow tree of $L\subset J^1(M)$} is a continuous map $\phi\colon\Gamma\to M$, where $\Gamma$ is a source-tree which satisfies the following conditions.
\begin{itemize}
\item[{\rm (a)}] If $e$ is an edge of $\Gamma$ then $\phi\colon e\to M$ is an injective parametrization of a flow line of $L$.
\item[{\rm (b)}] Let $v$ be a $k$-valent vertex with cyclically ordered adjacent edges $e_1,\dots, e_k$. Let $\{\bar\phi_j^1,\bar\phi_j^2\}$ be the cotangent lift corresponding to $e_j$, $1\le j\le k$. We require that there exists a pairing of lift components such that for every $1\le j\le k$ (with $k+1=1$)
$$
\bar\phi_j^2(v)=\bar\phi_{j+1}^1(v)=\bar m\in\Pi_\C(L)\subset T^\ast M,
$$
and such that the flow orientation of $\bar\phi_j^2$ at $\bar m$ is directed toward $\bar m$ if and only if the flow orientation of $\bar\phi_{j+1}^1$ at $\bar m$ is directed away from $\bar m$.
\item[{\rm (c)}] The cotangent lifts of the edges of $\Gamma$ fit together to an oriented curve $\bar\phi$ in $\Pi_\C(L)$. We require that this curve is closed.
\end{itemize}
\end{dfn}
For simpler notation we will often denote flow trees simply by $\Gamma$, suppressing the parametrization map $\phi$ from the notation. We will also write $\bar\Gamma$ and $\tilde\Gamma$ for the cotangent- and the $1$-jet lifts of $\Gamma$, respectively.

We next define punctures of a flow tree $\Gamma$. Let $v$ be a $k$-valent vertex of $\Gamma$ with cyclically ordered edges $e_1,\dots,e_k$. Consider two paired cotangent lifts $\bar\phi_j^2$ and $\bar\phi_{j+1}^1$ and the corresponding $1$-jet lifts $\tilde\phi_j^2$ and $\tilde\phi_{j+1}^1$ at $v$. If $\tilde\phi_j^2(v)\ne \tilde\phi_{j+1}^1(v)$ then both must equal Reeb chord endpoints. If this is the case then we say that $v$ contains a {\em puncture} after $e_j$.

\begin{lma}\label{l:1pperv}
Any flow tree $\Gamma$ with a vertex $v$ which contains more than one puncture is a union of flow trees $\Gamma=\Gamma_1\cup\dots\cup\Gamma_m$ such that every vertex of each $\Gamma_j$ contains at most one puncture.
\end{lma}

\begin{pf}
It is enough to show how to split off one flow tree from $\Gamma$ at $v$. Consider the cyclically ordered vertices $e_1,\dots,e_k$ at $v$. Let $f_1$ and $f_2$ be the local functions such that $f_1-f_2$ has a critical point at $\phi(v)=m$. We will identify the local functions with the sheets of the Legendrian submanifold they define (or rather their Lagrangian projections). Assume that notation is chosen so that $\bar\phi_1^1$ maps to $f_1$ and assume for definiteness that $\bar\phi_1^1$ is oriented toward $\bar m=\bar\phi_1^1(v)$. If $\bar\phi_1^2$ maps to $f_2$ then $\bar\phi_1^2$ is oriented away from $\bar m$ and we can split $\Gamma$ at $e_1$. If not then $\bar\phi_1^2$ lies in some other sheet $f_3$ and is oriented away from $\bar\phi_1^2(v)$ therefore $\bar\phi_2^1$ maps to $f_3$ and is oriented toward $\bar\phi_2^1(v)$. Continuing in this manner we find a smallest $r<k$ such that $\bar\phi_r^2$ maps to $f_1$ or $f_2$ and is oriented away from $\bar m$. But then we can split off $e_1,\dots, e_r$.
\end{pf}

In view of Lemma \ref{l:1pperv} we will assume throughout that every vertex of the flow trees we study contains at most one puncture. Hence, to simplify notation, we will say that any vertex which contains a puncture is a puncture.

In some of the proofs below it will be useful to have a notion for parts of flow trees. We define a {\em partial flow tree} in the same way as a flow tree except that condition (c) in Definition \ref{d:ftree} is weakened. More precisely, a partial flow tree is allowed to have $1$-valent vertices $v$ which map to some point on the flow line of the adjacent edge $e$. We say that the partial flow tree has a {\em special puncture} at $v$ and call the vertical line segment connecting its $1$-jet lift at $v$ a {\em special chord}. Note that cutting an edge in a flow tree subdivides it into two partial flow trees with one special puncture each.

\begin{dfn}\label{d:tpunctsign}
Let $p$ be a puncture of a flow tree (a special puncture of a partial flow tree). Let $\tilde\phi^1$ and $\tilde\phi^2$ be the $1$-jet lifts which map to the Reeb chord (special chord) at $p$, with notation chosen so that $\tilde\phi^1$ is oriented toward $\tilde\phi^1(p)$ and $\tilde\phi^2$ oriented away from $\tilde\phi_2(p)$. Then we say that $p$ is a {\em positive} (special) puncture if
$$
z(\tilde\phi^1(p))<z(\tilde\phi^2(p)),
$$
and we say it is {\em negative} is the opposite inequality holds. (Recall, $z$ is the coordinate in the $\R$-direction of $J^1(M)=T^\ast M\times\R$.)
\end{dfn}

\begin{lma}\label{l:ft1+}
Every (partial) flow tree $\Gamma$ in $M$ has at least one positive puncture.
\end{lma}

\begin{pf}
To see this, consider the $1$-jet lift $\tilde\Gamma$ of $\Gamma$. Since $L$ is Legendrian we have $dz-p\,dq=0$ along $L$. In particular, we find that
\begin{equation}\label{e:ftarea}
-\int_{\tilde\Gamma} p\,dq=-\int_{\tilde\Gamma} dz=\sum(z(p^+)-z(p^-))-\sum(z(q^+)-z(q^-)),
\end{equation}
where the first sum runs over the positive punctures of $\Gamma$, the second over the negative punctures, and where $p^+,p^-$ and $q^+,q^-$ denotes the upper and lower endpoints of the corresponding Reeb chords (and special chords) on $L$. On the other hand if $\alpha\colon[0,1]\to M$ is a parametrization of a flow line of $-\nabla(f_1-f_2)$ with $1$-jet lift $\tilde\alpha$ oriented according to the flow orientation then
\begin{equation}\label{e:ftarealoc}
-\int_{\tilde\alpha} dz= f_1(0)-f_1(1)-(f_2(1)-f_2(0))=(f_1(0)-f_2(0))-(f_1(1)-f_2(1))>0.
\end{equation}
It follows that the expression in \eqref{e:ftarea} must be positive, and hence the tree has at least one positive puncture.
\end{pf}

\begin{rmk}
We will refer to the integral in \eqref{e:ftarea} as the {\em symplectic area} of $\Gamma$.
\end{rmk}

Consider a flow tree $\Gamma$ with only one positive puncture. The proof of Lemma \ref{l:ft1+} implies that the function difference along an edge in $\Gamma$ cannot change sign: if it did, by cutting, we could produce a partial flow tree with only negative punctures.

We next discuss multiply covered vertices in a flow tree $\Gamma$. Let $v$ be a $kp$-valent, $k\ge 2$, vertex of $\Gamma$ and let $\Theta\subset\Gamma$ be the partial flow tree consisting of small pieces of the edges $e_1,\dots,e_{kp}$ adjacent to $v$. If there exists a partial flow tree $\Omega$ with a $p$ valent vertex $v'$ such that the $1$-jet lift $\tilde\theta$ of $\Theta$ factors through the $1$-jet lift $\tilde\omega$ of $\Omega$ then we say that the vertex $v$ is {\em multiply covered}.
\begin{dfn}\label{d:simple}
A flow tree without multiply covered vertices will be called {\em simple}.
\end{dfn}
Note that our assumption that each vertex contains only one puncture implies that no vertex which contains a puncture is multiply covered.
\begin{lma}\label{l:simple}
Any (partial) flow tree with exactly one positive puncture is simple.
\end{lma}
\begin{pf}
Assume that $\Gamma$ is a tree with a multiply covered vertex $v$ and with only one positive puncture. Cutting all edges adjacent to $v$ near $v$ we produce a partial flow tree $\Theta$ containing $v$ which covers a tree $\Theta'$ at least twice. The complement of $\Theta$ in $\Gamma$ is a disjoint collection of trees. At least one of the special vertices in $\Theta'$ must be positive by Lemma \ref{l:ft1+} and this means that at least two special vertices of $\Theta$ is positive. Thus in at least two trees in the disjoint collection the special puncture near $\Theta$ is negative. An application of Lemma \ref{l:ft1+} then implies that $\Gamma$ has at least two positive punctures.
\end{pf}

%% file: Sec/3dim.tex
\section{Dimension counts and transversality}\label{S:dimcount}
In Subsection \ref{s:tdimcount} we introduce some notation related to flow trees, discuss a preliminary transversality condition, and define the geometric- and the formal dimension of flow trees. In Subsection \ref{s:ttv} we prove a transversality result for flow trees and as a consequence we obtain a proof of Theorem \ref{t:main} (a). In Subsection \ref{s:ddimcount}, we relate the formal dimension of flow trees to the formal dimension of holomorphic disks.

\subsection{Dimension formulas for trees}\label{s:tdimcount}
Let $L\subset J^1(M)$ be an $n$-dimensional Legendrian submanifold and fix a metric $g$ on $M$.
\subsubsection{Notation for trees and a preliminary transversality condition}\label{ss:nottree}
Let $\Gamma$ be a flow tree of $L$. We subdivide the set of vertices of $\Gamma$ as follows.
\begin{itemize}
\item $P(\Gamma)$ is the set of positive punctures of $\Gamma$.
\item $Q(\Gamma)$ is the set of negative punctures of $\Gamma$.
\item $R(\Gamma)$ is the set of vertices of $\Gamma$ which are not punctures.
\end{itemize}
We will make use of several functions of trees and vertices which we introduce next.
\begin{itemize}
\item If $v$ is a vertex of $\Gamma$ then $\delta(v)$ denotes its valence.
\end{itemize}
An edge of $\Gamma$ is called {\em interior} if it connects two vertices $v_1$ and $v_2$ such that $\delta(v_1)>1$ and $\delta(v_2)>1$.
\begin{itemize}
\item
If $\Gamma$ has at least one vertex $v$ with $\delta(v)\ge 2$ then  $\iota(\Gamma)$ is the number of interior edges of $\Gamma$.
\item
If $\delta(v)=1$ for all vertices $v$ of $\Gamma$ then $\iota(\Gamma)=-1$.
\end{itemize}
Consider the stratification
$$
\Pi(\Sigma)=\Sigma_1\supset\Sigma_2\supset\dots\supset\Sigma_k,\quad \Sigma_j^\circ=\Sigma_j-\Sigma_{j+1}.
$$
\begin{itemize}
\item If $v$ is a vertex of $\Gamma$ then $\sigma(v)=\sigma$, where  $v\in\Sigma_\sigma^\circ$. (We take $\Sigma_0^\circ=M-\Sigma_1$.)
\end{itemize}
Note that $\sigma(p)=0$ for any puncture $p$ and that $\sigma(r)>0$ for any $r\in R(\Gamma)$ with $\delta(r)=1$.
\begin{rmk}
In the case that $\dim(L)=2$ we also define $\sigma(v)=2$ for vertices $v$
which map to $\Sigma_2^{\rm sw}$.
\end{rmk}

Assume that $t$ is a puncture of $\Gamma$ and consider the Reeb chord $c$ corresponding to $t$. Let $c^\pm\in L\subset J^1(M)$ be the endpoints of $c$, where $z(c^+)>z(c^-)$ as usual.
\begin{itemize}
\item
The {\em index} $I(t)$ of $t$ equals the Morse index of the critical point $\Pi(c)$ of $f^+-f^-$, where $f^\pm$ is the local function of the sheet of $L$ containing $c^{\pm}$.
\item
If $t$ is a special puncture then $I(t)=n+1$ if $t$ is positive and $I(t)=-1$ if $t$ is negative.
\end{itemize}

For the next two definitions, we will impose preliminary transversality conditions on $L$ and $g$. (These will also be used in Subsection \ref{s:ttv}.)

To this end we first look at the transversality conditions in two purely local situations. Consider first interactions between a smooth and a cusp sheet. Let $m\in\Pi(\Sigma)\subset     M$ and let $m_1,m_2\in \Pi^{-1}(m)\subset L$. Assume that $m_1\in\Sigma$ and $m_2\notin\Sigma$. Let $U_1$ be a small neighborhood of $m_1$ and let $U_2$ be a small neighborhood of $m_2$ which does not intersect $\Sigma$. (The local sheet of $m_j$ is $S_j=\Pi_\C(U_j)$.) Choose coordinates $(q,s)\in\R^{n-1}\times\R$ around $m\in M$ such that $\Pi(\Sigma\cap U_1)$ corresponds to the subset $\R^{n-1}=\{s=0\}$. If $f_1$ is the local function along $\R^{n-1}$ determined by $S_1$ and $f_2$ that determined by $S_2$ then $\nabla(f_1-f_2)$ is a section of $T^\ast M$ along $\R^{n-1}$. We will require that the subset $W=\{q\colon \nabla(f_1-f_2)\in T(\Pi(\Sigma\cap U_1))=T\R^{n-1}\}$ is a transversally cut out submanifold of $\R^{n-1}$. Picking a normal $\nu$ of $T\R^{n-1}$ we can write any vector field $v$ along $\R^{n-1}$ as $v=v^T+v^\perp\nu$, where $v^T$ is a vector field parallel to $T\R^{n-1}$ and $v^\perp$ is a function. Let $v=\nabla(f_1-f_2)$. The requirement on $W$ is then equivalent to $0$ being a regular value of $v^\perp$, which clearly can be achieved by small Legendrian isotopy of the sheet $S_2$. We will call $W$ the {\em tangency locus of $S_1$ and $S_2$}. Consider second the vector field $v^T$. We require that this vector field is maximally transverse to $W$. In particular there is a stratification $W=W_1\cup W_2\cup\dots \cup W_{n-1}$ such that $v^T$ is transverse to $W$ along $W_1$, $v^T$ has a first order tangency with $W$ along the codimension $1$ subset $W_2$ of $W$, $v^T$ has a second order tangency along the codimension $2$ subset $W_3$, etc. It is not hard to see that each order of tangency increases the codimension by $1$ and that this can be achieved by a small Legendrian isotopy of $S_2$. Furthermore we note that these conditions are open.

Consider second interactions between two cusp sheets. Let $m\in\Pi(\Sigma)\subset M$ and let $m_1,m_2\in \Pi^{-1}(m)\subset L$. Assume that $m_j\in\Sigma$, $j=1,2$. Let $U_j$ be a small neighborhood of $m_j$. Again we denote the local sheet of $m_j$ by $S_j$. Choose coordinates $(q,s_1,s_2)\in\R^{n-2}\times\R^2$ around $m\in M$ such that $\Pi(\Sigma\cap U_j)$ corresponds to the subset $\R^{n-1}_j=\{s_j=0\}$. If $f_j$ is the local function along $\R^{n-2}$ determined by $S_j$ then $\nabla(f_1-f_2)$ is a section of $T^\ast M$ along $\R^{n-2}$. We require that the subsets $A^j=\{q\colon \nabla(f_1-f_2)\in T(\Pi(\Sigma\cap U_j))=T\R^{n-1}_j\}$ are transversally cut out submanifolds of $\R^{n-2}$, which meet transversally. To see that this is possible we pick normals $\nu_j$ of $T\R^{n-1}_j$ and write vector fields $v$ as $v=v^T+v^\perp_1\nu_1+v^\perp_2\nu_2$, where $v^T$ is a vector field parallel to $T\R^{n-2}$ and $v^\perp_j$ are functions. We also require that there are transverse stratifications $A^j=A^j_1\cup \dots \cup A^j_{n-2}$ such that $v^T$ is transverse to $A^j$ along $A_1^j$, $v^T$ has a first order tangency with $A^j$ along the codimension $1$ subset $A_2^j$ of $A^j$, etc.

Consider next the subset $\R^{n-1}_{1+}=\{s_1=0, s_2\le 0\}\subset \Pi(U_2)$. It is not hard to see (for example by making the function $v^\perp_1$ and the vector field $v^T$ independent of $s$ for $s$ very close to $0$) that if the above transversality conditions are satisfied then for any of the two newborn local sheets $S'_2$ of $S_2$ near $\{s_2=0\}$ the tangency locus $W$ of $S_2'$ and $S_1$ is transversally cut out and has $A^1$ as its boundary. Furthermore, the stratification of $W$ is transverse to the boundary and has the stratification of $A^1$ as its boundary $\pa W=\pa A^1$. In particular, the deepest tangency stratum $W_{n-1}$ is empty in a small neighborhood of $\{s_2=0\}$. We point out that it is in general not possible to make the two stratifications corresponding to tangencies between the two newborn sheets of $S_2$ and $S_1$ transverse for $s_2$ close to $0$.

With the local requirements described we next globalize the situation.
Consider the stratification
$$
\Pi(\Sigma)=\Sigma_1\supset\Sigma_2\supset\dots\supset\Sigma_k.
$$
We start over $\Sigma_k$. Here there are $\frac{k(k-1)}{2}$ pairs of cusp sheets and $m$, say, pairs of smooth and cusp sheets. We require that the transversality conditions above hold for any pair and that furthermore all distinct stratified tangency loci are mutually stratum-wise transverse as well as stratum-wise transverse to the stratified self intersection of $\Pi(\Sigma)$. Using the arguments above it is not hard to see that this can be achieved by a small Legendrian isotopy. We next extend the construction stratum by stratum (over $\Sigma_{j}$), requiring in each step that the local conditions above hold as well as stratum-wise transversality between distinct tangency loci and to the self intersection of $\Pi(\Sigma)$, until we reach $\Sigma_1$.

A Legendrian submanifold which satisfies the local transversality conditions as well as the stratum-wise transversality between distinct pairs of sheets and to the self intersections of $\Pi(\Sigma)$ is said to satisfy the {\em preliminary transversality conditions}. The argument above shows that any Legendrian submanifold satisfies these conditions after small Legendrian isotopy. Also, it is not hard to see that the preliminary transversality condition is an open condition. Finally, we remark that the preliminary transversality condition involves the definition of a sheet. More precisely, we must choose a small $\epsilon>0$ such that any connected subset of $\Pi_\C(L)$ of distance no larger than $\epsilon$ from $\Pi_\C(\Sigma)\subset\Pi_\C(L)$  must be considered one sheet. This has to do with the local transversality problems for the two newborn sheets near a cusp edge mentioned above.

It is a consequence of the preliminary transversality condition that no flow line has order of contact with $\Pi(\Sigma)$ which is larger than $n$. (Here the order of contact of a transverse intersection is $1$, the order of contact of an intersection where the distance between points on the flow line and $\Pi(\Sigma)$ approaches $0$ quadratically as the points on the flow line tends to $v$ is $2$, etc.)

\begin{rmk}
In the case when $\dim(L)=2$ and in the presence of swallow tail singularities we assume that all flow lines are transverse to the image of the differential of $\Pi$ at the swallow tail point.
\end{rmk}

Let $v$ be a vertex of $\Gamma$ with $\sigma(v)\ge 0$ and consider a flow line of an edge adjacent to $v$. Such a flow line determines a pair of sheets $S_1$ and $S_2$. If $S_1$ and $S_2$ are the same sheet then we define the order of tangency of $v$ with respect to $S_1$ and $S_2$ to be $0$. (This corresponds to the fact that any edge of a flow tree which ends at a $1$-valent vertex in $\Pi(\Sigma)$ is transverse to $\Pi(\Sigma)$.) If $v$ lies in the image of the cusp edge of exactly one of these sheets, $S_1$ say, then let $W=W_1\cup\dots\cup W_{n-1}$ be the stratification of the tangency locus of $S_1$ and $S_2$, let $W_0=\Pi(\Sigma\cap S_1)-W$, and define the order of tangency of $v$ with respect to $S_1$ and $S_2$ as the number $k$ such that $v\in W_k$. If $v$ lies in the image of the cusp edges of both $S_1$ and $S_2$ then let $A^j=A^j_1\cup\dots A^j_{n-2}$, $j=1,2$, be the two tangency loci described above, let $A^j_0=\Pi(\Sigma\cap S_1)\cap\Pi(\Sigma\cap S_2)-A^j$, and define the order of tangency of $v$ with respect to $S_1$ and $S_2$ as the sum $k_1+k_2$, where $v\in A^1_{k_1}\cap A^2_{k_2}$.

\begin{itemize}
\item If $v$ is a vertex with $\sigma(v)>0$ then $\tau(v)$ equals the sum, over all {\em distinct} pairs of local sheets $S_1$ and $S_2$ determined by edges adjacent to $v$, of the order of tangency of $v$ with respect to $S_1$ and $S_2$.
\item If $\sigma(v)=0$ then $\tau(v)=0$.
\end{itemize}
Note that if $\delta(r)=\sigma(r)=1$ then $\tau(r)=0$.

If $r\in R(\Gamma)$ and if $x\in\Sigma$ is a cusp point over $r$ which lies in the $1$-jet lift of $\Gamma$ then we define $\tilde\mu(x)=+1$ ($\tilde\mu(x)=-1$) if the incoming arc of the $1$-jet lift at $x$ lies in the upper (lower) of the newborn sheets of $L$ and the outgoing arc lies in the lower (upper). In other cases we define $\tilde\mu(x)=0$.
\begin{itemize}
\item
For $r\in R(\Gamma)$, the {\em Maslov content} of $r$ is
$$
\mu(r)=\sum_x\tilde\mu(x),
$$
where the sum runs over all cusp points $x$ in the $1$-jet lift of $\Gamma$ which lies over $r$.
\end{itemize}
Note that if $\sigma(r)=0$ then $\mu(r)=0$ and that if $r\in R(\Gamma)$ is $1$-valent then $\mu(r)=1$.
\begin{rmk}
In the case that $\dim(L)=2$ we define $\mu(r)=0$ if $r\in R(\Gamma)$ maps to $\Sigma_2^{\rm sw}$.
\end{rmk}

\subsubsection{Geometric- and formal dimension of a tree}\label{ss:geofordim}
If $X$ is a finite set then let $|X|$ denote the number of elements in $X$. Let $\Gamma$ be a (partial) flow tree of $L\subset J^1(M)$ where $\dim(L)=\dim(M)=n$.

\begin{dfn}\label{d:tfdim}
The {\em formal dimension} of $\Gamma$ is
\begin{align*}
\dim(\Gamma) &=(n-3) + \sum_{p\in P(\Gamma)}(I(p)-(n-1))
- \sum_{q\in Q(\Gamma)} (I(q)-1)
+ \sum_{r\in R(\Gamma)} \mu(r)\\
&= \sum_{p\in P(\Gamma)}I(p) + \sum_{q\in Q(\Gamma)} (n-I(q))+ \sum_{r\in R(\Gamma)} \mu(r)\\
&-\bigl(|P(\Gamma)|+|Q(\Gamma)|-1\bigr)n + \bigl(|P(\Gamma)|+|Q(\Gamma)|-3\bigr)
\end{align*}
\end{dfn}

\begin{dfn}\label{d:tgdim}
The geometric dimension of $\Gamma$ is
\begin{align*}
\gdim(\Gamma)&=\sum_{\{p\in P(\Gamma)\colon \delta(p)=1\}}I(p)\quad\quad
+ \sum_{\{q\in Q(\Gamma)\colon \delta(q)=1\}}(n-I(q))\\
&+\sum_{\{r\in R(\Gamma)\colon \delta(r)=1\}} (n-(\sigma(r)-1))\quad\quad
-\sum_{\{t\in P(\Gamma)\cup Q(\Gamma)\colon \delta(t)>1\}}\delta(t)n \\
&
-\sum_{\{r\in R(\Gamma)\colon \delta(r)>1\}}((\delta(r)-1)n + \sigma(r)+\tau(r))
\quad+\quad\iota(\Gamma)(n+1).
\end{align*}
\end{dfn}
The reason for the absence of a $\tau(r)$-term in the third sum in the formula for the geometric dimension is that $\delta(r)=1$ implies $\tau(r)=0$. We next relate the two dimension concepts. In order to do so we will first present a lemma, the proof of which utilizes the following observation: $\Sigma\subset L$ is dual to an integer cohomology class. (We may thus define intersection numbers with $\Sigma$.) To see this we argue as follows. The normal bundle of $\Sigma$ is trivial since along the cusp edge in the front we can distinguish locally between the upper and the lower sheets (in the $z$-direction). Using this trivialization we may discuss intersection numbers of oriented paths, transverse to $\Sigma$ in $L$, and $\Sigma$: at an intersection point we compare the trivialization of the normal bundle with the tangent vector of the curve.

\begin{lma}\label{l:m(r)d(r)t(r)}
If $r\in R(\Gamma)$, then $\mu(r)\ge -\delta(r)-\tau(r)+2$.
\end{lma}
\begin{pf}
If $\sigma(r)=0$ the lemma trivial. If $\sigma(r)=1$ and $\delta(r)=1$ then $\tau(r)=0$ and $\mu(r)=1$, and the lemma holds. If $\sigma(r)=1$ and $\delta(r)=2$ then $\mu(r)<0$ only if the flow line on which $r$ lies is tangent to $\Pi(\Sigma)$ at $r$ (if not then the $1$-jet lift does not meet the orientation conditions at the point in the $1$-jet lift over $r$ which lies in a smooth sheet). Thus $\mu(r)\ge-\tau(r)=-\delta(r)-\tau(r)+2$.

Consider the case $\delta(r)=3$, $\sigma(r)=1$ and let $\Theta$ be the patrial flow tree which arises if all edges adjacent to $r$ are cut off close to $r$. Note that the cotangent lift of $\Theta$ consists of $6$ arcs. Let $a_1$ and $a_2$ be two arcs which contributes $-1$ to $\mu(v)$ and consider the arcs $b_1$ and $b_2$ which are the other components of the cotangent lifts of gradient lines corresponding to $a_1$ and $a_2$, respectively. Note that the third edge adjacent to $r$ must be a gradient line between the sheets of $b_1$ and $b_2$ which, since $\sigma(r)=1$, contributes $0$ to $\mu(r)$. Thus, $\mu(r)\ge -1$.

Consider the general case. Again let $\Theta$ denote the partial flow tree around $r$ as above. Since $\mu(r)$ equals the intersection number between the $1$-jet lift of $\Theta$ and $\Sigma$ it does not change under sufficiently small perturbation of $\Theta$. It is easy to see that there exists arbitrarily small perturbations of $\Theta$ into a flow tree $\Theta'$ with vertices only of the following kinds: special $1$-valent vertices, $1$-valent vertices with $\sigma(r)=1$, $2$-valent vertices with $\sigma(r)=1$ and $\tau(r)=1$, and $3$-valent vertices with $\sigma(r)=1$ or $\sigma(r)=0$. The number of special $1$-valent vertices equals $\delta=\delta(r)$. Let $\epsilon$ denote the number of non-special $1$-valent vertices and let $\tau$ denote the number of $2$-valent vertices. An elementary argument using the fact that the Euler characteristic of a tree equals $1$ shows that the number of $3$-valent vertices equals $\delta+\epsilon-2$. By the above, each non-special $1$-valent vertex contributes $1$ to $\mu(r)$, each $2$-valent vertex contributes  $-1$, and each $3$-valent vertex of $\Theta'$ contributes at least $-1$. Hence $\mu(r)\ge-\delta-\tau+2$. Consider a $2$-valent vertex as $\Theta'\to\Theta$, since the condition of a flow line being tangent to a branch of $\Pi(\Sigma)$ is closed it follows that each $2$-valent vertex of $\Theta$ contributes at least $1$ to $\tau(r)$ in the limit. Thus, $\tau(r)\ge \tau$ and $\mu(r)\ge -\delta(r)-\tau(r)+2$.
\end{pf}

\begin{lma}\label{l:fdim>gdim}
If $\Gamma$ is a (partial) flow tree then
\begin{equation}\label{e:fdim>gdim}
\gdim(\Gamma)\le \dim(\Gamma).
\end{equation}
Equality in \eqref{e:fdim>gdim} holds if and only if $\Gamma$ has the following properties.
\begin{itemize}
\item[{\rm (a)}] If $p\in P(\Gamma)\cup Q(\Gamma)$ then it satisfies one of the following three conditions.
\begin{itemize}
\item[{\rm (1)}] $\delta(p)=1$, or
\item[{\rm (2)}] $p\in P(\Gamma)$, $\delta(p)=2$ and $I(p)=n$.
\item[{\rm (3)}] $p\in Q(\Gamma)$, $\delta(p)=2$ and $I(p)=0$.
\end{itemize}
\item[{\rm (b)}] If $r\in R(\Gamma)$ then it satisfies one of the following four conditions.
\begin{itemize}
\item[{\rm (1)}] $\delta(r)=\sigma(r)=1$, or
\item[{\rm (2)}] $\delta(r)=2$, $\sigma(r)=\tau(r)=1$, and $\mu(r)=-1$, or
\item[{\rm (3)}] $\delta(r)=3$ and $\sigma(r)=0$, or
\item[{\rm (4)}] $\delta(r)=3$, $\sigma(r)=\tau(r)=1$, and $\mu(r)=-1$.
\end{itemize}
\end{itemize}
\end{lma}

\begin{rmk}\label{r:vertnot}
We give names to the vertices described in Lemma \ref{l:fdim>gdim}. The vertices in (a) will be called $1$- and {\em $2$-valent punctures}, respectively. The vertices in  in (b) will be called (1) {\em end}, (2) {\em switch}, (3) {\em $Y_0$-vertex}, and (4) {\em $Y_1$-vertex}.
\end{rmk}

\begin{pf}
We prove this in the case that $L$ has simple front singularities. The modifications needed to include swallow-tail points in the $2$-dimensional case are straightforward. We use induction on $\iota(\Gamma)$. If $\iota(\Gamma)=-1$ then the only vertices of the tree are two $1$-valent vertices and it follows immediately that $\dim(\Gamma)\ge\gdim(\Gamma)$.

Consider the case $\iota(\Gamma)=0$. In this case $\Gamma$ has only one vertex $v$ with $\delta(v)\ge 2$. Consider first the case when $v\in R(\Gamma)$. If $\delta(v)=2$ then $\sigma(v)\ge 1$, $\tau(v)\ge 1$ for orientation reasons, and $\mu(v)+\tau(v)\ge 0$ by Lemma \ref{l:m(r)d(r)t(r)}. Thus
$$
\dim(\Gamma)-\gdim(\Gamma)\ge \mu(v)-1+\sigma(v)+\tau(v)\ge 0.
$$
Since $\delta(v)=2$, $|\mu(v)|\le 2$ and $|\mu(v)|\le\sigma(v)$. Hence,
equality holds if and only if $\mu(v)=-1$ and $\sigma(v)=\tau(v)=1$, as claimed. If $\delta(v)\ge 3$ then
\begin{align*}
\dim(\Gamma)-\gdim(\Gamma) &= \mu(v)+\sigma(v)+\tau(v)
+\sum_{\{r\in R(\Gamma)\colon \delta(r)=1\}} (\sigma(r)-1+\mu(r))\\
&+ (|P(\Gamma)|+|Q(\Gamma)|-3)\\
&=
\mu(v)+\sigma(v)+\tau(v)+(\delta(v)-3)+ \sum_{\{r\in R(\Gamma)\colon\delta(r)=1\}} (\sigma(r)-1),
\end{align*}
since $\mu(r)=1$ for each $r\in R(\Gamma)$ with $\delta(r)=1$, and since $\delta(v)=|P(\Gamma)|+|Q(\Gamma)|+|R(\Gamma)-\{v\}|$. If $\sigma(v)=0$ then $\mu(v)=\tau(v)=0$ and we conclude that $\dim(\Gamma)-\gdim(\Gamma)\ge 0$ with equality if and only if $\delta(v)=3$ and each $r\in R(\Gamma)$ with $\delta(r)=1$ satisfies $\sigma(r)=1$. If $\sigma(v)\ge 1$ then, by Lemma \ref{l:m(r)d(r)t(r)}, $\mu(v)+\sigma(v)+\tau(v)+\delta(v)-3\ge \sigma(v)-1$ and again we conclude $\dim(\Gamma)-\gdim(\Gamma)\ge 0$. In this case equality holds if and only if $\sigma(v)=1$, $\mu(v)=-1$, $\tau(v)=0$, and each $r\in R(\Gamma)$ with $\delta(r)=1$ satisfies $\sigma(r)=1$.

Next we consider the case when $v\in P(\Gamma)\cup Q(\Gamma)$. Write $J(v)=I(v)$ if $v\in P(\Gamma)$ and $J(v)=n-I(v)$ if $v\in Q(\Gamma)$. Then
\begin{align*}
\dim(\Gamma)-\gdim(\Gamma)&=
J(v)+\sum_{\{r\in R(\Gamma)\colon \delta(r)=1\}}(\sigma(r)-1+\mu(r))\quad+\quad (|P(\Gamma)|+|Q(\Gamma)|-3)\\
&=J(v)+\sum_{\{r\in R(\Gamma)\colon \delta(r)=1\}}\sigma(r)\quad+\quad (|P(\Gamma)|+|Q(\Gamma)|-3).
\end{align*}
This is non-negative and equals $0$ only if the two second terms sum to $0$ and the first one equals $0$. That is, if and only if $\delta(v)=2$, $J(v)=0$, and each $r\in R(\Gamma)$ with $\delta(r)=1$ satisfies $\sigma(r)=1$.

Assume now that the Lemma holds for all partial trees $\Delta$ with $\iota(\Delta)< N$. Consider a tree $\Gamma$ with $\iota(\Gamma)=N$. Cutting $\Gamma$ along one of its interior edges we split it into two partial flow trees $\Gamma_1$ and $\Gamma_2$, such that $\iota(\Gamma_j)<N$, $j=1,2$. Moreover,
\begin{align*}
\dim(\Gamma)&=\dim(\Gamma_1)+\dim(\Gamma_2)-(n+1),\\
\gdim(\Gamma)&=\gdim(\Gamma_1)+\gdim(\Gamma_2)-(n+1).
\end{align*}
This finishes the proof.
\end{pf}

\subsection{Transversality for trees}\label{s:ttv}
Let $L\subset J^1(M)$ be a Legendrian submanifold with simple front singularities and fix a metric on $M$. Assume that the preliminary transversality conditions are met.
\begin{lma}\label{l:noshort}
There exists $\eta_0>0$ such that if $r\in R(\Gamma)$ is a vertex of a flow tree $\Gamma$ with $\delta(r)=1$ and if $r$ is connected by an edge $E$ to a vertex $v\in R(\Gamma)$ with $\delta(v)=2$ and $\mu(v)< 0$ then the length of $e$ is at least $\eta_0$.
\end{lma}

\begin{pf}
Noting that the gradient flow of the two newborn functions at $m\in\Pi(\Sigma)$ is transverse to the branch of $\Pi(\Sigma)$ corresponding to the lift of the vertex $r$ and that $\Sigma\subset L$ is an embedded submanifold, a short edge $E$ connecting to $v$ with $\delta(v)=2$ and $\mu(v)<0$ can be ruled out.
\end{pf}

The next lemma limits the density of $2$- and certain $3$-valent vertices with $\mu<0$ in a simple flow tree $\Gamma$ such that $\dim(\Gamma)=\gdim(\Gamma)$. Recall the restrictions on the vertices on such flow trees, see Lemma \ref{l:fdim>gdim}, let $\Gamma$ be such a tree, and let $\eta_1\le\eta_0$, where $\eta_0$ is as in Lemma \ref{l:noshort}. Consider a vertex $v$ which is connected by an edge of length at most $\eta_1$ to a vertex $r\in R(\Gamma)$ with $\delta(r)=1$. Note that $\delta(v)=3$ for any such vertex. We say that such a vertex which is a $Y_1$-vertex is an {\em $\eta_1$-close $Y_1$-vertex}. Define the distance between two points in $\Gamma$ to be the length of the shortest path in $\Gamma$ connecting them and recall that $\dim(L)=n$.
\begin{lma}\label{l:dens2}
There exists $\eta_1>0$ such that any subset $A$ of $\Gamma$ of diameter $< \eta_1$ contains at most $n$ switches and  $\eta_1$-close $Y_1$-vertices.
\end{lma}

\begin{pf}
Let $s$ be a switch. Then the local gradient difference $\nabla(f_1-f_2)$, where $f_1$ is the local function of the cusp sheet and where $f_2$ is the local function of the smooth sheet, is tangent to the cusp edge at $s$.

Consider next a sequence of $\eta_1$-close $Y_1$-vertices as $\eta_1\to 0$. Let $y$ be the $Y_1$-vertex and let $e$ be the nearby end. Let the newborn functions of the sheet of $e$ be $f_1'$ and $f_1''$ and the new born functions of the sheet of $y$ be $f_2'$ and $f_2''$. Let $f_1$ and $f_2$ denote the common limits of the corresponding  primed functions along the respective branches $S_1$ and $S_2$ of $\Pi(\Sigma)$ and let $S_2^+$ be the image of the local sheet of $y$ under $\Pi$.
The three edges adjacent to $y$ are flow lines of $\nabla(f_1-f_1')$, $\nabla(f_1'-f_2')$, and $\nabla(f_2''-f_1'')$, respectively. In particular, at $y$, both $\nabla(f_1'-f_2)$ and $\nabla(f_2-f_1'')$ must point into $S_2^+$. As $\eta_1\to 0$, we conclude that the component of $\nabla(f_1-f_2)$   perpendicular to $S_2$ at the limit $y_0$ of $y$ must equal $0$. In other words $\nabla(f_1-f_2)$ is tangent to $S_1$ at $y_0$.

With these local pictures at hand we prove the lemma by contradiction: assume that no $\eta_1>0$ exists and consider $A\subset\Gamma$ with at least $n$ switches and $\eta_1$-close $Y_1$-vertices as $\eta_1\to 0$. After passing to a subsequence, $A$ converges to some point $p\in\Pi(\Sigma)$. However, the local study of the limit above shows that, as $\eta_1\to 0$, each switch and each $\eta_1$-close $Y_1$-vertex gives an independent tangency condition and more than $n$ such tangency conditions contradict the preliminary transversality condition. (In other words, as the switches and $\eta_1$-close $Y_1$-vertices collide we find that some stratum of the intersection of tangency loci and self intersections of $\Pi(\Sigma)$ at $p$, which must be empty by the preliminary transversality condition, would contain $p$.) We conclude that there exists $\eta_1>0$ with properties as claimed.
\end{pf}

\begin{rmk}
The bound $n$ in Lemma \ref{l:dens2} is not optimal. For example, the limit of an $\eta_1$-close $Y_1$-vertex lies in $\Sigma_2$ and in the tangency locus. Also, if there are only switches in $A$ a local study of the limit shows that there can be no more than $\frac{n}{2}$ of them within $A$.
\end{rmk}

\begin{lma}\label{l:tfinite}
There exists a constant $C>0$ such that any flow tree with at most $P$ positive punctures and of formal dimension at most $D>0$ has no more than $C$ edges and vertices.
\end{lma}

\begin{pf}
The case of multiply covered trees follows from the case of simple trees by subdivision (see \S\ref{ss:many+}). Let $\Gamma$ be a simple flow tree with at most $P$ positive punctures.  The symplectic area of the tree is positive and by \eqref{e:ftarea} it equals the sum of the lengths of the Reeb chords at positive punctures minus the sum of the lengths of the Reeb chords at its negative punctures. Since there is only a finite number of Reeb chords all of non-zero length, it follows that the number of negative punctures of the tree is bounded by $k_0P$ for some $k_0>0$.

We will use a slight perturbation $\hat\Gamma$ of the flow tree. Take $\hat\Gamma$ to be a tree close to $\Gamma$ with vertices only of the types of a tree for which the geometric and formal dimensions agree. More precisely, we take $\hat\Gamma$ to be a true flow tree in some neighborhood of $\Pi(\Sigma)$ and near its punctures by perturbing the partial flow trees of $\Gamma$ near any puncture or vertex which is not of the right form. Outside of these regions we construct $\hat\Gamma$ using gradient-like curves with $1$-jet lifts. (That is we relax condition (a) in Definition \ref{d:ftree} but keep the other conditions. A gradient-like curve $\gamma$ with respect to the function $f$ is such that $g(\dot\gamma,-\nabla f)>0$, where $g$ is the metric.) In fact we can choose these curves to lie very close to gradient flow lines. Since $\hat\Gamma$ is a true flow tree close to $\Pi(\Sigma)$, Lemmas \ref{l:noshort} and \ref{l:dens2} apply also to $\hat\Gamma$. Note that the number of edges and vertices of $\hat \Gamma$ is at least as large as the corresponding numbers for $\Gamma$ and that as $\hat\Gamma$ is gradient-like we can still bound its length using the symplectic area.

Pick $\rho>0$ so that Lemma \ref{l:dens2} holds for $\eta_1=\rho$ (in particular, Lemma \ref{l:noshort} holds with $\eta_0=\rho$). Subdivide the $1$-valent vertices of $\hat\Gamma$ as follows. Let ${\mathbf 1}_{\rm l}$ be the set of all punctures and all $1$-valent vertices with adjacent edge of length at least $\rho$, and let ${\mathbf 1}_{\rm s}$ be the set of all other $1$-valent vertices. Outside a neighborhood $W$ of the Reeb chord projections and $\Pi(\Sigma)$ the length of the gradient of any local function difference is bounded from below. Using \eqref{e:ftarealoc} we estimate the symplectic area contribution from a flow line outside $W$ from below by a constant times its length. Thus the global bound on the symplectic area implies that the number of elements $|{\mathbf 1}_{\rm l}|$ in ${\mathbf 1}_{\rm l}$ is bounded by $k_1 P$ for some $k_1>0$.

Let $\hat\Gamma_0$ be the tree obtained by first erasing all edges in ${\mathbf 1}_{\rm s}$ from $\hat\Gamma$ and then forgetting all $2$-valent punctures in the resulting tree. Then the number of $3$-valent vertices of $\hat\Gamma_0$ is
$$
|{\mathbf 1}_{\rm l}|-2\le k_1 P
$$
and its number of edges is
$$
2|{\mathbf 1}_{\rm l}|-3\le 2k_1P.
$$
Note that the fiber difference between any two points in the cotangent lift of a point $p\in\hat\Gamma_0$ which lies outside $W$ and at distance larger than $\rho$ from an end in $\hat\Gamma_0$ is bounded from below. In particular, the total length of that part of the tree is bounded by the symplectic area which in turn is bounded by $k_0P$. If $E$ is an edge of $\hat\Gamma_0$ which contains a switch or a $\rho$-close $Y_1$-vertex, let $E'$ be the part of that edge which lies outside $W$ and outside a $\rho$-neighborhood of the ends of $\hat\Gamma_0$. If $L(E')$ is the length of $E'$ then it follows from Lemma \ref{l:dens2} that the number of switches and $\rho$-close $Y_1$-vertices in $E'$ (and therefore in $E$) is bounded by $(\rho^{-1}L+2)n$. In particular the total number of switches and $\rho$-close $Y_1$-vertices of $\hat\Gamma$ is bounded by $k_2 P$ for some $k_2>0$. Since any $Y_1$-vertex of $\hat\Gamma$ is either a $3$-valent vertex of $\hat\Gamma_0$ or a $\rho$-close $Y_1$-vertex it follows that the total number $N$ of switches and $Y_1$-vertices of $\hat\Gamma$ is bounded by $k_3 P$ for some $k_3>0$.

Since $\dim(\Gamma)=\dim(\hat\Gamma)$ (the formal dimension is determined by the homotopy class of the $1$-jet lift), since each vertex in ${\mathbf 1}_s$ contributes $1$ to $\dim(\hat\Gamma)$, since each vertex in ${\mathbf 1}_{\rm l}$ contributes at least $-1$, since each switch and $Y_1$-vertex contributes $-1$, and since $Y_0$-vertices contributes $0$, we find
$$
D\ge\dim(\Gamma)\ge (n-3)+|{\mathbf 1}_{\rm s}|-|{\mathbf 1}_{\rm l}|-N.
$$
Hence $|{\mathbf 1}_{\rm s}|\le k_4 P+D+3$ for some $k_4>0$. Finally
the number of $3$-valent vertices of $\hat\Gamma$ equals
$$
|{\mathbf 1}_{\rm s}|+|{\mathbf 1}_{\rm l}|-2
$$
and thus the number of $Y_0$-vertices of $\hat\Gamma$ is bounded by $k_6(P+D+3)$. In conclusion the total number of vertices (or edges) of $\hat\Gamma$ is bounded by $C=k(P+D+3)$ for some $k>0$. The lemma follows.
\end{pf}

\begin{rmk}
The modifications needed to prove Lemma \ref{l:tfinite} in the $2$-dimensional case in presence of swallow tail singularities are straightforward. Our preliminary transversality condition implies that that any tree can be changed slightly into a gradient-like tree which avoids an $\epsilon$-neighborhood of the swallow tail points and which is a gradient tree near $\Pi(\Sigma)$. The argument in the above proof can then be used to estimate the number of vertices and edges of this gradient-like tree. Looking at intersections between the $1$-jet lift and $\Sigma$ near the swallow tail point, again using the preliminary transversality condition, it is straightforward to check that the number of vertices added in going back to the original tree is bounded by $2$ times the number of edges with $1$-jet lift passing near the swallow tail points. This number is in turn bounded by twice the finite number of edges in the approximate tree.
\end{rmk}

It follows from Lemma \ref{l:maxintvl} that if $m\in M$ is a critical point of some difference of local functions (i.e. $m$ is the image under $\Pi$ of a Reeb chord of $L$) then we may talk about the stable and unstable manifolds of the gradient flow determined by $L$ at $m$. The main difference from standard Morse theory is that in the present case stable and unstable manifolds may end at $\Pi(\Sigma)$, where the corresponding function difference seize to exist. Also, there might exist manifolds of flow lines which are not related to any critical points, since flow lines can begin and end at $\Pi(\Sigma)$. The preliminary transversality conditions imply that these flow lines meet $\Pi(\Sigma)$ in a locally stable fashion. In particular, the tangency locus is stratified by the order of tangency. Except for these differences the standard properties of flow-manifolds in Morse theory hold. For example, a manifold of flow lines emanating from some compact submanifold and ending at another is in general not compact but it has a natural compactification consisting of broken flow lines connecting the two. If $K\subset M$ is a subset contained in a region where some local function difference $f_1-f_2$ of $L$ is defined then we define the {\em $(f_1-f_2)$-flow-out} as the union of all maximal flow lines of $\nabla(f_1-f_2)$ which pass through $K$.

\begin{prp}\label{p:ttv}
After small perturbation of $L\subset J^1(M)$ and small perturbation of a given metric $g$ on $M$ the following hold. Any simple flow tree $\Gamma$ with at most one special vertex, with less than $P$ positive punctures, and of formal dimension at most $D>0$, is transversely cut out and the space of flow trees with the same geometric properties in a neighborhood of $\Gamma$ is a manifold of dimension $\gdim(\Gamma)$. Moreover, this is an open condition on $L$ and $g$.
\end{prp}

\begin{pf}
We will impose conditions inductively on $L$ and $g$. Note first that the assumption that $\Pi_\C(L)$ has only transverse double points translates into that all critical points of function differences are non-degenerate. We first impose the condition that all stable and unstable manifolds of all critical points meet transversely, that they meet the stratified space $\Pi_\C(\Sigma)=\Sigma_1\supset\Sigma_2\supset\dots\supset\Sigma_k$ transversely, and that they are transverse to the stratified tangency locus in $\Pi(\Sigma)$ discussed above. We next consider the intersections of such manifolds. As mentioned above such intersections have natural compactifications which are stratified. We require that all flow-outs of such intersections and of $\Sigma_j$, $j\ge 2$, meet stable- and unstable manifolds, $\Pi(\Sigma)$ and its stratified subspaces transversely, and also that two such flow outs meet transversely. The inductive definition now continues in the obvious way: we require transversality between flow-objects considered in previous steps and flow-outs of their intersections, and also transversality between intersections of new flow-outs. We call a stratified subset which arises as the intersection of flow manifolds in this inductive construction a {\em flow intersection chain}.

It follows from Lemmas \ref{l:fdim>gdim} and \ref{l:tfinite} that in order to prove transversality properties for flow trees with a bounded number of positive punctures and with bounded geometric dimension one needs only consider a finite number of transversality conditions on flow intersection chains. It is clear that the intersection of a finite number of such conditions specifies an open subset of the space of Legendrian submanifolds and metrics. The fact that the condition is dense can be proved by arguments from ordinary finite dimensional Morse theory, see e.g. \cite{Sm}. Here we only sketch the proof.  Let $m\in M$ be a critical point of a local function difference $f_1-f_2$. Consider a small ball $B_r(m)$ of radius $r>0$ around $m$ and let $\bar W_s(m, f_1-f_2)$ and $\bar W_u(m,f_1-f_2)$ be the intersections of the stable and unstable manifolds of $f_1-f_2$ with $\pa B_r(m)$. If the index of $m$ equals $k$ then $\bar W_s(m,f_1-f_2)$ is an $(n-k-1)$-dimensional sphere and $\bar W_u(m,f_1-f_2)$ a $(k-1)$-dimensional sphere. Let $b\colon (B_{\frac12r}(m)-B_{\frac14r})\to\R$ be a function with critical point at $m$. It is easy to see that for any normal vector field $n$ of $\bar W_\ast(m,f_1-f_2)$ ($\ast\in\{u,s\}$) there exists $b$ such that $\bar W_\ast(m, f_1-f_2+\epsilon b)=W_\ast(m)+\epsilon n$ up to first order in $\epsilon$, where the right hand side is defined using a suitable exponential map. Thus varying $f_1-f_2$ near the critical point we span the normal bundle of its stable and unstable manifolds. It is easy to see from this that manifolds of flow lines can be made transversely cut out by small perturbation and that their normal bundles are spanned by variations of the Legendrian submanifolds near its Reeb chords. Consider next the intersection of two manifolds of flow lines. Note first that we can make the intersection transverse by small perturbation around all Reeb chords involved and then that the normal bundle of the transverse intersection is spanned by the sum of the normal bundles of the intersecting flow manifolds and hence is spanned by variations of the Legendrian submanifold near its Reeb chords. In order to show that the intersection of two distinct flow outs in a flow intersection chain, $K_1$ and $K_2$ say, are transverse we use induction over the strata in combination with the above observations. An inductive argument shows that the normal bundles of the strata of $K_1$ and $K_2$ are spanned by variations of the Legendrian submanifold near its critical points, and a small perturbation of the Legendrian submanifold and the metric over the flow out of $K_1$ ($K_2$) makes the gradient field defining the flow out of $K_2$ ($K_1$) a generic section in $TM|K_1$ with respect to $TK_1$ (of $TM|K_2$ with respect to $TK_2$). It follows from this that we can achieve the transversality properties described above by an arbitrarily small perturbation. The properties of flow-outs with respect to $\Pi(\Sigma)$ can be achieved in a similar way. Since the perturbations in each step can be made arbitrarily small they need not affect the transversality properties achieved in earlier steps of the construction.

With these finitely many transversality conditions satisfied, in order to finish the proof, we need only check that the formula for the geometric dimension of a tree is correct for simple flow trees and that the spaces of flow trees locally are manifolds. We check this by induction on $\iota(\Gamma)$. We start by showing that the formula holds for arbitrary partial flow trees with $\iota(\Gamma)\le 0$.

If $\iota(\Gamma)=-1$ then $\Gamma$ has only two vertices which both have valence $1$. If both are punctures then
$$
\gdim(\Gamma)= I(p)+(n-I(q))-(n+1)=I(p)-I(q)-1.
$$
If $p$ and $q$ are both non-special punctures then $\gdim(\Gamma)$ equals the dimension of the intersection of the unstable manifold at $p$ and the stable manifold at $q$ with $1$ subtracted (for the dimension for the flow line itself). Hence the dimension formula holds in that case. If $p$ ($q$) is a special puncture then the  manifold of flow lines near $p$ ($q$) has dimension $n$. This case is analogous to the case when $p$ is a maximum ($q$ is a minimum) except that there is one more degree of freedom for each special puncture since it is free to move along its flow line. The dimension formula thus holds also in this case because of our index conventions for special punctures.

If we replace $q$ by a vertex in $r\in R(\Gamma)$ above, then
$$
\gdim(\Gamma)= I(p) +(n-(\sigma(r)-1)) - (n+1)=I(p)-\sigma(r),
$$
which equals the dimension of the intersection of the unstable manifold of $p$ with $\Sigma_{\sigma(r)}^\circ$. The modification needed when $p$ is a special puncture is similar to those discussed above. We conclude that the proposition holds for all (partial) flow trees $\Gamma$ with $\iota(\Gamma)=-1$.

We make one more observation regarding special punctures: if $s$ is a special puncture in a flow tree $\Gamma$ (we think of $s$ as a point in the source of the tree) then there is a natural evaluation map $\ev_s\colon\Omega\to M$, where $\Omega$ is a small neighborhood in the space of trees near $\Gamma$: $\ev_s(\Gamma)$ equals the image of $s$ under the map parameterizing $\Gamma$. We note that for trees $\Gamma$ with special puncture $s$ and $\iota(\Gamma)=-1$ a neighborhood $\Omega$ of $\Gamma$ has the form of a product $M^0(s)\times M^1(s)$, where $M^0(s)=\ev_s(\Omega)$ and $M^1(s)$ are disks. (For trees with $\iota(\Gamma)=-1$, $M^1(s)$ is $0$-dimensional if the other puncture of $\Gamma$ is non-special, otherwise it is $1$-dimensional.)

Consider the case of simple partial flow trees with $\iota(\Gamma)=0$. Such trees have exactly one vertex $v$ with $\delta(v)>1$. Assume first that $v$ is a puncture and that all punctures of $\Gamma$ are non-special. Then
\begin{align*}
\gdim(\Gamma)&=\sum_{\{p\in P(\Gamma)\colon \delta(p)=1\}} I(p)
+\sum_{\{q\in Q(\Gamma)\colon \delta(q)=1\}}(n-I(q))\\
&+\sum_{\{r\in R(\Gamma)\colon \delta(r)=1\}}(n-(\sigma(r)-1))-n\delta(v).
\end{align*}
Note that this measures the expected dimension of the space of flow lines from all $1$-valent vertices to the critical point $v$ and that it is smaller than $0$ unless all $r\in R(\Gamma)$ with $\delta(r)=1$ has $\sigma(r)=1$, all $q\in Q(\Gamma)$ with $\delta(q)=1$ has $I(q)=0$, and all $p\in P(\Gamma)$ with $\delta(p)=1$ has $I(p)=n$. Also our transversality assumptions guarantee that there are no such trees unless these conditions are met and that if the conditions are met then the set of such trees for a $0$-manifold. Consider next the case of a partial flow tree with special punctures. In this case our transversality conditions imply that the tree is empty unless all conditions above are met at non-special punctures and that if the conditions are met and if the set of such trees is non-empty then any tree in a neighborhood is determined by the locations of the special punctures on their respective flow lines. Thus the dimension count is correct and a neighborhood $\Omega$ of a tree with special puncture $s$ can again be written as $M^0(s)\times M^1(s)$, where $M^0(s)=\ev_s(\Omega)$ is a $1$-disk and where $M^1(s)$ is a product of $1$-disks with one factor for each special puncture different from $s$.

Assume second that $v$ does not contain a puncture. Then
\begin{align*}
\gdim(\Gamma) &=\sum_{\{p\in P(\Gamma)\colon \delta(p)=1\}} I(p)
+\sum_{\{q\in Q(\Gamma)\colon \delta(q)=1\}}(n-I(q))\\
&+ \sum_{\{r\in R(\Gamma)\colon \delta(r)=1\}}(n-(\sigma(r)-1))
-((\delta(v)-1)n + \sigma(v)+\tau(v)).
\end{align*}
To see that the dimension formula is correct in this case, we start in the case when there are no special punctures. Note that the dimension $d_w$ of the flow manifold emanating from a $1$-valent vertex $w$ is as follows. If $p\in P(\Gamma)$ then $d_p=I(p)$, if $q\in Q(\Gamma)$ then $d_q=(n-I(q))$, and if $r\in R(\Gamma)$ then $d_r=(n-(\sigma(r)-1))$. A neighborhood of trees with the same geometric properties as $\Gamma$ can be parameterized as follows. If $A$ denotes the intersection of the flow manifolds intersected with the  $(n-\sigma(r))$-dimensional manifold $\Sigma_{\sigma(r)}^\circ$ and if $B$ denotes the intersection of the strata of the tangency loci of distinct sheet pairs corresponding to the edges adjacent to $v$ with $\Sigma_{\sigma(r)}^\circ$ then the trees are parameterized by $A\cap B$. Our transversality conditions implies that $A\cap B$ is a manifold of dimension
$$
\sum_p d_p + \sum_q d_q +\sum_r d_r -(\delta(v)-1)n-\sigma(v)-\tau(v)=
\gdim(\Gamma).
$$
The case when $\Gamma$ is a partial flow tree with special vertices follows by the usual modification: the extra degrees of freedom come from moving the special vertices along the flow lines. We conclude that the dimension formula holds for all (partial) flow trees with $\iota(\Gamma)\le 0$. Moreover, we find again that if $s$ is a special puncture of a tree $\Gamma$ as considered above then a neighborhood $\Omega$ of $\Gamma$ in the space of nearby trees has the form $M^0(s)\times M^1(s)$, where $M^0(s)=\ev_s(\Omega)$ and $M^1(s)$ are disks. In this case, $M^0(s)$ is obtained as the product of the flow out along the edge ending at $s$ of (a disk in) the intersection manifold near the vertex $v$ (i.e. the vertex with $\delta>0$) times a small interval in the flow line itself.

Assume inductively that the dimension count is correct for any simple partial tree $\Delta$ with $\iota(\Delta)<N$ and with one special puncture $p$ and that moreover a neighborhood of such a tree can be written as $M^0(p)\times M^1(p)$, where $M^0(p)$ and $M^1(p)$ are disks with properties as above. Consider a partial flow tree $\Gamma$ with $\iota(\Gamma)=N$ and with one special puncture. Let $v$ be the vertex which is connected to the special puncture $s$ of $\Gamma$. Assume first that $v$ is not a puncture. Cut all edges connecting $v$ to a vertex with $\delta>1$ at a point $p_j$ and fix a small $(n-1)$-disk $D(p_j)$ at $p_j$ which is transverse to the flow line on which $p_j$ lies. We obtain a subdivision $\Gamma=\Gamma_0\cup\Gamma_1\cup\dots\cup\Gamma_r$ where $\Gamma_0$ is a partial tree with $\iota(\Gamma_0)=0$ and where each $\Gamma_j$, $j>0$, is a partial flow tree with one special puncture and $\iota(\Gamma_j)<N$. Note that $\delta(v)=r+1$. Let $\Gamma'_j$ be the partial flow tree with special puncture corresponding to the cut points $p_j$ constrained to lie in $D(p_j)$. Then a neighborhood $\Omega_j$ of trees near $\Gamma_j'$ has the form $M^0(p_j)'\times M^1(p_j)$ where a neighborhood of $\Gamma_j$ has the form $M^0(p_j)\times M^1(p_j)$, where $M^k(p_j)$ are disks and $M^0(p_j)'=M^0(p_j)\cap D(p_j)=\ev_{p_j}(\Omega_j)$. Note also that $D(p_j)$ is transverse to $M^0(p_j)$. If $\sigma(v)=\sigma$ and $\tau(v)=\tau$ it is then clear that a neighborhood of the tree $\Gamma$ is given by a product
$$
W\times (-\epsilon,\epsilon)\times M^1(p_1)\times\dots M^1(p_r),
$$
where $W$ is the flow out along the gradient difference determined by the edge of $\Gamma$ ending at $s$ of the intersection of the flow outs of all the $M^0(p_j)'$, $j=1,\dots,r$, with a codimension $\tau$ subset of $\Sigma_\sigma^\circ$ determined by the tangency condition at $v$ and where $(-\epsilon,\epsilon)$ is a small interval in the flow line ending at $s$. Our transversality condition then implies that this is a manifold of dimension
\begin{align*}
&\sum_{j=1}^r M^1(p_j)+\sum_{j=1}^r \dim(M^0(p_j)'+1)-(r-1)n - \sigma -\tau + 1\\
=&\sum_{j=1}^r \gdim(\Gamma_j) -\sigma(r)-\tau(r) -(r-1)n +1\\
=&\sum_{j=0}^r \gdim(\Gamma_j) -(r-1)n +1 + rn - (r+1)(n+1)\\
=&\sum_{j=0}^r \gdim(\Gamma_j) - r(n+1)=\gdim(\Gamma).
\end{align*}
Here the last equality holds since each edge that was cut contributes $1$ to $\iota(\Gamma)$. We conclude that the dimension formula holds. Moreover, the evaluation at the special puncture $s$ of $\Gamma$ satisfies $\ev_s(W\times M^1(p_1)\times\dots M^1(p_r))\approx W\times(-\epsilon,\epsilon)$. We thus find that the desired product decomposition of a neighborhood holds with $M^0(s)\approx D_W\times(-\epsilon,\epsilon)$ and $M^1(p)\approx\Pi_{j=1}^r M^1(p_j)$, where $D_W\subset W$ is a disk. The case when $v$ is a puncture can be treated in a similar way using a computation similar to the one used for high valence punctures in the case $\iota(\Gamma)=0$ above. We conclude that the dimension formula and product decomposition of a neighborhood hold for all partial flow trees with one special puncture.

Finally, consider a flow tree $\Gamma$ with $\iota(\Gamma)>0$. Cutting $\Gamma$ in an edge connecting two vertices with $\delta>1$ at a point $s$, we obtain two partial flow trees $\Gamma_1$ and $\Gamma_2$ with  special punctures $s_1$ and $s_2$, respectively. Fix an $(n-1)$ disk $D(s)$ transverse to the edge and containing $s$. If $\Omega_j=M^0(s_j)\times M^1(s_j)$, $j=1,2$ are neighborhoods of the trees $\Gamma_1$ and $\Gamma_2$, respectively, then a neighborhood of $\Gamma$ is parameterized by
$$
\bigl(M^0(s_1)\cap M^0(s_2)\cap D(s)\bigr)\times M^1(s_1)\times M^1(s_2).
$$
Thus, the transversality conditions implies that a neighborhood of the tree $\Gamma$ is a manifold of dimension
$$
\gdim(\Gamma_1)+\gdim(\Gamma_2)-n-1=\gdim(\Gamma),
$$
where the last equality follows since
$\iota(\Gamma)=\iota(\Gamma_1)+\iota(\Gamma_2)+1$ and since the special punctures of $\Gamma_1$ and $\Gamma_2$ contributes $n+1$ to the respective geometric dimensions. This finishes the proof.
\end{pf}

\begin{rmk}\label{r:simpleimp}
Note that the assumption that the tree is simple is essential in the proof of Lemma \ref{p:ttv}. For a multiply covered vertex the transversality conditions which are needed for the dimension formula to hold cannot be achieved by perturbations of the metric and the Legendrian submanifold since they would involve making a flow manifold transverse to itself. To achieve transversality in this more general setting one would have to perturb the gradient equation itself.
\end{rmk}

\begin{rmk}\label{r:dimrigidsubtree}
Note that if $\Gamma$ is a rigid tree (i.e. $\dim(\Gamma)=\gdim(\Gamma)=0$) and if $\Gamma'$ is a partial flow tree obtained by cutting one of the edges of $\Gamma$ then also its complement $\Gamma''$ is a partial flow tree and
$$
\dim(\Gamma')+\dim(\Gamma'')-(n+1)=0.
$$
Now $\dim(\Gamma')$ and $\dim(\Gamma'')$ both have dimension at least $1$ by our transversality conditions. Therefore $\dim(\Gamma')\le n$. Moreover, if the special puncture of $\Gamma'$ is $s'$ then the transversality condition implies that $\ev_{s'}\colon\Omega'\to M$ is a local embedding. That is,  $\Omega'=M^0(s')$ and $M^1(s)$ is $0$-dimensional. (If this was not the case then $\ev_{s'}(\Omega')\cap \ev_{s''}(\Omega'')=\emptyset$ for dimensional reasons, which contradicts the existence of $\Gamma$.)
\end{rmk}

\begin{pf}[Proof of Theorem \ref{t:main} (a)]
Theorem \ref{t:main} (a) is a consequence of Lemma \ref{l:fdim>gdim} and Proposition \ref{p:ttv}.
\end{pf}

\begin{rmk}\label{r:perpsi}
It will be convenient in later sections to have a metric such that $\Pi(\Sigma)$ self intersects orthogonally. The proof of Proposition \ref{p:ttv} carries over to show that one can achieve the transversality properties for trees as stated there also for large classes of restricted metrics. For example one could require that the transverse self intersections of $\Pi(\Sigma)$ are orthogonal in the metric used. In fact, for a generic metric no rigid flow tree passes through $\Sigma_2$ which has codimension $2$ in $M$. Thus, by the transversality conditions satisfied by a generic metric, changing the metric so that $\Pi(\Sigma)$ self intersects orthogonally by a change supported only near $\Sigma_2$ affects flow outs only in a small neighborhood of a subset of codimension $2$. Hence the space of rigid flow trees need not be affected by such a change.
\end{rmk}

\subsection{Dimension formula for disks}\label{s:ddimcount}
If $\Gamma$ is a flow tree with $p$ positive- and $q$ negative punctures then its $1$-jet lift $\tilde\Gamma$ is a collection of oriented curves which can be thought of as a boundary condition for a holomorphic disk $u\colon D_{p+q}\to T^\ast M$ with boundary on $L$ with positive (negative) punctures mapping to the Reeb chords corresponding to positive (negative) punctures of $\Gamma$. The tangent spaces of $\Pi_\C(L)$ along $\tilde\Gamma$ give boundary conditions for the linearized $\bar\pa_J$-operator. We denote the operator with these boundary conditions $\bar\pa_\Gamma$. For details on such operators in the present set up, see \cite{EES1} and \cite{EES4}. Here we just mention that the operator $\bar\pa_\Gamma$ is Fredholm and note that together with the dimension of the space of conformal structures its index give the expected dimension of the moduli space $\M_\Gamma$ of holomorphic disks with boundaries homotopic to $\tilde\Gamma$ according to the formula
$$
\dim(\M_\Gamma)=\ix(\bar\pa_\Gamma)+(p+q-3),
$$
see e.g. \cite{EES4}.
\begin{prp}\label{p:tfdim=dfdim}
The index of $\bar\pa_\Gamma$ equals
$$
\ind(\bar\pa_\Gamma)= n+\sum_{p\in P(\Gamma)}(I(p)-n)-\sum_{q\in Q(\Gamma)} I(q)+ \sum_{r\in R(\Gamma)}\mu(r).
$$
In particular, the expected dimension of the space of holomorphic
disks with boundary conditions homotopic to $\tilde\Gamma$ equals
$$
\dim(\M_\Gamma)=
n-3 + \sum_{p\in P(\Gamma)}(I(p)-(n-1))-\sum_{q\in Q(\Gamma)}(I(q)-1)+\sum_{r\in R(\Gamma)}\mu(r),
$$
and thus agrees with $\dim(\Gamma)$.
\end{prp}

\begin{pf}
As in \cite{EES2, EES4} the index of $\bar\pa_\Gamma$ equals
$$
n+ m(\hat\gamma),
$$
where $\hat\gamma$ is the path obtained by closing the path $\gamma$ of Lagrangian tangent planes along $\tilde\Gamma$, rotating the incoming tangent plane to the outgoing one in the negative direction (see \cite{EES2, EES4} for details) and where $m$ is the Maslov index. The Fredholm index is unchanged under continuous deformations of the Fredholm operator $\bar\pa_\Gamma$. In particular, it remains constant under the scaling procedure, $L\to L_\lambda$, discussed above. The only intersections of $\hat\gamma$ with the tangent spaces of the fibers happen when $\tilde\Gamma$ passes $\Sigma$ and during the close up connecting the tangent spaces at double points. An easy local check, see \cite{EES1}, shows that for $\lambda>0$ small enough the contribution of each positive puncture $p$ is $-(n-I(p))$ and that of each negative puncture $-I(q)$. Furthermore, each passage of $\Sigma$ in the positive $z$-direction contributes $1$, and of each passage in the negative $z$-direction, $-1$. The proposition then follows by the definition of $\mu$.
\end{pf}

%% file: Sec/4metric.tex
\section{Metrics and perturbations}\label{S:metrpert}
In Subsection \ref{s:Rchandce} we introduce deformations of $L_\lambda\subset J^1(M)$ in a two-step process. In the first step we Legendrian isotope $L_\lambda$ to $\tilde L_\lambda$ and we change the metric near $\Pi(\Sigma)$ into a metric which splits as a product metric in a neighborhood of $\Pi(\Sigma)$. The main difference between $\tilde L_\lambda$ and $L_\lambda$ is that the cusps of $\tilde L_\lambda$ are "rounded". In the second step we deform $\tilde L_\lambda$ further so that it respects the product structure introduced near $\Pi(\Sigma)$ in the first step. This is in general not possible keeping $\tilde L_\lambda$ Legendrian. The result of the second perturbation is therefore an immersed totally real submanifold $\hat L_\lambda\subset T^\ast M$ which lies very close to $\Pi_\C(\tilde L_\lambda)$. The specific features of $\hat L_\lambda$ will be used in Section \ref{S:disktotree} when proving that holomorphic disks limit to flow trees. In Subsection \ref{s:metmapfltr}, we introduce further deformations of $\hat L_\lambda$ near its rigid flow trees. These are designed to guarantee the existence of a "piecewise" holomorphic disk near each rigid flow tree which will be used in Section \ref{S:treetodisk} when proving existence of holomorphic disks near rigid flow trees. The constructions mentioned above are presented for Legendrian submanifolds with simple front singularities, modifications needed for the $2$-dimensional case in the presence of swallow tail singularities are given in \S\ref{ss:metmap2D}. In Subsection \ref{s:acs} we briefly review properties of the almost complex structure on a cotangent bundle induced by a Riemannian metric on the base manifold.

\subsection{The metric and the map near Reeb chords and cusp edges}\label{s:Rchandce}
Let $L\subset J^1(M)$ and let $h$ be a $1$-regular metric on $M$.

\subsubsection{Legendrian isotopies near Reeb chords}\label{ss:LegisoRC}
If $c$ is a Reeb chord of $L$ then $m=\Pi(c)$ is a critical point of some difference of local functions $f_2-f_1$ of $L$. We require that there exists coordinates  $x=(x_1,\dots,x_n)$ in some $\epsilon_0$-neighborhood $N_{\epsilon_0}$ of $m$ such that
\begin{equation}\label{e:morseform}
(f_2-f_1)(x)= c + \sum_{j=1}^n\sigma_jx_j^2,\quad c>0
\end{equation}
where $\sigma_j\ne 0$ for all $j$, and such that the metric $g$ is given by  $\sum_{j=1} dx_j^2$. To see that this can be achieved by small changes of the metric we argue as follows. Use the Morse lemma to find coordinates so that \eqref{e:morseform} holds in some neighborhood of $m$. In fact, applying the standard proof of the Morse lemma, see \cite{Mi} Lemma 2.2, with initial coordinates chosen as Riemann normal coordinates around $m$ it is not hard to see that the change of the metric can be taken $C^1$-small.

We next discuss Legendrian isotopies of the family $L_\lambda$ supported in $J^1(N_{\epsilon_0})$. We require that the new family $\tilde L_\lambda$ has the property that all sheets of $\Pi_\C(\tilde L_\lambda)$, except the sheet corresponding to $\lambda\cdot f_2$, are given by (covariantly) constant sections of $T^\ast M$ in this neighborhood. We achieve this as follows. Let $(x,y)$ be standard coordinates on $T^\ast N_{\epsilon_0}$ and consider the Hamiltonian flow of the function
$$
H(x,y)=\lambda(f_1(x)-df_1(0)x),
$$
where we use the flat coordinates to view $x$ as a tangent vector to $M$ at $0$,
in some neighborhood of $0$ and act with the time $1$-map on all the local sheets of $L_\lambda$ over $M$. Using a cut-off function this can be realized as a small contact isotopy. Note that it makes the differential of the lower sheet covariantly constant and that it preserves all gradients of function differences. With this done we make all sheets over $m$ (except the one corresponding to $\lambda f_2$) covariantly constant near $m$ as follows. In the local coordinates $x$ around $m$ a sheet of $L_\lambda$ corresponding to some local function $\lambda f\ne \lambda f_j$, $j=1,2$, is given by $x\mapsto (x,\lambda df(x),\lambda f(x))$. We replace it in some neighborhood of $0$ by its first order Taylor expansion $x\mapsto (x,\lambda df(0), \lambda(f(0)+df(0)x))$ and interpolate to the original Legendrian submanifold. Since all these deformations can be expressed as isotopies of the front they correspond to Legendrian isotopies.

We will refer to the parameter $\epsilon_0$ determining the size of the neighborhood as a {\em deformation parameter}. We denote the Legendrian submanifold obtained by deforming $L_\lambda$ by $\tilde L_\lambda$. We will consider further deformations of $\tilde L_\lambda$, for simplicity we will keep the notation $\tilde L_\lambda$ through all steps in the deformation.

\subsubsection{Legendrian isotopies near cusp edges - rounding cusps}\label{ss:LegisoCE}
Assume that the metric $h$ on $M$ is such that the self intersections of $\Pi(\Sigma)$ are orthogonal, see Remark \ref{r:perpsi}. Consider the stratification $\Sigma_1\supset\Sigma_2\supset\dots\supset\Sigma_k$ of $\Pi(\Sigma)$. Note that the normal bundle of the self-transverse immersion $\Pi\colon \Sigma\to M$ is trivial, see the discussion preceding Lemma \ref{l:dens2}. Let $U(k,\epsilon_1)\subset M$ be an $\epsilon_1$-neighborhood of $\Sigma_k$ ($\epsilon_1$ is another deformation parameter). We use the vector bundle $N\to\Sigma_k$ (the geometric normal bundle of $\Sigma_k$ in $M$) as coordinates on $U_k$ in the following way. The triviality of the normal bundle of $\Pi$ and the orthogonality of the self intersections of $\Pi(\Sigma)$ imply that there exists a field of $k$ (unordered) everywhere linearly independent sections  $\{v_1(q),\dots,v_k(q)\}$ of $N$ with the property that $\Pi(\Sigma)\cap U_k$ corresponds to
$$
\Spa\{v_2,\dots,v_k\}\cup\Spa\{v_1, v_3,\dots,v_k\}\cup\dots
\cup\Spa\{v_1,\dots,v_{k-1}\},
$$
where $\Spa$ denotes the fiberwise linear span.

Let $b$ be the restriction of the original metric $h$ to the submanifold $\Sigma^k$. Using these two metrics we get a metric $g$ on $U(k,\epsilon_1)$ which can be made arbitrarily close to $h$ by choosing $\epsilon_1$ sufficiently small. Moreover, if $m\in\Sigma^k$ then there are coordinates $(q,s)\in\R^{n-k}\times\R^k$, $s=(s_1,\dots,s_k)$, around $m$ such that in these coordinates
\begin{itemize}
\item $\Sigma_k$ corresponds to the subset $\{s=0\}$.
\item $\Sigma_1$ corresponds to the subset $\cup_{j=1}^k\{s_j=0\}$.
\item The metric $g$ is given by
$$
g(q,s)=b_{ij}(q)\,dq_i\otimes dq_j + ds_j\otimes ds_j,
$$
where the summation convention: sum over repeated indices, is used.
\end{itemize}

Assume inductively that we have defined a metric on a neighborhood $U(r,\frac{\epsilon_1}{2^{k-r}})$ of $\Sigma_r$, $k\ge r\ge j$, such that for any point in $\Sigma_r- U(r+1,\frac{\epsilon_1}{2^{k-r}})$ we have coordinates $(q,s)\in\R^{n-r}\times\R^r$ such that
\begin{itemize}
\item $\Sigma_r$ corresponds to the subset $\{s=0\}$.
\item $\Sigma_1$ corresponds to the subset $\cup_{j=1}^r\{s_j=0\}$.
\item The metric $g$ is given by
$$
g(q,s)=b_{ij}(q)\,dq_i\otimes dq_j + ds_j\otimes ds_j.
$$
\end{itemize}
Then the above construction can be repeated using a tubular neighborhood $U(j-1,\frac{\epsilon_1}{2^{k-j+1}})$ of $\Sigma_{j-1}- U(j,\frac{\epsilon_1}{2^{k-j+1}})$ extending the coordinates and the metric defined in previous steps. We thus have coordinates and a metric $g$ on a $K\epsilon_1$-neighborhood, where $K>0$ is some constant, $U$ of $\Sigma_1$ such that on subsets $U(j,K\epsilon_1)$ (around $\Sigma_j$) of $U$ the metric has the special features above and it can be made arbitrarily $C^1$-close to the original metric by choosing $\epsilon_1$ sufficiently small.

With this product metric and these product coordinates introduced we next Legendrian isotope $L_\lambda$ inside $U\supset\Sigma_1$. (Note that $\tilde L_\lambda$ constructed in \S\ref{ss:LegisoRC} agrees with $L_\lambda$ in $U$ if $\epsilon_0$ and $\epsilon_1$ are sufficiently small).

Let $m\in\Sigma^r$ and let $V$ be a small neighborhood around $m$ in which the product coordinates $(q,s)\in\R^{n-r}\times\R^r$ are defined. There are two types of components in $J^1(U)\cap L$: those along which $\Pi$ is an immersion and those along which $\Pi$ has a fold singularity. We first describe deformations of the former sheets, then the deformations of the later (which are more involved).

Any smooth sheet of $L_\lambda$ is locally parameterized by the $1$-jet extension of some function. In local coordinates $(q,s,p,\sigma,z)$ on $J^1(M)\approx T^\ast M\times\R$, where $(q,s)\in\R^{n-k}\times\R^k$ is as above and where $1$-forms are written as $\sum_j p_j\,dq_j+\sum_k \sigma_k\,ds_k$, it appears as follows
$$
(q,s)\mapsto \Bigl(q,s_1,\dots,s_r,\lambda d_q f,\lambda d_{s_1}f,\dots, \lambda d_{s_r} f, \lambda f(q,s)\Bigr).
$$
The initial step of the deformation is a Legendrian isotopy, supported in $J^1(U)$ which replaces the function $\lambda f(q,s)$, for small $s$, $|s|\ll\epsilon_1$, with its Taylor polynomial of degree $1$ in $s$. That is by,
$$
\lambda a(q)+\sum_{j=1}^r s_j \lambda h_j(q),
$$
where $h_j(q)=\frac{\pa f}{\pa s_j}(q,0)$ and $a(q)=f(q,0)$. The corresponding Legendrian submanifold is then locally parameterized by
\begin{equation}\label{e:Legsmooth}
(q,s)\mapsto \Bigl(q,s_1,\dots,s_k,\lambda d_q a + \sum_j
s_j\,\lambda d_qh_j,\lambda h_1,\dots,\lambda h_r,\lambda(a+\sum_j s_j\,h_j)\Bigr)\in J^1(M).
\end{equation}
Clearly, such an isotopy need not introduce any new Reeb chords of $L$. Note also that this procedure can be extended inductively over the strata of $\Sigma_1$, much like in the construction of the metric above: we start at the deepest stratum $\Sigma_k$ and work our way up to $\Sigma_1$.

We next consider a similar construction for sheets which project with folds. Let $m\in\Sigma^r$ and let $U$ be as above. Consider now a sheet with cusp edge over $U$. For simplicity we assume that the cusp edge is given by $\{s_1=0\}$ in the standard coordinates and that the projection of the sheet lies in $\{s_1\ge 0\}$. We write $s=(u,\hat s)\in\R\times\R^{r-1}$, the sheet then admits a parametrization of the form
$$
(q,u,\hat s)\mapsto \Bigl(q,\frac12 u^2,\hat s,\lambda d_q
f,\lambda \gamma(q,u,\hat s),\lambda d_{\hat s} f, \lambda f(q,s)\Bigr),
$$
where $\frac{\pa f}{\pa u}(q,0,\hat s)=0$, $\gamma=\frac{1}{u}\frac{\pa f}{\pa u}(q,u,\hat s)$, and $\frac{\pa^3 f}{\pa u^3}(q,0,\hat s)\ne 0$. We next make an isotopy normalizing this, replacing all terms with suitable Taylor polynomials. The resulting Legendrian submanifold is locally,
\begin{align*}
(q,u,\hat s)\mapsto\Bigl(& q,\frac12 u^2,\hat s,\lambda (d a+\sum_{j=2}^rs_j\,dh_j
+ \frac12u^2\,d\beta
+\frac13u^3\,d\alpha),\lambda(\beta+u\alpha),\\
&\lambda h_2,\dots,\lambda h_r,
\lambda (a+\sum_{j=r}^rs_j\,h_j+\frac12u^2\beta+\frac13u^3\alpha)\Bigr),
\end{align*}
where $\alpha$, $\beta$, $a$, and $h_2,\dots,h_r$ are functions of $q$ only. Observe that $\alpha(q)\ne 0$.

We will Legendrian isotope this sheet in order to "round" its cusps. The deformation will involve only the terms $\beta+u\alpha$, $\frac13 u^3\,d\alpha$, and $\frac13 u^3\alpha$ in the above local coordinate expression. To explain this isotopy we look at the family $\tilde\theta_\lambda$ of Legendrian space curves
$$
u\mapsto \Bigl(\frac12 u^2,\lambda u, \lambda \frac13u^3\Bigr)
$$
for $-\delta<u<\delta$, where $(x,y,z)$ are coordinates in $3$-space and the contact form is $dz-y\,dx$. Let $\theta_\lambda=\pi\tilde\theta_\lambda$, where $\pi$ is the projection to the $xy$-plane. We will Legendrian isotope the given curve in a small neighborhood of $0$ to a curve with $xy$-projection which consists of a half-circle with straight line segments attached. To see that such an isotopy exists we argue as follows.

Fix $\eta\ll\delta$. Consider the curve $c_\lambda$ which is the union of the half circle $(x_1-\lambda\sqrt{\eta})^2+y^2=\lambda^2\eta$, $x\le \lambda\sqrt{\eta}$, and the two line segments $y=\lambda\sqrt{\eta}$, $\lambda\sqrt{\eta}\le x\le l$. Note that the area bounded by the curve $c_\lambda$ and a vertical line segment connecting its endpoints equals $\frac12\pi\eta\lambda^2+2l\lambda\sqrt{\eta}$ and that it intersects the curve $\theta_\lambda$ at $x=\frac{\eta}{2}$. On the other hand, the curve $\theta_\lambda$ together with a vertical line segment at $x=l$ bounds an area of magnitude $\frac{2^{\frac52}}{3}\lambda l^{\frac32}$. Thus we can find a Hamiltonian isotopy supported in $-\eta< x < 10\eta$ which deforms the curve $\theta_\lambda$ so that it agrees with (a smoothened version of) $c_\lambda$ for $0\le x\le \frac{\eta}{2}$ and so that the first and second derivatives of the image curve are bounded by $K\lambda$ for some constant $K>0$. It follows that we can Legendrian isotope $\tilde\theta_\lambda$ to $\tilde c_\lambda$ with $xy$-projection as desired. Moreover, $\tilde c_\lambda$ can be parameterized as
$$
u\mapsto \Bigl(\frac12 u^2,\gamma_\lambda(u),\psi_\lambda(u)\Bigr),
$$
where $\phi_\lambda(0)=\psi_\lambda(0)=0$, $\gamma_\lambda(u)=\lambda u$ and $\psi_\lambda(u)=\lambda \frac13 u^3$ for $|u|\ge 100\eta$.

We define the Legendrian submanifold $\tilde L_\lambda$ in local coordinates as
\begin{align}\label{e:Legcusp}
(q,u,\hat s)\mapsto\Bigl(& q,\frac12 u^2,\hat s,\lambda (d a+\sum_{j=2}^rs_j\,dh_j
+ \frac12u^2\,d\beta)
+\psi_\lambda(u)\,d\alpha,\lambda\beta+\alpha\gamma_\lambda(u),\\
&\lambda h_2,\dots,\lambda h_r,
\lambda (a+\sum_{j=r}^rs_j\,h_j+\frac12u^2\beta)+\psi_\lambda(u)\alpha\Bigr),
\end{align}
and again using an inductive procedure starting over $\Sigma_k$ and ending over $\Sigma_1$ we can achieve this local form along all of $\Sigma_1$. We denote the resulting deformed Legendrian submanifold $\tilde L_\lambda$. Note that the isotopy from $\tilde\theta_\lambda$ to $\tilde c_\lambda$ can be used to produce an isotopy from $L_\lambda$ to $\tilde L_\lambda$. Moreover, if  the parameter $\eta>0$ is chosen small enough the isotopy between $L_\lambda$ and $\tilde L_\lambda$ need not introduce any new Reeb chords (since the distance between any two sheets of $\Pi_\C(L_\lambda)$ over $\Sigma_1$ are bounded from below by $k\lambda$ for some $k>0$). The parameter $\epsilon_2=\eta>0$ just mentioned is another deformation parameter.

\begin{rmk}
We denote the family of Legendrian submanifolds obtained through the deformations described above $\tilde L_\lambda(\epsilon_0,\epsilon_1,\epsilon_2)$, where $\epsilon_j$, $j=0,1,2$ are the deformation parameters discussed above. Note that as $\epsilon_j\to 0$, the fiber scaling $s_{\lambda^{-1}}(\tilde L_\lambda)$ $C^0$-converges to $L$.
\end{rmk}

\subsubsection{Further isotopies near cusp edges - local product structure}\label{ss:TRisoCE}
We deform the Legendrian submanifolds $\tilde L_\lambda$ in order to make them compatible with the product structure near the singularity set. We use the notation introduced in \S\ref{ss:LegisoCE}. Fix a function $\phi\colon \R\to\R$ such that $\phi(s)=0$ for $|s|<\delta$ and $\phi(s)=s$ for $|s|>2\delta$, $\delta\ll\epsilon_1$.

We first describe initial deformations of the smooth sheets over $\Sigma_1$. The Lagrangian submanifold in \eqref{e:Legsmooth} (the projection of the Legendrian submanifold) is isotoped to the totally real submanifold $\hat L_\lambda^0$ given by
\begin{equation}\label{e:TR0smooth}
(q,s)\mapsto \Bigl(q,s_1,\dots,s_k,\lambda (d_q a + \sum_j
\phi(s_j)\,d_qh_j),\lambda h_1,\dots,\lambda h_r\Bigr)\in T^\ast M.
\end{equation}
We second describe initial deformations of the cusped sheets which are entirely similar. The Lagrangian submanifold in \eqref{e:Legcusp} is isotoped to the totally real submanifold $\hat L_\lambda^0$ given by
\begin{align}\notag
(q,s)\mapsto\Bigl(& q,\frac12 s_1^2,\hat s, \lambda (da
+\sum_{j=2}^r\phi(s_j)\,dh_j + \phi(s_1)^2\,d\beta)+\psi_\lambda(\phi(s_1))d\alpha,\\\label{e:TR0cusp}
&\lambda\beta+\alpha\gamma_\lambda(s_1),\lambda h_2,\dots,\lambda
h_r\Bigr)
\end{align}
Note that it is easy to find local diffeomorphisms $F_\lambda$ of $T^\ast M$ defined in a small neighborhood of the sheet discussed such that $F_\lambda(\hat L_\lambda^0)=\tilde L_\lambda$, and such that (for fixed deformation parameters) $d_{C^1}(F_\lambda,\id)=\Ordo(\lambda)$.

\subsubsection{Further isotopies compensating cusp rounding}\label{ss:TRisocomp}
The perturbations described in \S\ref{ss:LegisoRC} and \S\ref{ss:TRisoCE} do not affect flow lines between two sheets much as long as the corresponding local gradient is transverse to the cusp edge. However, near the tangency locus, the change in the flow is more substantial. The last step in the construction of $\hat L_\lambda$ is a deformation of $\hat L_\lambda^0$ designed to compensate this change.

Let $S_0$ be a sheet of $\hat L_\lambda^0$ intersecting the cusp edge. Consider the tangency locus along the fold of $S_0$ of the local gradient difference of $S_0$ and another sheet $S_1$. We must consider two cases: when $S_1$ is a smooth sheet and when $S_1$ is itself a cusped sheet.

Consider first the case when $S_1$ is smooth and let $m$ be a tangency point of $S_1$ and $S_0$ which lies in $\Sigma_j^\circ$, $j\ge 1$. In a neighborhood of $m$ (where $\phi(s_j)=0$) we have the following local coordinate expressions for $S_1$,
\begin{align*}
(q,s)\mapsto \Bigl(q,s_1,\dots,s_j,\lambda da^0(q),\lambda h_1^1(q),\dots,\lambda h_r^1(q)\Bigr),
\end{align*}
where $q\in\Sigma_j$ and $s_r\in\R$, see \eqref{e:TR0smooth}, and for $S_0$
\begin{align*}
(q,s)\mapsto\Bigl(& q,\frac12 s_1^2,s_2,\dots,s_j \lambda da^0(q)
,\\
&\lambda\beta^0(q)+\alpha^0(q)\gamma_\lambda(s_1),\lambda h_2^0,\dots,\lambda
h_j^0\Bigr),
\end{align*}
respectively. In these coordinates the tangency locus is given by $\{s_1=0\}\cap\{h_1^1(q)=\beta^0(q)\}$. The preliminary transversality condition implies that $0$ is a regular value of $h_1^1(q)-\beta^0(q)$ and therefore $V=\{h_1^1(q)=\beta^0(q)\}$ is a smooth $(n-j-1)$-dimensional submanifold of $\Sigma_j$.

Let $\xi_\eta(q)$ be the vertical vector in the $(s_1,\sigma_1)$-plane with positive $\sigma_1$-component such that the horizontal lines which are part of the curve $(\frac12 s_1^2,\alpha^0(q)\gamma_\lambda(s_1))$, see \eqref{e:TR0cusp}, are given by displacing the $s_1$-axis along $\pm \lambda\xi_\eta(q)$. (We use the subscript $\eta$ to indicate that this vector depends on the deformation parameter $\epsilon_2=\eta$. In fact, $\xi_\eta(q)=(0,\alpha^0(q)\sqrt{\eta})$.)

Let $\alpha$ be a number larger than $\alpha(q)$, for all $q\in\Sigma_1$ and all functions $\alpha$ in local coordinate expressions as above. (Such a number exists since $\Sigma_1$ is compact.) Since $0$ is a regular value of $h_1^1(q)-\beta(q)$ there exists for sufficiently small $\eta>0$ a collar neighborhood $V\times[-\delta,\delta]$ of $V$ such that $h_1^1(v,t)-\beta^0(v,t)$ is monotone in $t$ for all $v$ and such that $|h_1^1(v,\pm\frac12\delta)-\beta^0(v,\pm\frac12\delta)|\ge\alpha\sqrt{\eta}$. In fact there exists constants $0<k_-<k_+$ such that $k_-\sqrt{\eta}\le \delta\le k_+\sqrt{\eta}$. We assume for definiteness that $h^1_1(v,t)-\beta^0(v,t)$ is decreasing in $t$.

\begin{figure}[htbp]
\begin{center}
\includegraphics[angle=0, width=8cm]{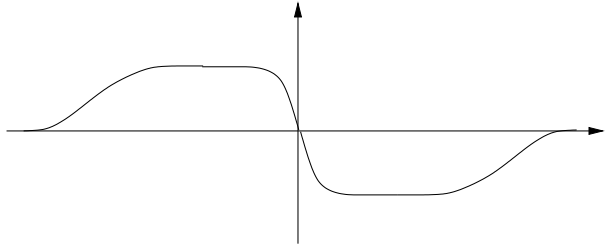}
\end{center}
\caption{The function $\phi_\lambda$.}
\label{f:cutoff}
\end{figure}

Fix a small $\rho>0$. (Here $\rho=\epsilon_3$ is another deformation parameter, which will sometimes be called a {\em speed parameter} to distinguish it from the other deformation parameters. It will be important later on that it is sufficiently small). Let $\phi_\lambda\colon [-\delta,\delta]\to[-1,1]$ be a function which equals $0$ in a neighborhood of $\pm\delta$, which increases from $0$ to $1$ on $[-\delta,-\frac12\delta]$, which equals $1$ on $[-\frac12\delta,-\rho^{-1}\lambda]$, which decreases from $1$ to $-1$ on $[-\rho^{-1}\lambda,\rho^{-1}\lambda]$, which equals $-1$ on $[\rho^{-1}\lambda,\frac12\delta]$, and which increases from $-1$ to $0$ on $[\frac12\delta,\delta]$, see Figure \ref{f:cutoff}. Moreover, we chose this function so that $|\frac{d\phi_\lambda}{dt}|\le 10\delta^{-1}$ for $\frac12\delta\le|t|\le\delta$ and so that $|\frac{d\phi_\lambda}{dt}|\le 2\rho\lambda^{-1}$.

Also, fix a small $\theta>0$ and let $\psi\colon[-\theta,\theta]\to[0,1]$ be a function which equals $0$ in a neighborhood of $|t|=0$, and which equals $1$ on $[-\frac12\theta,\frac12\theta]$.

We redefine the sheet $S_1$ by replacing its local coordinate function $h_1^1(q,s)=h^1_1(q)$ in $V\times[-\delta,\delta]\times[-\theta,\theta]\times\R^{j-1}$ by the function
$$
{\hat h}_1^1(v,t,s_1,\hat s)=h_1^1(v,t)+\lambda\phi_\lambda(t)\psi(s_1)\xi_\lambda(v,t).
$$

Consider next the case when $S_1$ is a cusp sheet as well. In this case we repeat the above construction with the following trivial modification. In the local coordinate expression of $S_1$ replace one coordinate $s_r$, $r\ne 1$, by $\frac12 s_r^2$ and the corresponding function $h_r(q)$, by a cusp function $\beta^1(q)+\alpha^1(q)\gamma_\lambda(s_r)$.

The main reason for this construction is the following. For $|s_1|\le\frac12\theta$, the region where the
horizontal line $\sigma_1={\hat h}_1^1(q,s)$ in the $(s_1,\sigma_1)$-plane intersects the half circular arc in the curve $(\frac12 s_1^2,\alpha^0(q)\gamma_\lambda(s_1))$ lies in a $(\rho^{-1}\lambda)$-neighborhood of the codimension one submanifold $V\times\R^j$. The corresponding region for $\hat L_\lambda^0$ is not as concentrated. It lies only in a uniformly finite neighborhood of $V\times\R^j$.

Note that there exists a fiber preserving diffeomorphism $F_\lambda\colon T^\ast M\to T^\ast M$ which takes the new sheet to the old one such that the $C^0$-distance from $F_\lambda$ to $\id$ is $\Ordo(\lambda)$, and such that the $C^1$-distance is $\Ordo(\rho)$.

Finally, we note that the construction described above is well adapted to the stratification of $\Sigma_1$ and that the deformations defined extend from a deeper stratum to less deep ones. The local deformations of $\hat L_\lambda^0$ described thus match up to give a global deformation. We call the totally real immersed submanifold resulting from this global deformation $\hat L_\lambda$.

We also discuss certain genericity conditions related to the preliminary transversality conditions for $L$. We assume that with notation as above, various tangency loci $V_1$, $V_2$,\dots, $V_m$ corresponding to different sheets over $\Sigma_j^\circ$ intersect transversely in $\Sigma_j^\circ$ for each $j$. Moreover we assume that the flow of any function difference $a_1-a_0\colon \Sigma_j^\circ\to \R$, where $a_0(q)$ and $a_1(q)$ are the local coordinate functions of two sheets over $\Sigma_j^\circ$ as in \eqref{e:TR0smooth} or \eqref{e:TR0cusp} are maximally transverse to all $V_r$. This means in particular that no flow line has order of tangency higher than $(n-j-1)$ with any $V_r$.

\subsubsection{Fiber preserving diffeomorphisms}\label{ss:fibpresdif}
Using the above notation, for small enough deformation parameters, there exists a fiber preserving diffeomorphism $\Phi_\lambda\colon T^\ast M\to T^\ast M$ such that $\Phi_\lambda(\hat L_\lambda)=\tilde L_\lambda$. Moreover, this map can be chosen so that the $C^0$-distance between $\Phi_\lambda$ and $\id$ is $\Ordo(\lambda)$ and the $C^1$-distance is $\epsilon$, for all sufficiently small $\lambda>0$, where $\epsilon\to 0$ with the deformation parameters $\epsilon_j$, $j=0,\dots,3$.

We note that the $C^1$-distance cannot in general be assumed to go to $0$ as $\lambda\to 0$. The reason for this is on the one hand the uniformly finite size of the derivative of the diffeomorphism $F_\lambda$ taking $\hat L_\lambda$ to $\hat L_\lambda^0$ in \S\ref{ss:TRisocomp}. On the other hand the $C^1$-distance for a diffeomorphism taking $\hat L_\lambda^0$ to $\tilde L_\lambda$ can neither be assumed to go to $0$ with $\lambda$: although the local pieces of the diffeomorphism (denoted $F_\lambda$ in \S\ref{ss:TRisoCE}) are of distance $\Ordo(\lambda)$ from the identity the total diffeomorphism must in general remain a finite distance from the identity. Namely, it must change on the $\lambda$-scale over distances of order of magnitude $\lambda$. The fact that the diffeomorphism stays an arbitrarily small finite $C^1$-distance from the identity follows from the fact that the $C^0$-changes on the $\lambda$-scale can be controlled by $k\lambda$ for some $k$ which goes to $0$ with the deformation parameters and that the distances are bounded below by $K\lambda$ for a constant $K$, where $K\lambda$ is the minimum distance between two points in $\Pi_\C(L_\lambda)\cap T^\ast_m M$ for $m\in\Pi(\Sigma)$. Here $K>0$ since there are no Reeb chords near $\Sigma$.

Below we will always assume that the $C^1$-distance $\epsilon$ is sufficiently small. One of the reasons for that is related to certain taming conditions: the immersed totally real submanifold $\hat L_\lambda$ is Lagrangian with respect to the symplectic form $\Phi_\lambda^\ast\omega$, see Lemma \ref{l:Stoke}.

\begin{rmk}
The natural fiber preserving diffeomorphism $\Psi_\lambda\colon T^\ast M\to T^\ast M$ which takes $L_\lambda$ to $\tilde L_\lambda$ is not $C^1$-bounded as $\lambda\to 0$. Close to the cusp edge $\Psi_\lambda$ takes an interval of length $\lambda^{\frac32}$ to an interval of length $\lambda$. Hence its derivative is at least $\lambda^{-\frac12}\to\infty$ as $\lambda\to 0$. Thus, rounding the cusps is in a sense a large deformation. The reason for introducing it will become clear in later sections, see for example Lemma \ref{l:blowup}.
\end{rmk}

\begin{rmk}
As we shall see below the flow trees determined by $\hat L_\lambda$ as $\lambda\to 0$ are not as easily described as flow trees of $L_\lambda$. The reason for changing the nice manifold $L_\lambda$ to $\hat L_\lambda$ is that it is easier to control holomorphic disks with boundary on $\hat L_\lambda$ than the corresponding objects for $L_\lambda$. The difference between trees and disks arises for the following reason. The (local or global) study of trees for a family of Legendrian submanifolds converging to the $0$-section at rate $\lambda\to 0$ is simplified by fiber scaling of the family by $\lambda^{-1}$. On the other hand, the local study of holomorphic disks with boundary on the corresponding family is simplified if both the fiber and the base coordinates are scaled by $\lambda^{-1}$ (this makes the almost complex structure converge to the standard complex structure on $\C^n$).

For example the fiber scaling of the deformation function $\phi_\lambda$ in \S\ref{ss:TRisocomp} converges to a function with a jump-discontinuity at $0$ as $\lambda\to 0$, whereas after scaling both the base- and the fiber coordinate, it converges to a smooth function with derivative bounded by $\rho$ as $\lambda\to 0$.
\end{rmk}

\subsubsection{Limits of flow trees}\label{ss:totrealfltr}
Let $\hat L_\lambda$ denote the totally real immersed submanifold which is the result of the deformations of $L_\lambda$ as discussed above. In later sections, we will study flow trees determined by the totally real (non-Lagrangian) $\hat L_\lambda$. We employ the same definition (Definition \ref{d:ftree}) as in the ordinary case just replacing the local gradient differences at a point $m\in M$ by the metric duals of differences of the intersections $T^\ast_m M\cap\hat L_\lambda$. Note that the flow trees of $\hat L_\lambda$ are not entirely independent of $\lambda$, they change with $\lambda$ in a small neighborhood of $\Sigma_1$. However the flow trees of $\hat L_\lambda$ converges as $\lambda\to 0$. As in Section \ref{S:intr}, we call the limits {\em flow trees of $\hat L_0$}. To describe the limits of flow trees of $\hat L_\lambda$ we fix the deformation parameters and take $\lambda\to 0$. As flow trees are not affected by fiber scaling, fiber scale $\hat L_\lambda$ by $\lambda^{-1}$. The limit of $\hat L_0=\lim_{\lambda\to 0}s_{\lambda^{-1}}(\hat L_\lambda)$ can be described as follows. Outside a neighborhood of $\Sigma_1$ it agrees with the projection of a small Legendrian deformation of $L$. Inside the neighborhood its behavior is less standard. Here $\hat L_0$ is a piecewise smooth submanifold. Local models for it are the following.

Consider the limit of a  smooth sheet $S_\lambda\subset\hat L_\lambda$ lying over a neighborhood of $p\in \Sigma_j^\circ$. Let $V_1,\dots,V_r$ be intersecting tangency loci between $S$ and cusp sheets over $p$. Note that our genericity assumption implies that these sheets intersects in general position. In particular, $r\le\max\{j,n-j\}$. To describe the local model we introduce local coordinates $q=(q_1,\dots,q_{n-j})\in\R^{n-j}$ on $\Sigma_j^\circ$ in a neighborhood of $q$ so that $V_m$ corresponds to the subset $\{q_m=0\}$. Local coordinates on $M$ around $p$ are then $(q,s)=(q_1,\dots,q_{n-j},s_1,\dots,s_j)\in\R^{n-j}\times\R^j$. Assume for simplicity of notation that the projection of the sheet with which $S$ has a tangency along $V_m$ lies in the region $\{s_m\ge 0\}$. In local coordinates the re-scaled sheet $s_{\lambda^{-1}}(S_\lambda)$ is then given by
$$
(q,s)\mapsto (q,s,a(q),\hat h_1(s_1,q),\dots,\hat h_r(s_r,q),h_{r+1}(q),\dots,h_j(q)),
$$
where
$$
\hat h_m(q_1,\dots,q_{n-j},s_m)= h_m(q)-\phi_\lambda(q_m)\psi(s_m)\alpha^m(q)\xi^m_\eta(q).
$$
The limit is only piecewise smooth since $\phi_\lambda$ converges to a function with jump discontinuity. More precisely, it can be described as follows. Identify the fibers of $T^\ast M$  with $\R^n$ equipped with coordinates $(\xi,\sigma)=(\xi_1,\dots,\xi_{n-j},\sigma_1,\dots,\sigma_j)$ dual to $(q,s)$ in the base. Let $\phi_0=\lim_{\lambda\to 0}\phi_\lambda$ be the function with jump discontinuity. (Then $\phi_0(t)=1$ for $t<0$ and $\phi_0(t)=-1$ for $t>0$ if $|t|$ is sufficiently small, we consider $\phi_0(0)$ as undefined.)

We describe the intersection of $\hat L_0$ with a fiber $\R^n$ over a point $q$  with $q_m=0$, $1\le m\le s\le r$, and $q_m\ne 0$, $s+1\le m\le r$. (We choose this particular combination of coordinates which are zero respectively nonzero in order to keep the notation simple, other combinations are similar.)  The intersection lies in a copy of $\R^s$ given by the equations
\begin{align*}
\xi &=a(q),\\
\sigma_{m}&=h_m(q)-\phi_0(q_m)\psi(s_m)\alpha_m(q)\xi^m_\eta(q),\quad s+1<m\le r,\\
\sigma_m &=h_m(q),\quad r+1<m\le j.
\end{align*}
and consists of the $s$-dimensional parallelepiped with corners at the $2^s$ points $(\sigma_1,\dots,\sigma_s)$ with
\begin{align*}
\sigma_m &=h_m(q)+\alpha_m(q)\psi(s_m)\xi^m_\eta(q)\text{ or,}\\
\sigma_m &=h_m(q)-\alpha_m(q)\psi(s_m)\xi^m_\eta(q).
\end{align*}

\begin{figure}[htbp]
\begin{center}
\includegraphics[angle=0, width=8cm]{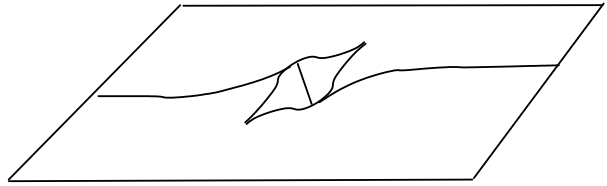}
\end{center}
\caption{The piecewise smooth limit of a sheet over the tangency locus.}
\label{f:slit}
\end{figure}

\begin{figure}[htbp]
\begin{center}
\includegraphics[angle=0, width=4cm]{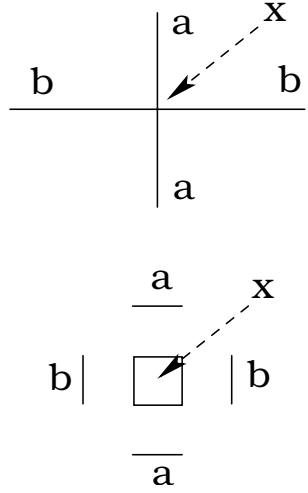}
\end{center}
\caption{The stratification of $\hat L_0$; the upper picture shows a part of $\Pi(\Sigma)$ and the lower depicts the intersections of $\hat L_0$ with the fibers over the corresponding points.}
\label{f:resol}
\end{figure}

Note that these local pieces over $V_1\times\R^{j}\cup \dots V_r\times\R^j$ glue together with the obvious limits of $s_\lambda(\hat L_\lambda)\cap T^\ast \Sigma_j$ outside this locus to an $n$-dimensional stratified space the strata of which are  smooth manifolds. We illustrate the limit in Figures \ref{f:slit} and \ref{f:resol}.

Consider next the limit of a cusp sheet $S_\lambda\subset\hat L_\lambda$ lying over a neighborhood of $p\in \Sigma_j^\circ$. Let $V_2,\dots,V_r$ be intersecting tangency loci between $S$ and other cusp sheets over $p$. Note that our genericity assumption implies that these sheets intersects in general position. In particular, $r\le\max\{j,n-j\}$. To describe the local model we introduce local coordinates $q=(q_1,\dots,q_{n-j})\in\R^{n-j}$ on $\Sigma_j$ in a neighborhood of $q$ so that $V_m$ corresponds to the subset $\{q_m=0\}$. Local coordinates on $M$ around $p$ are then $(q,s)=(q_1,\dots,q_{n-j},s_1,\dots,s_j)\in\R^{n-j}\times\R^j$, where the sheet $S_\lambda$ itself projects to $\{s_1\ge 0\}$. Assume for simplicity of notation the projection of the sheet with which $S$ as a tangency along $V_m$ lies in the region $\{s_m\ge 0\}$. In local coordinates the re-scaled sheet $s_{\lambda^{-1}}(S_\lambda)$ is then given by
$$
(q,s)\mapsto (q,\frac12 s_1^2,\hat s,a(q),\hat \beta(q)+\lambda^{-1}\alpha(q)\gamma_\lambda(s_1),\dots,\hat h_r(s_r,q),h_{r+1}(q),\dots,h_j(q)),
$$
with $\hat h_m$ as above.

\begin{figure}[htbp]
\begin{center}
\includegraphics[angle=0, width=8cm]{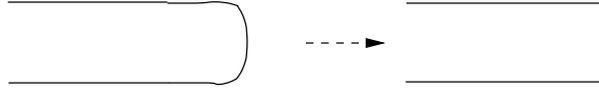}
\end{center}
\caption{The limit of the re-scaling of $\hat L_\lambda$ near the cusp edge.}
\label{f:bendlim}
\end{figure}

Again the limit is only piecewise smooth and we must also take into account that the curve $(\frac12 s_1^2,\lambda^{-1}\gamma_\lambda(s_1))$ converges to a curve with corners, see Figure \ref{f:bendlim}.

More precisely it can be described as follows. Identify the fibers of $T^\ast M$  with $\R^n$ equipped with coordinates $(\xi,\sigma)=(\xi_1,\dots,\xi_{n-j},\sigma_1,\dots,\sigma_j)$ dual to $(q,s)$ in the base. Let $\phi_0=\lim_{\lambda\to 0}\phi_\lambda$ be as above.

We describe the intersection of $\hat L_0$ with a fiber $\R^n$ over a point $q$  with $q_m=0$, $2\le m\le s\le r$, and $q_m\ne 0$, $s+1\le m\le r$ and $s_1=0$. The intersection lies in a copy of $\R^s$ given by the equations
\begin{align*}
\xi &=a(q),\\
\sigma_{m}&=h_m(q)-\phi_0(q_m)\psi(s_m)\alpha_m(q)\xi^m_\eta(q),\quad s+1<m\le r,\\
\sigma_m &=h_m(q),\quad r+1<m\le j.
\end{align*}
and consists of the $s$-dimensional parallelepiped with corners at the $2^s$ points $(\sigma_1,\dots,\sigma_s)$ with
\begin{align*}
\sigma_1 &=\beta(q)+\xi_\eta(q)\text{ or,}\\
\sigma_1 &=\beta(q)-\xi_\eta(q),\\
\sigma_m &=h_m(q)+\alpha_m(q)\psi(s_m)\xi^m_\eta(q)\text{ or,}\\
\sigma_m &=h_m(q)-\alpha_m(q)\psi(s_m)\xi^m_\eta(q).
\end{align*}
At a similar point where $s_1\ne 0$ the image consists of two parallelepipeds in two copies of $\R^{s-1}$ which are defined by the equations for $\R^s$ above and one of the two equations
\begin{align*}
\sigma_1 &=\beta(q)+\xi_\eta(q)\text{ or,}\\
\sigma_1 &=\beta(q)-\xi_\eta(q).\\
\end{align*}

Again, these local pieces over $\{s_1=0\}\cup V_2\times\R^{j}\cup \dots V_r\times\R^j$ glue together with the obvious limits of $s_\lambda(\hat L_\lambda)\cap T^\ast \Sigma_j$ outside this locus to an $n$-dimensional stratified space the strata of which are  smooth manifolds.

In conclusion, $\hat L_0$ is a piecewise smooth submanifold of $T^\ast M$. A local sheet of $\hat L_0$ projects in a $1-1$ fashion to $M$ outside a stratified codimension one subset. The intersection of the local sheet of $\hat L_0$ with $T^\ast_p M$ where $p$ is a point in the codimension $k$ stratum of the stratified simplex is a $k$-dimensional parallelepiped. We call the union of parallelepipeds over this stratified subset a {\em vertical subset} of $\hat L_0$ and we call the preimage of the subset itself the {\em discontinuity locus} of $\hat L_0$.

We next define flow trees determined by $\hat L_0$. First we must define the local gradients. Initially these are well defined for sheets outside the discontinuity locus. We extend the definition over the discontinuity locus by associating to a point $m$ in the discontinuity locus the center of mass in the parallelepiped $\hat L_0\cap T^\ast M$. The local gradient of a sheet then has jump-discontinuities at all strata of the discontinuity locus. We then define local gradient differences and flow lines. Again we note that any flow line admits a cotangent lift into $\hat L_0$, where we take the cotangent lift to be defined by the local gradients except at endpoints of an interval where we define the lift by taking the limit. The cotangent lift of a flow line is then a pair of curves with jump discontinuities at the discontinuity locus of a sheet.

To define flow trees determined by $\hat L_0$ we employ Definition \ref{d:ftree} with the following modifications. We keep condition (a). We change condition (b) as follows. Instead of requiring equality of $\bar \phi^2_j(v)$ and $\bar \phi^1_{j+1}(v)$ for $\phi(v)$ such that $T^\ast_{\phi(v)}M$ intersects the vertical subset of $\hat L_0$, we require that the two cotangent lifts can be joined by a line segment in $T^\ast_{\phi(v)} M\cap \hat L_0$. We also keep the condition that the flow orientation of one arc in each pair is toward its boundary point in the fiber and the other one away from it. We define the cotangent lift of the tree to include these fiber line segments. Condition (c) then remains unchanged.

\begin{lma}\label{l:ftrel}
Fix a speed parameter $\epsilon_3>0$. Any sequence of flow trees determined by $\hat L_0(\epsilon_0,\epsilon_1,\epsilon_2)$ has, as the deformation parameters $\epsilon_j\to 0$, $j=0,1,2$, a subsequence which converges to a flow tree of $L$. Moreover, near any rigid flow tree of $L$ there exists a unique rigid flow tree of $\hat L_\lambda(\epsilon_0,\epsilon_1,\epsilon_2)$ for all sufficiently small $\epsilon_j$ and $\lambda$. In particular for small enough deformation parameters there is a 1-1 correspondence between rigid flow trees of $L$ and rigid flow trees of $\hat L_0$.
\end{lma}

\begin{pf}
The flow lines of $\hat L_0$ and the flow lines of $L$ agree outside a neighborhood of $\Pi(\Sigma)$ and the Reeb chords. The size of this neighborhood goes to $0$ with the deformation parameters. Moreover, outside a neighborhood of the Reeb chords, the length of the local gradient differences is uniformly bounded.

We first consider the situation near a Reeb chord. Since the metric $g$ $C^1$-converges to $h$ and since $\hat L_0$ $C^0$-converges to $L$ it follows that the stable and unstable manifolds of $\hat L_0$ $C^1$-converges to the corresponding objects for $L$ in a neighborhood of each Reeb chord. Since all local flows of $\hat L_0$ and $L$ agree outside a neighborhood of the Reeb chords and $\Pi(\Sigma)$ we find that parts of a trees of $\hat L_0$ outside a neighborhood of $\Pi(\Sigma)$ $C^1$-converges to parts of trees of $L$.

Consider the situation near $\Pi(\Sigma)$. Here we must consider two cases separately. First there are flow lines of $\hat L_0$ between two sheets corresponding to a newborn function difference inside the region where the deformation is supported. Such flow lines may lie entirely in the discontinuity locus. However, as the deformation parameters go to $0$ the length of such lines go to $0$. For other parts of a flow tree, the preliminary transversality condition implies that the time it spends in the deformation region goes to $0$ as the deformation parameters go to $0$. (In local coordinates, the preliminary transversality condition translates into the flow lines along $\Sigma_r^\circ$ of the functions $a(q)$ in \eqref{e:TR0smooth} and \eqref{e:TR0cusp} being as transverse as is possible to the tangency locus $h(q)=\beta(q)$. Thus, the time that a flow line spends near any "cut" in $\hat L_0$ goes to $0$ with the deformation parameters.) Thus the length of the part of a tree of $\hat L_0$ where the tree is not $C^1$-converging to a tree of $L$ goes to $0$. The first statement follows.

Consider the second statement. To construct rigid flow trees of $\hat L_\lambda$ from those of $L$ we cut the latter, subdividing them into small partial flow trees in the deformation regions and larger pieces outside. If the deformation parameters are small enough, the vertices of such a partial tree have one of the following forms.
\begin{itemize}
\item[{\rm (1)}] Two special punctures on the boundary of the deformation region.
\item[{\rm (2)}] One puncture at a Reeb chord $p$ and one special puncture on the boundary of the deformation region.
\item[{\rm (3)}] One $3$-valent puncture in $\Sigma_1$ with $\mu=-1$ and three special punctures on the boundary of the deformation region.
\item[{\rm (4)}] One $2$-valent puncture in the tangency locus in $\Sigma_1$ and two special punctures on the boundary of the deformation region.
\item[{\rm (5)}] One $1$-valent puncture in $\Sigma_1$ and one special puncture on the boundary of the deformation region.
\end{itemize}
Note that the space of such local partial trees for both $L$ and $\hat L_\lambda$ have dimensions given by the formula in Definition \ref{d:tfdim}. (One must modify the formula slightly, taking into account the property that the special punctures are constrained to lie on the boundary of the deformation region.) It then follows from the transversality of the defining conditions of the rigid flow trees for $L$ and the fact that outside the deformation region $\Pi_\C(L_\lambda)$ and $\hat L_\lambda$ agree that for sufficiently small deformation parameters there is exactly one rigid flow tree of $\hat L_\lambda$ near any rigid flow tree of $L$.
\end{pf}

\subsection{The metric and the map near rigid flow trees.}\label{s:metmapfltr}
It follows from the transversality conditions in Proposition \ref{p:ttv} that the following holds. If $\dim(M)\ge 3$ and if $e$ and $e'$ are distinct edges of some rigid flow trees determined by $L\subset J^1(M)$ then either $e$ and $e'$ are disjoint or they are both subsets of the same gradient line. If $\dim(M)=2$ then we will also have transverse intersection points between distinct edges of flow trees. In this subsection we will deform $\hat L_\lambda$ and the metric constructed in Subsection \ref{s:Rchandce} further. Lemma \ref{l:ftrel}, implies that for small enough deformation parameters all rigid flow trees of $\hat L_\lambda$ lie very close to flow trees of $L$. In particular every rigid flow tree $\Gamma$ of $\hat L_\lambda$ has the properties listed in Lemma \ref{l:fdim>gdim} and its edges has the properties discussed above. We use the notation introduced in Remark \ref{r:vertnot} for the vertices of a rigid flow tree.

\subsubsection{The metric and the map near $Y_0$-vertices}\label{ss:metmapY0}
Let $m\in M$ be a $Y_0$-vertex of some rigid flow tree. Define the metric $g$ to be flat in some neighborhood of $m$. Note that this can be done without affecting the gradient flow much: if $x\in\R^n$ are local coordinates around $m$ take $\hat g_{ij}(x)= g_{ij}(0)$ in some neighborhood of $0$.

As in \S\ref{ss:LegisoRC}, we then isotope $\hat L_\lambda$ so that all the local gradient differences determined by $\hat L_\lambda$ are (covariantly) constant in some small neighborhood around $m$.

\subsubsection{The metric and the map near $Y_1$-vertices}\label{ss:metmapY1}
Let $m\in\Pi(\Sigma)$ be a $Y_1$-vertex of some rigid flow tree. In a small neighborhood of $m$, we keep the local product structure of the metric but make it flat in the direction tangential to $\Sigma_1$. Furthermore we make all sheets except the local sheet which intersect $\Sigma$ covariantly constant. The sheet intersecting $\Sigma$ is made constant in the directions parallel to $\Sigma_1$.

\subsubsection{The map over flow lines close to punctures}\label{ss:mapfl}
Let $\gamma$ be a flow line, some part of which is part of a rigid flow tree ending at, beginning at, or (in the $2$-valent case) passing a critical point $m$. Pick point(s) $m$ on $\gamma$ near the critical point (inside the coordinate region where all sheets are covariantly constant) and make the differential of the local function difference of which $m$ is a critical point constant there. In local coordinates: if
$$
(f_1-f_2)(x)=c+\sum_j \sigma_jx_j^2
$$
and we flow along the $x_1$-axis. Isotope the Legendrian slightly so that it looks like
$$
(f_1-f_2)(x_1^0+h_1,h_2,\dots,h_n) = c + \sigma_1 (x_1^0)^2 + 2\sigma_1 x_1^0 h_1
$$
for small $h$, where $x_1^0$ is small.

\subsubsection{The metric and the map near switches and ends}\label{ss:metmapswitch}
\begin{figure}[htbp]
\begin{center}
\includegraphics[angle=0, width=8cm]{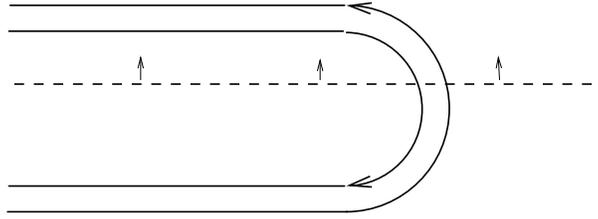}
\end{center}
\caption{The image of the holomorphic function $h$, the inner curve is the bend. The dashed line represents the moving line $\{\sigma=\hat h(\tau)\}$}.
\label{f:imageh}
\end{figure}

Near a switch, the metric is a product metric: $g(q,s)=b(q)+ds^2$, see \S\ref{ss:LegisoCE}. As in \S\ref{ss:metmapY1}, we make the tangential metric $b(q)$ flat, and make the two $\Sigma_1$-components of $\hat L_\lambda$ of the flow line covariantly constant, as well as the functions $\beta(q)$ and $\alpha(q)$. (For simplicity we assume $\beta(q)=0$ and $\alpha(q)=1$ below.) We also perturb the function $h(q)$ giving rise to the tangency so that the tangency locus is locally an affine subspace in the flat coordinates just discussed. Finally, we impose restrictions on the "jump function" $\hat h(q)$. Let $X$ be the tangency locus in $\Sigma_1$, let $k$ denote $\lambda^{-1}\cdot\nabla(a_1-a_0)$, and consider the speed parameter $\rho$ and the radius $\eta^{\frac12}$ of the re-scaled half-circle circle in the bend. We take $\hat h(q)$ to be constantly equal to $\pm \frac{1001}{1000}\sqrt{\eta}$ in a $\rho^{-1}\lambda$ neighborhood of $X$. Consider next a holomorphic function $s\colon\R\times[0,k]\to\C$ which takes $0$ to the extremum of the bend and has image as shown in Figure \ref{f:imageh}. (Note that this function converges exponentially to the function  $\sqrt{\eta}+\frac{1}{1000}\sqrt{\eta}z$ as $z\to\pm\infty$.)
We define the function $\hat h(q)$ by defining it along gradient lines of $\nabla(a_1-a_0)$: we take $\hat h(\tau)$ to equal the imaginary part of $s(\tau+ik)$. Note that this function has the properties required of the functions $\hat h$, in particular, it is constant outside a neighborhood of size $c\lambda$ for some $c>0$ of $X$.

Near an end we make the sheet covariantly constant in the $\Sigma_1$-direction and we make the functions $\alpha(q)$ and $\beta(q)$ constant.

\subsubsection{Plateau points and inflection points}\label{ss:ppip}
The following points, where the metric and $\hat L_\lambda$ were prescribed above will be called {\em plateau points}. Fix $K>0$.
\begin{itemize}
\item Points on the tree at small finite distance from $Y_0$-vertices and $Y_1$-vertices, see \S\ref{ss:metmapY0} and \S\ref{ss:metmapY1}.
\item Points in the deformation regions at small distance $\rho>0$ from Reeb chords, see \S\ref{ss:mapfl}.
\item Points at small finite distance (inside the deformation region) from end vertices and switches, where the difference of the local functions have constant differential, see \S\ref{ss:metmapswitch}.
\end{itemize}
By the construction of $\hat L_\lambda$ and $g$ above, all local functions have (covariantly) constant differentials in some neighborhood of each plateau point.

Let $e$ be an edge in some rigid flow tree. Then $e$ is a part of a gradient flow line $\gamma$. There is a maximal subset of $\gamma$, consisting of subintervals $I$ of $\gamma$, such that all points in this subset lies in some rigid flow tree. Moreover, there are a finite number of plateau points on $\gamma$. Choose two distinct points between any two plateau points on $\gamma$. Call these points {\em inflection points}.

\subsubsection{The metric and the map near inflection points}\label{ss:metmapip}
Let $p$ be an inflection point on a gradient flow line of some local function difference $f_1-f_2$. We use the flow of the vector field $\nabla(f_1-f_2)/|\nabla(f_1-f_2)|^2$ to change the metric in some small neighborhood of $p$ in $M$ and construct flat coordinates $(x_1,x_2,\dots,x_n)$, so that in these coordinates
$$
(f_1-f_2)(x)=k_p x_1 + c_p,
$$
where $k_p$ and $c_p$ are constants. (We take the metric to be the standard flat metric in these coordinates and note that the gradient flows in the original metric and the flat one agree (geometrically) in a neighborhood of $\gamma$. Note also that in these coordinates $d(f_1-f_2)=k_pdx_1$ is (covariantly) constant. Also, after scaling the flat metric we are free to choose $k_p$. We choose $k_p=k$, where $k dx_1=d(f_1-f_2)$ in the standard coordinates at the plateau point closest to $p$ along $\gamma$.

Let $q$ be the inflection point immediately following $p$ along the flow line $\gamma$ oriented by $\nabla(f_1-f_2)$. Fix a sufficiently small $\delta>0$. We isotope $\hat L_\lambda$ slightly (without changing gradient lines of rigid flow trees geometrically) so that for $0<x_1<\delta$
\begin{equation}\label{eqepnf1}
(f_1-f_2)(x)=k_p x_1 +c_p -\alpha x_1^2,
\end{equation}
in the local coordinates $x$ at $p$, and so that for $-\delta<y_1<0$,
\begin{equation}\label{eqepnf2}
(f_1-f_2)(y)=k_q y_1 + c_q -\alpha y_1^2,
\end{equation}
in the local coordinates $y$ at $q$, where $\alpha=\frac{k_q^2-k_p^2}{c_q-c_p}$. (Note that $c_q-c_p\ne 0$ since we follow a gradient line and that we can (and will) take $f_1-f_2$ so that it depends only on $x_1$ throughout the neighborhood.) We require this to hold for $0<x_1<\delta$ ($-\delta<y_1<0$) and that the differential of $f_1-f_2$ is constant for $x_1<-\delta$ ($y_1>\delta$). As $\lambda\to 0$ we can (and will) make this alteration in the region $-\delta\lambda\le x_1<0$.

\subsubsection{Two Morse theory lemmas}
To continue the construction of the metric we will use two general lemmas about gradient flows which we present next. Let $M$ be an $n$-manifold and let $f\colon M\to \R$ be a smooth function on $M$. Let $g$ be a metric on $M$ and let $\gamma\colon[0,T]\to M$ be a non-constant solution to the gradient equation
$$
\dot\gamma(t)=\nabla f(\gamma(t)).
$$
Assume that there are coordinates $x=(x_1,\dots,x_n)$ around $\gamma(0)$ such that in these coordinates
$$
f(x)=c_1+\mu x_1,\quad\text{ and }g(x)=\sum_j dx_j^2,
$$
and coordinates $y=(y_1,\dots,y_n)$ around $\gamma(T)$ such that in these coordinates
$$
f(x)=c_2+\mu y_1,\quad\text{ and }g(y)=\sum_j dy_j^2.
$$

\begin{lma}\label{l:linear}
There exists a metric $h$ on $M$ which agrees with $g$ in some neighborhoods of $\gamma(0)$ and $\gamma(T)$ and in the complement of some small neighborhood of $\gamma([0,T])$, and with the following property. There are coordinates $(u,u')\in [0,T]\times D^{n-1}(\epsilon)$ such that
$$
f(u,u')=\mu u, \quad h(u,u')=du^2 + (du')^2,
$$
where $D^{n-1}(\epsilon)$ denotes the $(n-1)$-disk of radius $\epsilon$, and where these coordinates agree with the given coordinates at the ends. Moreover, the Morse flow of $f$ in the metric $h$ is conjugated via a $C^1$-small diffeomorphism supported near $\gamma([0,T])$ to the Morse flow of $g$.
\end{lma}

\begin{pf}
Let $I\times D^{n-1}$ be a flow tube around the gradient flow $\gamma$. Then $I\times\pa D^{n-1}$ consists of flow lines and $\nabla f$ is transverse to $\pa I\times D^{n-1}$. We will alter the metric in the interior of the flow tube and it is obvious that this cannot change the dynamics much.

We define coordinates as follows. Let $E_1,\dots,E_{n-1}$ be an orthonormal frame of normal vectors along $\gamma$, which agrees with standard frames in a neighborhood of the end points. Let $\beta$ be the re-parametrization of $\gamma$ which satisfies the equation
$$
\dot\beta=\lambda\frac{1}{g(\nabla f,\nabla f)}\nabla f.
$$

For each $u\in[c_1,\frac{1}{\lambda}(c_2-c_1)]$ let $N^u$ denote the level surface of $f$ which contains $\beta(u)$. Let $U_r$ be a neighborhood of $N^u$ of points of distance at most $r$ from $N^u$. Then there exists $\epsilon>0$ such that the Fermi coordinate map
$$
\Pi^u\colon U_\epsilon\to N^u,
$$
is well defined for each $u$. (The Fermi coordinate map projects along geodesics perpendicular to $N^u$.) Define the map
$$
\phi(u,u')=\Pi^u(\exp_{\beta(u)}(u_j E_j(u)).
$$
It is easy to see that $d\phi(u,0)=\id$. Thus $\phi$ gives local coordinates on some neighborhood of $\gamma$. Also in these coordinates we have
$$
f(u,u')=u.
$$
We let $h$ be a suitable scaling of the flat metric in these coordinates $(u,u')$ and interpolate to the given metric outside.
\end{pf}

The second lemma is similar, we have the same assumptions as above but change the local form of $f$ near the ends as follows. Assume that there are coordinates $x=(x_1,\dots,x_n)$ around $\gamma(0)$ such that in these coordinates
$$
f(x)=c_1+k_1 x_1-\alpha x_1^2,\quad\text{ and }g(x)=\sum_j dx_j^2,
$$
and coordinates $y=(y_1,\dots,y_n)$ around $\gamma(T)$ such that in these coordinates
$$
f(x)=c_2+ k_2 y_1 -\alpha y_1^2,\quad\text{ and }g(y)=\sum_j dy_j^2,
$$
where $\alpha$ satisfies
\begin{equation}\label{eqmatch}
c_2-c_1=\frac{k_2^2-k_1^2}{2\alpha}.
\end{equation}
\begin{lma}\label{l:quadratic}
There exists a metric $h$ on $M$ which agrees with $g$ in some neighborhoods of $\gamma(0)$ and $\gamma(T)$ and in the complement of some small neighborhood of $\gamma([0,T])$, and with the following property. There are coordinates $(u,u')\in [0,T]\times D^{n-1}(\epsilon)$ such that
$$
f(u_1,u)=k_1 u-\alpha u^2, \quad h(u,u')=du^2 + (du')^2,
$$
where these coordinates agree with the coordinates at the ends (after a suitable translation in the $u$-direction). Moreover the Morse flow of $f$ in the metric $h$ is conjugated via a $C^1$-small diffeomorphism supported near $\gamma([0,T])$ to the Morse flow of $g$.
\end{lma}

\begin{pf}
The proof is similar to the proof of Lemma \ref{l:linear}.
\end{pf}

\subsubsection{Connecting the pieces and edge points}
We use Lemma \ref{l:linear} to connect the flat metrics at plateau- and inflection points, and Lemma \ref{l:quadratic} to connect the flat metrics between inflection points. This defines the metric in some small neighborhood of each (part of) a gradient line which is part of some rigid flow tree. Our transversality assumptions imply that the flow manifolds of function differences along an edge in a rigid flow tree (which are not the difference of the flow line of the edge) are transverse to the gradient of the flow line. Therefore, by making the flow tube in which the metric is altered sufficiently small, we change all stable and unstable manifolds and all flow outs of intersections by a small diffeomorphism.

We will consider one further deformation of $\hat L_\lambda$ as follows. Consider a part of a flow line of $f_1-f_2$ between two inflection points or between an inflection- and a plateau point. We will make both of the differentials of $f_1$ and $f_2$ constant except for a part of the flow line which shrink with $\lambda$. More precisely fix a small $\epsilon>0$ and subdivide the part of the flow line under consideration in finitely many pieces of length $\epsilon$. Use the Hamiltonian $H(x,y)= \lambda(f(x)-df(0)x)$ to make $df_1$ covariantly constant in a neighborhood of length $(1-\lambda)\epsilon$. Note that under this deformation $df_2$ becomes constant or linear. As in \S\ref{ss:metmapip}, we see that the interpolation regions may shrink as $\lambda\epsilon$ as $\lambda\to 0$. Moreover, note that such deformations do not affect gradient differences and hence do not affect flow lines.

For simpler notation below we will use the collective name {\em edge points} to denote plateau points, inflection points, and the points where we interpolate between the pieces where $df_1$ is constant as above. Also, we keep the notation $\hat L_\lambda$ for the totally real submanifold with special properties near its rigid flow trees constructed above.

\begin{rmk}\label{r:fltrel}
Lemma \ref{l:ftrel} still holds for (the further deformed) $\hat L_\lambda$ and the metric. This follows from a straightforward extension of its proof to incorporate the new deformation regions as well.
\end{rmk}

For future convenience we note that in coordinates $x=(x_1,\dots,x_n)$, in which the metric has the form $g(x)=\sum dx_j^2$, the differentials of the local functions along a flow line in a region between two edge points has one of the following forms:
\begin{align}\notag
df_1 &= c_1\,dx_1,\\\label{e:norform1}
df_2 &=c_2\,dx_1,
\end{align}
where $c_1\ne c_2$ are real constants, or
\begin{align}\notag
df_1 &= c_1\,dx_1,\\\label{e:norform2}
df_2 &= (c_2 + k x_1)\,dx_1,
\end{align}
where $c_1\ne c_2 + kx_1$.
Near critical points of some function difference $f_1-f_2$ it has the form
\begin{align}\notag
df_1 &= \sum c_j\,dx_j,\\\label{e:norform3}
df_2 &=\sum_j(c_j+\sigma_j x_j)\,dx_j,\quad \sigma_j=\pm1.
\end{align}

\subsubsection{The $2$-dimensional case}\label{ss:metmap2D}
We describe the modifications of the above construction in the case when $\dim(L)=2$ and when there are swallow tail points. Consider the normal form of a swallow tail point in Remark \ref{r:2Dswt}. After a small perturbation we may assume that the $y_2$-coordinates of all sheets near $p$ are distinct. We first isotope $L$ so that the $y_1$-coordinates of all sheets near $p$ agree. (Such an isotopy, although it is not necessarily small can have support in an arbitrarily small neighborhood of $p$ and need not introduce any new Reeb chords. Moreover, if the support of the isotopy is sufficiently small, there are no flow trees of formal dimension $0$ which pass through the neighborhood of the swallow tail points. To see this we note that even though the isotopy is not small, as the region of change approaches $0$ the flow lines approach flow lines of $L$. But by
Proposition \ref{p:ttv} none of the finitely many rigid flow trees of $L$ passes through the swallow-tail points.)

In local coordinates the above means that around a swallow-tail point the totally real submanifold has the following appearance. The swallow tail sheet looks like
\begin{align*}
x_1(u,q)&=u\\
x_2(u,q)&=q(\frac13 q^2 +\frac12 u),\\
y_1(u,q)&=\alpha,\\
y_2(u,q)&=q+\beta,
\end{align*}
any other sheet looks like
\begin{align*}
x_1(u,v)&=u\\
x_2(u,v)&=v\\
y_1(u,v)&=\alpha,\\
y_2(u,v)&=\gamma,
\end{align*}
where all $\gamma\ne \beta$ and the constants $\gamma$ corresponding to different smooth sheets are distinct.

After normalizing the situation near swallow-tail points, the construction proceeds exactly as described above with the only difference that in this dimension rigid flow trees may intersect and self-intersect transversely. We treat such intersection and self-intersection points exactly as the $Y_0$-vertices were treated in \S\ref{ss:metmapY0}.

\subsection{The almost complex structure}\label{s:acs}
Let $(M,g)$ be a Riemannian $n$-manifold. The metric $g$ determines the Levi-Civita connection $D$ on $TM$. In local coordinates $x=(q_1,\dots,q_n)$ on $M$ with $\pa_j=\frac{\pa}{\pa q_j}$,
$$
D(v_i\pa_i)=\Bigl(\frac{\pa v_k}{\pa q_j} + \Gamma_{ji}^k\, v_i\Bigr)
\,\pa_k\otimes dq_j,
$$
where $\Gamma_{ij}^k$ are the Christofel symbols of $D$ and where we use the summation convention. Moreover, $D$ induces a connection $\nabla$ on $T^\ast M$ by
$$
[\nabla_X\alpha](Y)=X(\alpha(Y))-\alpha(D_X Y),
$$
where $\alpha$ is a $1$-form, $X$ is a tangent vector, and $Y$ is vector fields on $M$. Then $\nabla$ is a metric connection with respect to the metric naturally induced on $T^\ast M$ (if $g=g_{i\!j}$ then the metric on $T^\ast M$ is $g^{i\!j}$, the inverse of $g_{i\!j}$). In local coordinates we have
$$
\nabla(a_i dq_i)=\Bigl(\frac{\pa a_k}{\pa q_j} -
\Gamma_{jk}^ia_i\Bigr)dq_k\otimes d q_j.
$$

The connection $\nabla$ gives a direct sum decomposition $T(T^\ast M)=V\oplus H$, where the vertical bundle $V$ equals the kernel of $d\pi\colon T(T^\ast M)\to TM$, the differential of the projection $\pi\colon T^\ast M\to M$. The fiber $H_\alpha$ of the horizontal bundle $H$ at $\alpha\in T^\ast M$ is defined as the velocity vectors of covariantly constant lifts of curves through $\pi(\alpha)$ with initial value $\alpha$. Thus, in local coordinates $(q,p)$ on $T^\ast M$, where $q=(q_1,\dots,q_n)$ are coordinates on $M$, where $p=(p_1,\dots,p_n)$, and where $(q,p)$ corresponds to the cotangent vector $p_jdq_j$ at the point $q$, if  $\pa_j=\frac{\pa}{\pa q_j}$ and $\pa_{j^\ast}=\frac{\pa}{\pa p_j}$ we have
\begin{align*}
V_{(q,p)}&=\Spa\left\{\pa_{j^\ast}\right\}_{j=1}^n,\\
H_{(q,p)}&=\Spa\left\{\pa_j + p_r\Gamma^r_{js}\pa_{s^\ast}\right\}_{j=1}^n.
\end{align*}
The natural almost complex structure $J$ on $T^\ast M$ is required to satisfy $J(V)=H$ and defined as follows on vertical vectors $\xi\in V_\alpha$. Translate $\xi$ to the origin in $T^\ast_{\pi(\alpha)} M$. Identify this translate with a tangent vector to $M$ and let $J(\xi)$ be the negative of its horizontal lift. In local coordinates this gives
\begin{align*}
J(\pa_{j^\ast})&=g^{kj}(-\pa_k+\Gamma^r_{ks}\pa_{s^\ast}),\\
J(\pa_j)&=g_{jk}\pa_{j^\ast}-g^{ts}p_r\Gamma^r_{ks}\pa_t+
p_rp_m\Gamma^r_{ks}\Gamma^m_{tl}\pa_{l^\ast},
\end{align*}
where $q=\pi(\alpha)$ and $\alpha=p_jdq_j$. Note that $\Pi_H\circ J=J\circ\Pi_V$ and $\Pi_V\circ  J=J\circ\Pi_H$, where $\Pi_H$ and $\Pi_V$ denotes the projections to $H$ and $V$, respectively.

\begin{rmk}\label{r:notation}
In Sections \ref{S:disktotree} and \ref{S:treetodisk} we will use the following notation. The immersed totally real submanifold constructed in Subsections \ref{s:Rchandce} and \ref{s:metmapfltr} will be denoted $\hat L_\lambda$. We will consider $J$-holomorphic disks $u\colon \Delta_m\to T^\ast M$ with boundary on $\hat L_\lambda$, where $J$ is the almost complex structure associated to the metric constructed. Let $\Phi_\lambda\colon T^\ast M\to T^\ast M$ denote the diffeomorphism $C^1$-close to identity which maps $\Pi_\C(\tilde L_\lambda)$ to $\hat L_\lambda$. As in Definition \ref{d:phol}, we impose the condition that if $u\colon D_m\to T^\ast M$ is a $J$-holomorphic disk with boundary on $\hat L_\lambda$ then the restriction of $\Phi_\lambda^{-1}\circ u$ to the boundary has a continuous lift to $\tilde L_\lambda\subset J^1(M)$. As mentioned in Section \ref{S:intr}, the complex structure $J_\lambda=d\Phi_\lambda^{-1}\circ J\circ\Phi_\lambda$ is the one used in Theorem \ref{t:main}. Thus, composition with $\Phi_\lambda^{-1}$ gives a 1-1 correspondence between $J$-holomorphic disks with boundary on $\hat L_\lambda$ and $J_\lambda$-holomorphic disks with boundary on $\tilde L_\lambda$.
\end{rmk}

%% file: Sec/5disktotree.tex
\section{From disks to trees}\label{S:disktotree}
In Subsection \ref{s:areamonbst} we present standard analytical estimates for holomorphic disks with boundary on $\hat L_\lambda$. In Subsection \ref{s:subdiv}, the domains in a sequence of holomorphic disks with boundary on $\hat L_\lambda$ are subdivided into pieces which map close to the singularity strata $\Sigma_1\supset\dots\supset\Sigma_k$ and pieces which map far from them. We construct the subdivision by adding a uniformly finite (as $\lambda\to 0$) number of punctures to the domains. In Subsection \ref{s:projftconv}, we show that pieces of the disks which map far from $\Sigma_1$ converge to flow trees and, moreover, that pieces which map near $\Sigma_j$, $j\ge 1$, when composed with a local projection, converge to projected flow trees. In Subsection \ref{s:totavlkbu} we study the pieces of the disk that map close to $\Sigma_j$ in more detail. Using the result about convergence to projected flow trees, we associate an average linking number to a sequence of holomorphic disks and prove that it is finite. This finiteness in combination with  (re-scaled) blow-up arguments allow us to prove Theorem \ref{t:disktotree}.

\begin{rmk}
Throughout the rest of the paper we will make use of the following notation. If $f(x)$ and $g(x)$ are functions then we write $f(x)=\Ordo(g(x))$ if there exists a constant $K$ such that $|f(x)|\le K|g(x)|$ as $x\to 0$, or sometimes, as $x\to\infty$. Also we write $f(x)=\ordo(g(x))$ if $\lim \frac{f(x)}{|g(x)|}=0$, where the limit is taken either as $x\to 0$ or as $x\to\infty$.
\end{rmk}

\subsection{Area, monotonicity, and bootstrapping}\label{s:areamonbst}
It is a consequence of Stokes theorem that if $u$ is a holomorphic disk with boundary on the Legendrian submanifold $\tilde L_\lambda$ which has at most $P$ positive punctures then
$$
\area(u)\le C\lambda,
$$
where the constant $C>0$ depends on the length of the Reeb chords of $L$ and on the constant in the taming condition on $\omega$. We show that a similar estimate holds for holomorphic disks with boundary on $\hat L_\lambda$, provided $\hat L_\lambda$ is obtained from $\tilde L_\lambda$ by a  sufficiently small perturbation.
\begin{lma}\label{l:Stoke}
If the deformation parameters of $\hat L_\lambda$ are sufficiently small then there exists $K>0$ such that for any $J$-holomorphic disk $u\colon \Delta_m\to T^\ast M$ with boundary on $\hat L_\lambda$ and less than $P$ positive punctures
$$
\area(u)\le K\lambda.
$$
\end{lma}

\begin{pf}
Since the standard complex structure $J$ on $T^\ast M$ is tamed by the standard symplectic form $\omega$ there exists $\epsilon >0$ such that any complex structure $\tilde J$ of distance less than $\epsilon$ from $J$ is tamed by $\omega$. Recall that there exists a fiber preserving diffeomorphism $\Phi_\lambda\colon T^\ast M\to T^\ast M$ such that $\Phi_\lambda(\Pi_\C(\tilde L_\lambda))= \hat L_{\lambda}$ and such that $d_{C^1}(\Phi_\lambda,\id)<\epsilon$, see Remark \ref{r:notation}. In particular, the complex structures $J_\lambda=d\Phi_\lambda^{-1}\circ J\circ d\Phi_\lambda$ are uniformly tamed by $\omega$. If $u$ is $J$-holomorphic then $\Phi_\lambda^{-1} \circ u$ is a $J_\lambda$-holomorphic map with boundary on $L_\lambda$. It follows that
$$
\area(\Phi_\lambda^{-1}\circ u)\le
C'\int_{D_m}(\Phi^{-1}\circ u)^\ast\omega\le C'K'\lambda,
$$
for some constants $C'>0$ and $K'>0$. But since $d_{C^1}(\Phi_\lambda,\id)\le \epsilon$ we find that $\area(u)\le K\lambda$.
\end{pf}

\subsubsection{Consequences of monotonicity}
Let $u\colon S\to T^\ast M$ be a $J$-holomorphic map, where $S$ is some Riemann surface. Let $p\in T^\ast M$ and let $B(p,r)$ be the ball of radius $r$ around $p$. Assume that $u(\pa S)\subset \pa B(p,r)$ and that $u(0)=p$ then the monotonicity property for holomorphic curves, see e.g. Proposition 4.3.1 \cite{Si}, asserts that there exists a constant $C>0$ such that
\begin{equation}\label{e:monint}
\area(u)\ge C r^2.
\end{equation}

The boundary version of ordinary monotonicity is as follows. Let $v\colon S\to T^\ast M$ be a holomorphic map, where $S$ is some Riemann surface. Let $p\in T^\ast M$ and let $B(p,r)$ be the ball of radius $r$ around $p$. Assume that $v(\pa S)\subset \pa B(p,r)\cup L$ and that $v(0)=p$ where $L$ is some Lagrangian submanifold. Then the monotonicity property for $J$-holomorphic disks asserts that there exists a constant $C'>0$ such
\begin{equation}\label{e:monbdry}
\area(u)\ge C' r^2.
\end{equation}

The estimate \eqref{e:monbdry} depends on the comparison between some intrinsic metric on $L$ and the metric induced by the embedding (outside neighborhoods of the double points). We will apply boundary monotonicity to $\hat L_\lambda$, which converges to the $0$-section, and hence no such estimate can hold uniformly. However, if we restrict attention to balls of radius smaller than $\delta\lambda$ for some small fixed $\delta>0$, such that $\delta\lambda$ is small compared to the minimal radius of the half circles in $\hat L_\lambda$ over $\Sigma$ then the estimate holds. (These radii have the form $\alpha(q)\lambda\sqrt{\eta}$, see \S\ref{ss:LegisoCE}.) Also, for \eqref{e:monbdry} to hold in our case it is essential that our totally real submanifold $\hat L_\lambda$ is Lagrangian with respect to some symplectic form, we take for example $(\Phi_\lambda^{-1})^\ast\omega$.
The next result is one of the reasons for rounding all cusp edges.
\begin{lma}\label{l:finitelength}
There exists a constant $K>0$ such that for a $J$-holomorphic disk
$u\colon \Delta_m\to T^\ast M$ with boundary on $\hat L_\lambda$ and with less than $P$ positive punctures, the length of  $u_\lambda(\pa \Delta_m)$ is smaller than $K$.
\end{lma}
\begin{pf}
Assume that the length of $u(\pa\Delta_m)=L$. For sufficiently small $\delta>0$, we can find $\frac{L}{\lambda}$ disjoint balls of radius $\delta\lambda$ centered on points of the path. Inequality \eqref{e:monbdry} then implies that the area of the disk is at least
$$
C'\frac{L}{\lambda}\delta^2\lambda^2=C' L\delta^2\lambda.
$$
Together with Lemma \ref{l:Stoke} this implies the lemma.
\end{pf}

\subsubsection{A subharmonic function}
Consider the operator $\bar\pa_J u=du + J\circ du\circ i$ for functions $u\colon \Delta_m\to T^\ast M$. In coordinates $\tau+it$ in $\Delta_m$, $\bar\pa_J u$ is determined by $J$ and $\bar\pa_J(\pa_\tau)$, and $\bar\pa_Ju=0$ is equivalent to
$$
\frac{\pa u}{\pa \tau}+ J\frac{\pa u}{\pa t}=0.
$$
Writing $u(z)=(q(z),p(z))$ where $q(z)\in M$ and $p(z)\in T_{q(z)}^\ast M$ and using the projections $\Pi_V$ and $\Pi_H$, see \S\ref{s:acs}, we find that the above equation is equivalent to the system of equations
\begin{align*}
&\frac{\pa q}{\pa\tau}+(\nabla_t p)^\ast=0,\\
&\frac{\pa q}{\pa t}-(\nabla_\tau p)^\ast=0.
\end{align*}
Here we view $p(z)$ as a section of $T^\ast M$ over $q(z)$ and we think of the covectors $\nabla_t p$ and $\nabla_\tau p$ as vertical tangent vectors of $T^\ast M$ at $p(z)$, and $\alpha^\ast$ denotes the tangent vector which corresponds to the cotangent vector $\alpha$ via the given metric.

\begin{lma}\label{l:subh}
Let $u=(q,p)\colon \Delta_m\to T^\ast M$ be a holomorphic disk with boundary on $\hat L_\lambda$. Then the function $z\mapsto |p_\lambda(z)|^2$ is subharmonic. In particular it achieves its maximum on the boundary and thus there exists a constant $K>0$ such that
$$
|p_\lambda(z)|\le K\lambda.
$$
\end{lma}
\begin{pf}
Recall that the Levi-Civita connections on $TM$ and on $T^\ast M$ were denoted $D$ and $\nabla$, respectively. Since $u=(q,p)$ is $J$-holomorphic we have
$$
(\nabla_\tau\nabla_\tau p)^\ast=D_\tau(\nabla_\tau p)^\ast=D_\tau\frac{\pa q}{\pa t}=D_t\frac{\pa q}{\pa \tau},
$$
and analogously,
$$
(\nabla_t\nabla_t p)^\ast=-D_t\frac{\pa q}{\pa \tau}.
$$
Thus,
\begin{align*}
&\Delta\la p,p\ra=\left(\frac{\pa^2}{\pa\tau^2}+\frac{\pa^2}{\pa t^2}\right)\la p,p\ra\\
&=2\Bigl(\left\la \nabla_\tau p,\nabla_\tau p\ra + \la \nabla_t p,\nabla_t p\right\ra
+\left\la\bigl(\nabla_\tau\nabla_\tau p + \nabla_t\nabla_t p),p \right\ra\Bigr)\\
&=2\Bigl(\left\la \nabla_\tau p,\nabla_\tau p\ra + \la \nabla_t p,\nabla_t p\right\ra
+\left\la\bigl(\nabla_\tau\nabla_\tau p + \nabla_t\nabla_t p)^\ast,p^\ast \right\ra\Bigr)\\
&=2\Bigl(\left\la \nabla_\tau p,\nabla_\tau p\ra + \la \nabla_t p,\nabla_t p\right\ra\Bigl)\ge 0.
\end{align*}
\end{pf}

\subsubsection{Bootstrap estimates}
Consider the singularity set $\Sigma\subset\hat L_\lambda$ and the stratification
$$
\Pi(\Sigma)=\Sigma_1\supset\Sigma_2\supset\dots\supset\Sigma_k
$$
of its image. Note that there exists $d>0$ such that the distance between any two points $p,q\in\Sigma\cap \Pi^{-1}(\Sigma_2)$ with $p\ne q$ and with $\Pi(p)=\Pi(q)$ is larger than $d$. Let $S\subset\Sigma$ have diameter less than $d$.
For small $\delta>0$ (smaller than the size of the neighborhood of $\Sigma$ in which $\hat L_\lambda$ is compatible with the local product structure), let $U(\Pi(S),\delta)$ be the restriction to $S$ of a $\delta$-tubular neighborhood in the product neighborhood, see \S\ref{ss:TRisocomp}, of the branch of $\Pi(\Sigma)$ containing $S$. Let $V(S,\delta)\subset\hat L_\lambda$ be the connected component of $\Pi^{-1}(U(\Pi(S),\delta))$ containing $S$.

By definition of the product metric there exists a projection $\pi\colon T^\ast U(\Pi(S),\delta)\to T^\ast S$.
(In the the local coordinates $(q,p,s,\sigma)$ of \S\ref{ss:TRisocomp}, $\pi_j(q,p,s,\sigma)=(q,p)$.) If $u$ is a map with image in $T^\ast U(\Pi(S),\delta)$, then let $u^T=\pi\circ u$. Since the complex structure is split near $\Pi(S)$ and since $\hat L_\lambda$ is a product in a neighborhood of $S$ the following holds. If $u$ is $J$-holomorphic then $u^T$ is $J_T$-holomorphic, where $J_T$ is the almost complex structure corresponding to the induced metric on $S$.
Also, if $u$ maps the boundary to the branch of $\hat L_\lambda$ containing $S$, then $u^T$ maps the boundary to the submanifold of $T^\ast S$ the general form of which in local coordinates $(q,p)=(q_0,p_0,s,\sigma)\in\R^{2(n-j-1)}\times\R^{2j}$ on $T^\ast S$ is
$$
(q_0,p_0,s,\sigma)\mapsto \bigl(q_0,s,a(q_0),\lambda\hat h_1(q_0,s_1),\dots,\lambda\hat h_r(q_0,s_r),\lambda h_{j-r}(q_0),\dots,\lambda h_j(q_0)\bigr),
$$
see \S\ref{ss:TRisocomp}.
For $r>0$, define $D_r=\{\tau+it\colon \tau^2+t^2< r^2\}$ and $E_r=D_r\cap\{t\ge 0\}$.
\begin{lma}\label{l:bootstr}
The following interior estimate holds: there exists $C>0$ such that if $u\colon D_2\to T^\ast M$ is a $J$-holomorphic map then
$$
\sup_{D_1}|Du|\le C\|Du\|_{L^2,D_2}.
$$
(By the taming condition $\|Du\|_{L^2,D_2}^2\le C\area(u(D_2))$.)

The following boundary estimates hold: for all sufficiently large $K>0$ and for all sufficiently small $\delta>0$ there exists $C>0$ with the following properties.
\begin{itemize}
\item[{\rm (a)}]
If $u\colon (E_2,\pa E_2)\to (T^\ast M,\hat L_\lambda)$ is $J$-holomorphic and such that $u(\pa E_2)$ lies outside a $K\lambda$-neighborhood of $\Sigma\subset\hat L_\lambda$ then
$$
\sup_{E_1}|Du|\le C\|Du\|_{L^2,E_2},
$$
\item[{\rm (b)}]
If $u\colon (E_2,\pa E_2)\to (T^\ast M, V(S,\delta))$ is a $J$-holomorphic map such that $u(\pa E_2)\subset T^\ast(U(\Pi(S),\delta))$, then
$$
\sup_{E_1}|Du^T|\le C \|Du^T\|_{L^2,E_2}.
$$
\end{itemize}
\end{lma}

\begin{pf}
This follows from a standard bootstrap argument, see e.g. \cite{Fl:un-reg, Oh, EES1}. Note that the condition on the image in (b) guarantees that the
diffeomorphisms $(T^\ast M,\hat L_\lambda)\to(\C^{n},\R^{n})$, used to translate to the standard situation are uniformly bounded.
\end{pf}

Let $u\colon \Delta_m\to T^\ast M$ be a $J$-holomorphic disk with boundary on $\hat L_\lambda$. For a point $p\in\Delta_m$ which lies at distance at least $r>0$ from $\pa\Delta_m$, let $B_r(p)\subset\Delta_m$ be an $r$-ball around $p$. For $q\in\pa\Delta_m$ let $A_r(q)$ denote the subset of points in $\Delta_m$ which can be connected to $q$ by a path in $\Delta_m$ of length at most $r$. Note that if $r>0$ is small enough then $A_r(q)$ is a topological half disk.

\begin{lma}\label{l:estgood}
If the speed parameter $\epsilon_3>0$, see \S\ref{ss:TRisocomp}, $\delta>0$, and $\epsilon>0$ are sufficiently small, and if $K>0$ is sufficiently large. Then there exists a constant $C=C(\epsilon)>0$ such that the following estimates hold  for all sufficiently small $\lambda>0$.
\begin{itemize}
\item[{\rm (a)}]
Let $p\in\Delta_m$ be any point of distance at least $4\epsilon$ from $\pa\Delta_m$. Then
$$
\sup_{B_{\epsilon}(p)}|Du|\le C\lambda.
$$
\item[{\rm (b)}]
Let $p\in\pa\Delta_m$. If $u(A_{4\epsilon}(p)\cap\pa\Delta_m)$ lies outside a $K\lambda$-neighborhood of $\Sigma\subset\hat L_\lambda$ then
$$
\sup_{A_{\epsilon}(p)}|Du|\le C\lambda.
$$
\item[{\rm (c)}]
If $u(A_{4\epsilon}\cap\pa\Delta_m)\subset V(S,\delta)$ and $u(A_{4\epsilon})\subset T^\ast U(\Pi(S),\delta)$ then
$$
\sup_{A_{\epsilon}(p)}|Du^T|\le C\lambda.
$$
\end{itemize}
\end{lma}

\begin{pf}
Case (a) follows from a simpler version of the argument we give for (b) and case (c) is proved in exactly the same way as case (b) just suppressing one of the coordinates. We therefore consider case (b).

Note that any half-disk $A_{4\epsilon}(p)\subset \Delta_m$ can be mapped to $E_\epsilon$ by a biholomorphic map with uniformly bounded derivatives. (In fact we can take the inclusion map except near boundary minima.) It is then a consequence of Lemma \ref{l:bootstr} that there exists a constant $K'$ such that the weaker estimate $|Du(z)|\le K'\lambda^{\frac12}$, $z\in A_{2\epsilon}$ holds. For simpler notation, we think of $u$ as defined on the half-disk $E$ of radius $1$ corresponding to $B_{2\epsilon}$ where the coordinates are chosen so that
\begin{equation}\label{e:sqrtest}
|Du(z)|\le K\lambda^{\frac12}.
\end{equation}

Equation \eqref{e:sqrtest} implies that $u(z)\subset B(\pi(u(0)),K_0\lambda^{\frac12})\subset T^\ast M$ for some $K_0>0$, where $B(x,r)$ denotes the $r$-ball around $x$. Let $(\C^n,\R^n)$ be a standard coordinate system in $T^\ast M$ around $\pi(u(0))$ such that the standard complex structure $J_0$ on $\C^n$ corresponds to the complex structure $J(\pi(u(0)))$. Then for any $q\in B(0,K'\lambda^{\frac12})$ we have $|J(q)-J_0|=\Ordo(\lambda^{\frac12})$.

Let the sheet of $\hat L_\lambda$ in which $u(0)$ lies be the graph $\Gamma_{\lambda\alpha}$ of the $1$-form $\lambda\alpha$. We claim that for all sufficiently small $\lambda>0$, there exist a diffeomorphism $\Theta$ of $\C^n$ which maps $\Gamma_{\lambda\alpha}$ to the $0$-section, such that $d\Theta+J_0\circ d\Theta\circ J=0$ along $\Gamma_{\lambda\alpha}$, and such that $d_{C^1}(\Theta, \id)\le\eta$, for an arbitrarily small $\eta>0$. To see this we construct the inverse $\Psi$ of $\Theta$. Let $z=x+iy$ be standard coordinates on $\C^n$ where the $0$-section $\R^n$ corresponds to $y=0$. Define $\psi(x)=(x,\lambda \alpha(x))$. Then let $\Psi(x,y)=\psi(x)+y\cdot J(\psi(x)) \frac{\pa\psi}{\pa x}$ for very small $y$ and extend the map to the rest of the chart. It is then clear that $\Psi$ is $(J,J_0)$-holomorphic along the $0$-section and that it can be chosen to lie no further than $\eta$ from the identity in $C^1$-distance. (In fact the $C^1$-distance is $\Ordo(\lambda)$ except for sheets near the tangency locus in $\Pi(\Sigma)$, see \S\ref{ss:TRisocomp}, where $\eta\to 0$ with the $C^1$-distance between $\Phi_\lambda$ and $\id$ which in turn goes to $0$ with the speed parameter $\epsilon_3$.) Note that $\Theta\circ u$ is $J_\Theta$-holomorphic where $J_\Theta=d\Theta\circ J\circ d\Theta^{-1}$ and that $|J_\Theta-J_0|=\Ordo(\eta)$.

Consider the scaling by $\lambda^{-1}$, $z\mapsto\lambda^{-1} z$, of the chart $\C^n$. Let $\hat u=\lambda^{-1}\Theta_\lambda\circ u$ and note that $\hat u$ is $\hat J$-holomorphic, where $\hat J(x,y)=J_\Theta(\lambda x,\lambda y)$ so that still $|\hat J-J_0|=\Ordo(\eta)$. Moreover, by Lemma \ref{l:subh} the imaginary part of $\hat u$ is bounded. Since $\Theta$ is $(J,J_0)$-holomorphic on $\Gamma_{\lambda\alpha}$ it follows that $\hat u$ is $J_0$-holomorphic on $\pa E$. Thus, doubling $\hat u$ over $\R^n$ we obtain a map $\tilde u\colon D\to \C^n$, where $D$ is the unit disk, which is $C^1$ and which satisfies the differential equation
$$
d\tilde u + \tilde J(\zeta)\circ d\tilde u\circ i =0,
$$
where $\tilde J(\zeta)=\hat J(\hat u(\zeta))$ for $\zeta\in E$ and $\tilde J(\zeta)=\hat J(\hat u(\zeta^\ast)^\ast)$ with $a^\ast$ denoting the complex conjugate of $a$.

Let $F\colon\C^n\to\C^n$ be the map
$$
F(z_1,\dots,z_n)=(e^{iz_1},\dots,e^{iz_n})
$$
and let $f=F\circ\tilde u$. Then $|f|$ is uniformly bounded and the derivatives of $F$ are uniformly bounded in a neighborhood of the image of $\tilde u$ (since the imaginary part of $\hat u$ and therefore of $\tilde u$ is bounded). Moreover, $f$ satisfies the equation
$$
df + J_F(z)\circ df\circ i=0,
$$
where $J_F(z)=dF\circ \tilde J(z)\circ dF^{-1}$. Since $F$ is $J_0$-holomorphic we have $|J_F(z)-J_0|\le |dF||\tilde J(z)-J_0||dF^{-1}|=\Ordo(\eta)$. We conclude (writing $\bar\pa_{J_0}=\bar\pa$ and $\pa_{J_0}=\pa$) that
\begin{equation}\label{e:skew}
\bar\pa f + q(z)\pa f=0,
\end{equation}
where $q(z)=(J_0+ J_F(z))^{-1}(J_0-\hat J_F(z))=\Ordo(\eta)$.

We next use an argument similar to the proof of Proposition 2.3.6 in \cite{Si}: let $\beta\colon D\to\R$ be a cut-off function which equals $1$ on the disk of radius $\frac12$ and equals $0$ outside the disk of radius $\frac34$, and let $f_1=\beta f$. Then
$$
\bar\pa f_1 + q\pa f_1 = g_1,
$$
where $g_1=(\pa\pa\beta +\pa\beta q) f$. In particular, $|g_1|\le K|f|$. We conclude that
$$
|\bar\pa f_1|\le K_0\eta|\pa f_1|+ K_1|f_1|
$$
which, provided $\eta$ is sufficiently small, in turn implies that
$$
|\pa f_1|^2+|\bar\pa f_1|^2\le K'_0(|\pa f_1|^2-|\bar\pa f_1|^2) + K_1'|f_1|^2.
$$
Now $\int_D(|\pa f_1|^2-|\bar\pa f_1|^2)dA=\int_Df_1^\ast\omega=0$ by exactness of $\omega$. Therefore
$$
\|Df\|_{L^2}^2\le K\sup_D|f|^2\le K'
$$
since $f$ is bounded. Note that
$$
\|D\hat u\|_{L^2}\le K\|Df\|_{L^2},
$$
since the derivatives of $F$ are uniformly bounded on a neighborhood of the image of $\tilde u$. It follows that $\|D\hat u\|_{L^2}=\Ordo(1)$ and therefore by scaling that $\|Du\|_{L^2}^2=\Ordo(\lambda^2)$. The lemma now follows from Lemma \ref{l:bootstr}. As mentioned, case (a) is is simpler: no doubling and no diffeomorphism $\Theta$ is needed.
\end{pf}

\begin{rmk}\label{r:highstrap}
The bootstrapping argument in Lemma \ref{l:bootstr} can also be used to bound higher derivatives of $u$. Consequently, the estimates in Lemma \ref{l:estgood} holds for higher derivatives of $u$ as well.
\end{rmk}

\subsection{Domain subdivision}\label{s:subdiv}
\begin{figure}[htbp]
\begin{center}
\includegraphics[angle=0, width=8cm]{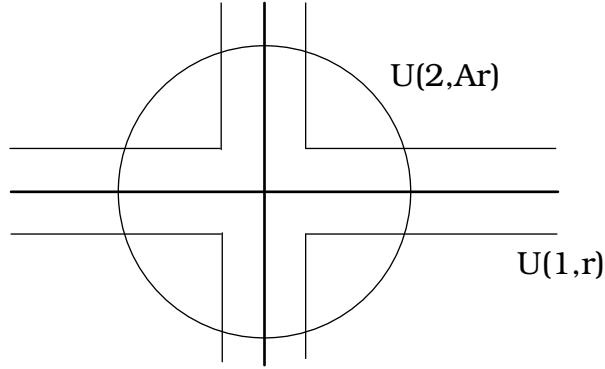}
\end{center}
\caption{Neighborhoods of strata.}
\label{f:projwd}
\end{figure}

We describe a procedure for adding punctures in the domain of a $J$-holomorphic disk with boundary on $\hat L_\lambda$ in such a way that the new domain can be subdivided into parts which are mapped into neighborhoods of the strata of the image of the singularity set $\Sigma$ of $\hat L_\lambda$ provided $\lambda>0$ is sufficiently small. As usual, $\Delta_m$ will denote a standard domain with $m$ punctures.
Consider the stratification
$$
\Pi(\Sigma)=\Sigma_1\supset\dots\supset\Sigma_k.
$$
For $r>0$ and $1\le s\le k$, let $U(s,r)$ denote an $r$-neighborhood of $\Sigma_s$ in $M$. If $A$ is sufficiently large and $r$ sufficiently small then the projection $\pi_j\colon T^\ast (U(j,r)-U(j+1,Ar))\to T^\ast\Sigma_j^\circ$ which in local coordinates, see \S\ref{ss:TRisoCE} is $\pi_j(q,p,s,\sigma)=(q,p)$ is well defined, see Figure \ref{f:projwd}. (Here $A$ depends on the dimension and $j$, we will not need the detailed information.) If $u$ maps into $T^\ast (U(j,r)-U(j+1,Ar))$ we will write $\pi_j\circ u=u^{T_j}$.

Let $u\colon\Delta_m\to T^\ast M$ be a $J$-holomorphic disk with boundary on $\hat L_\lambda$ and less than $P>0$ positive punctures.  Fix a small $\delta>0$ such that $u|\pa\Delta_m$ is transverse to $\pa(T^\ast U(1,c\delta))$ for $c\in\{1,2,3\}$. Let $I\approx\R$ be a boundary component of $\Delta_m$. Define
$$
b_1^c< b_2^c<\dots<b_{n(c)}^c,\quad c=1,2,3,4
$$
to be the points in $I$ such that $u(b^c_j)\in \pa(T^\ast U(1,c\delta))$. For $2\le c\le 4$, puncture the disk $\Delta_m$ at every $b_j^c$ such that there exists $b_k^{c-1}$ with $b_j^c<b_k^{c-1}<b_{j+1}^c$. We obtain in this way a new domain $\Delta_r=\Delta_{m+m_1}$ with a holomorphic map (which we continue to denote $u$) $u\colon \Delta_r\to T^\ast M$.
\begin{lma}\label{l:finpunct1}
There exists a constant $C=C(\delta)>0$ such that the number $m_1$ of added punctures satisfies $m_1\le C$.
\end{lma}
\begin{pf}
Note that to each puncture corresponds a segment in $u(\pa \Delta_m)$ connecting $\pa(T^\ast U(1,c\delta))$ to $\pa(T^\ast U(1,(c-1)\delta))$, $c=2,3,4$. Since the lengths of such segments are bounded from below, the lemma follows from Lemma \ref{l:finitelength}.
\end{pf}
The boundary components $I$ of $\Delta_{r}$ are of three different types.
\begin{itemize}
\item[$\out$:]$u(I)\subset T^\ast (M-U(1,3\delta))$,
\item[$\0$:]$u(I)\subset T^\ast (U(1,4\delta)-U(1,\delta))$, and
\item[$\inn$:]$u(I)\subset T^\ast U(1,2\delta)$.
\end{itemize}

For $\frac14>d>0$ and $I$ a boundary component of $\Delta_{r}$, let $N_d(I)$ be a $d$-neighborhood of $I$ in $\Delta_r$ and let $\Omega_d=\Delta_{r}-\bigcup_{I\subset \pa \Delta_{r}} N_d(I)$. Fix a small $\epsilon>0$.
\begin{lma}\label{l:estgood1}
There exists a constant $C$ such that if $\lambda>0$ is sufficiently small and if $\Theta_\epsilon=\Omega_\epsilon\cup \bigcup_{I\in \out\cup\0}N_\epsilon(I)$, then
$$
\sup_{z\in\Theta_\epsilon}|Du(z)|\le C\lambda.
$$
\end{lma}

\begin{pf}
This is immediate from Lemma \ref{l:estgood}.
\end{pf}

Let $D'\subset \Delta_r$ be the subset containing all vertical line segments in $\Delta_r$ which connects a point in a boundary component of type $\inn$ to some other boundary point of $\Delta_r$. Note that $\pa D'-\pa\Delta_r$ is a collection of vertical line segments.
\begin{lma}\label{l:boundary1}
For any $0<a<1$ and for sufficiently small $\lambda>0$ the distance from any point $p\in I$, where $I$ is a boundary segment of type $\out$, to $D'$ is larger than $\lambda^{-a}$. In particular if
$l$ is a vertical line segment in $\pa D'-\pa\Delta_r$ and if $q\in\pa l$ then $q$ is either a boundary minimum on a boundary segment of type $\inn$ or it lies in a boundary segment of type $\0$.
\end{lma}

\begin{pf}
We argue by contradiction. Assume that there exists $p\in I$ where $I$ is a boundary component of type $\out$ with properties as stated. Then there exists a path in $\Omega_\epsilon\cup N_\epsilon(I)\subset \Delta_r$ of length smaller than $\lambda^{-a}+5r$ ($r=m+m_1\le m+C$ see Lemma \ref{l:finpunct1}) connecting $p$ to a point midway between two horizontal boundary segments of length $1$ at least one of which, $E$ say, is of type $\inn$. Since $|Du|=\Ordo(\lambda)$ along the path and since $u(p)\in T^\ast(M-U(1,3\delta))$ we find that for $\lambda>0$ small enough and for a horizontal segment $A$ of length $1$ parallel to $E$, $u(A)\subset T^\ast(M-U(1,\frac52\delta))$. By definition $u(E)\subset T^\ast U(1,2\delta)$. The region $R$ between $E$ and $A$ is of the form $[0,1]\times [0,b]$, $\frac12\le b\le \frac{r}{2}$ and every path $t\mapsto u(\tau+it)$ has length at least $\frac12\delta$. Therefore
$$
\int_R \left|\frac{\pa u}{\pa t}\right|^2\,dA\ge \frac{1}{4r\delta}.
$$
Thus, the $L^2$-norm of $|Du|$ is bounded from below by a constant which contradicts Lemma \ref{l:Stoke} and the taming condition for $\lambda$ sufficiently small.
\end{pf}

Let $F_d$ be the region consisting of all points of distance $\le d$ from $D'$ and where we chose $\log(\lambda^{-1})\le d\le 2\log(\lambda^{-1})$ in such a way that $\pa F_d-\pa \Delta_r$ and $\pa F_{\frac12 d}-\pa \Delta_r$ consists of vertical line segments disjoint from all boundary minima. Define $D_1'=F_d$ and $D_0=\Delta_r-(F_{\frac12 d})$. Note that if $p\in\pa D_0\cap \pa\Delta_r$ then $p$ lies in a boundary component of type $\0$ or $\out$, that if $q\in\pa D'_1\cap\pa \Delta_r$ then $q$ lies in a boundary segment of type $\0$ or $\inn$ and that Lemma \ref{l:estgood1} implies  $\sup_{z\in D_0}|Du(z)|\le K\lambda$.

\begin{lma}\label{l:localize1}
For sufficiently small $\lambda$, $u(D_1')\subset T^\ast U(1,\frac92\delta)$ and $u(D_0)\subset T^\ast(M-U(1,\frac12\delta))$.
\end{lma}
\begin{pf}
Any point $q\in D_0$ can be joined by a path $\gamma\subset D_0$ of length at most $5r$ to a point $p\in \pa D_0\cap \pa\Delta_r$. Since $p\in T^\ast(M-U(1,\delta))$ and since $|Du|\le K\lambda$ along $\gamma$ the first statement follows. The boundary of $D_1'$ consist of (parts of) boundary segments of type $\0$ and $\inn$ and vertical lines $l\subset D_0$ which ends at boundary components of type $\0$. The estimate in $D_0$ and the definition of the subdivision of the boundary then imply that $u(\pa D_1)\subset T^\ast U(1,\frac{17}{4}\delta)$. Thus, the second statement follows from a monotonicity argument: if $u_\lambda(D_1')$ does not stay inside $T^\ast U(1,\frac92\delta)$ its area contribution is bounded from below, see \eqref{e:monint}.
\end{pf}

We continue this construction in $k-1$ steps, where $\Sigma_k$ is the deepest stratum of the stratification of $\Sigma_1$. The next step is the following. Consider $D_1'$ above and consider neighborhoods $U(2,c\delta_2)$ such that $u|\pa D_1'\cap \pa\Delta_r$ is transverse to $\pa(T^\ast U(2,c\delta_2))$, $c=1,2,3,4$, where $\delta_2\ge 16\delta=16\delta_1$ is small but sufficiently large in order for the projection $\pi_1$ mapping into $T^\ast\Sigma_1^\circ$ to be well defined on $T^\ast (U(1,5\delta_1)-U(2,\frac{3}{4}\delta_2))$, see Figure \ref{f:projwd}. Let $I$ be a component of $\pa D_1'\cap\pa\Delta_r$ and let
$$
b_1^c< b_2^c<\dots<b_{m(c)}^c,\quad c=1,2,3,4
$$
be the points such that $u(b_j^c)\in \pa(T^\ast U(2,c\delta_1))$. We puncture $I$ exactly as described above and the argument in Lemma \ref{l:finpunct1} again gives a uniform bound on the number of new punctures. As above, the punctures give rise to a new domain which we continue to call $\Delta_r$ with boundary components corresponding to the previous $\pa D_1'$ of three types $\inn$, $\out$, and $\0$. Consider a $\frac12$-neighborhood of a boundary segment of type $\out$ or $\0$. Note that $|Du|=\Ordo(\lambda)$ along the inner boundary component of this region (the one which does not intersect the boundary of $\Delta_r$) and that the outer boundary of the region stays inside the region $T^\ast (U(1,\frac92\delta_1)-U(2,\delta_2))$. Hence the energy argument in the proof of Lemma \ref{l:boundary1} implies that the inner boundary component stays inside $T^\ast(U(1,\frac{19}{4}\delta_1)-U(2,\frac78\delta_2))$ and a monotonicity argument then shows that the entire region maps into $T^\ast(U(1,5\delta_1)-U(2,\frac34\delta_2))$ where $\pi_1$ is defined.  Lemma \ref{l:estgood} then shows that
\begin{equation}\label{e:estgood2}
|Du^{T_1}(z)|=\Ordo(\lambda),
\end{equation}
for $z$ in a neighborhood of the boundary components of types $\out$ and $\0$. We then take $F_d$ with $d$ as before to be a neighborhood of the region which is the union of all vertical line segments emanating from $\inn$-segments. We take $D_2'=F_d$ and $D_1=D_1'-F_{\frac12d}$. The arguments of the proofs of Lemmas \ref{l:boundary1} and \ref{l:localize1} carries over to show that the new subdivision has the properties of the subdivision in the previous step.

To finish the inductive construction, we assume that we have constructed a subdivision of the domain $\Delta_r=D_0\cup D_1\cup\dots\cup D_s'$ where $r\le m+Ks$ for some integer $K$ depending only on the area of the disk and the finite neighborhood sizes $\delta_1,\delta_2,\dots,\delta_s$, with the following properties.
\begin{itemize}
\item $\pa D_j-\pa\Delta_r$, $1\le j<s$ and $D_s'-\pa\Delta_m$ consists of vertical line segments disjoint from boundary minima.
\item $u_\lambda(D_j)\subset T^\ast\bigl(U(j,\frac92\delta_j)-U(j+1,\frac34\delta_{j+1})\bigr)$,
\item $|Du^{T_j}|=\Ordo(\lambda)$ on $D_j$, and
\item $u(D_s')\subset T^\ast U(s,\frac52\delta_{s})$.
\end{itemize}
We continue by constructing $D_s$ and $D'_{s+1}$. The construction is completely parallel to the construction of $D_2'$ above so we only sketch the argument. First, consider intersections with $\pa (T^\ast U(s+1,c\delta_{s+1}))$, $c=1,2,3,4$ and add punctures. The number of punctures is bounded in exactly the same way as in Lemma \ref{l:finpunct1}. The new boundary segments come in three types $\inn$, $\out$, and $\0$. Again it follows from Lemma \ref{l:estgood} that
\begin{equation}\label{e:estgoods}
|Du^{T_s}(z)|\le K\lambda,
\end{equation}
for $z$ in a neighborhood of the new born boundary segments of type ${\mathbf 0}$. We then take $D_{s+1}'$ as a neighborhood of the region which is the union of all vertical line segments emanating from $\inn$-segments. The arguments of the proofs of Lemmas \ref{l:boundary1} and \ref{l:localize1} carries over to show that the new subdivision has the properties of the subdivision in the previous step. After the $(k-1)^{\rm th}$ step, where $\Sigma_k$ is the deepest stratum of $\Pi(\Sigma)$, we take $D_k=D'_k$ and we have obtained a domain $\Delta_r$ with finitely many punctures and a subdivision of this domain with properties as listed above.

\subsubsection{Domain subdivision for sequences of disks}\label{ss:subdivseq}
Consider an infinite sequence of $J$-holomorphic disks $u_\lambda\colon \Delta_m\to T^\ast M$, $\lambda\to 0$, such that $u_\lambda$ has boundary on $\hat L_\lambda$ and such that the Reeb chords at positive and negative punctures of the disks in the sequence remain constant.
We apply the subdivision procedure described above to the maps in the sequence. (In general the transversality conditions used in the subdivision are not met for all $u_\lambda$ if $\delta_1,\dots,\delta_k$ are fixed. We deal with this by letting $\delta_j$ be functions of $\lambda$ which varies in an arbitrarily small neighborhood of a fixed $\delta_j$. This way the transversality conditions can be assumed to hold for all $u_\lambda$ in the sequence.)

We obtain a sequence of domains $\Delta_{r(\lambda)}$ with subdivisions
$\Delta_{r(\lambda)}=D_0(\lambda)\cup\dots\cup D_k(\lambda)$. Here $u_\lambda$ takes each $D_j(\lambda)$ into $T^\ast(U(j,a)-U(j+1,b))$ where the projection $\pi_j$ mapping to $T^\ast \Sigma_j^\circ$ is defined and $|Du^{T_j}_\lambda|=\Ordo(\lambda)$ on $D_j(\lambda)$. By construction, the number of punctures $r(\lambda)$ is uniformly bounded. Hence by passing to a subsequence we may assume that $r=r(\lambda)$ is constant. Furthermore, since $D_j(\lambda)$ is defined as a (neighborhood of) a collection of vertical line segments of $\Delta_r$ emanating from boundary components of $\Delta_r$ of a certain type we may assume by passing to a further subsequence that also the topology of each component of each $D_j(\lambda)$ remains constant. We will assume without further mentioning of passing to subsequences that the sequences considered below have these properties.

\subsection{Convergence to gradient equation and projected flow trees}\label{s:projftconv}
Consider a sequence $u_\lambda\colon \Delta_r^\lambda\to T^\ast M$ of $J$-holomorphic disks with boundary on $\hat L_\lambda$ and with associated subdivisions $\Delta_r^\lambda=D_0(\lambda)\cup D_1(\lambda)\cup\dots\cup D_k(\lambda)$ as described in \S\ref{ss:subdivseq}. For $1\le j\le k$, let $W_j(\lambda)$ be a neighborhood of the boundary minima of $D_j(\lambda)$ such that $\pa W_j(\lambda)$ is a union of arcs in $\pa D_j(\lambda)$ and of vertical line segments and such that each component of $W$ contains at least one boundary minimum. Then $D_j(\lambda)-W_j(\lambda)$ is a finite collection of strip regions. Let $l\approx\{0\}\times[0,1]$ be a vertical segment $l\subset D_j(\lambda)-W_j(\lambda)$ with $\pa l\subset\pa D_j(\lambda)$. Let $[-c,c]\times[0,1]\subset D_j(\lambda)$ be a neighborhood of $l\approx\{0\}\times[0,1]$. With $(\tau,t)\in[-c,c]\times[0,1]$, we write $u^{T_j}_\lambda(\tau,t)=(q_\lambda(\tau,t),p_\lambda(\tau,t))\in T^\ast \Sigma_j^\circ$ (using the convention that $\Sigma_0^\circ=M-\Sigma_1$). Let $\lambda b_\sigma\colon \Sigma_j^\circ\to T^\ast\Sigma_j^\circ$, $\sigma=0,1$ denote the section of the sheet which contain $u^{T_j}_\lambda(0,\sigma)$. For a connected component $W\subset W_j(\lambda)$, we define the {\em width} of $W$ as the maximum distance from a vertical line segment in the boundary of $W$ to a boundary minimum inside $W$.

\begin{lma}\label{l:gradconv}
For $0\le j\le k$ and all sufficiently small $\lambda>0$, there exist neighborhoods $W_j(\lambda)$ as described above with each connected component of width at most $2\log(\lambda^{-1})$ such that along any vertical line segment $l\subset D_j(\lambda)-W_j(\lambda)$
\begin{align*}
&\frac{1}{\lambda}\nabla_t p_\lambda(0,t)
-\bigl(b_1(q_\lambda(0,0))-b_0(q_\lambda(0,0))\bigr)=\Ordo(\lambda),\\
&\frac{1}{\lambda}\nabla_\tau p_\lambda(0,t)=\Ordo(\lambda).
\end{align*}
\end{lma}

\begin{pf}
We prove this for $D_0(\lambda)$. The proof for $D_j(\lambda)$, $j\ge 0$ is completely analogous, just replace $u_\lambda$ in the proof below by $u^{T_j}_\lambda$ and $\hat L_\lambda$ by its projection to $T^\ast \Sigma_j^\circ$.

Let $\Theta_c=[-c,c]\times[0,1]\subset D_j(\lambda)$ denote a neighborhood of $l$. Note that the condition on the width of the components of $W_j(\lambda)$ allows us to choose $c=\log(\lambda^{-1})$, so let $c=\log(\lambda^{-1})$. The derivative bound $|Du|=\Ordo(\lambda)$, see Lemma \ref{l:estgood1}, implies that $u_\lambda(\Theta_c)$ is contained in a ball of radius $\Ordo(\lambda\log(\lambda^{-1}))$

Consider the geodesic disk $D_r$ around $\pi(u_\lambda(0,0))$ with normal coordinates and radius $r$ such that $u_\lambda(\Theta_c)\subset  T^\ast D_{R\lambda}$. Let $g_0$ be the flat metric on $D_r$ and let $J_0$ be the corresponding (standard) complex structure and let $J$ be that induced from $T^\ast M$. We think of $T^\ast D_r\subset\C^n$ (in standard coordinates $q+ip$ on $\C^n$, $\{q+ip\colon |q|\le r\}$). Consider the scaling of $\C^n$ by $\lambda^{-1}$ note that the map $\hat u=\lambda^{-1}u\colon\Theta_c\to\C^n$ is $J_\lambda$ holomorphic where $J_\lambda(q,p)=J(\lambda q,\lambda p)$ and that
\begin{equation}\label{e:Jclose}
|J_\lambda-J_0|_{C^1}=\Ordo(\lambda).
\end{equation}
The two sheets of $\hat L_\lambda$ over $D_r$ are given by $q\mapsto (q,\lambda b_\sigma(q))$, where $b_\sigma\colon D_r\to T^\ast$ are sections $j=0,1$. After scaling these become $q\mapsto (q, b_\sigma(\lambda q))$. Let $L_\sigma$ be the Lagrangian subspaces $\{x+iy\colon y=b_\sigma(0)\}$ and note that there exists a function $f\colon\Theta_c\to\C^n$ such that $f(0,0)=0$
\begin{equation}\label{e:bdryclose}
\sup_{\Theta_c}|D^k f|=\Ordo(\lambda),\quad k=1,2,3
\end{equation}
and with the following properties
\begin{align}
&\hat u(\tau+0i) + f(\tau+0i)\in L_0,\\
&\hat u(\tau+i) + f(\tau+i)\in L_1,
\end{align}
and $\hat u+f$ is $J_0$-holomorphic to first order on the boundary.
Let $u_1=\hat u+f$ and consider the standard solution $u_0\colon\R\times[0,1]\to\C^n$ to the $\bar\pa_{J_0}$-equation which satisfies the boundary conditions of $u_1$: $u_0(\tau+it)=[b_1(0)-b_1(0)](\tau+it)$.

If $v=u_1-u_0\colon\Theta_c\to\C^n$ then $v(\pa\Theta_c)\subset\R^n$, $v$ is holomorphic to first order on the boundary, and $v(0,0)=0$. Moreover, \eqref{e:Jclose} and \eqref{e:bdryclose} imply that $\sup_{\Theta_c}|D^k\bar\pa_{J_0} v|=\Ordo(\lambda)$, $k=0,1,2$.

Consider the $\bar\pa=\bar\pa_{J_0}$-operator with $\R^n$ boundary conditions, with a negative exponential weight $-\gamma=-3<0$ at both ends
$$
\bar\pa\colon {\mathcal H}_{3,-\gamma}\bigl(\R\times[0,1],\C^n;\R^n,0_2\bigr)
\to{\mathcal H}_{2,-\gamma}\bigl(\C^n; 0_1\bigr).
$$
Here ${\mathcal H}_{3,-\gamma}\bigl(\R\times[0,1],\C^n;\R^n,0_2\bigr)$ denotes the space of $\C^n$-valued functions $w$ on the strip with three derivatives in $L^2$ with respect to the negative exponential weight at both ends, such that the restriction of $\bar\pa w$ to the boundary and the restriction of its derivatives to the boundary vanishes and ${\mathcal H}_{2,-\gamma}\bigl(\C^n; 0_1\bigr)$ denotes the space of $\C^n$-valued functions with $2$ derivatives in the weighted Sobolev space and which vanish to first order along the boundary. To be more specific we take the weight function to equal $e$ for $|\tau|\le 1$ and $e^{-\gamma|\tau|}$ for $|\tau|\ge 1$.  As in \cite{EES1}, we use the vanishing on the boundary to double the functions and still remain in a Sobolev space with $3$-derivatives in $L^2$. The calculations given in the proof of Proposition 5.3 \cite{EES1} then shows that $\bar\pa$ is a Fredholm operator of index $n$ and kernel spanned by the constant functions. In particular there is a constant $C$ such that for any function in the closed codimension $n$ subspace $W\subset {\mathcal H}_{3,-\gamma}\bigl(\R\times[0,1],\C^n;\R^n,0_2\bigr)$ of functions $w$ with $w(0)=0$,
\begin{equation}\label{e:ellest}
\|w\|_{3,-\gamma}\le C\|\bar\pa w\|_{2,-\gamma}.
\end{equation}
Let $\beta\colon\Theta_c\to\C$ be a cut off function which is real valued and holomorphic to first order on the boundary, which equals $1$ on $[-(c-2),(c-2)]\times[0,1]$, and equals $0$ outside $[-(c-1),(c-1)]\times[0,1]$. Then
\begin{align}\notag
&\|\bar\pa(\beta v)\|_{2,-\gamma}^2\le\\\notag
& \|\bar\pa v|[-1,1]\times[0,1]\|_{2}^2+
2K_0\int_1^\infty\lambda^2 e^{-2\gamma|\tau|}\,d\tau + 2K_1|c|^2\int_{c-1}^\infty e^{-2\gamma|\tau|}\,d\tau\\\label{e:smartweight}
&=\Ordo(\lambda^2) + \frac{K'|c|^2}{2\gamma} e^{-2\gamma(c-1)}
=\Ordo(\lambda^2)+K''(\log(\lambda^{-1}))^2\lambda^6=\Ordo(\lambda^2).
\end{align}
The first integral estimates the terms in $\bar\pa(\beta v)$ where all derivative falls on $v$ and the second estimates the terms where some derivative fall on $\beta$. Note that by the derivative bound $|D^k u|=\Ordo(\lambda)$, $k=1,2$, see Remark \ref{r:highstrap}, we get $|D^kv|$ bounded and hence $|D^{k-1}v|=\Ordo(|c|)$.

Thus \eqref{e:ellest} implies $\|\beta v\|_{3,-\gamma}=\Ordo(\lambda)$.
Since the Sobolev $(2,3)$-norm controls the supremum norm of the derivative we have
\begin{equation}\label{e:vsmall}
|D^k v|=\Ordo(\lambda),\quad k=0,1.
\end{equation}
In particular, writing $\hat u(\tau,t)=(\hat q_\lambda(\tau,t),\hat p_\lambda(\tau,t))$ we find
\begin{align*}
\nabla_t \hat p_\lambda(0,t)&=\frac{\pa \hat p_\lambda}{\pa t}(0,t)+\Gamma(\hat q_\lambda(0,t))\left(\frac{\pa \hat q_\lambda}{\pa t},\hat p_\lambda\right),\\
\nabla_\tau \hat p_\lambda &=\frac{\pa \hat p_\lambda}{\pa \tau}(0,t)+\Gamma(\hat q_\lambda(0,t))\left(\frac{\pa \hat q_\lambda}{\pa \tau},\hat p_\lambda\right),
\end{align*}
writing $u_0(\tau+it)=(q_0(\tau,t),p_0(\tau,t))$ we have $\frac{\pa q_0}{\pa t}=0=\frac{\pa p_0}{\pa \tau}$ and $\frac{\pa p_0}{\pa t}=[b_1(0)-b_0(0)]$. This together with \eqref{e:bdryclose}, \eqref{e:vsmall} and $\hat u=\lambda^{-1} u$ this implies
\begin{align*}
\lambda^{-1}\nabla_t p_\lambda(0,t)&=(b_1(0)-b_0(0))+\Ordo(\lambda),\\
\lambda^{-1}\nabla_\tau p_\lambda &=\Ordo(\lambda).
\end{align*}
\end{pf}

\begin{rmk}
Lemma \ref{l:gradconv} shows that the holomorphic disk are $C^1$-close to a map associated to a gradient flow. It is not hard to show that also higher derivatives are close, one just needs to control more derivatives in \eqref{e:bdryclose}, see Remark \ref{r:highstrap}, and use the estimate \eqref{e:ellest} for Sobolev spaces with more derivatives.
\end{rmk}

\begin{rmk}\label{r:gradconv}
A result similar to Lemma \ref{l:gradconv}, with less detailed convergence information, appears as Proposition 9.8 in \cite{FO}. The proof presented there does however contain some misstatements. In particular, the estimate $|\nabla_sp|\le |\nabla p||\frac{\pa q}{\pa s}|$ appearing in that proof is not true as can be checked by applying the construction given there to the holomorphic map $\sigma+is\mapsto e^{-\theta(\sigma+is)}$, $\sigma+is\in [0,\infty)\times[0,1]$ with boundary components mapping to $\R$ and $e^{-i\theta}\R$, where $q+ip$ are coordinates on $\C$.
\end{rmk}

\subsubsection{Projected flow trees}\label{ss:projft}
We next make a more detailed study of the parts of the disks in a sequence which map close $\Pi(\Sigma)$. That is we restrict attention to $D_j(\lambda)\subset\Delta_r$, $j\ge 1$. Let $W_j(\lambda)$ be neighborhoods as in Lemma \ref{l:gradconv} and let $\Theta\subset D_j(\lambda)-W_j(\lambda)$ be a strip region. Since $|Du^{T_j}|=\Ordo(\lambda)$ in $\Theta$ it follows that if $l\subset\Theta$ is any vertical line segment then $\pi(u^{T_j}(l))$ is contained in a $\Ordo(\lambda)$-ball around a point. Note that there exists $d>0$ such that over any $d$-ball in $\Sigma_j^\circ$ the projection $\pi_j(\hat L_\lambda\cap U(j,a)-U(j+1,b))$ consists of a number of distinct sheets. We say that $\Theta$ is {\em special} if it contains a vertical line segment $l$ such that $u^{T_j}(\pa l)$ lies in the same sheet. Otherwise we say that $\Theta$ is non-special.

\begin{lma}\label{l:specstrip}
Let $\Theta\subset D_j(\lambda)$ be a special strip region and let $\epsilon>0$. Then, after passing to a subsequence, there exists a point $p\in\Sigma_j^\circ$ such that $u^{T_j}_\lambda(\Theta)$ is contained in an $\epsilon$-ball around $p$.
\end{lma}

\begin{pf}
After passing to a subsequence we may assume that $\pi_j(u_\lambda(l))$ converges to $p\in \Sigma_j^\circ$. Assume that for all small $\lambda>0$ $\pi_j(u_\lambda)$ leaves a $\epsilon$-ball around $p$. Consider the preimage of the $\frac12\epsilon$-ball around $p$. Note that for $\lambda$ large enough we can find a strip region $[-c,c]\times[0,1]$ containing this preimage. Note next that the projection of $\hat L_\lambda$ is an exact Lagrangian submanifold. Letting $a$ denote the generating function of the sheet in which the endpoints lie we have,
\begin{align*}
\area(u_\lambda^{T_j})&\le K\int_{\pa[-c,c]\times[0,1]} p\,dq\\
&=-\int_{\{-c\}\times[0,1]}p\,dq+\int_{\{c\}\times[0,1]}p\,dq\\
&+a(u^{T_j}(-c,1))-a(u^{T_j}(c,1))-a(u^{T_j}(-c,0))+a(u^{T_j}(c,0)).
\end{align*}
Since $D u^{T_j}_\lambda=\Ordo(\lambda)$ and $|p|=\Ordo(\lambda)$ the first two terms are $\Ordo(\lambda^2)$. Moreover, since $a=\lambda \hat a$ for some smooth function $\hat a$ the derivative bound also implies that the remaining terms are $\Ordo(\lambda^2)$. Thus $\area(u_\lambda^{T_j}([-c,c]\times[0,1]))=\Ordo(\lambda^2)$ and therefore, by monotonicity \eqref{e:monbdry} $u_\lambda^{T_j}$ cannot leave a $K\lambda$-ball around $u^{T_j}_\lambda(l)$. This contradicts the strip region mapping to points at distance $\frac12\epsilon$ from $p$ for $\lambda$ small enough.
\end{pf}

We next turn our attention to non-special strip regions. The sheets of $\pi_j\hat L_\lambda$ determines local gradients in $\Sigma_j^\circ$. We define flow lines of these local gradient differences as usual and call them {\em projected flow lines}. A connected finite union of projected flow lines joined at zeros of the gradient differences will be called a {\em broken projected flow line}.

\begin{lma}\label{l:nspecstrip}
Let $\Theta(\lambda)\subset D_j(\lambda)$ be a non-special strip region and let $\epsilon>0$. Then, after passing to a subsequence, there exists a possibly broken projected flow line $\gamma$ such that $u^{T_j}_\lambda(\Theta(\lambda))$ is contained in an $\epsilon$-neighborhood of $\gamma$.
\end{lma}

\begin{pf}
We consider two cases. Assume first that $\Theta=[-c_\lambda,c_\lambda]\times[0,1]$ is such that $\lambda c_\lambda\le K$ for some $K$. Write $u_\lambda^{T_j}=(q_\lambda,p_\lambda)$. Since $u^{T_j}_\lambda$ is $J$-holomorphic we have
\begin{align*}
&\frac{\pa q_\lambda}{\pa\tau}+(\nabla_t p_\lambda)^\ast=0,\\
&\frac{\pa q_\lambda}{\pa t}-(\nabla_\tau p_\lambda)^\ast=0.
\end{align*}
Consider the map $(q_\lambda,p_\lambda)(\lambda^{-1}\tau,\lambda^{-1}t)$, $(\tau, t)\in [-\lambda c_\lambda,\lambda c_\lambda]\times[0,\lambda]$. Lemma \ref{l:gradconv} together with the above equations then imply
\begin{align}\label{e:conveq}
&\frac{\pa q_\lambda}{\pa \tau}-Y=\Ordo(\lambda),\\
&\frac{\pa q_\lambda}{\pa t}=\Ordo(\lambda),
\end{align}
where $\lambda Y$ is the local gradient difference determined by the (distinct) local sheets to which $u^{T_j}_\lambda$ maps the boundary. The statement now follows easily: after passing to a subsequence for which both the re-scaled lengths $\lambda^{-1}c_\lambda$ and the points $\pi(u^{T_j}_\lambda(-c_\lambda,0))$ converge it follows by standard ODE-results from \eqref{e:conveq} that the image of the strip region in this case lies in a small neighborhood of a flow line.

We next consider the case when $\lambda c_\lambda$ is unbounded.  First we note that the length $L$ of the part of $\Theta$ which map outside an $\epsilon$-neighborhood of the points where the local gradient difference vanishes must satisfy $\lambda L< K$ for some $K$. If not then we could use the argument above to show that the length of $u^{T_j}_\lambda(\pa\Theta)$ is unbounded. This however contradicts Lemma \ref{l:finitelength}. To finish the proof we write $\Theta=\Theta'\cup\Theta''$ where $\Theta''$ is the preimage of $\epsilon$-neighborhoods of the points where the local gradient difference vanishes. If we know that the number of components of $\Theta''$ is bounded then the lemma follows using the above argument on each component of $\Theta'$ since the rest of the image is contained in small balls around the vanishing points. However, it follows from the local form of a non-degenerate gradient zero that the length of any flow line in the exterior of an $\epsilon$-neighborhood of two vanishing points (possibly non-distinct) has length bounded from below and therefore the number of components is finite.
\end{pf}

Consider the local gradients of the projection of $\hat L_\lambda$ to $T^\ast\Sigma_j^\circ$. A {\em projected flow tree} is a partial flow tree defined as usual with respect to these local gradients. We also allow degenerate projected flow trees which are just points.

\begin{lma}\label{l:projft}
For any $\epsilon>0$, after passing to a subsequence, there exists for $1\le j\le k$, a finite collection $\Gamma_j^1,\dots\Gamma_j^{s_j}$ of projected flow trees such that for all sufficiently small $\lambda$, any connected component of $u^{T_j}_\lambda(D_j(\lambda))$ lies in an $\epsilon$-neighborhood of some $\Gamma_j^{s_r}$.
\end{lma}

\begin{pf}
After Lemmas \ref{l:specstrip} and \ref{l:nspecstrip} we need only note that the finitely many components of $W_j(\lambda)$ each have diameter $\Ordo(\log(\lambda^{-1}))$ therefore the derivative estimate implies that the image under $u^{T_j}_\lambda$ of such a component lies in an $\Ordo(\lambda\log(\lambda^{-1}))$-ball around some point.
\end{pf}

\begin{rmk}
As a consequence of Lemma \ref{l:gradconv} and the choice of the overlaps between $D_j(\lambda)$ and $D_s(\lambda)$, $s\ge j$ we find that the limiting (projected) flow trees of $u^{T_j}$, $j<s$ give endpoint conditions for the projected flow trees of $u^{T_s}$.
\end{rmk}

\subsection{Total average linking number and blow-up}
\label{s:totavlkbu}
In order to limit the amount of derivative blow up on the $\lambda$-scale we will define the total average linking number of a holomorphic disk with boundary on $\hat L_\lambda$ and prove that it is uniformly bounded. The definition is presented in \S\ref{ss:deftotavlk} (and is rather long and technical). It depends on a lemma which we consider next.

Let $u_\lambda\colon \Delta_r^\lambda\to T^\ast M$ be a sequence of holomorphic maps as in Lemma \ref{l:projft}. Assume that $\delta_1<\dots<\delta_k$ in the construction of the subdivision of the domain of $u_\lambda$ into $D_j(\lambda)$, $j=0,\dots,k$ are sufficiently small so that the local coordinate description $(q,s)=(q,s_1,\dots,s_j)$ of a neighborhood of $\Sigma_j^\circ$, see \S\ref{ss:TRisoCE}, is valid for $|s_i|\le 100\delta_k$. Assume that Lemma \ref{l:projft} holds with $\epsilon_0\ll\delta_1$ and fix $\epsilon$ with $0<\epsilon_0\ll\epsilon\ll \delta_1$ . For $p\in\Sigma_j^\circ$, define $V(j,\epsilon;p)$ as the set of all $(q,s)\in\Sigma_j^\circ\times\R^{n-j}$ in standard local coordinates on $M$ around $p$ such that $q\in B(p;\epsilon)$, the $\epsilon$-ball around $p$, and $|s|\le 10\delta_j$. Also define $\tilde V(j,2\epsilon;p)$ as the subset of all $(q,s)$ such that $q\in B(p;2\epsilon)$, $|s|\le 10(\delta_j+\epsilon)$.

Let $p\in \Sigma_j- U(j+1,\frac78\delta_{j+1})$ be any point in a projected flow tree $\Gamma$ in $\Sigma_j^\circ$ as in Lemma \ref{l:projft} (i.e. some component of $u_\lambda^{T_j}(D_j(\lambda))$ converges to $\Gamma$). Let $\Gamma_r(p)$ denote the connected component of the projected flow tree $\Gamma$ which lies in $B(p;r)$ and let $A_\lambda(p;r)\subset D_j(\lambda)$ be a connected component of the preimage of an $\epsilon_0$-neighborhood of $\Gamma_r$ under $u_\lambda^{T_j}$.

\begin{lma}\label{l:smallaction}
There exists $M>0$ and $C>0$ such that for all sufficiently small $\lambda>0$ there exists a connected subset $E_\lambda(p;\epsilon)\subset \Delta_r^\lambda$ with $E_\lambda(p;\epsilon)\supset A_\lambda(p;\epsilon)$ such that $E_\lambda(p;\epsilon)$ is bounded by at most $C$ vertical segments $\approx\{\tau\}\times[0,1]$ of $\Delta_r^\lambda$ with
$$
\int_0^1\left|\frac{\pa u_\lambda}{\pa t}\right|^2\, dt\le M\lambda^2,
$$
and such that $u_\lambda(E_\lambda(p;\epsilon))\subset \tilde V(j,2\epsilon;p)$.
\end{lma}

\begin{pf}
Assume first that $p\in\Sigma_1-U(2,\frac78\delta_2)$. Note that there exists $\epsilon'$, slightly larger than $\epsilon$, such that the component of the projected flow tree containing $p$ intersects $B(p;\epsilon')$ and $B(p;2\epsilon')$ transversely in its edges. The boundary of the projected tree $\Gamma_{\epsilon'}(p)$ consists of two types of points: $1$-valent vertices inside $B(p;\epsilon')$ and intersection points with $\pa B(p;\epsilon')$. The first type of points correspond to an asymptotic condition from a tree in $\Sigma_0^\circ$ (i.e. a point where the disk $u_\lambda$ enters $U(1,\delta_1)$). In particular, it follows from Lemma \ref{l:estgood1} that $|D u_\lambda|=\Ordo(\lambda)$ in this region and we find a vertical segment of the desired type. For the second kind of boundary point $q$ consider the connected component of the tree $\Gamma_{2\epsilon'}(p)-\Gamma_{\epsilon'}(p)$ which contains $q$. Note that other boundary points of this tree lie either on $\pa B(p;2\epsilon')$ or at vertices inside $B(p;2\epsilon')-B(p;\epsilon')$. (Note that $\epsilon'\ll\delta_1\le\frac18\delta_2$ so that the flow tree is still defined in this region.) At the inside vertices we find vertical segments as above. For the other vertices, Lemmas \ref{l:projft} and \ref{l:gradconv} imply that the preimage of the neighborhood of any tree of length at least $\epsilon'$ contains a collection of strip regions $R$ of total length at least $\epsilon'(c\lambda)^{-1}$, for some $c>0$. Moreover, the number of strip regions is bounded by a constant since the topology of $D_j(\lambda)$ is fixed. We have
$$
K\lambda\ge\area(u_\lambda)\ge K'\|Du_\lambda\|_{L^2}\ge K'\int\!\int_R\left|\frac{\pa u_\lambda}{\pa t}\right|^2\,d\tau\wedge dt,
$$
It follows that the total length of the regions where $\int\left|\frac{\pa u_\lambda}{\pa t}\right|^2\,dt\ge M\lambda^2$ is bounded by $\frac{K}{K' M}\lambda^{-1}$. Thus, choosing $M>0$ large enough the statement follows.

The general case follows in a similar way using induction. If $p\in\Sigma_j-U(j+1,\frac12\delta_{j+1})$ then the partial tree $\Gamma_{\epsilon'}(p)$ has the same type of boundary components as before. Boundary points which are not vertices but intersections with $\pa B(p;\epsilon')$ connected to $\pa B(p;2\epsilon')$ are treated exactly as above. The boundary points which are true vertices of $\Gamma_{\epsilon'}(p)$ come from asymptotic conditions of projected flow trees in $\Sigma_r^\circ$, $r< j$ by construction of the decomposition $D_j(\lambda)$. However, there are $\epsilon$-neighborhoods of the parts of the flow tree which gives the asymptotics. Using these neighborhoods and the inductive assumption it is easy to produce $E_\lambda(p;\epsilon)$.
\end{pf}

\subsubsection{Definition of total average linking number}\label{ss:deftotavlk}
Let $u_\lambda\colon\Delta_r^\lambda\to T^\ast M$, $\lambda\to 0$ be a sequence of holomorphic disks. The finite collection of finite length projected flow trees of Lemma \ref{l:projft} constitute a compact set. Let $\epsilon>0$ be small and consider a finite cover of the projected flow trees by subsets $V(j,\epsilon;p)$ as above. More precisely, we first cover the projected trees in $\Sigma_k$ with $\epsilon$-balls. Then we cover the projected flow trees in $\Sigma_{k-1}-U(k,\frac12\delta_k)$ in a similar way. Continuing the construction we cover, in the last step the projected flow trees in $\Sigma_1-U(2,\frac12\delta_2)$. By compactness, the number of $V(j,\epsilon;p)$-neighborhoods may be chosen finite.

Let $u_\lambda$ and $E_\lambda(p;\epsilon)$ be as in in Lemma \ref{l:smallaction}. We will define an average linking number for the restriction of $u_\lambda\colon E_\lambda(p;\epsilon)\to V(j,p;2\epsilon)$ as follows. Inside the region where the special coordinates $(q,s)\in\Sigma_j^\circ\times\R^j$ around $\Sigma_j^\circ$ are defined, $\hat L_\lambda$ splits, respecting the local product structure $T^\ast \Sigma_j^\circ\times\C^{n-j}$, as a fibration over a Lagrangian submanifold in $T^\ast \Sigma_j^\circ$ and {\em line-curves} and {\em bend-curves} in each orthogonal $\C$-direction. As in \S\ref{ss:TRisoCE} we will write the Lagrangian submanifold as a graph of a $1$-form
$$
q\mapsto \bigl(q,\lambda b(q)\bigr)\in T^\ast\Sigma_j^\circ.
$$
We will denote the product of the line- and bend-curves in the orthogonal $\C$-directions as follows,
$$
q\mapsto
\hat l_1(q,s_1)\times \hat l_2(q,s_2)\times\dots\times \hat l_j(q,s_r)\subset\C^j,
\quad,\quad (s_1,\dots,s_j)\in\R^j, q\in\Sigma_j^\circ.
$$
In the notation of \S\ref{ss:TRisocomp} there are three different forms for $\hat l_r$. If the sheet under consideration is a cusp sheet then exactly one of the $\hat l_r$ is a bend-curve and has the form $s_r\mapsto \bigl(\frac12 s_r^2,\beta(q)+\alpha(q)\gamma_\lambda(s_r)\bigr)\in\R^2=\C$. If the sheet is smooth then it contains no bend curve. Other components are line-curves which have either the form $s_r\mapsto \bigl(s_r,\lambda h(q)\bigr)$ or, close to the tangency locus, $s_r\mapsto \bigl(s_r,\hat h_\lambda(q,s_r)\bigr)$.

Assume that $\epsilon>0$ is small and consider the $r^{\rm th}$ complex plane factor $\C\subset\C^j$ perpendicular to $T^\ast\Sigma_j^\circ$ in the fiber over the point $p\in\Sigma_j^\circ$. We will denote the map $u_\lambda$ composed with the projection to this copy of $\C$, $u_\lambda^{\perp_r}$. Note that by definition of the complex structure on $T^\ast M$, $u_\lambda^{\perp_r}$ is holomorphic with respect to the standard complex structure on $\C$.

Consider the line and bend curves, $\hat l^0_r(q),\dots,\hat l^m_r(q)$ of all the sheets over $p$. Note that exactly one of these is a bend-curve, we take it to be $\hat l^0_r(q)$, and that all others are line-curves. Let $l^s_r(q)=\lambda^{-1}\hat l^s_r(q)$. Since $\epsilon>0$ is small there exists a rectangle $R_d$ of width $1$ and height $d>c\sqrt{\epsilon_2}$ (where $\epsilon_2$ is the deformation parameter mentioned in \S\ref{ss:LegisoCE}) for some $c>0$ (independent of $\epsilon>0$ provided it is sufficiently small) such that the distance form any bend-curve $l^0_r(q)$ to $R_d$ is at least $2d$, for $q$ varying in an $\epsilon$-neighborhood of $p$.
\begin{lma}
For almost every $\zeta\in R_d$ the curve $\lambda^{-1}u^{\perp_r}_\lambda(\pa E_\lambda(p;\epsilon))\subset \C$ is disjoint from $\zeta$.
\end{lma}

\begin{pf}
By Lemma \ref{l:smallaction} the images of the vertical boundary components of $E_\lambda(p;\epsilon)$ under $\lambda^{-1}u_\lambda^{\perp_r}$ are finite length curves. Hence they meet $R_d$ in a subset of measure $0$. Since no bend meets $R_d$ the boundary components mapping to the sheet with bend never intersects $R_d$. Finally, by Lemma \ref{l:estgood} the derivatives of the re-scaled map of other boundary components are uniformly bounded. Hence also these curves meet $R_d$ in a subset of measure $0$.
\end{pf}

Thus for almost every $\zeta\in R_d$ the linking number $\lk\bigl(\lambda^{-1}u_\lambda^{\perp_r}(\pa E_\lambda(p;\epsilon)),\zeta-\infty\bigr)$ is well defined. (Here we think of $\zeta-\infty$ as a $0$-cycle in the one point compactification of $\C$. This linking number is also known as the winding number with respect to $\zeta$). We define the {\em average linking number} of $u_\lambda$ as
$$
\bar\lk(u_\lambda^\perp)=\sum_{r=1}^j\frac{1}{\area(R_d)}\int_{R_d}
\lk\bigl(\lambda^{-1}u_\lambda^{\perp_r}(\pa E_\lambda(p;\epsilon)),\zeta-\infty\bigr)\,dA.
$$

We define the {\em total average linking number} of a member of a sequence of holomorphic maps $u_\lambda\colon \Delta_r(\lambda)\to T^\ast M$, with respect to a fixed $\epsilon$-cover of the projected flow tree as described above, to be the sum of the average linking numbers of $u_\lambda|E_\lambda(p;\epsilon)$, where $E_\lambda(p;\epsilon)$ are as in Lemma \ref{l:smallaction}, with respect to the chosen $\epsilon$-cover of the projected trees to which $u_\lambda^{T_j}$ limit. Note that this is a sum with a finite number of terms (by finiteness of the cover).

\begin{lma}\label{l:totavlkfin}
There exists a constant $C$ such that for sufficiently small $\epsilon>0$, the average linking number of $u_\lambda^{\perp_r}|E_\lambda(p;\epsilon)$ in any $\C$-direction is bounded by $C$. In particular, the total average linking number of $u_\lambda$ is uniformly bounded for fixed $\epsilon>0$ small enough.
\end{lma}

\begin{pf}
For simpler notation, let $\lambda^{-1}u_\lambda^{\perp_r}=v$, $E_\lambda(p;\epsilon)=X$ and let
$\gamma=\lambda^{-1}u_\lambda^\perp(\pa E_\lambda(p;\epsilon))$. Also, let $z=x+iy$ be standard coordinates on $\C$ and assume that the bend $l_0(p)$ has its minimum at $0$ (i.e. $\beta(q)=0$). To compute
$$
\lk\bigl(\gamma,\zeta-\infty\bigr),
$$
let $\zeta=\xi + i\eta$ and consider the vertical line $l_\zeta=\{x+iy\colon x=\xi\}$ subdivided into two parts $l_\zeta^+=\{x+iy\colon x=\xi, y>\eta \}$ and $l_\zeta^-=\{x+iy\colon x=\xi, y<\eta\}$. Note that $\lk\bigl(\gamma,\zeta-\infty\bigr)$ equals the algebraic sum of intersections of $\gamma$ and $l_\zeta^+$. Consider first the parts of $\gamma$ which are the images under $v$ of the components of $\pa X$ which map to the bends. Since every such component is homotopic, in the complement of $\zeta$, to a curve which lies in a constant bend it follows that the contribution to $\lk\bigl(\gamma,\zeta-\infty\bigr)$ from such a component is $0$ or $\pm 1$.

To estimate the contribution from other components note that the algebraic number of intersections with $l_\zeta^+$ is estimated from above by the number of pairs of intersection points $(a,b)$, $a\in\gamma\cap l_\zeta^\pm$ and $b\in \gamma\cap l_\zeta^\mp$ such that there is no intersection point $c\in \gamma\cap l_\zeta$ between $a$ and $b$.

We next subdivide $X$. First we remove $\log(\lambda^{-1})$-neighborhoods of the vertical boundary components and of the boundary minima in $X$. The remaining parts of $X$ are strip regions and by Lemma \ref{l:projft} the restriction of $u_\lambda^{T_j}$ to such a region converges to a flow line at rate $\Ordo(\lambda)$.

Consider first such a strip region the length of which grows faster than  $\lambda^{-1}$. As in the proof of Lemma \ref{l:projft} we see that the image under $u^{T_j}$ of such a region is confined to a small neighborhood of a point where the local gradient difference vanishes. Since such points lie far from the tangency locus we see that there can be no pairs $(a,b)$ in such a region since no $l^s(q)$ passes $y=\eta$.

Consider next a strip region of length $\Ordo(\lambda^{-1})$. Fix $\delta>0$ such that $\delta$ equals half of the support of the cut-off function $\psi$, see \S\ref{ss:TRisocomp}. We subdivide the pairs $(a,b)$ into long pairs and short pairs. A long pair is a pair connected by an arc intersecting $\{|x|=\lambda^{-1}\delta\}$. In particular, such an arc has length bounded from below and by Lemma \ref{l:finitelength} the number of such pairs is bounded. Consider next the short arcs. Note that between any two points in a short arc $h(q)-\eta$ must change sign and since we are in the region $\{|x|\le \lambda^{-1}\delta\}$, $u_\lambda^{T_j}$ must map the point of sign change into the region near the tangency locus. Let $\sigma(\eta)$ denote the number of solutions to $h(q)-\eta=0$ along the short arcs. Then the contribution to the average linking number can be estimated by the integral
$$
\frac{1}{\area(R_d)}\int_{R_d} \sigma(\zeta)\, dA.
$$
On the other hand this integral measures the length of the part of the curve $h(q(\tau))$ which lies inside $\{|y|\le d\}$. Consider the corresponding integral for the cotangent lift of the flow line to which $u_\lambda^{T_j}$ converges. Note that since flow lines are maximally transverse to the level surfaces of $h$ near the tangency locus there exists $\epsilon'$ so that for any flow line of length at most $\epsilon'$ there are no more than $j$ intersections with a level surface of $h$ (compare the proof of Lemma \ref{l:dens2}). Thus for the flow tree the corresponding integral is bounded. On the other hand the length of the curve can be written as
$$
\int_{I'}\left|\frac{dh(q(\tau))}{d\tau}\right|\,d\tau,
$$
where $I'$ is the part of the finite length interval corresponding to the re-scaled domain in the proof of Lemma \ref{l:gradconv} which map into the region $\{|y|\le d\}$. As the curve $\lambda^{-1}u_\lambda^{T_j}$ on this domain $C^1$-converges at an $\Ordo(\lambda)$-rate to the cotangent lift we find that the contribution also from this part is bounded.

Finally, we consider the contribution of the remaining pieces. We bound it using the projected length as above. The vertical segments have bounded length hence also their projections have bounded length and they give a finite contribution. We subdivide the remaining pieces into a length $K>0$ piece near the vertical segment and a complementary piece. Using the argument in Lemma \ref{l:gradconv} we conclude that on the complementary piece we can bound the $C^1$-distance between the cotangent lift of a gradient tree and the map $\lambda^{-1}u^{T_j}$ by $K'$. (The last term in equation \eqref{e:smartweight} is bounded by a constant when $c=K$.) Since the re-scaled length of the domain of definition of this curve is $\Ordo(\lambda\log(\lambda^{-1}))$ we find that the length of the projection and therefore the linking contribution goes to $0$.

Any remaining piece $I$ of the boundary have finite length $\le 10K$. The derivative of $h(q)$ is bounded by $\epsilon_3\lambda^{-1}$, where $\epsilon_3$ is the speed parameter, and the derivative of $q(\tau)$ is $\Ordo(\lambda)$. Thus
$$
\int_I\left|\frac{dh(q(\tau))}{d\tau}\right|\,d\tau=\Ordo(1).
$$
The lemma follows.
\end{pf}

\begin{rmk}
Note that the linking number $\lk\bigl(\lambda^{-1}u_\lambda^{\perp_r}(\pa E_\lambda(p;\epsilon)),\zeta\bigr)$ can also be computed by counting the number of points in $({\lambda^{-1}u_\lambda^{\perp_r}})^{-1}(\zeta)$ with signs. As $\lambda^{-1}u_\lambda^{\perp_r}$ is holomorphic, all signs are positive. \end{rmk}

\subsubsection{Blow up analysis}
Consider a sequence of disks $u_\lambda\colon \Delta_{r}^\lambda\to T^\ast M$ with subdivision $\{D_j(\lambda)\}_{j=0}^k$. We will show that by adding a finite number of punctures to $\Delta_{r}^\lambda$ and considering the maps $u_\lambda$ as defined on the new domains, we obtain maps which satisfies an $\Ordo(\lambda)$-derivative bound everywhere for all sufficiently small $\lambda>0$. The first step in this process is to understand the geometry at points in $\Delta_r^\lambda$ where such a derivative bound does not hold for the original maps $u_\lambda$. To this end let $M(\lambda)=\sup_{\Delta_{r}^\lambda}|Du_\lambda|$ and assume that $\lambda [M(\lambda)]^{-1}\to 0$ as $\lambda\to 0$. By the asymptotics of $u_\lambda$ near its punctures it follows that there exists points $p_\lambda\in\Delta_{r}(\lambda)$ such that $|Du_\lambda(p_\lambda)|=M(\lambda)$. Moreover, Lemma \ref{l:estgood1} implies that $p_\lambda\in D_j(\lambda)$, $j>0$.
\begin{lma}\label{l:blowup}
There exists $\lambda_0>0$ such that for all $\lambda<\lambda_0$ there exists exactly one point in a $\frac{\lambda}{M(\lambda)}$-neighborhood of $p_\lambda$ which maps to $\Sigma\subset \hat L_\lambda$. Moreover,  $u_\lambda|B(p_\lambda,\frac{\lambda}{M(\lambda)})$ contributes $1$ to the total average linking number of $u_\lambda$.
\end{lma}

\begin{pf}
By \eqref{e:estgoods}, the image of $u_\lambda^{T_j}$ restricted to $B(p_\lambda,\frac{\lambda}{M(\lambda)})$ approaches a point $q\in \Sigma_j^\circ$.
Consider local coordinates $\C$ in one of the directions orthogonal to $T^\ast \Sigma_j^\circ$ near $q$ and consider the sequence of maps
$$
\tilde u_\lambda=\lambda^{-1} u^\perp \left(p_\lambda+\frac{\lambda}{M(\lambda)}z\right),
$$
defined for all $z\in \C$ such that $p_\lambda+\frac{\lambda}{M(\lambda)}z\in B(p_\lambda,1)\cap \Delta_r^\lambda$. We use a standard blow-up argument, see e.g. \cite{Si}, and so we only sketch the details.  Note that the domain of these maps approach either $\C$ if $p_\lambda$ stays away from the boundary or the upper half plane if $p_\lambda$ goes to the boundary sufficiently fast.
Since the derivative of $\tilde u_\lambda$ is uniformly bounded, we conclude that some subsequence converges to a non-constant holomorphic map with derivative $1$ at $0$. Moreover, by Lemma \ref{l:subh}, the imaginary part of $\tilde u_\lambda$ is uniformly bounded. If follows that $\tilde u_\lambda$ converges to a non-constant map with boundary on the re-scaling of $\C\cap\hat L_\lambda$ and with bounded imaginary part. Thus, the domain of the limit map must be the upper half-plane and the limit map itself must be a map in the $1$-parameter family of maps taking the upper half plane to the bend, which are asymptotic to a linear map at infinity, and which have derivative of magnitude $1$ at $0$. In particular, for all $\lambda$ small enough $\lambda^{-1}u_\lambda^{\perp}$ covers the region $R_d$ exactly once. The statement about the total average linking number follows.
\end{pf}

\begin{cor}\label{c:blowup}
After adding a uniformly finite number $r'$ of punctures to $\Delta_r^\lambda$ there exists $K>0$ such that the maps $u_\lambda\colon\Delta_{r+r'}^\lambda\to T^\ast M$ satisfies
$$
|Du_\lambda|\le K\lambda.
$$
\end{cor}

\begin{pf}
If $p_\lambda$ is as above, note that after adding a puncture at the point near $p_\lambda$ which maps to $\Sigma$, the derivative satisfies the bound in a neighborhood of the puncture. We may now repeat the argument for a new maximizing sequence $p_\lambda'$ on the new domains. Lemma \ref{l:blowup} implies that the total average liking number contribution of each blow up sequence is $1$. The bound in Lemma \ref{l:totavlkfin} then implies that we obtain the desired estimate by adding a finite number of punctures.
\end{pf}

\subsubsection{Reconstructing the domain}
Note that the domain $\Delta_r^\lambda$ has three kinds of punctures: Reeb chord punctures, blow-up punctures, and intersection punctures (i.e. the punctures added in \S\ref{ss:subdivseq}). We next show how to change them into only two types (in a sense removing unnecessary  intersection punctures).

For two integers $m$ and $n$, let $\Delta_{m,n}$ be a standard domain $\Delta_m$ with $n$ marked points on the boundary. If $\Delta_{m+n}$ is a standard domain then note that there exists a unique $\Delta_{m,n}$ that admits a conformal map to $\Delta_{m+n}$ taking the marked points of $\Delta_{m+n}$ to punctures of $\Delta_{m+n}$. Moreover, in a neighborhood of any puncture of this map has the form $z\mapsto z +\sum_{n\le 0}c_n e^{n\pi z}$ and in a neighborhood of a marked point $z\mapsto \frac{1}{\pi}\log z+\sum_{n\ge 0} c_n z^n$, where $c_n$ are real constants.

Let the intersection punctures of $\Delta_r^\lambda$ be $p_1,\dots,p_s$. Choose any $k$ of these. Then there is a unique holomorphic map $\alpha\colon\Delta_{r-k,k}^\lambda\to\Delta_r^\lambda$.

\begin{lma}\label{l:unnecpun}
There exists $k$, with $0\le k\le s$ such that $u_\lambda\circ\alpha\colon \Delta_{r-k}^\lambda\to T^\ast M$ satisfies
$$
|Du_\lambda\circ\alpha|=\Ordo(\lambda),
$$
and such that near any puncture of $\Delta_{r-k}^\lambda$ the map looks like a blow-up puncture or a Reeb chord puncture.
\end{lma}

\begin{pf}
Start by removing all intersection punctures. If the $\Ordo(\lambda)$ derivative bound holds we are done. Thus, assume it does not hold. Cutting out any $\epsilon$-neighborhood of the marked points we claim that the conformal map $\alpha\colon \Delta_{r-k,k}^\lambda\to\Delta_r^\lambda$ satisfies $|D\alpha|=\Ordo(1)$ on the complement of the cut out neighborhoods. To see this we view the map as a map into $\C$ with bounded imaginary part. Assuming that no such bound for $\alpha$ holds, we find a sequence of points where the derivative is maximal. Using similar (but simpler) blow-up arguments as above we derive a contradiction.

It follows from this claim that any blow up sequence $p_\lambda$ of the type considered in Lemma \ref{l:blowup} on the domain $\Delta_{r-k,k}^\lambda$ must converge to a marked point in $\Delta_{r-k,k}^\lambda$. Repeating the argument given in the proof of Lemma \ref{l:blowup} we conclude that in a neighborhood of such a marked point the maps looks like a blow-up puncture. We then repeat the argument with $\Delta_{r-k+1,k-1}^\lambda$ not making the intersection puncture corresponding to the limit of the blow up sequence a marked point. This way, we we find the desired domain in a finite number of steps.
\end{pf}

We call the new domains $\Delta_t^\lambda$ and its two types of punctures Reeb chord punctures and blow-up punctures.

\begin{cor}\label{c:leave}
Let $p$ be a blow-up puncture of $u_\lambda\colon \Delta_t^\lambda\to T^\ast M$ and let $[0,\infty)\times[0,1]\subset\Delta_t^\lambda$ be its neighborhood. Then there exists $\rho_0>0$ such that for all $\lambda>0$ large enough $u_\lambda([\rho_0,\infty)\times[0,1])$ lies inside a $K\lambda$-neighborhood of $\Sigma$ and the region $u_\lambda([\rho_0-10,\rho_0]\times[0,1])$ lies outside of a $\frac{K}{2}\lambda$ neighborhood of $\Sigma$. Moreover, there exists $a>0$ such that
$$
\area(u_\lambda|[\rho_0,\infty)\times[0,1])\ge a\lambda^2,
$$
for all sufficiently small $\lambda>0$.
\end{cor}

\begin{pf}
Since the map limits to a standard map from the upper half-plane to the bend we conclude that there is $\rho_0\in[0,\infty]$ such that for all $\lambda>\lambda_0$, $[\rho_0,\infty)\times[0,1]$ maps to the part of the bend at distance $K\lambda$ from $\Sigma_j$. The second statement follows from the asymptotics of the limit map and the third follows since the standard map which covers the half disk bounded by the half-circle in the bend.
\end{pf}

\begin{rmk}\label{r:areaend}
Since the diffeomorphism taking $\Pi_\C(\tilde L_\lambda)$ to $\hat L_\lambda$ is $C^0$ $\Ordo(\lambda)$-close to $\id$ and $C^1$ $\epsilon$-close. It follows from Corollary \ref{c:leave} that an area bound similar to the area bound for $u_\lambda$ holds also for $\Phi_\lambda^{-1}\circ u_\lambda$.
\end{rmk}

Let $W_0(\lambda)$ be a $\log(\lambda^{-1})$-neighborhood of all boundary minima of $\Delta_t^\lambda$ and let $W_1(\lambda)$ be the union of the neighborhoods of all blow-up punctures as in Corollary \ref{c:leave}.

\begin{lma}\label{l:onesheet}
Let $u_\lambda\colon\Delta_t^\lambda\to T^\ast M$ be a sequence of holomorphic disks with boundary on $\hat L_\lambda$, with only one positive puncture, and let $l\subset \Delta_t^\lambda -\bigl(W_0(\lambda)\cup W_1(\lambda)\bigr) $ be a vertical line segment such that $\{p_1,p_2\}=\pa l\subset \pa\Delta_t^\lambda$ then $p_1$ and $p_2$ map to different sheets of $\hat L_\lambda$.
\end{lma}

\begin{pf}
Note that $l$ lies in a strip region of $\Delta_t^\lambda$ and that by definition of $W_0(\lambda)$ and $W_1(\lambda)$ we can find a strip region $\Theta_\lambda
\approx[-\log(\lambda^{-1}),\log(\lambda^{-1})]\times[0,1]\subset\Delta_t^\lambda$  containing $l\approx\{0\}\times[0,1]$. Passing to a subsequence we may assume that $\Pi(u_\lambda((0,0)))$ converges to $m$ as $\lambda\to 0$ and by the derivative estimate $|Du_\lambda|=\Ordo(\lambda)$ on $\Delta_t^\lambda$, $u_\lambda(\Theta_\lambda)\subset T^\ast B(m_\lambda,K\lambda\log(\lambda^{-1}))$ for some $K>0$ and all sufficiently small $\lambda$, where $m_\lambda\to m$ as $\lambda\to 0$ and where $B(p,r)$ is an $r$-ball around $p$.

Consider the $\lambda^{-1}$-scaling $\lambda^{-1}u_\lambda$ of $u_\lambda|\Theta_\lambda$ in local coordinates around $m$. Note that this map has uniformly bounded derivative and takes the boundary to the scaling of $\hat L_\lambda$. We can thus extract a subsequence which converges to a standard holomorphic map from $\R\times[0,1]$ with boundary on the limit of the scaling of $\hat L_\lambda$ and with bounded imaginary part by Lemma \ref{l:subh}. However, the scaling approaches a totally real submanifold which is a product of a totally real submanifold $\Lambda\subset \C^{n-1}$ which is $C^1$-close (the distance is controlled by the speed parameter $\epsilon_3$) to an affine Lagrangian subspace parallel to the $0$-section, and a bend in the orthogonal $\C$-direction. Hence the $n-1$ first components of the limit map must be constants and the other component is either a standard map to the bend or constant. It then follows by scaling that for any $b>0$ there exists vertical line segments $l'\subset\Theta_\lambda$ with endpoints $p_1'$ and $p_2'$, such that $\left|\int_{l'} p\,dq\right|\le b\lambda^2$ and such that the length of $u_\lambda(l')$ is smaller than $b\lambda$. We conclude from this and the fact that $\Phi_\lambda^{-1}$ is $C^0$ $\Ordo(\lambda)$-close and $C^1$ $\epsilon$-close to $\id$ that
$$
\left|\int_{l'} (\Phi_\lambda^{-1}\circ u_\lambda)^{\ast}p\, dq\right| + |f(\Phi_\lambda^{-1}(u_\lambda(p_1')))-f(\Phi^{-1}(u_\lambda(p_2')))|\le K b\lambda^2,
$$
for some $K$ depending only on $\Phi_\lambda$, where $f$ is the local function generating $\tilde L_\lambda$.

Note next that $l'$ subdivides $\Delta_t^\lambda$ in two components: one component containing the positive puncture and another component which we call the {\em non-positive} component. The punctures of the non-positive component are either negative punctures or blow-up punctures and there exists at least one such puncture. Since the area contribution at each negative puncture can be bounded from below by $c\lambda$ for some constant $c>0$ we find that the area of the image of the non-positive component of $l'$ under $\Phi_\lambda^{-1}\circ u_\lambda$ is bounded from below by $a\lambda^2$, see Remark \ref{r:areaend}. However, by exactness of $\tilde L_\lambda$, the integral of $p\,dq$ along $\Phi_\lambda^{-1}\circ u_\lambda$ restricted to the boundary of the non-positive component is bounded below by $Kb\lambda^2$. Taking $b>0$ small enough we get a contradiction. Thus, $p_1$ and $p_2$ map to different sheets.
\end{pf}

Let $p$ be a blow up puncture and let $[0,\infty)\times[0,1]$ be a neighborhood of it where $\{0\}\times[0,1]$ lies at finite distance from the nearest boundary minimum. Let $T$ be the the minimum of $\log(\lambda^{-1})$ and the $\tau$-coordinate of the first point of $[0,\infty)\times[0,1]$, where $u_\lambda$ leaves a small $\delta$-neighborhood of $\Sigma$. Here we take $\delta>\delta_k$ where $\delta_k$ is the largest distance used to subdivide $\Delta_m(\lambda)$. In particular, all the cut-off functions $\psi(s)$, see \S\ref{ss:TRisocomp} have derivative with support in $|s|<\delta$.

\begin{lma}\label{l:glineblup}
The restriction of $u_\lambda$ to $[T,\infty)\times[0,1]$ converges to a flow line.
\end{lma}

\begin{rmk}
The argument of Lemma \ref{l:gradconv} gives a weaker result in this case since the diffeomorphism $\Phi_\lambda$ near the tangency locus satisfies a weaker bound. More precisely, in  \eqref{e:Jclose} and \eqref{e:bdryclose} the $\Ordo(\lambda)$-bound is replaced with an $\Ordo(\epsilon_3)$-bound, where $\epsilon_3$ is the speed parameter.
\end{rmk}

\begin{pf}
Let $S$ be the local sheet of $\hat L_\lambda$ into which the blow up puncture maps. The projection $\pi_T\colon T^\ast M\to T^\ast(\Pi (S\cap\Sigma))$ maps $S$ to the graph of a $1$-form on $\Pi(S\cap\Sigma)$. Moreover, since the metric on $M$ is locally a product metric around $\Pi(S\cap\Sigma)$ the map $\pi_T\circ u_\lambda$ is $J_T$-holomorphic where $J_T$ is the complex structure associated to the restriction of the metric on $M$ to $\Pi(S\cap\Sigma)$. The local projections just described fit together to give a map which we still denote $\pi_T\circ u_\lambda$, where $\pi_T\circ u_\lambda\colon [T,\infty)\times[0,1]\to T^\ast\Sigma$ and where $\pi_T\circ u_\lambda$ takes $\pa([T,\infty)\times[0,1])$ to the graph $\hat \Lambda_\lambda$ of a $1$-form. Furthermore, the map $\pi_T\circ u_\lambda$ is $J_T$-holomorphic where $J_T$ is the complex structure associated to the metric on $\Sigma$ induced from pulling back the metric on $M$ along local branches. As for $\hat L_\lambda$ itself it is not hard to see that there exists a fiber preserving diffeomorphism $\Psi_\lambda\colon T^\ast\Sigma\to T^\ast\Sigma$, which is at $C^0$-distance $\Ordo(\lambda)$ and $C^1$-distance $\epsilon$ from the identity, and which takes $\hat \Lambda_\lambda$ to an exact Lagrangian submanifold $\tilde\Lambda_\lambda$. Thus arguing as in the proof of Lemma \ref{l:onesheet} we find that $\int_{\{T\}\times[0,1]} (\pi_T\circ u_\lambda)^\ast p\,dq=\ordo(\lambda^2)$ and that the area of $\Psi_\lambda\circ \pi_T\circ u_\lambda$ is $\ordo(\lambda^2)$. Thus the same is true for the area of $\pi_T\circ u_\lambda$ and it follows from the monotonicity property of  $J_T$-holomorphic disks, see \eqref{e:monint} and \eqref{e:monbdry}, that the image of $\pi_T\circ u_\lambda$ lies in an $\ordo(\lambda)$-ball around some point $p_\lambda\in\Pi(\Sigma)$ which must lie in $\Pi(S\cap\Sigma)$.

Locally we have $T^\ast M= T^\ast(\Pi(S\cap\Sigma))\times\C$ and the intersection $\hat L_\lambda\cap \{(q,p)\}\times\C$ is a bend curve $\hat l(q)$ which depends on $q$. Note that the variation of $\hat l(q)$ is $C^0$-bounded by $\Ordo(\lambda)$ and $C^1$-bounded by $\Ordo(\epsilon_3)$, where $\epsilon_3$ is the speed parameter. Let $u_\lambda^\perp$ denote the composition of $u_\lambda$ with the projection to $\C$. Using the variation bounds in combination with the fact that image of $\pi_T\circ u_\lambda$ lies in an $\ordo(\lambda)$-neighborhood of a point and Corollary \ref{c:leave} which controls $u^\perp_\lambda$ in a neighborhood of $\infty$ it is straightforward to see that the argument in Lemma \ref{l:gradconv} applies to show that $u^\perp_\lambda$ converges to a flow line. This finishes the proof since flow lines of the gradient difference of the two newborn sheets of $\hat L_\lambda$ corresponding to $S$ are constant in the directions of $\Pi(S\cap\Sigma)$ in the neighborhood we consider.
\end{pf}

We enlarge $W_1(\lambda)$ including in it also the neighborhoods $[T,\infty)\times[0,1]$ of all boundary punctures with properties as in Lemma \ref{l:glineblup}. We next remove neighborhoods also of points mapping close to the tangency locus in $\Sigma_1$. Consider a strip region $\Theta_\lambda\subset \Delta_t^\lambda-(W_0(\lambda)\cup W_1(\lambda))$ and let $l$ be a vertical line segment in $\Theta_\lambda$. By Lemma \ref{l:onesheet} we know that the two points in $\pa l$ map to two distinct sheets $S_0$ and $S_1$ of $\hat L_\lambda$. We consider the case when one or both of these are cusp sheets. Assume that $l$ maps into the $(\rho^{-1}\lambda)$-neighborhood of the tangency locus of these two sheets, where $\rho=\epsilon_3$ is the speed parameter, see \S\ref{ss:TRisocomp}. Then $u_\lambda(l)$ lies in $U(j,\delta_j)$ for some $j$ and Lemma \ref{l:projft} implies that the projection $u^{T_j}$ converges to a projected flow line at rate $\Ordo(\lambda)$. We make two observations. First, projected flow lines are maximally transverse to the tangency loci in $\Sigma_j^\circ$ an thus the order of tangency between the projected flow line and the tangency locus is bounded. (For a uniform treatment of all tangency loci we use the upper bound $n-2$, $n=\dim(L)$, but the actual bound in $\Sigma_j^\circ$ is $n-j-1$.) Second, the points where the projected gradient difference vanishes lie outside the tangency locus therefore the length of any projected gradient in a neighborhood of the tangency locus is bounded from below by $c\lambda$ for some $c>0$. It follows form these two observations that
the time a projected flow line spends in a $(\rho^{-1}\lambda)$-neighborhood of the tangency locus can be bounded by $C\lambda^{\frac{1}{n-2}}$ for some constant $C>0$ and all sufficiently small $\lambda>0$.

\begin{rmk}\label{r:rigtan}
If we know a priori that the intersection between the tangency locus and the projected flow line is transverse then the estimate on the flow time can be strengthened to $\Ordo(\lambda)$.
\end{rmk}

Assume next that $\Theta_\lambda=[-c_\lambda,c_\lambda]\times[0,1]$ with $l=\{0\}\times[0,1]$ where  $c_\lambda\ge \lambda^{-(1-\frac{1}{2(n-2)})}$. Then it follows from projected flow tree convergence that there are  vertical segments $l_-$ and $l_+$ at distance at most $C\lambda^{-(1-\frac{1}{2(n-2)})}$ which map outside a $(2\rho^{-1}\lambda)$-neighborhood of the tangency locus. Let
$$
b_1^k<b_2^k<\dots< b_m^k, \quad k=1,2,
$$
be the points in $[-c_\lambda,c_\lambda]\times\{0\}$ such that $u_\lambda^{T_j}(b_j^k)$ lies in the boundary of the $(k\rho^{-1}\lambda)$-neighborhood of the tangency locus.
Mark each $b_j^2$ such that there exists $b_m^1$ with $b_j^2<b_m^1<b_{j+1}^2$. Call the union of all regions bounded by vertical segments at marked points and containing points mapping to the $(\rho^{-1}\lambda)$-neighborhood $W_3(\lambda)$. Note that each component in $W_3(\lambda)$ has diameter bounded by $C\lambda^{-(1-\frac{1}{2(n-2)})}$.

\begin{lma}
The number of components of $W_3(\lambda)$ is uniformly finite.
\end{lma}

\begin{pf}
We first note that there exists $\eta_1>0$ such that for small $\lambda>0$ any projected flow line segment of length $<\eta_1$ can contain at most $n$ consecutive intersections with a $(2\rho^{-1}\lambda)$- and a $(\rho^{-1}\lambda)$-neighborhood of the tangency locus. (If no such $\eta_1$ exists it is easy to construct a flow line with too high order of contact with the tangency locus by taking the limit as $\eta_1\to 0$, see Lemma \ref{l:dens2}.) Since the total length of all limiting projected flow lines is finite as is the number of strip regions the lemma follows.
\end{pf}

\begin{lma}\label{l:leaveall}
For any $K>0$, $u_\lambda\Bigl(\Delta_t^\lambda-\bigl(W_0(\lambda)\cup W_1(\lambda)\cup W_2(\lambda)\bigr)\Bigr)$ lies outside a $K\lambda$ neighborhood of $\Sigma$ for all sufficiently small $\lambda$.
\end{lma}

\begin{pf}
By definition of $W_j(\lambda)$, any point in $\Delta_t^\lambda-(W_0(\lambda)\cup W_1(\lambda)\cup W_2(\lambda))$ lies in the middle of a strip region $\Theta_\lambda$ which has length at least $\log(\lambda^{-1})$ and which maps outside a $(\rho^{-1}\lambda)$-neighborhood of the tangency locus. Since the endpoints of a vertical strip maps to different sheets $S_1$ and $S_2$ of $\hat L_\lambda$ the composition of $u_\lambda$ with the projection to the orthogonal complement of $T^\ast(\Pi(\Sigma\cap S_j))$, $j=1,2$, takes at most one boundary component to a bend and at least one to a line. Since we are not in the tangency locus the line does not intersect the bend. Noting that $\lambda^{-1}u_\lambda^\perp$ has bounded derivative, we extract a subsequence which converges to a holomorphic map from $\R\times[0,1]$ with boundary on the re-scaling of the bend and the line and with bounded imaginary part. Since there are no such maps passing the extremum of the bend, the extremum of the bend is not part of the limit. Therefore the image of $\Theta_\lambda$ lies outside a $K\lambda$-neighborhood of $\Sigma$.
\end{pf}

\subsubsection{Flow tree convergence}
We next associate a source tree to $\Delta_t^\lambda$: all components of $W_j(\lambda)$, $j=0,2$ are vertices and the mid-lines of remaining strip regions are edges. We denote this tree $\tilde\Gamma$. There is a natural map $\pi\colon\Delta_t^\lambda\to\tilde\Gamma$. Let $\phi_\lambda\colon\tilde\Gamma\to M$ be a parametrization of a flow tree $\Gamma$. We say that the sequence $u_\lambda\colon\Delta_t^\lambda\to T^\ast M$ converges to the flow tree $\Gamma$ if the $C^0$-distance between the maps $u_\lambda$ and $\pi^\ast\phi_\lambda$ goes to $0$ for some parametrization $\phi_\lambda\colon\tilde\Gamma\to\Gamma$ (see also Remark \ref{r:domaincontrolI}).

\begin{pf}[Proof of Theorem \ref{t:disktotree} in the case of simple front singularities]
Let $u_\lambda\colon \Delta_m\to T^\ast M$ be a sequence of $J$-holomorphic disks with boundary on $\hat L_\lambda$. As described above we construct $\Delta_t^\lambda$ and we note that our finiteness results imply that we may, after passing to a subsequence, assume that $t$ is constant.

The finitely many subsets of $\Delta_t^\lambda$ which correspond to vertices in the tree $\tilde\Gamma$ have diameter bounded by $M\lambda^{-(1-\frac{1}{2(n-2)})}$ for some $M>0$. The $\Ordo(\lambda)$-derivative bound for $|Du_\lambda|$ thus implies that the images of each connected component in such a set converges to a point at a $\Ordo(\lambda^{\frac{1}{2(n-2)}})$-rate. Moreover, on the complement of these vertex regions (this complement is a finite collection of strip regions), the map takes the boundaries of each strip to different sheets of $\hat L_\lambda$. As in the proof of Lemma \ref{l:projft} we study two cases for the remaining strip regions depending on their lengths $L$. We separate the strip regions into short, for which  $L=\Ordo(\lambda^{-1})$, and long, for which $L\lambda\to\infty$ as $\lambda\to 0$. As in the proof of Lemma \ref{l:projft} we note that the long strips must map close to Reeb chords and that their number is finite by the finite length condition. The limit of each such piece is thus a broken flow line. The short regions converge at rate $\Ordo(\lambda)$ to flow lines, as follows from Lemma \ref{l:gradconv}.

Since the dimension of a tree and a disk is determined by the homotopy class of its boundary data the theorem follows from Definition \ref{d:tfdim} and Proposition \ref{p:tfdim=dfdim}.
\end{pf}

\begin{rmk}\label{r:domaincontrolI}
It follows from our analysis that the $\R^{t-3}$-coordinates of the conformal structure on the domains $\Delta_t^\lambda$ of the maps $u_\lambda$ in a sequence converging to a tree goes to infinity at least as fast as $\lambda^{-1}$ as $\lambda\to 0$. More precisely, if the the limiting tree $\Gamma$ is unbroken then it determines the coordinates of the conformal structure on $\Delta_t^\lambda$ up to $\Ordo(\lambda^{-(1-\frac{1}{2(n-2)})})$ in the following sense. The flow time lengths of the edges of $\Gamma$ determines a $c\in\R^{t-3}$ such that the coordinates of the conformal structure on $\Delta^\lambda_t$ equals $\lambda^{-1} c +\Ordo(\lambda^{-(1-\frac{1}{2(n-2)})})$.

Also, if we know a priori (for example for dimension reasons) that the limit tree $\Gamma$ is rigid. Then it follows from Remark \ref{r:rigtan} that the diameter of the regions $W_3(\lambda)$ can be taken as small as $\log(\lambda^{-1})$, so in this case the tree determines the coordinates of the conformal structures on $\Delta_t^\lambda$ to within $\Ordo(\log(\lambda^{-1}))$, i.e. the coordinates in $\R^{t-3}$ of the conformal structure on $\Delta_t^\lambda$ are $\lambda^{-1} c+\Ordo(\log(\lambda^{-1}))$.

If the limit tree $\Gamma=\Gamma_1\cup\dots\cup\Gamma_s$ is broken then the coordinates of the conformal structure on $\Delta_t^\lambda$ are not determined in the same way by the tree. In this case they depend also on rate at which the long edges grow.
\end{rmk}

\begin{pf}[Proof of Theorem \ref{t:disktotree} in the $2$-dimensional case]
The proof is a modification of the proof just given so we only sketch it pointing out the main differences. We have
$$
\Pi(\Sigma)=\Sigma_1\supset\Sigma_2\supset\Sigma^{\rm sw}_2.
$$
Note that all results in Subsections \ref{s:areamonbst} still hold. We modify the results in Subsection \ref{s:subdiv} slightly by including also a neighborhood $U^{\rm sw}(2,\delta_2)$ of the swallow tail points. We add punctures to the disk as in \S\ref{ss:subdivseq} taking also intersections with $\pa(T^\ast U^{\rm sw}(2,c\delta_2))$, $c=1,2,3,4$, into account and obtain the corresponding subdivision $\Delta_r^\lambda= D_0(\lambda)\cup D_1(\lambda)\cup D_2(\lambda)\cup D_2^{\rm sw}(\lambda)$. As in Lemma \ref{l:projft} we prove that $u^{T_1}_\lambda$ lies near a projected flow tree. In $U^{\rm sw}(2,\delta_2)$, we consider the projection $\pi^{\rm sw}\colon T^\ast M\to T^\ast Y$, where $Y$ is a line segment in the direction of the swallow tail. We apply the argument in Lemma \ref{l:projft} to show that $\pi^{\rm sw}\circ u$ converges to a projected flow tree at rate $\Ordo(\lambda)$. However, by definition of $\hat L_\lambda$ all local projected gradient differences vanishes around $Y$ and hence the restriction of the map $\pi^{\rm sw}\circ u_\lambda$ to any one of the components of $D_2^{\rm sw}(\lambda)$ converges to a constant at rate $\Ordo(\lambda^2)$ by the argument in Lemma \ref{l:onesheet}. In particular, if $U_\lambda$ is a component of $D_2^{\rm sw}(\lambda)$ then $u_\lambda(U_\lambda)$ lies in an $\Ordo(\lambda^2)$-neighborhood of a copy of $\C$ perpendicular to $T^\ast Y$ in $T_p^\ast M$ for some $p\in U^{\rm sw}(2,\delta_2)$. Thus the limit of the map can be understood by $1$-dimensional complex analysis or by the $1$-dimensional version of Theorem \ref{t:disktotree}. In either case it follows that the limit is a flow tree in the limiting complex line. In particular, we find vertical line segments subdividing  $D_2^{\rm sw}(\lambda)$ into finitely many components of two kinds: those with image inside an $\frac{1}{2}\delta_2$-ball around the swallow tail and those with image outside a $\frac{1}{4}\delta_2$-ball. Moreover, there are $\log(\lambda^{-1})$-strip neighborhoods of the vertical line segments where the $\Ordo(\lambda)$ derivative estimate holds. With this established we may use the total average linking number and Lemma \ref{l:blowup} exactly as above in the complement of the inside region. We then apply the reasoning in Lemma \ref{l:unnecpun} to decide which punctures to keep and which to throw away. The proof can then be finished in much the same way as above.
\end{pf}

\begin{rmk}\label{r:domaincontrolII}
The main difference between the $2$-dimensional case with swallow tail points and the case of simple front singularities is that in the former case we lose some control of the conformal structure on the source for holomorphic disks with images near the swallow tail points. However, if we know a priori that the limit tree is rigid then it follows that it stays a finite distance away from the swallow tail points and all results from the case of simple front singularities hold.
\end{rmk}

\subsubsection{Many positive punctures}\label{ss:many+}
\begin{figure}[htbp]
\begin{center}
\includegraphics[angle=0, width=8cm]{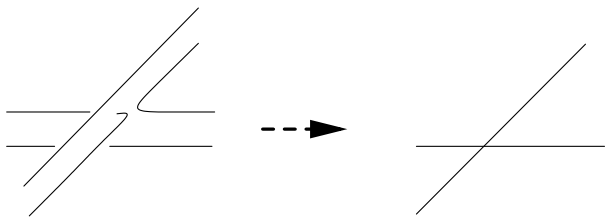}
\end{center}
\caption{A boundary bubble forming.}
\label{f:bdrybub}
\end{figure}

Theorem \ref{t:disktotree} shows that any sequence of rigid holomorphic disks $u_\lambda\colon \Delta_m\to T^\ast M$ with boundary on $\hat L_\lambda$ and only one positive puncture converges to a flow tree of formal dimension $0$. The counterpart of this result for disks with many positive punctures holds as well. The only point in the proof of Theorem \ref{t:disktotree} where the assumption on only one positive puncture was used was in Lemma \ref{l:onesheet}. In the case that the disk has many positive punctures the corresponding result is not true. This has the consequence that in order to prove flow tree convergence one must also deal with strip regions with boundaries which map to the same sheet of $\hat L_\lambda$. The existence of such regions are connected with  "boundary bubbling". We depict the situation in Figure \ref{f:bdrybub}. The argument in Lemma \ref{l:onesheet} can be used to show that the image of such a strip region converges to a point. If $W_3(\lambda)$ denotes the union of all such strip regions then flow tree convergence on $\Delta_t^\lambda-(\cup_{j=0}^3 W_j(\lambda))$ follows from the same arguments as above. Thus also disks with many positive punctures converge to flow trees. (As in the case of broken flow trees the appearance of "boundary bubbling" also affects the control of the coordinates of the conformal structures.)

There are also more profound difficulties associated to disks with many positive punctures: when there are more than one positive puncture, one cannot rule out the possibility that the limiting tree is multiply covered. In fact, using the construction in Section \ref{S:treetodisk} we can produce disks near any rigid flow tree and composing any such disk with a multiple cover we get disks which limit to multiply covered trees.

The dimension of multiply covered trees (disks) depend in a crucial way on the dimension of the Legendrian submanifold. Consider the formal dimension of a multiply covered tree. Note that by cutting the tree at one of its multiple vertices and choosing edges there we construct out of the multiply covered tree, trees with fewer multiple vertices, and eventually many simple trees. We consider the formal dimension. Let $\Gamma$ be the multiple tree broken up into $p$ trees $\Gamma_1,\dots,\Gamma_p$. We have
$$
\dim(\Gamma)=\sum_{j=1}^p\dim(\Gamma_j) - (p-1)(n-3).
$$
Repeating this argument we find that the same formula holds true when all the trees $\Gamma_j$ are simple. Our genericity assumptions shows that $\dim(\Gamma_j)\ge 0$ and it follows in particular that if $n\le 2$, that in order for $\dim(\Gamma)$ to equal $0$, $p=1$. That is $\Gamma$ itself must be simple. Note however, that in the higher dimensional case it may well be that $\dim(\Gamma_j)>0$ and still $\dim(\Gamma)=0$. In particular, it follows that in order to study rigid disks with many positive punctures in dimensions higher than $3$ we must study also flow trees of high formal dimension. Moreover, in order to count such disks one should study their obstruction bundles.

The author decided not to pursue the study of obstruction bundles in this context any further since at the time of writing he was not aware of any applications of holomorphic disks with boundaries on a Legendrian submanifold and with more than one positive puncture.

%% file: Sec/6treetodisk.tex
\section{From trees to disks}\label{S:treetodisk}
In Subsection \ref{s:locsol} we discuss local solutions of the $\bar\pa_J$-equation and its linearization near parts of a rigid flow tree. In Subsection \ref{s:almohol} we use the local solutions to associate an almost holomorphic map to a rigid flow tree. In Subsection \ref{s:barpabund} we describe a neighborhood of this map in a suitable space of maps and conformal structures on the domain. In Subsection \ref{s:unifinv} we study the $\bar\pa_J$-equation and its linearization in these function spaces. This leads to a proof Theorem \ref{t:treetodisk}.

\subsection{Local solutions}\label{s:locsol}
In Subsection \ref{s:metmapfltr} we constructed specific local forms of $\hat L_\lambda$ near its rigid flow trees. These specific forms allow us to write down explicit solutions to the $\bar\pa_J$-equation around any vertex of a rigid flow tree, see Remark \ref{r:vertnot}, and around any edge of it, outside an $\Ordo(\lambda)$-neighborhood of the edge points. We will also prove results about the linearized $\bar\pa_J$-operator acting on function spaces which are direct sums of, on the one hand, a Sobolev space on a standard domain with $2$ or $3$ punctures and positive exponential weights at the ends, and, on the other, a finite dimensional space of cut-off solutions to the $\bar\pa$-equation. (It is spaces of this form we use in Subsection \ref{s:almohol}.)

\begin{rmk}
Many of the lemmas presented in this subsection concern Fredholm operators. That the operators studied have the Fredholm property  follows in all cases from \cite{EES1} Proposition 6.8. The Fredholm index is easily computed using the argument in \cite{EES1} Proposition 6.13. We will not repeat this information in the proofs below.
\end{rmk}

\subsubsection{Edges}\label{ss:soledge}
Consider a flow line segment $\sigma$ of length $l>0$ which connects two edge points on a rigid flow tree $\Gamma$. There are two local forms for $\hat L_\lambda$ in a neighborhood of the cotangent lift of $\sigma$, see \eqref{e:norform1} and \eqref{e:norform2}. We will consider a slightly larger piece $\sigma'\supset\sigma$ of the flow line which includes the edge points and small finite parts of the neighboring flow line segments which are attached at the other end of the edge points bounding $\sigma$. By construction of the metric on $M$ it is flat in a neighborhood of $\sigma'$. It will be convenient to use complex coordinates $(z_1,\dots,z_n)$ on $T^\ast M$. Here $z_j=x_j+iy_j$, where $(x_1,\dots,x_n)$ are the local flat coordinates on $M$ and where elements in the fiber $T^\ast_x M$ are represented as $\sum y_j\,dx_j$. Note that the complex structure $J$ induced by the metric corresponds to the standard complex structure $i$ on $\C^n$ in these coordinates.

The local form of the two sheets of $\hat L_\lambda$ along the cotangent lift of $\sigma'$ is
\begin{align*}
df_1&= \lambda c_1(x_1)\,dx_1,\\
df_2&= \lambda c_2(x_1)\,dx_1,
\end{align*}
Here the function $c_1(x)$ has one of the following forms outside an $\Ordo(\lambda)$-neighborhood of the edge points
\begin{align*}
df_1&= \lambda c_1\,dx_1,\\
df_2&= \lambda c_2\,dx_1,
\end{align*}
or
\begin{align*}
df_1&= \lambda c_1\,dx_1,\\
df_2&= \lambda (c_2+kx_1)\,dx_1.
\end{align*}
Consider the scaling by $\lambda^{-1}$ of these two sheets in the $\C^n$-coordinates. Note that the sheets can be continued to a twice as big region in such a way that the functions $c_1(x)$ and $c_2(x)$ become constant outside a compact set and that they converge to constant functions everywhere. Fix a disk $D_p\subset L$ in one of the sheets of $L$ passing through the midpoint of $\sigma$ and orthogonal to its cotangent lift at this point. It follows from the Riemann mapping theorem that there exists an $D_p$-family of holomorphic maps $\tilde s_\lambda\colon\R\times[0,1]\to\C$ satisfying these boundary conditions and mapping $0$ to $D_p$.

Define $s_\lambda\colon [T_1,T_2]\times[0,1]\to T^\ast M$ by
$$
s_\lambda(\tau+it)=\lambda\tilde s_\lambda,
$$
where $s_\lambda$ is chosen so that $s_\lambda(0,0)$ maps to the intersection point of the tree and $D_p$ and where $T_j$, $j=1,2$, is chosen so that the distance to the nearest edge point is $\Ordo(\lambda)$. Then $\bar\pa_J s_\lambda=0$, and the argument in Lemma \ref{l:gradconv} shows that $T_j$ increases at rate $\lambda^{-1}$ as $\lambda\to 0$. We call $s_\lambda$ as defined above {\em a local solution at an edge} and we call the disk $D_p$ a {\em marking disk}.

We next turn our attention to the linearized $\bar\pa$-equation. Along the cotangent lift of an edge of a rigid flow tree, the boundary condition for the linearized $\bar\pa$-equation in the special coordinates of Subsection \ref{s:metmapfltr}, approaches constant $\R^n$ boundary conditions at rate $\Ordo(\lambda)$. We will consider such boundary conditions for the $\bar\pa$-equation on the strip $\R\times[0,1]$.

If $V\subset\R^n$ is a subspace let $V^\perp$ denote its orthogonal complement.
Fix a small $\delta>0$ and fix subspaces $F(-\infty)$ and $F(+\infty)$ of $\R^n$.
Let ${\mathbf e}\colon \R\times[0,1]\to GL(n)$ be a function such that in some neighborhood of $\{\pm\infty\}\times[0,1]$
\begin{align*}
{\mathbf e}(\tau+it)|F(\pm\infty)=e^{-\delta|\tau|}\id,\\
{\mathbf e}(\tau+it)|F(\pm\infty)^\perp=e^{\delta|\tau|}\id.
\end{align*}
Let $\sblv_{2,\delta}\bigl(F(-\infty),F(+\infty)\bigr)$ denote the weighted Sobolev space of functions $u\colon \R\times[0,1]\to \C^n$ such that
$$
u(\tau+0i)\in\R^n,\quad u(\tau+i)\in\R^n,\quad \bar\pa u|\pa(\R\times[0,1])=0, \text{ and } \|{\mathbf e}u\|_2<\infty,
$$
where $\|\bullet\|_k$ denotes the Sobolev $k$-norm. Let $\sblv_{1,\delta}\bigl(F(-\infty),F(+\infty)\bigr)$ be the weighted Sobolev space  of functions $u$ such that
$$
u|\pa(\R\times[0,1])=0 \text{ and } \|{\mathbf e}u\|_1<\infty.
$$
(Note that the boundary condition for $\bar\pa u$ and the boundary condition on $u\in\sblv_{1,\delta}\bigl(F(-\infty),F(\infty)\bigr)$ should be understood in terms of the trace. For future reference we call these conditions {\em auxiliary boundary conditions}.)
\begin{lma}\label{l:Iw}
The operator
$$
\bar\pa\colon\sblv_{2,\delta}\bigl(F(-\infty),F(+\infty)\bigr)
\to\sblv_{1,\delta}\bigl(F(-\infty),F(+\infty)\bigr)
$$
is Fredholm. Its index satisfies
$$
\ind(\bar\pa)=\dim(F(-\infty))+\dim(F(+\infty))-n.
$$
Its kernel has dimension
$$
\dim(\krn(\bar\pa))=\dim\bigl(F(-\infty)\cap F(+\infty)\bigr),
$$
and it is spanned by constant functions with values in $F(-\infty)\cap F(+\infty)$.
\end{lma}

\begin{pf}
After doubling the functions over $\R+i$, elementary Fourier analysis allows us to write down the local $L^2$ solutions to the equation explicitly, see \cite{EES2} Lemma 5.2:
$$
v(\zeta)=\sum_{k\in\Z} c_k e^{i k\pi\zeta},\quad \zeta=\tau+it,\quad c_k\in\R^n.
$$
Hence, for any subspace $A\subset\R^n$, where the weights at both infinities are negative, constant functions with values in this subspace give solutions which lie in $\sblv_{2,\delta}\bigl(F(-\infty), F(\infty)\bigr)$. Moreover, if $\alpha\in(\R^n)^\ast$ is a linear functional such that at least one of $F(-\infty)$ and $F(\infty)$ is not a subspace of $\krn(\alpha)$ and if $v\in\sblv_{2,\delta}\bigl(F(-\infty),F(\infty)\bigr)$ is a solution (which must be a constant vector) then $\alpha\circ v$ equals $0$.
\end{pf}

We will apply Lemma \ref{l:Iw} to prove a result for the direct sum function spaces mentioned above. Consider the following finite dimensional space of functions on $\R\times[0,1]$ with $\R^n$ boundary conditions. Let $F(\pm\infty)\subset\R^n$ be subspaces. Fix functions $\beta_\pm\colon \R\times[0,1]\to\C$ such that for some fixed $M>1$, $\beta_+=1$ on $[M,\infty)\times[0,1]$, $\beta_+=0$ on $(-\infty,M-1]\times[0,1]$, and such that  $\beta_+$ is real valued and holomorphic on $\pa(\R\times[0,1])$. Define $F_{\rm sol}(\infty)$ to be the linear space spanned by the functions
$\beta_+ c$, $c\in F(\infty)$. We define $\beta_-$ similarly with support in a neighborhood of $-\infty$ and take $F_{\rm sol}(-\infty)$ as the space spanned by $\beta_-c$, $c\in F(-\infty)$. Write $F_{\rm sol}=F_{\rm sol}(-\infty)\oplus F_{\rm sol}(\infty)$ and endow this space with the supremum norm.

We  consider the $\bar\pa$-operator
$$
\bar\pa\colon H_{2,\delta}\oplus F_{\rm sol}\to H_{1,\delta},
$$
where $H_{k,\delta}=\sblv_{k,\delta}(\{0\},\{0\})$ in the notation of Lemma \ref{l:Iw}. Note that any function $v\in H_{2,\delta}\oplus F_{\rm sol}$ has well defined limits $v(\pm\infty)=\lim_{z\to\pm\infty}v(z)$. Let $\ev_{\pm\infty}\colon H_{2,\delta}\oplus F_{\rm sol}\to\R^n$ be the map $\ev_{\pm\infty}(v)=v(\pm\infty)$.

If $A\subset(\R^n)^\ast$ then we let the {\em vanishing condition $A$ at $\infty$} be the following linear constraint for functions $v\in H_{2,\delta}\oplus F_{\rm sol}$:
$$
\alpha\bigl(\ev_\infty(v)\bigr)=0,\quad\text{for all }\alpha\in A.
$$
Note that these constraints specify a closed subspace $H_A$ of $H_{2,\delta}\oplus F_{\rm sol}$. (A vanishing condition  at $-\infty$ is defined in an analogous way.)

Also, if $Q\subset \R^n$ is a subspace, we write $A\lr Q\subset Q\subset\R^n$ for the subspace of $v\in Q$ such that
$$
\alpha(v)=0,\quad\text{ all }\alpha\in A.
$$

\begin{lma}\label{l:Is}
The operator
$$
\bar\pa\colon H_{2,\delta}\oplus F_{\rm sol}\to H_{1,\delta}
$$
is Fredholm of index
$$
\ind(\bar\pa)=\dim(F_{\rm sol})-n.
$$
Its kernel has dimension
$$
\dim(\krn(\bar\pa))=\dim\bigl(F(-\infty))\cap F(\infty)\bigr).
$$
The restrictions $\ev_{\pm\infty}|\krn(\bar\pa)$ are isomorphisms onto $F(-\infty)\cap F(\infty)$. Moreover, if $H_A\subset H_{2,\delta}\oplus F_{\rm sol}$ is the subspace determined by a vanishing condition $A$ at $\infty$ or at $-\infty$ such that
\begin{equation*}
A\lr (F(-\infty)\cap F(+\infty))=\{0\},
\end{equation*}
then there exists a constant $C_A$ such
that
$$
\|w\|\le C_A\|\bar\pa w\|_{H_{1,\delta}},\quad w\in H_A.
$$
\end{lma}

\begin{pf}
The operator under study is obtained from a Fredholm operator by adding a finite dimensional space to the source. Hence it is Fredholm and the index is as claimed. Moreover, the statement about the kernel is an easy consequence of the formula in the proof of Lemma \ref{l:Iw} for solutions. We turn to the last statement. The boundary conditions are split so we may consider one direction at the time. Let $\alpha\in(\R^n)^\ast$ and consider $\alpha(v)$, where $v\in H_A$. The assumption on the vanishing condition ensures that $\alpha(v)$ lies in a weighted Sobolev space with a small positive exponential weight at one end and a small positive or negative exponential weight at the other. In case the weights are positive at both ends the estimate is immediate. Thus the lemma follows once we prove it in the $1$-dimensional case with a vanishing condition at only one end. To show this we argue by contradiction.

Assume that the estimate does not hold. Then there exists a sequence of functions $v_j$ with
\begin{equation}\label{e:contr}
\|\bar\pa v_j\|=1,\quad \|\bar\pa v_j\|_{H_{1,\delta}}\to 0.
\end{equation}
The un-decorated norm $\|\bullet\|$ refers to the norm in the direct sum $H_{2,\delta}\oplus F_{\rm sol}$. Write $v_j=\hat v_j + \eta_j$, where $\eta_j$ is  the component along the added cut-off solution. Then by definition of the norm in the direct sum
$$
\|v_j\|=\|\hat v_j\|_{H_{2,\delta}}+\|\eta_j\|_{F_{\rm sol}}.
$$
Let $\|\cdot\|_{k,\delta,-\delta}$, $k=1,2$ denote the weighted Sobolev $k$-norm with a negative weight $e^{-\delta|\tau|}$ on the end where the cut-off solution is supported and a positive weight $e^{\delta|\tau|}$ on the other end. For definiteness we assume the weight is negative at positive infinity. Equation \eqref{e:contr} implies
$$
\|\bar\pa(\hat v_j+\eta_j)\|_{1,\delta,-\delta}\to 0.
$$
We conclude from the elliptic estimate in this weighted norm, see Proposition 5.5 \cite{EES2}, that
\begin{equation}\label{e:mixest}
\|\hat v_j+\eta_j\|_{2,\delta,-\delta}\to 0.
\end{equation}
We show that $\eta_j\to 0$. Assume not then there exists an infinite sequence of $\hat v_j+\eta_j$ such that $|\eta_j|>\epsilon>0$ for all $j$. Fix $T_0>0$ such that
$$
\frac{\epsilon^2}{4\delta}(1-e^{-2\delta})e^{2\delta T_0} >100.
$$
We conclude from \eqref{e:mixest} that
$$
\|\hat v_j+\eta_j|[T_0,T_0+1]\times[0,1]\|_{2,\delta,-\delta}\to 0.
$$
The triangle inequality then implies that for all sufficiently large $j$
$$
\|\hat v_j|[T_0,T_0+1]\times[0,1]\|^2_{2,\delta,-\delta}>
\frac12\|\eta_j|[T_0,T_0+1]\times[0,1]\|^2_{2,\delta,-\delta}\ge
\frac{\epsilon^2}{4\delta}(1-e^{-2\delta})e^{-2\delta T_0}.
$$
A simple calculation using the fact that the minimum of the weight function of the norm $\|\bullet\|_{H_{2,\delta}}$ on $[T_0,T_0+1]\times[0,1]$ equals $e^{\delta T_0}$ then shows that
$$
\|\hat v_j\|_{H_{2,\delta}}^2\ge
\frac{\epsilon^2}{4\delta}(1-e^{-2\delta})e^{2\delta T_0}>100.
$$
This contradicts $\|v_j\|=1$. It follows that
\begin{equation}\label{e:contr1}
\eta_j\to 0.
\end{equation}
Thus
$$
\|\bar\pa\hat v_j\|_{H_{1,\delta}}\to 0
$$
and we conclude
\begin{equation}\label{e:contr2}
\|\hat v_j\|_{H_{2,\delta}}\to 0
\end{equation}
from the usual elliptic estimate. However, \eqref{e:contr1} and \eqref{e:contr2} contradicts $\|v_j\|=1$. The estimate follows.
\end{pf}

\subsubsection{$1$-valent punctures}\label{ss:sol1punct}
Consider a $1$-valent puncture. We define local solutions along flow line segments near $1$-valent punctures. The local normal form of $\hat L_\lambda$ near a puncture is, see \eqref{e:norform3},
\begin{align*}
df_1&=\sum_j \lambda c_j\,dx_j,\\
df_2&=\sum_j \lambda (c_j+\sigma_j)\,dx_j,
\end{align*}
where $\sigma_j\ne 0$. If the flow line of the rigid tree $\Gamma$ is tangent to the $x_1$-direction then we define $s_\lambda\colon [0,\infty)\times[0,1]\to T^\ast M$ (in local coordinates as above) as
$$
z_j(\tau+it)=\begin{cases}
\lambda i c_1 + c_\lambda e^{i\theta_\lambda}\cdot e^{-\theta_\lambda(\tau+it)}, &\text{ if }j=1,\\
c_j &\text{ otherwise,}
\end{cases}
$$
if $\sigma_1>0$, and
$$
z_j(\tau+it)=\begin{cases}
\lambda i c_1 + c_\lambda e^{-\theta_\lambda(\tau+it)}, &\text{ if }j=1,\\
c_j &\text{ otherwise,}
\end{cases}
$$
if $\sigma_1<0$. Here $\theta_\lambda$ is the complex angle between $\lambda df_1$ and $\lambda df_2$ and $c_\lambda$ is chosen so that the distance from $s_\lambda(0+it)$ to the edge point closest to the puncture is $\Ordo(\lambda)$. Note that $\bar\pa_J s_\lambda=0$ and that $s_\lambda$ takes the boundary to $\hat L_\lambda$. We call $s_\lambda$ an {\em local solution at a $1$-valent puncture}.

The results needed for the linearized $\bar\pa$-equation near
$1$-valent punctures follow from Lemma \ref{l:Is}.

\subsubsection{Ends}\label{ss:solend}
By construction of $\hat L_\lambda$, in a neighborhood of an end of a rigid flow tree, there are coordinates $\C\times\C^{n-1}$ on $T^\ast M$ with the standard complex structure corresponding to $J$, in which $\hat L_\lambda$ is a product of an affine Lagrangian subspace $\Lambda_\lambda\subset\C^{n-1}$ and a curve $\gamma_\lambda$ consisting of a half-circle of radius $\Ordo(\lambda)$ and the boundary of a strip of finite length and of width
$\Ordo(\lambda)$ attached to it ($\gamma_\lambda$ is smoothened at the junctions), see \S\ref{ss:LegisoCE}. We construct a holomorphic map $\hat s_\lambda\colon [0,\infty)\times[0,1]\to T^\ast M$ with boundary in the arc in the first coordinate plane as follows. The map
$$
z\mapsto r_\lambda\frac{2i e^{-\pi z}}{1-ie^{-\pi z}}
$$
maps the half-strip $[0,\infty)\times[0,1]$ to the half of the disk with radius
$r_\lambda$ and center $-r_\lambda+0i$, which lies to the right of the vertical line through its center. Also, it maps $\infty$ to the origin. We find a holomorphic map from $[0,a_\lambda]\times[0,1]$ into the finite strip region between the lines continuing the half circle, starting in an $\Ordo(\lambda)$-neighborhood of the nearest edge point and extended slightly beyond the junction points of the half-circle and the lines. Using these two maps and their inverses we construct, by patching the pieces via transition functions, a holomorphic map $\hat s_\lambda$ on the domain $[0,\infty)\times[0,1]$ with the required boundary conditions. (Standard quasi-conformal estimates shows that
the order of magnitude of $a_\lambda$ satisfies $a_\lambda=\Ordo(\lambda^{-1})$.)
We take $s_\lambda=(\hat s_\lambda, c)$, where $c\in\Lambda_\lambda$ is the coordinate of the cotangent lift of the end. We call $s_\lambda$ a {\em local solution at an end}.

Consider next the linearized boundary condition along a local solution as above. More precisely, consider the strip $\R\times[0,1]$ with boundary conditions which rotates a $1$-dimensional subspace $W\subset\R^n$, inside $W\otimes\C$, an angle $\frac{\pi}{2}$ along $[0,\infty)\times\{0\}$ in the positive direction and $\frac{\pi}{2}$ in the negative direction along $[0,\infty)\times\{1\}$, and which keeps $W^\perp$ fixed. More specifically, for our applications we consider the following split Lagrangian boundary condition in $(\C\otimes W)\times(\C\otimes W^\perp)$. The first component equals the tangent lines of the map $[0,\infty)\times[0,1]\to\C\otimes W\approx\C$, $\zeta\mapsto\frac{2i e^{-\pi \zeta}}{1-ie^{-\pi \zeta}}$, $\zeta\in[0,\infty)\times[0,1]$, and continues constantly equal to $W$ along the boundary of $(-\infty,0]\times[0,1]$. The boundary condition in the second component are constantly equal to $W^\perp\subset\C\otimes W^\perp\approx \R^{n-1}\subset\C^{n-1}$.

Let $F(-\infty)\subset\R^n$ be a subspace and let
${\mathbf e}\colon \R\times[0,1]\to GL(n)$ be a function such that in a
neighborhood of $-\infty$
\begin{align*}
&{\mathbf e}(\tau+it)|F(-\infty)^\perp=e^{\delta|\tau|}\id,\\
&{\mathbf e}(\tau+it)|F(-\infty)=e^{-\delta|\tau|}\id,
\end{align*}
and such that in a neighborhood of $\infty$
\begin{align*}
&{\mathbf e}(\tau+it)|W=e^{\delta|\tau|}\id,\\
&{\mathbf e}(\tau+it)|W^\perp=e^{-\delta|\tau|}\id.
\end{align*}

Define $\sblv_{k,\delta}\bigl(F(-\infty),W\bigr)$, $k=1,2$ as the Sobolev space of functions $u$, weighted by ${\mathbf e}$, which satisfy the boundary conditions (for $k=2$) and which satisfies auxiliary boundary conditions analogous to those in \S\ref{ss:soledge}.

\begin{lma}\label{l:Iw+}
The operator
$$
\bar\pa\colon \sblv_{2,\delta}\bigl(F(-\infty),W\bigr)\to
\sblv_{1,\delta}\bigl(F(-\infty),W\bigr)
$$
is Fredholm and
$$
\ind(\bar\pa)=\dim(F(-\infty))=\dim(\krn(\bar\pa)).
$$
\end{lma}

\begin{pf}
As in \cite{EES1} Proposition 5.13, since the boundary conditions are parallel at the puncture, we may compactify the problem to a problem on the unit disk $D$ in the complex plane with punctures at $-1$ and $1$. The change of coordinates from the non-compact version of the disk to the compact one near punctures can (up to bounded bi-holomorphic map) be taken as  $z\mapsto e^{-\pi z}$ with inverse $w\mapsto -\frac{1}{\pi}\log(w)$. Thus, the pull-backs of holomorphic functions which do not vanish at $1$ live in a Sobolev space with small negative exponential weight but not in a Sobolev space with small positive weight. However, the pull backs of solutions which do vanish at $1$ lie in the Sobolev spaces with small positive weight.

Consider first the problem where $F(-\infty)=\R^n$ and $F(\infty)=W^\perp$. Then the kernel (after small perturbation of the boundary condition in the first coordinate) is $n$-dimensional, spanned by constant functions in $\R^{n-1}$ and real linear combinations of two linearly independent solutions to the Maslov index $1$-problem in $W\otimes\C$ which vanish at $z=1$. (For the standard uniform $\pi$-rotation the two solutions are given by  $z\mapsto i(z-1)$ and $z\mapsto z+1$.)  We next change $F(-\infty)$ by adding directions to $F(-\infty)^\perp$. Each new direction in $F(\infty)^\perp$ should be regarded as a vanishing condition at $-1$ in the compactified problem. The lemma follows.
\end{pf}

As in \S\ref{ss:soledge}, we give corresponding results for positively weighted spaces augmented by cut-off solutions. We use notation as there. Consider boundary conditions as above and let $H_{k,\delta}=\sblv_{k,\delta}(\{0\},\{0\})$, $k=1,2$. Let $F_{\rm sol}(-\infty)$ be a space of cut-off constant solutions supported in a neighborhood of $-\infty$ taking values in $F(-\infty)\subset\R^n$ and let $F_{\rm sol}(\infty)$ to be the corresponding space of cut off constant solutions with support in a neighborhood of $\infty$ and values in $W^\perp$. Let $F_{\rm sol}=F_{\rm sol}(-\infty)\oplus F_{\rm sol}(\infty)$.
\begin{lma}\label{l:Is+}
The operator
$$
\bar\pa\colon H_{2,\delta}\oplus F_{\rm sol}\to H_{1,\delta}
$$
is Fredholm of index
$$
\ind(\bar\pa)=\dim(F_{\rm sol}(-\infty)).
$$
Its kernel has dimension
$$
\dim(\krn(\bar\pa))=\dim\bigl(F(-\infty)\bigr).
$$
The restriction $\ev_{-\infty}|\krn(\bar\pa)$ is an isomorphism onto $F(-\infty)$. Moreover, if $H_A$ is the subspace determined by
a vanishing condition $A$ at $-\infty$ such that
$$
A\lr F(-\infty)=\{0\}
$$
then there exists a constant $C_A$ such that
$$
\|w\|\le C_A\|\bar\pa w\|_{H_{1,\delta}},\quad w\in H_A.
$$
\end{lma}

\begin{pf}
The proof is analogous to the proof of Lemma \ref{l:Is}. (Here Lemma \ref{l:Iw+} plays the role of Lemma \ref{l:Iw} in the proof of Lemma \ref{l:Is}.)
\end{pf}

\subsubsection{Switches}\label{ss:solswitch}
Consider the $\lambda^{-1}$ re-scaling of a finite neighborhood of a switch vertex and continue the re-scaling of $\hat L_\lambda$ constantly to infinity. Note that by the choice of function $\hat h(q)$, see \S\ref{ss:metmapswitch}, the function $\tilde s\colon\R\times[0,1]\to\C^n=\C^{n-1}\times\C$,
$$
\tilde s(z)=((a_1-a_0)z,s(z)),
$$
is holomorphic and has its boundary on the continuation of the scaling of $\hat L_\lambda$. Let $s_\lambda\colon [T_1,T_2]\times[0,1]\to T^\ast M$ equal $\lambda\tilde s$, where $T_j$ is chosen so that the distance of $s_\lambda(T_j+it)$ to the nearest edge points is $\Ordo(\lambda\log(\lambda^{-1}))$. (The existence of such $T_j$ follows from the exponential convergence of $s$ to standard solutions near infinity and from the argument in the proof of Theorem \ref{t:disktotree} which establishes flow line convergence at rate $\Ordo(\lambda\log(\lambda^{-1}))$.)

We also consider the uniqueness of the function $\tilde s$ (coordinates and notation for these near the switch are as in \S\ref{ss:metmapswitch}). Note that there is an obvious $(n-1)$-parameter family of such functions arising from translations along $X$ and from pre-composition with translations of $\R\times[0,1]$. We claim that if $v$ is any holomorphic map with boundary mapping to the continuation of the re-scaling of $\hat L_\lambda$ which is asymptotic to a standard map at infinity then $v$ is an element in this family.

First it is elementary to see that the $\C^{n-1}$-components $s^T$ and $u^T$ of $s$ and $u$, respectively, only differs by translation along $X$ after precomposing with the right translation. Therefore the functions $u^\perp$ and $s^\perp$, taking values in the copy of $\C$ perpendicular to $\C^{n-1}$, have the same $y$-coordinate along one boundary component, since it is determined by $\hat h\circ s^T=\hat h\circ u^T$. Thus $u^\perp-s^\perp$ maps one boundary component to $\R$ and has bounded imaginary part and bounded derivative. Consider the double $f$ of $u^\perp-s^\perp$. Note that the real part $r$ of $f$ is periodic and can be expanded in a Fourier series. Form the fact that it is harmonic and that its mean value is $0$ we find
$$
r(z)=\sum_{n\ne 0} c_n e^{n\pi\tau}\cos(n\pi t).
$$
It is then easy to see that the $C^0$-bound implies that all $c_n=0$. Therefore the solution is unique. We conclude that the (local) moduli space of solutions is isomorphic to $X$ (i.e. the tangency locus: any solution admits an obvious translation by a vector in $X$.)

Finally, we pick two marking disks, $D_p$ and $D_q$ in $\hat L_\lambda$, orthogonal to the cotangent lift of the flow line of the rigid tree, one on each side of the switch vertex at finite distance from it, and we consider the preimages of these as marked points in $[T_1,T_2]\times[0,1]$.

We next turn to the linearized problem. Consider the strip $\R\times[0,1]$ with boundary conditions which rotates a $1$-dimensional subspace $W\subset\R^n$, in $W\otimes\C$, an angle $-\pi$ along $[-M,M]\times\{0\}$ in the negative direction and except that is constantly equal to $\R^n$. Let $F(\infty)\subset\R^n$ be a subspace. Let ${\mathbf e}\colon \R\times[0,1]\to GL(n)$ be a function such that in a
neighborhood of $\infty$
\begin{align*}
{\mathbf e}(\tau+it)|F(\infty)^\perp=e^{\delta|\tau|}\id,\\
{\mathbf e}(\tau+it)|F(\infty)=e^{-\delta|\tau|}\id,
\end{align*}
and such that in a neighborhood of $-\infty$,
$$
{\mathbf e}(\tau+it)=e^{-\delta|\tau|}\id.
$$
Define $\sblv_{k,\delta}\bigl(F(\infty),W\bigr)$, $k=1,2$ as the Sobolev space of functions $v$ which satisfy the boundary conditions (for $k=2$) and which satisfy auxiliary boundary conditions analogous to the auxiliary boundary conditions in \S\ref{ss:soledge}.

\begin{lma}\label{l:Iw-}
The operator
$$
\bar\pa\colon \sblv_{2,\delta}\bigl(F(\infty),W\bigr)\to
\sblv_{1,\delta}\bigl(F(\infty),W\bigr)
$$
is Fredholm. Its index satisfies
$$
\ind(\bar\pa)=\dim(F(\infty))-1.
$$
Its kernel has dimension
$$
\dim(\krn(\bar\pa))=\dim\bigl(F(\infty)\cap W^\perp\bigr),
$$
and is spanned by the constant functions with values in $F(\infty)\cap W^\perp$.
\end{lma}

\begin{pf}
As in the proof of Lemma \ref{l:Iw+} we compactify the problem to a problem on $D$ with punctures at $-1$ and $1$. The $1$-dimensional problem corresponding to the $W\otimes\C$-component has Maslov index $-1$. If there is no vanishing condition at $-1$ (or $1$) then the $\bar\pa$-operator with such a boundary condition is an isomorphism. If $F(\infty)\subset W^\perp$ then there is a vanishing condition at $-1$ and the $\bar\pa$-operator has $1$-dimensional cokernel. The lemma follows from these observations and Lemma \ref{l:Iw} applied to the orthogonal complement of $W$.
\end{pf}

Consider boundary condition as above and let $H_{k,\delta}=\sblv_{k,\delta}(\{0\},\{0\})$. Let $F_{\rm sol}(\pm\infty)$ be a space of cut-off constant solutions supported in a neighborhood of $\pm\infty$ with values in $F(\infty)\subset\R^n$ and $F(-\infty)=\R^n$, respectively. Write $F_{\rm sol}=F_{\rm sol}(-\infty)\oplus F_{\rm sol}(\infty)$.

We specify vanishing conditions as in \ref{ss:soledge} by $A$ at $\infty$ or at $-\infty$, where $A\subset(\R^n)^\ast$.

\begin{lma}\label{l:Is-}
The operator
$$
\bar\pa\colon H_{2,\delta}\oplus F_{\rm sol}\to H_{1,\delta}
$$
is Fredholm of index
$$
\ind(\bar\pa)=\dim(F_{\rm sol})-n-1.
$$
Its kernel has dimension
$$
\dim(\krn(\bar\pa))=\dim\bigl(F(\infty)\cap W^\perp\bigr)
$$
The restriction $\ev_{\pm\infty}|\krn(\bar\pa)$ is an isomorphism onto $F(\infty)\cap W^\perp$. Moreover, let $H_A$ is the subspace determined by a vanishing condition $A$ at $\pm\infty$ such that
$$
A\lr (F(\infty)\cap W^\perp)=\{0\}.
$$
then there exists a constant $C_A$ such that
$$
\|w\|\le C_A\|\bar\pa w\|_{H_{1,\delta}},\quad w\in H_A.
$$
\end{lma}

\begin{pf}
The proof is analogous to the proof of Lemma \ref{l:Is}. (Here Lemma \ref{l:Iw-} plays the role of Lemma \ref{l:Iw} in the proof of Lemma \ref{l:Is}.)
\end{pf}

\subsubsection{$Y_0$-vertices}\label{ss:solY0}
Consider a $Y_0$-vertex of a rigid flow tree. By construction of $\hat L_\lambda$ the three sheets involved in the $Y_0$-vertex are locally (covariantly) constant. In flat coordinates as above on $M$
$$
df_j = \lambda \alpha_j ,\quad j=1,2,3,
$$
where $\alpha_j\in(\R^n)^\ast$ is a covector. Using the flat metric we translate between covectors and vectors. Let $v_j\in\R^n$ be the vector corresponding to $\alpha_j$. Assume that coordinates are chosen so that the $Y_0$-vertex lies at the origin, where flow lines of $v_1-v_3$ and $v_3-v_2$ meet and bifurcate into a flow line of $v_1-v_2$.

Consider a $3$-punctured disk $D_3$ as a strip $\R\times[0,2]$ with a slit around the half line $[0,\infty)\times\{1\}$ removed. Denote its punctures  $p_1$, $p_2$, $p_3$. Map $D_3$ biholomorpically to a strip $\R\times[0,1]$ with maps $U_2$ and $U_3$, where $U_2(p_1)=-\infty$, $U_2(p_3)=i$, $U_2(p_2)=\infty$, and where $U_3(p_1)=\infty$, $U_3(p_2)=0$, and $U_3(p_3)=\infty$. Consider the maps
$a_j\colon\R\times[0,1]\to \C^n$, $a_2(z)=(v_3-v_2)z + iv_3$ and $a_3=(v_1-v_3)z$. Define
$$
\tilde s_\lambda = \lambda( a_2\circ U_2 + a_3\circ U_3).
$$
We take $s_\lambda$ to be the restriction of $\tilde s_\lambda$ to $D_3'\subset D_3$ such that the image of each end of $D_3'$ under $s_\lambda$ lies at distance $\Ordo(\lambda\log(\lambda^{-1}))$ from an $\epsilon$-sphere around $0$, for some fixed small $\epsilon>0$ (smaller than the distance to the nearest edge point). Note that the length of the ends of $D_3'$ then scale as $\lambda^{-1}$, that $\bar\pa_J s_\lambda=0$, and that $s_\lambda$ maps the boundary of $D_3$ to $\hat L_\lambda$. To simplify notation below we will consider the three points of intersection of  $\Gamma$ and the $\epsilon$-sphere discussed above as edge points of $\Gamma$.

We next discuss the uniqueness of the map $\tilde s_1$. Consider a holomorphic map $g\colon D_3\to\C^n$ which is asymptotic to $\tilde s_1$ at infinity in the sense that the $C^0$-distance between the maps goes to infinity as we go out the ends. Then the difference $\tilde s_1-g$ is a holomorphic map $D_3\to\C^n$ taking the boundary to $\R^n$. Also, it is easy to see from the explicit solutions of the $\bar\pa$-equation with $\R^n$-boundary conditions that the $C^0$-convergence actually implies exponential convergence near punctures. In particular, the map $s_1-g$ gives a holomorphic map of the closed unit disk which takes the boundary to $\R^n$. Such a map must be constant by the maximum principle. We conclude that $s_1$ is unique up to translation by a vector in $\R^n$.

Much as above, we pick three marking disks $D_{p_j}$, $j=1,2,3$ inside the $\epsilon$-sphere in $\hat L_\lambda$, orthogonal to the cotangent lift of the tree, one for each edge adjacent to the vertex at finite distance from the vertex, and we mark the points in the boundary of $D'_3$ mapping to these disks.

We turn our attention to the linearized $\bar\pa$-problem near the solution discussed above. Thus we consider a $3$-punctured disk with constant $\R^n$ boundary conditions. We think of the $3$-punctured disk as a standard domain and denote its punctures $p_j$, $j=1,2,3$.  Choose subspaces $F(p_j)\subset\R^n$, $j=2,3$, and let $F(p_1)=\R^n$. Let ${\mathbf e}\colon \Delta_3\to GL(n)$ be a function such that in a neighborhood of $p_j\approx \{\infty\}\times[0,1]$
\begin{align*}
&{\mathbf e}(\tau+it)|F(p_j)^\perp =e^{\delta|\tau|}\id,\\
&{\mathbf e}(\tau+it)|F(p_j)=e^{-\delta|\tau|}\id.
\end{align*}
Define the Sobolev spaces $\sblv_{k,\delta}\bigl(F(p_2),F(p_3)\bigr)$, $k=1,2$ as the Sobolev spaces, weighted by ${\mathbf e}$, of functions which satisfy boundary conditions, and auxiliary boundary conditions analogous to the auxiliary boundary conditions in \S\ref{ss:soledge}.
\begin{lma}\label{l:IIIw}
The operator
$$
\bar\pa\colon \sblv_{2,\delta}(F(p_2),F(p_3))
\to\sblv_{1,\delta}(F(p_2),F(p_3))
$$
is Fredholm. Its index satisfies
$$
\ind(\bar\pa)=\dim(F(p_2))+\dim(F(p_3))-n.
$$
Its kernel has dimension
$$
\dim(\krn(\bar\pa))=\dim(F(p_2)\cap F(p_3)),
$$
and is spanned by constant functions with values in $F(p_2)\cap F(p_3)$.
\end{lma}

\begin{pf}
As in the proof of Lemma \ref{l:Iw+} we compactify the problem and consider the corresponding problem on the unit disk $D$ with three punctures. If $F(p_j)=\R^n$, $j=1,2,3$ then the compactified problem has $n$-dimensional kernel spanned by constant functions. Each time a direction is added to some $F(p_j)^\perp$ the kernel dimension is reduced by one (unless the kernel is already trivial). This finishes the proof.
\end{pf}

We consider similar problems for positively weighted spaces augmented with cut-off solutions. There will be one further complication in the present set-up which arises from variations of conformal structures. We first deal with the cut-off solutions:
let $F_{\rm sol}(p_j)$ be spaces of cut-off constant solutions  with values in subspaces $F(p_j)\subset\R^n$, $j=1,2,3$, we assume that $\dim(F(p_2))+\dim(F(p_3))\le 2n-1$, see Remark \ref{r:dimrigidsubtree}, and let $F_{\rm sol}=\oplus_j F_{\rm sol}(p_j)$. Let $H_{k,\delta}=\sblv_{k,\delta}(\{0\},\{0\})$.

We consider also input from conformal variations. Let $v_2$ and $v_3$ be vectors in $F(p_2)$ and $F(p_3)$ respectively and consider the vector $v_2+v_3$ in $F(p_1)$ as above. (Note that if $f$ is the vector field $\pa_\tau$ along $\Delta_3$, used to vary the conformal structure on a larger domain containing a part of $\Delta_3$ as in \S\ref{ss:confmod}, and if $s_1$ is the standard solution with asymptotics $v_2$, $v_3$, and $v_2+v_3$ then $\pa s_1\cdot f$ is asymptotic to these vectors at $p_1$, $p_2$, and $p_3$, respectively.) Let $F_{\rm con}$ be the $1$-dimensional space spanned by $v_2(p_2)+ v_3(p_3)+ (v_2+v_3)(p_1)$, where the $w(p_j)$ denotes the cut-off solution at $p_j$ with value $w$.
We specify a vanishing condition at a puncture $p_j$ as in \S\ref{ss:soledge} by a subspace $A\subset(\R^n)^\ast$. We write $H_A$ for the closed subspace satisfying the vanishing conditions, we write $\R\la w\ra$ for the subspace generated by the vector $w$, and we write $A+B$ for the smallest subspace containing the subspaces $A$ and $B$.

\begin{lma}\label{l:IIIs}
The operator
$$
\bar\pa\colon H_{2,\delta}\oplus F_{\rm sol}\oplus F_{\rm con}\to H_{1,\delta}
$$
is Fredholm of index
$$
\ind(\bar\pa)=\dim(F(p_2))+\dim(F(p_3))-n+1.
$$
Its kernel has dimension
$$
\dim(\krn(\bar\pa))=\dim\bigl(F(p_2)\cap F(p_3))+1
$$
The restriction $\ev_{p_1}|\krn(\bar\pa)$ maps onto $F(p_2)\cap F(p_3) + \R\la v_1+v_2\ra$. Moreover, let $A$ be a vanishing condition at $p_1$ such that
$$
A\lr (F(p_2)\cap F(p_3)+\R\la v_1+v_2\ra)=0.
$$
Then there exists a constant $C_A$ such that
$$
\|w\|\le C_A\|\bar\pa w\|_{H_{1,\delta}},\quad w\in H_{A}.
$$
\end{lma}

\begin{pf}
The proof is analogous to the proof of Lemma \ref{l:Is}. (Here Lemma \ref{l:IIIw} plays the role of Lemma \ref{l:Iw} in the proof of Lemma \ref{l:Is}. The new conformal variation can be dealt with using finite dimensional linear algebra.)
\end{pf}

\subsubsection{$Y_1$-vertices}\label{ss:solY1}
\begin{figure}[htbp]
\begin{center}
\includegraphics[angle=0, width=8cm]{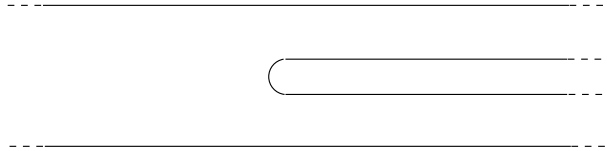}
\end{center}
\caption{The image of the obvious map in the first coordinate.}
\label{f:obv}
\end{figure}

Near a $Y_1$-vertex the three sheets of $\hat L_\lambda$ decompose as a product in local coordinates $\C\times\C^{n-1}$. In the first factor the three sheets look like a bend surrounded by two straight lines, one above and one below. In $\C^{n-1}$ $\hat L_\lambda$ looks just like near a $Y_0$-vertex. In this case we define $s_\lambda$, the {\em local solution at a $Y_1$-vertex}, as the product of the obvious map, see Figure \ref{f:obv}, in the first coordinate and the $Y_0$-solution discussed in \S\ref{ss:solY0} in $\C^{n-1}$.

Again we cut off the domain so that the end-segments lie at distance $\Ordo(\lambda\log(\lambda^{-1}))$ from a small $\epsilon$-sphere around the vertex and we consider the three intersections of $\Gamma$ with this sphere as edge points.

As in \S\ref{ss:solY0} we see that $s_1$ is unique up to translation by a constant vector in $\R^{n-1}$ (i.e. a vector parallel to the image  of the cusp edge under $\Pi$), and as there we also add one marking disk and marked point for each edge adjacent to the $Y_1$-vertex.

We next consider the corresponding linearized problem on a $3$-punctured standard domain $\Delta_3$ with boundary conditions as follows. Let $W\subset\R^n$ be a $1$-dimensional subspace. Along one edge in $\pa D_3$ the boundary condition rotates an angle $-\pi$ in $\C\otimes W$ and the orthogonal complement $W^\perp\subset\R^n$ remains fixed. Along other edges the boundary condition is constantly $\R^n$. Let the punctures be $p_1,p_2,p_3$. Choose subspaces $F(p_j)\subset\R^n$, $j=2,3$, and let $F(p_1)=\R^n$. Let
${\mathbf e}\colon D_3\to GL(n)$ be a function such that in neighborhoods of $p_j\approx\{\pm\infty\}\times[0,1]$
\begin{align*}
&{\mathbf e}(\tau+it)|F(p_j)^\perp=e^{\delta|\tau|}\id,\\
&{\mathbf e}(\tau+it)|F(p_j)=e^{-\delta|\tau|}\id.
\end{align*}
Define the Sobolev spaces $\sblv_{k,\delta}(F(p_2),F(p_3),W)$, $k=1,2$ as the Sobolev spaces, weighted by ${\mathbf e}$, of functions which satisfies boundary conditions as stated (for $k=2$) and auxiliary boundary conditions analogous to the auxiliary boundary conditions in \ref{ss:soledge}.

\begin{lma}\label{l:IIIw-}
The operator
$$
\bar\pa\colon \sblv_{2,\delta}(F(p_2),F(p_3),W)
\to\sblv_{1,\delta}(F(p_2),F(p_3),W)
$$
is Fredholm. Its index satisfies
$$
\ind(\bar\pa)=\dim\bigl(F(p_2)\bigr)+\dim\bigl(F(p_3)\bigr)-n-1.
$$
Its kernel has dimension
\begin{equation*}
\dim(\krn(\bar\pa))=
\dim\bigl(F(p_1)\cap F(p_2)\cap W^\perp\bigr),
\end{equation*}
and is spanned by constant vectors in this subspace.
\end{lma}

\begin{pf}
The boundary condition is split. As in the proof of Lemma \ref{l:Iw-} we note that the $W\otimes\C$ component of any solution must vanish. The lemma then follows from Lemma \ref{l:IIIw} applied to $W^\perp\otimes\C$.
\end{pf}

Let $H_{k,\delta}=H_{k,\delta}( \{0\},\{0\},W)$. Let $F_{\rm sol}(p_j)$ be spaces of cut-off constant solutions taking values in $F(p_j)\subset\R^n$, and let $F_{\rm con}$ be a one dimensional space of conformal variations spanned by $w(p_1)+w(p_2)+w(p_3)+v_2(p_2)+v_3(p_3)+(v_2+v_3)(p_1)$, where $w\in W$ and where $v_2\in F(p_2)\cap W^\perp$ and $v_3\in F(p_3)\cap W^\perp$.

\begin{lma}\label{l:IIIs-}
The operator
$$
\bar\pa\colon H_{2,\delta}\oplus F_{\rm sol}\oplus F_{\rm con}\to H_{1,\delta}
$$
is Fredholm of index
$$
\ind(\bar\pa)=\dim(F(p_2))+\dim(F(p_3))-n.
$$
Its kernel has dimension
$$
\dim(\krn(\bar\pa))=\dim\bigl(F(p_2)\cap F(p_3)\cap W^\perp)+1.
$$
The restriction $\ev_{p_1}|\krn(\bar\pa)$ maps onto $F(p_2)\cap F(p_3)\cap W^\perp + \R\la v_2+v_3+w\ra$. Moreover if $A$ is a vanishing condition at $p_1$ such that
$$
A\lr \bigl((F(p_1)\cap F(p_2)\cap W^\perp)+\R\la v_2+v_3+w\ra\bigr)=0
$$
and if $H_A$ is the subspace determined by $A$ then there exists a constant $C_A$ such that
$$
\|w\|\le C_A\|\bar\pa w\|_{H_{1,\delta}},\quad w\in H_A.
$$
\end{lma}

\begin{pf}
The proof is analogous to the proof of Lemma \ref{l:IIIs}.
\end{pf}

\subsubsection{$2$-valent punctures}\label{ss:sol2punct}
Consider next a $2$-valent puncture. Assume that the nearby flow lines are tangent to the $x_1$-direction then we consider the intersection of the three sheets of $\hat L_\lambda$ with the $(x_1,y_1)$-plane. It consists of three line-segments: $l_0$ and $l_1$ parallel to the $x_1$-axis and $l_2$ intersecting $l_0$ and $l_1$ at a small angle. The three lines containing $l_0$, $l_1$, and $l_2$, respectively, form a triangle with one corner at infinity. (The angles are $\alpha_0\approx \pi$, $\alpha_1=0$, and $\alpha_2\approx 0$ for small $\lambda$.) Consider a $3$-punctured disk $D_3$ with punctures at $1$, $i$, and $-1$ as a strip $\R\times[0,2]$ with a slit around the half line $[0,\infty)\times\{1\}$ removed. We define $\hat s_\lambda\colon D_3\to\C$ mapping into the triangle using the Schwarz-Christoffel formula
$$
z\mapsto A
\int_0^z(\zeta-1)^{-\frac{\alpha_0}{\pi}}
(\zeta+1)^{\frac{\alpha_2}{\pi}}\, d\zeta,
$$
where $A$ is a suitable constant. Here the $2$-valent puncture corresponds to the puncture $1\in D_3$. We define $s_\lambda$, the {\em local solution at a $2$-valent puncture}, by restricting to a subset of the strip with a slit obtained by cutting off neighborhoods of the other two punctures. We cut so that the distance from the image of the cuts to the edge points closest to the puncture is $\Ordo(\lambda)$. Note that $s_\lambda$ is the unique holomorphic map from the $3$-punctured disk which maps $1$ to the point described, which is asymptotic to a standard map at infinity, and which maps the third puncture to the intersection point. This uniqueness holds also in the limit which is a map of the strip with boundary conditions $\R^n$ and a parallel copy thereof mapping to the point which is the limit of the double point of the two intersecting Lagrangian sheets and with the limiting asymptotic conditions at $\pm\infty$.

We next consider the linearized problem. First we specify subspaces $F(p_2)\subset \R^n$ and $F(p_3)\subset\R^n$ with $\dim(F(p_2)\cap F(p_3))\ge 1$, we take $F(p_1)=0$, and let $F_{\rm sol}$ be the space of cut off solutions at $p_j$ with values in $F(p_j)$. We need again to take conformal variations into account. We let $v\in F(p_2)\cap F(p_3)$ be a vector and let $F_{\rm con}$ be the $1$-dimensional space spanned by the function $v(p_2)+v(p_3)$ with notation as above.  Let $F_{\rm sol}$ be the direct sum of cut-off constant solutions at $p_2$ and $p_3$ with arbitrary values in $\R^n$.

\begin{lma}\label{l:IIIs2}
The operator
$$
\bar\pa\colon H_{2,\delta}\oplus F_{\rm sol}\oplus F_{\rm con}\to H_{1,\delta}
$$
is Fredholm of index $1$. Its kernel has dimension $1$ and the restriction $\ev_{p_j}|\krn(\bar\pa)$ maps onto $\R\la v\ra$. Moreover, let $A$ be a vanishing condition at $p_2$ or $p_3$ such that
$$
A\lr \R\la v\ra =0.
$$
Then there exists a constant $C_A$ such that
$$
\|w\|\le C_A\|\bar\pa w\|_{H_{1,\delta}},\quad w\in H_{A}.
$$
\end{lma}

\begin{pf}
The proof is similar to the proof of Lemma \ref{l:Iw}. (Here Lemma \ref{l:IIIw} plays the role of Lemma \ref{l:Iw} in the proof of Lemma \ref{l:Is})
\end{pf}

\subsection{Approximately holomorphic maps}\label{s:almohol}
In this subsection we associate to each rigid tree $\Gamma$ of $\hat L_\lambda$ a standard domain $\Delta_{p,m}^{(0,0)}(\Gamma, \lambda)$, with $p$ punctures and $m$ marked points on the boundary and a map $u_\lambda^{(0,0)}\colon\Delta_{p,m}^{(0,0)}(\Gamma,\lambda)\to T^\ast M$ with boundary on $\hat L_\lambda$ which is holomorphic on most of its domain. The conformal structure on a domain such as $\Delta_{p,m}^{(0,0)}(\lambda)$ is determined by the map to $\R^{p-3}$ described in Lemma \ref{l:confmod} and by the location of the $m$ marked points. We will use a slightly different set of parameters to parameterize the conformal structures near $\Delta_{p,m}^{(0,0)}(\Gamma,\lambda)$: we use $p-2$ variations at the boundary minima (similar to those used in \S\ref{ss:confmod}) and $m-1$ variations which moves the marked points. In this way the conformal structures are parameterized by two parameters $(\kappa,\gamma)\in \R^{p-2}\times\R^{m-1}$. We then extend the family of maps $u_\lambda^{(0,0)}\colon \Delta_{p,m}^{(0,0)}(\Gamma,\lambda)\to T^\ast M$ to a family of maps $u_\lambda^{(\kappa,\gamma)}\colon \Delta_{p,m}^{(\kappa,\gamma)}(\Gamma,\lambda)\to T^\ast M$, where $(\kappa,\gamma)\in\R^{p-2}\times\R^{m-1}$ varies in an $\frac{\epsilon}{\lambda}$ neighborhood of $(0,0)$ in the space of conformal structures for a small $\epsilon>0$. With these maps defined we describe two bundles of function spaces over the space of conformal structures. One of these bundles will be used to parameterize a neighborhood of the maps $u_\lambda^{(\kappa,\gamma)}$ in the space of maps $\Delta^{(\kappa,\gamma)}_{p,m}(\Gamma,\lambda)\to T^\ast M$ and also as the source space of the $\bar\pa_J$-operator. The other bundle will be the target space of $\bar\pa_J$. Finally, we discuss this bundle  $\bar\pa_J$-operator and its linearization.

\subsubsection{The domain of a rigid tree}
Let $\Gamma$ be a rigid flow tree of $\hat L_\lambda\subset T^\ast M$. In order to define $\Delta_{p,m}^{(0,0)}(\Gamma,\lambda)$, we will associate pieces of this domain to various pieces of the tree. The following notation will be used: if $\Delta_3$ is a standard domain with punctures $p_j$, $j=1,2,3$ then $\Delta_3(T_1,T_2, T_3)$ denotes the subset of the domain obtained by cutting $\Delta_3$ along a vertical line segment in the neighborhood of $p_j$ at distance $T_j$ from the boundary minimum of $\Delta_3$. Let $\dot\Delta_2$ denote a strip with marked point at $(0,0)$ and let $\dot\Delta_2(T_1,T_2)$ be the subset of $\dot\Delta_2$ obtained by cutting at vertical segments in neighborhoods of $-\infty$ and $+\infty$ at distances $T_1$ and $T_2$ respectively from the marked point. We will allow the domains $\Delta_3$ and $\dot\Delta_2$ to have width between $1$ and $K$ for some $K$ and the location of the slit in $\Delta_3$ is allowed to vary between $2$ and $K-1$. The exact shapes of the regions will not matter much in our analytical estimates and so will not be discussed in any detail. The pieces are the following.

\begin{itemize}
\item
For each $1$-valent puncture or end, pick a domain of the form
$\dot\Delta_2(0,\infty)$.
\item
For each $2$-valent puncture $q$ pick a domain of the form $\Delta_3(\infty, T_2, T_3)$, where $T_2$ and $T_3$ are chosen so that the local solution at $q$ takes the vertical segments where $\Delta_3$ was cut to a $\Ordo(\lambda\log(\lambda^{-1}))$-neighborhood of the nearest edge point.
\item
For each $3$-valent vertex $v$. Pick a domain $\Delta_3(T_1,T_2,T_3)$ such that the local solution at $v$ maps the vertical segments at the cuts to an $\Ordo(\lambda\log(\lambda^{-1})$-neighborhood of the nearest edge point.
\item
For each switch $v$ pick a domain of the form $\dot\Delta_2(T_1,T_2)$ such that the restriction of the local solution at $v$ takes the vertical boundary line segment to an $\Ordo(\lambda\log(\lambda^{-1}))$-neighborhood of the nearest edge point.
\item For each edge point pick a domain of the form $\dot\Delta_2(a,a)$ so that the flow line passing through this point connects the two neighboring regions where the boundary conditions are standard in flow time $2\frac{a}{w}$, where $w$ is the width of $\dot\Delta_2(a,a)$.
\item
For each finite part of an edge $e$ between two edge points pick a domain $\dot\Delta_2(T_1,T_2)$ such that the restriction of the standard solution takes the vertical end points to an $\Ordo(\lambda)$-neighborhood of the nearest edge point.
\end{itemize}

\begin{rmk}\label{r:flowlength}
All distances labeled $T_j$ above will be called {\em flow length parameters}. Note that there exists $l_0>0$ such that $l_0\lambda^{-1}< T_j < l_0^{-1}\lambda^{-1}$ where $T_j$ is any flow length parameter. On the other hand by definition of $\hat L_\lambda$ near its rigid flow trees, the parameters labeled $a$ above satisfies $a_0<a<a_0^{-1}$ for some $a_0>0$. We call the neighborhoods $\dot\Delta_2(a,a)$ {\em edge point regions}.
\end{rmk}

We next glue the pieces above together to a standard domain. This involves specifying the shapes of the domains. To organize this we first orient the source tree of $\Gamma$ in the following way. Distinguish one of the positive punctures $p$ of $\Gamma$. If $p$ is $1$-valent then orient the edge ending at $p$ away from $p$. If $p$ is $2$-valent then orient the two edges adjacent to $p$ away from it. Orientations of the edges of $\Gamma$ are then uniquely specified by requiring that at each $2$-valent vertex $v\ne p$ one of the adjacent edges is oriented toward $v$ and the other away from $v$, and that at each $3$-valent vertex $w$ one adjacent edge is oriented toward $w$ and the other two away from it.

Consider a tree $\Gamma$ with distinguished positive puncture $p$ and equipped with an orientation as described. Define a function $W$ on its edges as follows. If $e$ is an edge adjacent to an end or a $1$-valent puncture $v\ne p$ then $W(e)=1$. If $v\ne p$ is a $2$-valent puncture with incoming edge $e_1$ and outgoing edge $e_2$ then $W(e_1)=W(e_2)+1$. If $v$ is a switch with adjacent edges $e_1$ and $e_2$ then $W(e_1)=W(e_2)$. If $v$ is a $3$-valent vertex with incoming edge $e_1$ and outgoing edges $e_2$ and $e_3$ then $W(e_1)=W(e_2)+W(e_3)$. This defines $W$ on all edges in $\Gamma$. If $v$ is a vertex of $\Gamma$ which is not a distinguished $2$-valent puncture then define $W(v)$ as the maximum of $W(e)$ over edges $e$ adjacent to $v$. If $p$ is a distinguished $2$-valent puncture then define $W(p)=W(e_1)+W(e_2)$, where $e_1$ and $e_2$ are the edges adjacent to $p$.

We are now in position to glue the pieces associated to pieces of $\Gamma$ together in the obvious way. More precisely, start with the regions $\dot\Delta_2(0,\infty)$ corresponding to punctures and ends with $W=1$ and let them have width slightly smaller than $1$ (the exact width is determined by the size of the slits in our standard domains). These regions correspond to a segment of an edge in the tree connecting a puncture to an edge point. We adjoin the corresponding edge point domains with the same widths. We then glue on regions $\dot\Delta_2(T_1,T_2)$ corresponding to flow lines between edge points, still with the same width. The construction continues until we come to a $3$-valent vertex or a $2$-valent puncture in which case we glue to $\Delta_3(T_1,T_2,T_3)$ the pieces adjacent to it respecting the cyclic order. (The complex orientation on $\Delta_3$ induces a cyclic order on its three end pieces.) The width of $\Delta_3$ should then be slightly smaller than $2$ and the location of its slit is determined by the width of the incoming edges (in this first step it equals $1$). This process can be continued upwards through the tree (against its orientation) and each time we meet a new region corresponding to a $3$-valent vertex or $2$-valent puncture the function $W$ dictates the shape of the region. The result is a standard domain with $p$ punctures.

\begin{figure}[htbp]
\begin{center}
\includegraphics[angle=0, width=8cm]{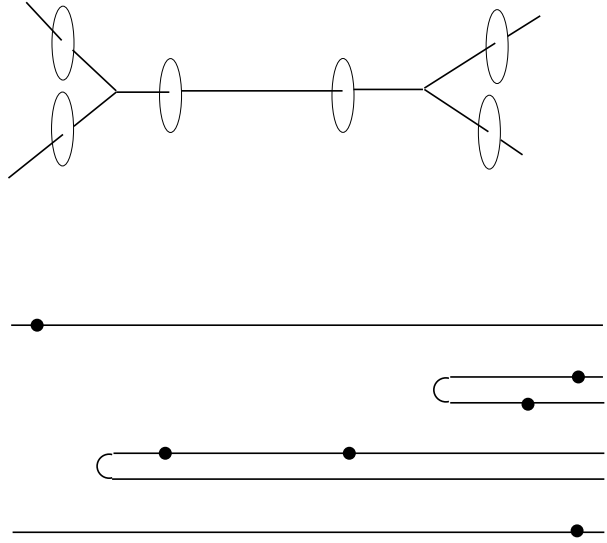}
\end{center}
\caption{Construction of $\Delta_{p,m}^{(0,0)}(\Gamma,\lambda)$. The top picture shows the tree and the marking disks, the bottom picture shows the domain.}
\label{f:cnstrdom}
\end{figure}

We equip this domain with additional marked points on the boundary. The marked points in all regions except those corresponding to edge points are considered marked points also in the domain constructed. The end result of this procedure is a family of standard domains with $p$ punctures and $m$ marked points which we denote $\Delta_{p,m}^{(0,0)}(\Gamma,\lambda)$, see Figure \ref{f:cnstrdom}.

We introduce notation for some subsets of $\Delta_{p,m}^{(0,0)}(\Gamma,\lambda)$. Fix a sufficiently large number $K>0$. If $v$ is a $Y_0$-vertex, a $Y_1$-vertex, or a $2$-valent puncture then let $B(v)$ be the region around the boundary minimum corresponding to $v$ bounded by vertical line segments at distance $K$ from it. If $p$ is a switch then let $B(p)$ be the region bounded by vertical line segments at distance $K$ from the point in $\dot\Delta_2(T_1,T_2)$ which map to $\Sigma$. If $q$ is an edge point then let $R(q)$ be the region bounded by vertical line segments of distance $K$ from the marked point in the region $\dot\Delta_2(a,a)$ corresponding to $q$.

\subsubsection{An approximately holomorphic map}\label{ss:apphol}
Using local solutions, we define maps
$$
u_\lambda^{(0,0)}\colon\Delta_{p,m}^{(0,0)}(\Gamma,\lambda)\to T^\ast M
$$
with boundary on $\hat L_\lambda$, such that $u_\lambda^{(0,0)}$ maps $\pa\Delta_{p,m}^{(0,0)}(\Gamma,\lambda)$ into a small neighborhood of the cotangent lift of $\Gamma$, such that $u_\lambda^{(0,0)}$ is $J$-holomorphic along most of $\Delta_{p,m}^{(0,0)}(\Gamma,\lambda)$ and along all of its boundary, and such that
$$
\sup_{\Delta_{p,m}^{(0,0)}(\Gamma,\lambda)}|\bar\pa_J u_\lambda^{(0,0)}|=\Ordo(\lambda\log(\lambda^{-1})).
$$
This map is defined in the natural way: in all regions of $\Delta_{p,m}^{(0,0)}(\Gamma,\lambda)$ except the in the edge point regions, take $u_\lambda^{(0,0)}$ to agree with the corresponding local solution. In particular, near the ends and $1$-valent punctures we take a local solution which map the vertical segment on the boundary of $\dot\Delta_2(0,\infty)$ to an $\Ordo(\lambda)$-neighborhood of the nearest edge point.

To complete the definition we extend the map $u_\lambda^{(0,0)}$ to the edge point regions ($\approx [-a,a]\times[0,1]$). Note that by definition of the local solutions and of $\Delta_{p,m}^{(0,0)}(\Gamma,\lambda)$ the restriction of $u_\lambda^{(0,0)}$ to an edge point region must interpolate between two maps of distance $\Ordo(\lambda\log(\lambda^{-1}))$. It follows from the flow line convergence of holomorphic disks, see Lemma \ref{l:gradconv} and Remark \ref{r:gradconv},  that we can choose $u_\lambda^{(0,0)}$ so that
$$
\sup_{\Delta_{p,m}^{(0,0)}(\Gamma,\lambda)}|D^k\bar\pa_J u_\lambda^{(0,0)}|=\Ordo(\lambda\log(\lambda^{-1})), \quad k=0,1
$$
and so that $\bar\pa u_\lambda^{(0,0)}$ is supported only in the edge point regions. Moreover, after a small perturbation (which does affect the above estimate) we may assume that the restriction of $\bar\pa_J u_\lambda^{(0,0)}$ to the boundary equals zero.

Note that the marked points on $\pa\Delta_{p,m}^{(0,0)}(\Gamma,\lambda)$ are related to the geometry of the map $u_\lambda^{(0,0)}$. If $p\in\pa\Delta_{p,m}^{(0,0)}(\Gamma,\lambda)$ then there is a corresponding marking $(n-1)$-disk $D_p\subset\hat L_\lambda$ perpendicular to the cotangent lift of $\Gamma$ such that $u_\lambda^{(0,0)}(p)\in D_p$.

\subsubsection{Varying the conformal structure}\label{ss:varyconf}
We will consider variations of the conformal structure of $\Delta_{m,p}^{(0,0)}(\Gamma,\lambda)$. We parameterize nearby conformal structures by fixing one marked point in a strip region corresponding to an edge segment connecting two edge points and by moving all boundary minima and all other marked points.

\begin{figure}[htbp]
\begin{center}
\includegraphics[angle=0, width=8cm]{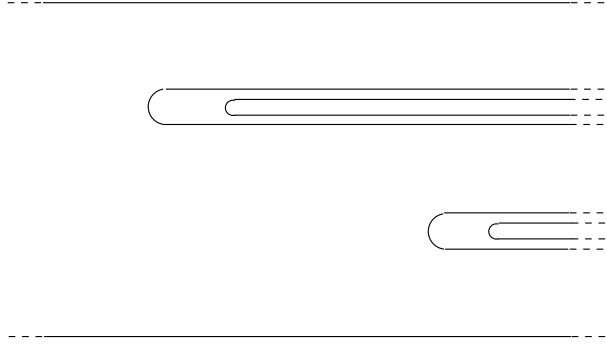}
\end{center}
\caption{The neighborhood $U(\Gamma,\lambda)$.}
\label{f:bigger}
\end{figure}

Fix $l\ll l_0$, where $l_0$ is as in Remark \ref{r:flowlength}. We first describe how to move the boundary minima. Consider $\Delta_{p,m}^{(0,0)}$ as a subset of $\C$ with coordinates $\zeta=\tau+it$ and with the $\tau$ coordinate of the fixed marked point equal to $0$. We use an auxiliary standard domain $U(\Gamma,\lambda)$ with $p$ punctures, which has slits of much smaller width than $\Delta_{p,m}^{(0,0)}(\Gamma,\lambda)$, and with one boundary minimum $q'$ for each boundary minimum $q$ of $\Delta_{p,m}^{(0,0)}(\Gamma,\lambda)$, such that if the $\tau$-coordinate of $q$ equals $\tau_q$ then the $\tau$-coordinate $\tau_{q'}$ of $q'$ equals $\tau_{q}+\frac14 l_0\lambda^{-1}$. See Figure \ref{f:bigger}. Note that $U(\Gamma,\lambda)$ is a neighborhood of $\Delta_{p,m}^{(0,0)}(\Gamma,\lambda)$ in $\C$.

Number the boundary minima of $\Delta_{p,m}^{(0,0)}(\Gamma,\lambda)$ from $1$ to $p-2$. Consider the $j^{\rm th}$ boundary minimum and assume that its $\tau$-coordinate equals $\tau_j$. Let the $\tau$-coordinates of the nearby marked points be $\tau_j^k$, $k=0,1,2$ if the boundary minimum corresponds to a $Y_0$- or $Y_1$-vertex, let them be $\tau_j^r$, $r=1,2$ if it corresponds to a $2$-valent puncture. We choose notation so that in the former case $\tau_j^0<\tau_j$, and $\tau_j^k>\tau_j$, $k=1,2$, and that in the latter $\tau_j^2>\tau_j$ and $\tau_j^1>\tau_j$ if the puncture is at $-\infty$, or $\tau_j^1<\tau_j$ if the puncture is at $+\infty$. By definition of $l_0$,  $|\tau_j-\tau_j^k|\ge \frac12l_0\lambda^{-1}$ and the distance from any  marked point to the nearest edge point region is bounded from below by $\frac12l_0\lambda^{-1}$ as well.

Consider the subdomains $W_{cl}$, $c=1,2$, of $U(\Gamma,\lambda)$ which in case the boundary minimum corresponds to a $Y_0$- or $Y_1$-vertex are bounded by vertical line segments with $\tau$-coordinates $\tau_j^0-cl\lambda^{-1}$,  and $\tau_j^k+cl\lambda^{-1}$, $k=1,2$, and in case it is a $2$-valent puncture are bounded by vertical line segments with $\tau$-coordinates $\tau_j^2+cl\lambda^{-1}$ and $\tau_j^1\pm cl\lambda^{-1}$ if the puncture is at $\mp\infty$. Pick a function $\alpha\colon U(\Gamma,\lambda)\to\C$ such that $\alpha=1$ on $W_l$, $\alpha=0$ outside $W_{2l}$, $\alpha$ is real valued and holomorphic on $\pa\Delta_{p,m}^{(0,0)}(\Gamma,\lambda)$, and such that $|D^k\alpha|\le 10^2 l^{-1}\lambda$, $k=1,2$. Let $b=10^{-3}l$ and let $f_j$ be the vector field on $U(\Gamma,\lambda)$ given by $f_j(\tau+it)= b\lambda^{-1}\alpha\pa_\tau$. If $\phi_j^s$, $-1\le t\le 1$ denotes the time $s$ flow of $f_j$ then $\phi_j^s(\Delta_{p,m}^{(0,0)})(\Gamma,\lambda)$ is a standard domain. Note also that
\begin{equation}\label{e:derssmall}
|d\phi_j^s-\id|=\Ordo(s)
\end{equation}
(here we use the canonical trivialization of $T\C$),
that the diffeomorphisms $\phi_j^{s'}$ and $\phi_r^{s''}$ commute, and that in a $\frac12 l_0\lambda^{-1}$ neighborhood of the boundary minimum $\phi_j^s(\tau+it)=\tau + sb\lambda^{-1} +it$. For $\kappa=(\kappa_1,\dots,\kappa_{p-2})\in(-1,1)^{p-2}$, let $\phi^\kappa=\phi_1^{\kappa_1}\circ\dots\circ\phi_{p-2}^{\kappa_{p-2}}$. Define
$$
\Delta_{p,m}^{(\kappa,0)}(\Gamma,\lambda)=
\phi^\kappa(\Delta_{p,m}^{(0,0)}(\Gamma,\lambda)),
$$
equipped with marked points which are the images of the marked points under $\phi^\kappa$.

We next define a map $u_\lambda^{(\kappa,0)}\colon\Delta_{p,m}^{(\kappa,0)}(\Gamma,\lambda)\to T^\ast M$. To this end we say that $\phi^\kappa$-images of the edge point regions of $\Delta_{p,m}^{(0,0)}(\Gamma,\lambda)$ are edge point regions in $\Delta_{p,m}^{(\kappa,0)}(\Gamma,\lambda)$. (Note that the restriction of $\phi_\kappa$ to any edge point region equals the identity.) Let $\dot\Delta_2(a_k,a_k)$ be the edge point region nearest the marked point $p_k$ and consider the neighborhood $R(p)=\dot\Delta_2(a_k+K,a_k+K)$ of it. Let $V_j=\Delta(T_1,T_2,T_3)\subset\Delta_{p,m}^{(0,0)}(\Gamma,\lambda)$ be the neighborhood of the boundary minimum where $u_\lambda$ equals the local solution $s_\lambda$. Let $\beta\colon\phi_\kappa(V_j)\to\C$ be a function which equals $1$ outside the edge point region neighborhoods $\dot\Delta_2(a_k+\frac12 K,a_k+\frac12K)$, $k=0,1,2$ or $k=1,2$, and which equals $0$ inside $\dot\Delta_2(a_k+\frac14 K,a_k+\frac14 K)$, $k=0,1,2$ or $k=1,2$. Also let $\beta$ be real valued and holomorphic on the boundary and let $|D^k\beta|=\Ordo(1)$, $k=1,2$. Define $u_\lambda^{(\kappa,0)}$ to be the map which equals
$$
\beta(\tau+it)\cdot s_\lambda(z-\kappa_jb\lambda^{-1}+it)
+(1-\beta(\tau+it))\cdot s_\lambda(z),
$$
on $\phi^\kappa(V_j)$ and equals $u_\lambda^{(0,0)}$ everywhere else. Then $u_\lambda^{(\kappa,0)}$ is holomorphic outside a finite neighborhood of the edge points and on the boundary. Moreover, since all sheets of $\hat L_\lambda$ are subsets of affine Lagrangian subspaces in a neighborhood of each $Y_0$- and $Y_1$-vertex and in a neighborhood of each puncture
$$
u_\lambda^{(\kappa,0)}(\pa\Delta_{p,m}^{(\kappa,0)}(\Gamma,\lambda))
\subset\hat L_\lambda.
$$

We consider in a similar way the conformal variations which arise from moving the punctures. To this end we first note that if $p$ is a marked point of $\Delta_{p,m}^{(0,0)}$ then $u_\lambda^{(\kappa,0)}(\phi_\kappa(p))=u_\lambda^{(0,0)}(p)$. It follows in particular that $u_\lambda^{(\kappa,0)}$ maps the marked point $p_\kappa\in\Delta_{p,m}^{(\kappa,0)}(\Gamma,\lambda)$ corresponding to $p\in\Delta_{p,m}^{(0,0)}(\Gamma,\lambda)$ to $D_p$.

Number the moving marked points in $\Delta_{p,m}^{(0,0)}(\Gamma,\lambda)$ from $1$ to $m-1$. This induces a numbering of the moving marked points of $\Delta_{p,m}^{(\kappa,0)}(\Gamma,\lambda)$ for each $\kappa\in(-1,1)^{p-2}$. Fix $\kappa$ and consider a the $j^{\rm th}$ marked point $p_j$ of $\Delta_\lambda^{(\kappa,0)}(\Gamma,\lambda)$.
Let $\dot\Delta_2(T_1,T_2)$ be a neighborhood of $p_\kappa^j$ which is bounded by two vertical line segments which either are of finite distance from a boundary minimum or from an edge point region, where the boundary minimum or the edge point region are the ones closest to $p_\kappa$. As above we consider a function $\alpha\colon\dot\Delta_2(T_1,T_2)\to\C$ which is real valued and holomorphic on the boundary of $\Delta_{p,m}^{(\kappa,0)}(\Gamma,\lambda)$, which equals $1$ on $\dot\Delta_2(l\lambda^{-1},l\lambda^{-1})$, which equals $0$ outside $\dot\Delta_2(2l\lambda^{-1},2l\lambda^{-1})$, and which has $|D^k\alpha|\le 10^2 l\lambda$. Let $\psi^s_j$  be the time $s$ flow of the vector field $\alpha b\lambda^{-1}\pa_\tau$ and note that $\psi^s$ is a diffeomorphism of $\Delta_{p,m}^{(\kappa,0)}(\Gamma,\lambda)$ with differential of distance $\Ordo(s)$ from $\id$. For $\nu=(\nu_1,\dots,\nu_{m-1})\in [-1,1]^{m-1}$ define $\Delta_{p,m}^{(\kappa,\nu)}(\Gamma,\lambda)$ to be the disk $\Delta_{p,m}^{(\kappa,0)}$ with marked points at the images of the marked points under the diffeomorphism $\psi_\nu=\psi_1^{\nu_1}\circ\dots\circ\psi_{m-1}^{\nu_{m-1}}$.

We next define a map $u_\lambda^{(\kappa,\nu)}$ which takes the marked points of $\Delta_{p,m}^{(\kappa,\nu)}(\Gamma,\lambda)$ to the corresponding marking disks.
In order to define the map $u_\lambda^{(\kappa,\nu)}$ we will use a similar construction as above, interpolating between the local solution and a shifted version of the local solution in a finite neighborhoods of the edge point regions. In order to define the interpolation for a marked point in a region connecting two edge points we must first discuss a certain exponential map. Consider a small $\delta$-neighborhood of an edge point in $\Gamma$. We scale this $\delta$-neighborhood by $\lambda$. Note that by definition $\lambda^{-1}\hat L_\lambda$ is flat outside a finite neighborhood of the center of these re-scaled coordinate system. Let $h'$ be a metric in these re-scaled coordinate neighborhoods such that $\lambda^{-1}\hat L_\lambda$ is totally geodesic and such that the space of Jacobi-fields along a geodesic on $\hat L_\lambda$ is closed under multiplication with $J_\lambda=J(\lambda x)$, and moreover so that the metric is flat outside a finite neighborhood and such that near $\Sigma_1$ the metric is flat in the directions orthogonal to the cusp. (The form of such a metric can be found in Proposition 4.3 \cite{EES1}.) Outside the $\delta$-neighborhoods we let $h''$ be the flat metric. Note that  $\exp^{h'}=\lambda\exp^{h''}\lambda^{-1}$ for sufficiently short tangent vectors in a region away from the center of the coordinate system where these two maps are both defined. Thus they patch together to a map $\exp\colon T^\ast M\to T^\ast M$, defined in a neighborhood of the tree. Note moreover that except in a $\lambda$-neighborhood of edge points and points in $\Sigma_1$ this exponential map is just addition.

To define the map $u^{(\kappa,\nu)}_\lambda$  note that there exists a vector field $\xi\colon\dot\Delta_2(a+K,a+K)\to T(T^\ast M)$ in a neighborhood of any edge point region $\dot\Delta_2(a,a)$ such that
$$
s_\lambda^{(\kappa,0)}(z-\nu_jb\lambda^{-1})=\exp_{s_\lambda(z)}(\xi(\zeta)),
$$
where $s_\lambda$ is the standard solution. Fix a cut-off function $\beta$ equal to $1$ on most of the strip region and becoming $0$ in a finite neighborhood of the edge point regions surrounding it. As above we require its first two derivatives to be bounded and that it is real valued and holomorphic on the boundary. In the strip region of the $j^{\rm th}$ marked point we take
$$
u_\lambda^{(\kappa,\nu)}=
\begin{cases}
u_\lambda^{(\kappa,0)}(z-\nu_j)=s_\lambda(z-\nu_jb\lambda^{-1}) &\text{ for }z\notin\dot\Delta_2(a+K,a+K),\\
\exp_{s_\lambda(z)}(\beta(z)\xi(z)), &\text{ for } z\in\Dot\Delta_2(a,a).
\end{cases}
$$

To define the map $u_\lambda^{(\kappa,\nu)}$ for marked points adjacent to a boundary minimum, or a switch region, or a $2$-valent puncture we use a different construction. Let $v$ be the vertex in question and let $B(v)$ denote a finite but large region around the boundary minimum or in the case of a switch a finite region outside which the boundary of the local solution maps to constant lines. Letting $\beta$ be a cut-off function with derivative which has support in $B(v)$ and in the nearest edge point region we define
$$
u_\lambda^{(\kappa,\nu)}=u_\lambda^{(\kappa,0)}+\beta\pa u_\lambda^{(\kappa,0)}\cdot \nu_j'\lambda^{-1}b\pa_\tau,
$$
where $\nu_j'$ is the unique number so that $u_\lambda^{(\kappa,\nu)}$ takes its marked point to $D_p$.

This completes the definition of $u_\lambda^{(\kappa,\nu)}$. It  takes the boundary of $\Delta_{p,m}^{(\kappa,\nu)}(\Gamma,\lambda)$ to $\hat L_\lambda$. It is holomorphic outside edge point regions and regions of the form $B(v)$ and it is holomorphic everywhere on the boundary. Finally, it takes the marked points $\psi_\nu(p)\in\Delta_{p,m}^{(\kappa,\nu)}(\lambda,\Gamma)$ to the marking disk $D_p$.

The construction above gives a local coordinate system $(-1,1)^{p+m-3}\to\conf_{p,m}$ on the space of conformal structures on the disk with $p$ punctures and $m$ marked points on the boundary around the conformal structure on $\Delta_{m,p}^{(0,0)}(\Gamma,\lambda)$. We denote the image of this coordinate system $\B_{p,m}$ and denote points in $\B_{p,m}$ by their coordinates $(\kappa,\nu)\in(-1,1)^{p-2}\times(-1,1)^{m-1}$.

\subsection{The $\bar\pa_J$-operator as a bundle map}\label{s:barpabund}
In this subsection we define a bundle of function spaces over $\B_{p,m}$. Its fibers consists of sections $s\colon\Delta_{p,m}^{(\kappa,\nu)}(\Gamma,\lambda)\to\bigl({u_\lambda^{(\kappa,\nu)}}\bigr)^\ast T(T^\ast M)$. More precisely, each fiber is a direct sum of an infinite dimensional weighted Sobolev space (see \S\ref{ss:wtsblv}) augmented by a finite dimensional function space (see \S\ref{ss:cutoffsol}). Once the fibers are described we give the bundle structure by defining a trivialization (see \S\ref{ss:triv}) and thereby equip its total space with a (Banach-) manifold structure.

\begin{rmk}
It seems to the author that some construction involving weighted spaces, similar to the one discussed below, is necessary also in the case when the singularity set of the Legendrian submanifold is empty. In fact, the construction below grew out of the authors attempts to understand a corresponding construction in \cite{FO} which is also used to prove existence of holomorphic disks near flow trees but which does not involve exponential weights. The details of the argument in \cite{FO} are presented only in the case of $3$-punctured disks and it contains misstatements. In particular, the crucial estimate $\|Q_0\|\le C_9$ in Lemma 6.1 is not true as stated: in the proof of Lemma 6.1 it is claimed that there exists a constant $C$ such that for any function $u\colon \Delta_3\to\C^n$ in $W^{1,p}$ such that $u(\pa \Delta_3)\subset\R^n$, $\|u\|_{1,p}\le C\|\bar\pa u\|_p$. Consider the sequence of functions $v_N\colon [0,\infty)\times[0,1]\to\C^n$,
$$
v_N(\tau,t)=
\begin{cases}
N^{-\frac{p}{p+1}}\tau, & 0\le \tau\le N,\\
N^{-\frac{p}{p+1}}(2N-\tau), & N\le \tau\le 2N,\\
0, &\text{otherwise}.
\end{cases}
$$
then $\|v_N\|_{1,p}\ge \frac{1}{1+p}$ and $\|\bar\pa v_N\|_p\to 0$. This shows that no such $C$ exists. The failure of this and similar estimates are among the reasons for the introduction of spaces weighted with somewhat complicated functions below.
\end{rmk}

\subsubsection{Weighted Sobolev spaces}\label{ss:wtsblv}
We specify weight functions, boundary conditions, and certain vanishing conditions. Consider a domain $\Delta_{p,m}^{(\kappa,\nu)}(\Gamma,\lambda)$. Note that the complement of all regions $B(v)$ and $R(q)$, where $v$ ranges over $2$- and $3$-valent vertices and where $q$ ranges over the edge points in $\Delta_{p,m}^{(\kappa,\nu)}(\Gamma,\lambda)$, is a collection of strip regions. For simpler notation below we re-scale any such region to have width $1$ and simply denote it $[-T_1,T_2]\times[0,1]$ if it has finite length and $[0,\infty)\times[0,1]$ otherwise. Here we choose $T_1$ and $T_2$ so that the marked point lies at $(0,0)$.

We define weight functions ${\mathbf w}^{(\kappa,\nu)}\colon \Delta_{p,m}^{(\kappa,\nu)}(\Gamma,\lambda)\to[1,\infty)$. This function will be the smoothing of an exponential function composed with a piecewise linear function. The details of the smoothing will not affect any of the arguments below and will therefore not be discussed. Fix a small $\delta>0$.
\begin{itemize}
\item
For $z\in\Delta_{p,m}^{(\kappa,\nu)}(\Gamma,\lambda)$ which lies in $B(v)$ or in $R(p)$, for some vertex $v$ or some edge point $p$, let ${\mathbf w}^{(\kappa,\nu)}(z)= 1$.
\item
For $z=\tau+it\in\Delta_{p,m}^{(\kappa,\nu)}(\Gamma,\lambda)$ which lies in a half infinite strip  $[0,\infty)\times[0,1]$ (attached to a region of the form $B(v)$ or $R(p)$ along $\{0\}\times[0,1]$) let
$$
{\mathbf w}^{(\kappa,\nu)}(\tau+it)=e^{\delta|\tau|}.
$$
\item
For $z=\tau+it\in\Delta_{p,m}^{(\kappa,\nu)}(\Gamma,\lambda)$ which lies in a finite strip $[-T_1,T_2]\times[0,1]$ (attached to a region of the form $B(v)$ or $R(p)$ along $\{-T_1\}\times[0,1]$ and along $\{T_2\}\times[0,1]$), let
$$
\begin{cases}
{\mathbf w}(\tau+it)=e^{\delta_1(T_1-|\tau|)}, & -T_1\le \tau\le 0,\\
{\mathbf w}(\tau+it)=e^{\delta_2(T_2-|\tau|)}, & 0\le \tau\le T_2,
\end{cases}
$$
where we choose $\delta_1$ and $\delta_2$ so that ${\mathbf w}(0,0)=e^{\delta T}$, where $T=\frac12(T_1+T_2)$.
\end{itemize}

For $k=1,2$, we let $S_{k,\delta}^{(\kappa,\nu)}(\Gamma,\lambda)$ denote the Sobolev space of sections
$$
s\colon \Delta_{p,m}^{(\kappa,\nu)}(\Gamma,\lambda)\to \bigl(u_\lambda^{(\kappa,\nu)}\bigr)^\ast T(T^\ast M)
$$
with $k$ derivatives in $L^2$ and weighted by the function ${\mathbf w}^{(\kappa,\nu)}$. Let $\hat V_{2,\delta}^{(\kappa,\nu)}(\Gamma,\lambda)\subset S_{2,\delta}^{(\kappa,\nu)}(\Gamma,\lambda)$ be the closed subspace of elements $\hat v$ such that $\hat v$ is tangent to $\hat L_\lambda$ along $\pa \Delta_{p,m}^{(\kappa,\nu)}(\Gamma,\lambda)$, and such that the restriction of $\bar \nabla_J \hat v=\nabla \hat v + J(\nabla \hat v) i$ to the boundary (the trace) equals $0$. (Here $\nabla$ is the connection of the metric on $T^\ast M$ and  $i$ is the complex structure on $\Delta_{p,m}^{(\kappa,\nu)}(\Gamma,\lambda)$ inherited from $\C$.) For more background on the requirements in this definition, we refer to Lemma 3.2 in \cite{EES4}, where the analogous conditions specifies the "tangent space of the space of candidate maps". Since our study in the present situation is more local than that in \cite{EES4} it suffices to work with the tangent space and an exponential map rather than working with the whole configuration space of candidate maps. The details of this definition will not be discussed here, we simply note that since the metric on $M$ is flat in a neighborhood of any rigid tree and we have local (holomorphic) $\C^n$-coordinates in a corresponding neighborhood in $T^\ast M$ we may study the $\bar\pa_J$-operator following \cite{EES2}. In particular, the boundary condition $\bar \nabla_J \hat v=0$ translates into the boundary condition $\bar\pa \hat v=0$ where $\bar\pa$ is the usual $\bar\pa$-operator on $\C^n$.

Define $V_{2,\delta}^{(\kappa,\nu)}(\Gamma,\lambda)\subset\hat V_{2,\delta}^{(\kappa,\nu)}(\Gamma,\lambda)$ to be the subspace of sections $v$ which vanishes at all marked points in $\Delta_{p,m}^{(\kappa,\lambda)}(\Gamma,\lambda)$. This subspace has finite codimension: each marked point increases the codimension by $n$.

Finally, define $W_{1,\delta}^{(\kappa,\nu)}\subset S_{1,\delta}^{(\kappa,\nu)}$ to be the closed subspace of $w$ such that the restriction of $w$ to $\pa \Delta_{p,m}^{(\kappa,\nu)}(\Gamma,\lambda)$ (the trace) vanishes.

\subsubsection{Cut-off solutions}\label{ss:cutoffsol}
We discuss cut-off solutions in the $\C^n$-coordinates corresponding to the flat coordinates on $M$ in neighborhoods of pieces of the rigid trees where these are defined. The boundary conditions in $V_{2,\delta}^{(\kappa,\nu)}(\Gamma,\lambda)$ in regions
$[-T_1,T_2]\times[0,1]$ between edge points are of one of the following forms (in standard coordinates $\tau+it$)
\begin{equation}\label{e:flat}
\begin{cases}
v(\tau)&\in \R^n,\\
v(\tau+i)&\in \R^n,
\end{cases},
\end{equation}
\begin{equation}\label{e:slope1}
\begin{cases}
v(\tau)&\in \R^n,\\
v(\tau+i)&\in \bigl(e^{i\theta_\lambda}\R\bigr)\times\R^{n-1},
\end{cases}
\end{equation}
or
\begin{equation}\label{e:slope2}
\begin{cases}
v(\tau)&\in \bigl(e^{i\theta_\lambda}\R\bigr)\times\R^{n-1},\\
v(\tau+i)&\in \R^n,
\end{cases}
\end{equation}
where $\theta_\lambda\to 0$ as $\lambda\to 0$. Let $p$ denote the marked point in $[-T_1,T_2]\times[0,1]$ and define $\hat V(p)$ to be the $n$-dimensional space of solutions of the $\bar\pa$-equation on $\R\times[0,1]$ which satisfies these boundary conditions and which converges to constant maps as $\lambda\to 0$. In case of \eqref{e:flat} these solutions have the form
$$
v(\tau+it)=c,\quad c\in\R^n,
$$
and in case of \eqref{e:slope1},
$$
v(\tau+it)=(c_1 e^{\theta_\lambda(\tau+it)},c'),\quad c_1\in\R,\,\,c'\in\R^{n-1}.
$$
The solutions in case of \eqref{e:slope2} are similar. Define $\hat V(p)$ to be the vector space spanned by these solutions. Pick a cut-off function
$\beta\colon[-T_1-1, T_2+1]\times[0,1]\to\C$
which equals $1$ on $[-T_1,T_2]\times[0,1]$ and which equals $0$ outside $[-\frac12- T_1,\frac12+ T_2]\times[0,1]$. Consider the functions $\beta v$, $v\in \hat V(p)$, and note that these functions extend by $0$ to the rest of $\Delta_{p,m}^{(\kappa,\nu)}(\Gamma,\lambda)$. Define $V_{\rm sol}^{(\kappa,\nu)}(p,\lambda)$ to be the $(n-1)$-dimensional space of functions of this form which takes values in $T D_p$, where $D_p$ is the marking disk at $p$. We endow these spaces with the norm induced by evaluation at $p$.

Similarly it is easy to find solutions of the (linearized) $\bar\pa$-equation which are bounded and converges to constants as $\lambda\to 0$, near each $1$-valent puncture of the tree. Consider a negative $1$-valent puncture $v$. Assume that the edge ending at $v$ is a flow line of some function difference $f_1-f_2$ and assume that $f_1(v)>f_2(v)$. Then the negative sign of the puncture implies that the orientation of the flow line given by the vector $-\nabla(f_1-f_2)$ points toward $v$. In particular, the flow line lies in the stable manifold of $f_1-f_2$. The  linearized boundary condition for $\tau+it\in[0,\infty)\times[0,1]$ in standard coordinates has the following  form:
$$
\begin{cases}
v(\tau)&\in \Bigl(e^{i\theta_\lambda}\R\,\times\stackrel{k}{\cdots}\times\, e^{i\theta_\lambda}\R\Bigr)
\times \Bigl(e^{-i\theta_\lambda}\R\,\times\stackrel{n-k}{\cdots}\times\, e^{-i\theta_\lambda}\R\Bigr),\\
v(\tau+i)&\in\R^n,\\
\end{cases}
$$
Here, the $k$ first factors correspond to stable directions at the critical point of $f_1-f_2$ and the $n-k$ last to unstable ones. The bounded solutions which converges to constants have the form $v(\tau+it)=(v_1,\dots,v_n)$ where
$$
v_j(\tau+it)=
\begin{cases}
ce^{i\theta_\lambda}e^{-\theta_\lambda(\tau+it)},\quad c\in\R&\text{ if }1\le j\le k,\\
0 & \text{ if } k+1\le j\le n,
\end{cases}
$$
In particular, the solutions form an $k$-dimensional vector space, where $k$ is the dimension of the stable manifold of the positive function difference at the critical point. The same argument gives an $l$ dimensional vector space of solutions at a positive puncture, where $l$ is the dimension of the unstable manifold of the positive function difference defining the critical point.

Pick a cut-off function $\beta\colon[0,\infty)\times[0,1]\to\C$ with properties analogous to the cut-off function discussed above and
such that $\beta=1$ on $[1,\infty)\times[0,1]$ and $\beta=0$ on $[0,\frac12]\times[0,1]$. Define the finite dimensional space $V_{\rm sol}^{(\kappa,\nu)}(v,\lambda)$ to be the space spanned by cut-off solutions at a puncture $v$ and endow it with the supremum norm. Note that if $v$ is a puncture ($1$- or $2$-valent) of a rigid flow tree mapping to a critical point $m$ of the function difference $f_1-f_2>0$ then
$$
\dim(V_{\rm sol}(v,\lambda))=
\begin{cases}
\dim\bigl(\text{stable manifold of $-\nabla(f_1-f_2)$ at $m$}\bigr)&\text{ if $v$ is negative,}\\
\dim\bigl(\text{unstable manifold of $-\nabla(f_1-f_2)$ at $m$}\bigr)&\text{ if $v$ is positive.}
\end{cases}
$$

In the region $[0,\infty)\times[0,1]\subset \Delta_{p,m}^{(\kappa,\nu)}(\Gamma,\lambda)$ corresponding to an end $e$, the boundary conditions in $V_{2,\delta}^{(\kappa,\nu)}(\Gamma,\lambda)$ split into a product of a boundary condition in $\C$, tangential to $u_\lambda(\pa \Delta_{p,m}^{(\kappa,\nu)}(\Gamma,\lambda))$, and a constant boundary condition in $\C^{n-1}$. If a function $v=(v_1,v')\in\C\times\C^{n-1}$ is decomposed into components then for $\tau+it\in[0,\infty)\times[0,1]$,
\begin{equation}\label{e:flatperp}
\begin{cases}
v'(\tau)&\in \R^{n-1},\\
v'(\tau+i)&\in \R^{n-1}.
\end{cases}
\end{equation}
In analogy with the above, we let $V_{\rm sol}^{(\kappa,\nu)}(e,\lambda)$ be the $(n-1)$-dimensional space of cut-off versions of the constant solutions with values in $\R^{n-1}$, satisfying these boundary conditions.

\begin{dfn}\label{d:Vsol}
Let $\Gamma$ be a rigid flow tree of $\hat L_\lambda$. Let $S(\Gamma)$ be the set of marked points, ends, and punctures of $\Gamma$. Define
$$
V_{\rm sol}^{(\kappa,\nu)}(\Gamma,\lambda)=\bigoplus_{a\in S(\Gamma)} V_{\rm sol}^{(\kappa,\nu)}(a,\lambda).
$$
\end{dfn}

\subsubsection{A trivialization}\label{ss:triv}
We define bundles $\V_{2,\delta}(\Gamma,\lambda)$, $\V_{\rm sol}(\Gamma,\lambda)$, and $\W_{1,\delta}(\Gamma,\lambda)$ over $\B_{p,m}$ such that the fibers over $(\kappa,\nu)$ are $V_{2,\delta}^{(\kappa,\nu)}(\Gamma,\lambda)$, $V_{\rm sol}^{(\kappa,\nu)}(\Gamma,\lambda)$, and $W_{1,\delta}^{(\kappa,\lambda)}(\Gamma,\lambda)$, respectively, simply by defining trivializations in the following way.

Consider first the bundle $\V_{\rm sol}(\Gamma,\lambda)$. To trivialize it, note that evaluation at the marked points or at points near a puncture give canonical identifications of the cut-off solutions in $V_{\rm sol}^{(0,0)}(\Gamma,\lambda)$ and $V_{\rm sol}^{(\kappa,\nu)}(\Gamma,\lambda)$ which we use to define the trivialization.

Consider second the bundle $\V_{2,\delta}(\Gamma,\lambda)$. To trivialize it we use the map
$$
\Psi[(0,0),(\kappa,\nu)]\colon V_{2,\delta}^{(0,0)}(\Gamma,\lambda)\to V_{2,\delta}^{(\kappa,\nu)}(\Gamma,\lambda)
$$
given by
\begin{equation}\label{e:wigthcomp}
\Psi[(0,0),(\kappa,\nu)]v(z)=
\left(\frac{{\mathbf w}^{(0,0)}(\phi_{-\kappa}\circ\psi_{-\nu}(z))}{{\mathbf w}^{(\kappa,\nu)}(z)}\right)^\dagger v\bigl(\phi_{-\kappa}\circ\psi_{-\nu}(z)\bigr)
\end{equation}
where $(f)^\dagger$ denotes a small perturbation of the function $f$ making it holomorphic on the boundary (and keeping it real valued) and where we use the canonical parallel translation in the flat coordinates to put the value of $\Psi[(0,0),(\kappa,\nu)]$ in the right fiber. Note that an inverse $\Psi[(0,0),(\kappa,\nu)]$ can be defined similarly and that these maps are uniformly (as $\lambda\to 0$) bounded isomorphisms. This trivializes our bundle.

We trivialize the bundle $\W_{1,\delta}(\Gamma,\lambda)$ in a completely analogous way.

\subsubsection{Source and target spaces of the $\bar\pa$-operator}
We consider the $\bar\pa_J$-operator in a neighborhood of $u_\lambda^{(0,0)}$ as a bundle map
$$
\bar\pa_J\colon\V_{2,\delta}(\Gamma,\lambda)\oplus \V_{\rm sol}(\Gamma,\lambda)\to \W_{1,\delta}(\Gamma,\lambda)
$$
We first define an exponential map associating maps $\Exp(v,f)\colon\Delta_{p,m}^{(\kappa,\nu)}(\Gamma,\lambda)\to T^\ast M$ to $(v,f)\in W_{2,\delta}^{(\kappa,\nu)}\oplus V_{\rm sol}^{(\kappa,\nu)}(\Gamma,\lambda)$. To define this map we use the map $\exp$ defined near edge points in \S\ref{ss:varyconf} and a map near ends defined in the same way: scale small $\delta$-neighborhoods of the ends by $\lambda^{-1}$. By definition $\lambda^{-1}\hat L_\lambda$ is flat outside a finite neighborhood of the center of these re-scaled coordinate system. Let $h'$ be a metric in these re-scaled coordinate neighborhoods with the properties of $h'$ in \S\ref{ss:varyconf} and such that $h'$ is flat in the directions orthogonal to the cusp edge. Outside the $\delta$-neighborhoods we let $h''$ be the flat metric. As before the two exponential maps patch together nicely and we get a map $\exp\colon T^\ast M\to T^\ast M$, defined in a neighborhood of the tree.  We define
$$
\Exp(v,f)=\exp_{u_\lambda^{(\kappa,\nu)}(z)}
(v(z)+f(z)),
$$
where $v+f$ denotes the sum in the fiber of $T(T^\ast M)$. The $\bar\pa_J$-operator on the bundle is then defined by
$$
\bar\pa_J\bigl((\kappa,\nu),v,f\bigr)=\bar\pa_J\Exp(v,f),
$$
where $\bar\pa_J h=dh+J\circ dh\circ i$, as usual. It is a $(J,i)$ complex anti linear map. Since a complex anti linear map $T_\zeta \Delta_{p,m}^{(\kappa,\nu)}(\Gamma,\lambda)\to T_{u_\lambda^{(\kappa,\nu)}(\zeta)}T^\ast M$ is determined by the image of one non-zero vector and since $T\Delta_{p,m}^{(\kappa,\lambda)}(\Gamma,\lambda)\subset\C$ is naturally trivialized (by $\pa_\tau$ for example), we view $\bar\pa_{J,j_\gamma}\Exp(v,f)$ as a vector field along $u_\lambda^{(\kappa,\nu)}$ and interpret the $\bar\pa_J$-operator as a bundle map
$$
\bar\pa_J\colon \V_{2,\delta}(\Gamma,\lambda)\oplus
\V_{\rm sol}(\Gamma,\lambda)\to
\W_{1,\delta}(\Gamma,\lambda).
$$

\subsubsection{The linearization of the $\bar\pa_J$-operator}
We will derive the local expression of the linearization of the $\bar\pa_J$-operator described above. More precisely we will compute the linearization $L\bar\pa_J$ at $u_\lambda^{(0,0)}$ in the coordinates on $\V_{2,\delta}(\Gamma,\kappa)$ around $((0,0),0,0)$ in the coordinate system $(-1,1)^{p+m-3}\times V_{2,\delta}^{(0,0)}(\Gamma,\lambda)\oplus V_{\rm sol}^{(0,0)}(\Gamma,\lambda)$ given by the trivialization above. A general expression for the the linearization of such a bundle map can be found in \cite{EES4}. However, as we have standard $\C^n$-coordinates around the tree we can use the results of \cite{EES2} to compute the linearization. We give a brief description.

The partial derivative of $\bar\pa_J$ with respect to $v\in V_{2,\delta}^{(0,0)}(\Gamma,\lambda)$ is
$$
L\bar\pa_J v=\bar\pa v,
$$
and if $f\in V_{\rm sol}^{(0,0)}(\Gamma,\lambda)$ then
$$
L\bar\pa_J f=\bar\pa f.
$$
We turn our attention to the conformal variations. We must compute
$$
\left(\frac{{\mathbf w}^{(\epsilon\kappa,\epsilon\nu)}(\phi_{\epsilon\kappa}\circ\psi_{\epsilon\nu}(z))}{{\mathbf w}^{(0,0)}(z)}\right)^\dagger \bar\pa u_\lambda^{(\epsilon\kappa,\epsilon\nu)}(z)-\bar\pa u_\lambda
$$
to first order in $\epsilon$. Note however that the first term is equal to $1$ on the support of $\bar\pa u_\lambda^{(0,0)}$ and $\bar\pa u_\lambda^{(\kappa,\nu)}$ and that these maps differ by a (cut-off) shift in the parameter domain along a holomorphic map, we find that the linearization is simply
$$
L\bar\pa_J (\kappa,\nu)=(\bar\pa\beta)\pa u_\lambda^{(0,0)}\cdot b\lambda^{-1}\pa_\tau.
$$
This expression agrees with $\bar\pa \bigl(\beta\pa u_\lambda\cdot b\lambda^{-1}\pa_\tau\bigr)$ since $u_\lambda^{(0,0)}$ is holomorphic on the support of $\beta$. We denote the space spanned by the vector fields $\beta\pa u_\lambda^{(0,0)}\cdot b\lambda^{-1}\pa_\tau$ corresponding to all boundary minima and all moving marked points endowed with the supremum norm by $V_{\rm con}^{(0,0)}(\Gamma,\lambda)$ and we will consider the linearization of $\bar\pa_J$ as a map
$$
V_{2,\delta}^{(0,0)}(\Gamma,\lambda)\oplus V_{\rm sol}^{(0,0)}(\Gamma,\lambda)\oplus V_{\rm con}^{(0,0)}(\Gamma,\lambda)\to W_{1,\delta}^{(0,0)}(\Gamma,\lambda).
$$
The advantage of the space $V_{\rm con}^{(0,0)}(\Gamma,\lambda)$ over the more abstract tangent space of conformal structures is that its elements take values in $T(T^\ast M)$ and thus can be acted on by evaluation maps.

\begin{rmk}\label{r:approxsol}
Let $\Gamma$ be a rigid flow tree determined by $\hat L_\lambda$. Note that the fact that $u_\lambda^{(0,0)}\colon \Delta_{p,m}^{(0,0)}(\Gamma,\lambda)\to T^\ast M$ is $J$-holomorphic in the regions where ${\mathbf  w}< 4$ implies
\begin{equation}\label{e:appnorm}
\|\bar\pa_J u_\lambda^{(0,0)}\|_{W_{1,\delta}}=\Ordo(\lambda\log(\lambda^{-1})).
\end{equation}
\end{rmk}

\subsection{Constructing all solutions}\label{s:unifinv}
In this subsection we construct a unique $J$-holomorphic disk near any rigid flow tree by applying the following result (often called Floer's Picard lemma).
\begin{lma}\label{l:FloerPicard}
Let $F\colon B_1\to B_2$ be a smooth Fredholm map of Banach spaces such that
$$
F(v)= F(0)+dF(v)+N(v),
$$
such that $dF$ has a bounded right inverse $Q$, and such that the non-linear term $N$ satisfies a quadratic estimate of the form
\begin{equation}\label{e:quadest}
\|N(u)-N(v)\|_{B_2}\le C(\|u\|_{B_1}+\|v\|_{B_1})\|u-v\|_{B_1}.
\end{equation}
If $\|Qf(0)\|\le \frac{1}{8C}$, then for $\epsilon<\frac{1}{4C}$, $F^{-1}(0)\cap B(0;\epsilon)$, where $B(0;\epsilon)$ is an $\epsilon$-ball around $0\in B_1$, is a smooth submanifold diffeomorphic to the $\epsilon$-ball in $\krn(dF(0))$.
\end{lma}
\begin{pf}
See \cite{Fl:mem}.
\end{pf}
In our application of Lemma \ref{l:FloerPicard}, $F=\bar\pa_J$, $B_1$ is a neighborhood of $((0,0),0,0)\in\V_{2,\delta}(\Gamma,\lambda)$, and $B_2$ is a neighborhood of $((0,0),0)\in\W_{1,\delta}(\Gamma,\lambda)$. We must therefore establish the existence of a uniformly bounded inverse of $L\bar\pa_J$ and produce a quadratic estimate for the remainder term in the Taylor expansion of $\bar\pa_J$.

We find the bounded inverse using an inductive procedure involving partial flow sub-trees of a given rigid flow tree. More precisely, we associate function spaces analogous to  $V_{2,\delta}^{(0,0)}(\Gamma,\lambda)\oplus V_{\rm sol}^{(0,0)}(\Gamma,\lambda)$ and $W_{1,\delta}^{(0,0)}(\Gamma,\lambda)$,  and an operator analogous to $L\bar\pa_J$ to partial flow trees. We then explain how to patch these operators and spaces. The quadratic estimate for the non-linear term follows from a slight modification of the argument given in \cite{EES1, EES4}.

\begin{rmk}\label{r:constbdrycond}
The boundary condition for $L\bar\pa_J$ converges to constant $\R^n$-boundary conditions as $\lambda\to 0$, except near ends, switches, and $Y_1$-vertices, where they converge to constant $\R^n$-boundary conditions with a $\pi$-rotation in one direction as in \S\ref{ss:solend}, \S\ref{ss:solswitch}, and \S\ref{ss:solY1}, respectively. Since we have uniformly positive exponential weights at the punctures the uniform invertibility of the standard $\bar\pa$-operator with constant $\R^n$-boundary conditions and $\pi$-rotations as just described for all sufficiently small $\lambda>0$ will imply the uniform invertibility of $L\bar\pa_J$ with its boundary condition for the same $\lambda>0$. For simplicity, we prove the invertibility in the former setting below.
\end{rmk}

\subsubsection{Auxiliary spaces of partial trees}\label{ss:fspprtree}
We consider partial trees which are subsets of a rigid flow tree $\Gamma$ and associate a subdomain of $\Delta_{p,m}^{(0,0)}(\Gamma,\lambda)$ to them. More precisely, the partial trees $\Gamma'\subset\Gamma$ are of the following form.
\begin{itemize}
\item[$(1)$]
$\Gamma'$ has exactly one special puncture in $D_p$ for some marked point $p\in\Delta_{p,m}^{(0,0)}(\Gamma,\lambda)$ and consists of one of the parts of $\Gamma$ obtained by cutting it at this point. The domain corresponding to $\Gamma'$ is the corresponding part of   $\Delta_{p,m}^{(0,0)}(\Gamma,\lambda)$ cut along the vertical line segment $\{p\}\times[0,1]$ where $p$ is the  marked point.
\item[{$(2_0)$}]
$\Gamma'$ has exactly two special punctures in $D_{p_1}$ and $D_{p_2}$, where $p_j$, $j=1,2$ are marked points in $\Delta_{p,m}^{(0,0)}(\Gamma,\lambda)$ such that the special punctures are connected by an edge segment with one edge point $q$ on it. The domain  corresponding to $\Gamma'$ is the region bounded by $\{p_1\}\times[0,1]$ and $\{p_2\}\times[0,1]$ which contains the edge point region $R(q)$.
\item[{$(2_1)$}]
$\Gamma'$ has exactly two special punctures in $D_{p_1}$ and $D_{p_2}$ and one switch $q$. The domain corresponding to $\Gamma'$ is the region bounded by $\{p_1\}\times[0,1]$ and $\{p_2\}\times[0,1]$, which contains the region $B(q)$ corresponding to the switch $q$.
\item[$(2_2)$]
$\Gamma'$ has exactly two special punctures in $D_{p_1}$ and $D_{p_2}$ and one $2$-valent puncture $q$. The domain corresponding to $\Gamma'$ is the region bounded by $\{p_1\}\times[0,1]$ and $\{p_2\}\times[0,1]$ containing the region $B(q)$ corresponding to $2$-valent puncture $q$.
\item[$(3_0)$]
$\Gamma'$ has exactly three special punctures in $D_{p_j}$, $j=1,2,3$ and one $Y_0$-vertex $v$. The domain corresponding to $\Gamma'$ is bounded by $\{p_j\}\times[0,1]$ where $p_j$, $j=1,2,3$ and contains the region $B(v)$.
\item[$(3_1)$]
$\Gamma'$ has exactly three special punctures in $D_{p_j}$, $j=1,2,3$ and one $Y_1$-vertex $v$. The domain corresponding to $\Gamma'$ is bounded by $\{p_j\}\times[0,1]$ where $p_j$, $j=1,2,3$ and contains the region $B(v)$.
\end{itemize}
If $\Gamma'\subset\Gamma$ is a tree of the form described above we associate a standard domain
$$
\bar\Delta_{p',m'}(\Gamma',\lambda)
$$
to $\Gamma'$. Here $\bar\Delta_{p',m'}(\Gamma',\lambda)$ is the part of $\Delta_{p,m}^{(0,0)}(\Gamma,\lambda)$ which correspond to $\Gamma'$ with half infinite strips of the form $[0,\infty)\times[0,1]$ or $(-\infty,0]\times[0,1]$ attached to all cuts $\{p\}\times[0,1]$ in $\Gamma'$ where $p$ is a marked point in $\Delta_{p,m}^{(0,0)}(\Gamma,\lambda)$. The marked points of $\bar \Delta_{p',m'}(\Gamma',\lambda)$ are the marked points it inherits from $\Delta_{p,m}^{(0,0)}(\Gamma,\lambda)$ except the marked points at cuts.

We define the space of functions $V_{2,\delta}(\Gamma',\lambda)$ on $\bar\Delta_{p',m'}(\Gamma',\lambda)$ much as $V_{2,\delta}^{(0,0)}(\Gamma,\lambda)$ was defined for a rigid tree above: the boundary conditions along the part of $\pa\bar\Delta_{p',m'}(\Gamma',\lambda)$ which is a subset of $\pa\Delta_{p,m}^{(0,0)}(\Gamma,\lambda)$ are left unchanged and along the added half infinite strips the constant boundary conditions (i.e. the boundary conditions in the limit $\lambda\to 0$) in a neighborhood of $\{p\}\times(\pa[0,1])$ are continued constantly. Also, we leave the weight function ${\mathbf w}^{(0,0)}$ unchanged on the subset of $\Delta_{p,m}^{(0,0)}(\Gamma,\lambda)$ and let it continue its exponential growth at the same rate in the infinite half strips added. This gives a weight function ${\mathbf w}$ on $\bar\Delta_{p',m'}(\Gamma',\lambda)$. We also define the space $W_{1,\delta}(\Gamma',\Lambda)$ in analogy with $W_{1,\delta}^{(0,0)}(\Gamma,\lambda)$ using this weight function.

If $p$ is a marked point where the tree $\Gamma$ was cut then we define $V_{\rm sol}(p,\lambda)$ as the finite dimensional space which is spanned by the $n$ cut-off solutions satisfying $\R^n$ boundary conditions that are defined in the infinite end corresponding to $p$. Let $S(\Gamma')$ be the set of all marked points, punctures, and ends of $\Gamma$ which lie in the part of $\Delta_{p,m}^{(0,0)}(\Gamma,\lambda)$ corresponding to $\Gamma'$ and define
$$
V_{\rm sol}(\Gamma',\lambda)=\sum_{a\in S(\Gamma')} V_{\rm sol}(a,\lambda).
$$
(Note that the marked points $p$ where $\Gamma$ was cut belongs to $S(\Gamma')$.)

The conformal variations of $\Gamma'$ are simply the conformal variations inherited from $\Gamma$. More precisely, let $C(\Gamma')$ denote the set of all marked points, $2$-valent punctures, and $3$-valent vertices which lie in $\Gamma'$, except the marked points where $\Gamma$ was cut. Define
$$
V_{\rm con}(\Gamma',\lambda)=\bigoplus_{a\in C(\Gamma)} V_{\rm con}(a,\lambda).
$$
Here we simply let the vector field $\pa u_\lambda\cdot b\lambda^{-1}\pa_\tau$ continue without being cut off along added half strips. As in the proof of Proposition \ref{p:ttv}, we use the convention that the index $I(r)$ of a special vertex $r$ equals $n+1$ if $r$ is positive and equals $-1$ if $r$ is negative.

\subsubsection{Building rigid flow trees}\label{ss:build}
Let $\Gamma$ be a rigid flow tree. We will study the $\bar\pa$-operator on $\Delta_{p,m}^{(0,0)}(\Gamma,\lambda)$ by breaking $\Gamma$ (and simultaneously $\Delta_{p,m}^{(0,0)}(\Gamma,\lambda)$) into simpler pieces. Let $\Gamma$ be a rigid flow tree. If $\Gamma'$ is a partial flow tree of type $(1)$ then we call $\Gamma'$ a {\em rigid flow sub-tree}. If $\Gamma'$ is a partial flow sub-tree then $1\le \dim(\Gamma')\le n$, see Remark \ref{r:dimrigidsubtree}.

\begin{lma}\label{l:build}
Any rigid flow sub-tree $\Gamma$ of a given rigid flow tree can be constructed using only the following operations
\begin{itemize}
\item[$(2_0)$] Join a partial flow tree of type $(2_0)$ to a rigid flow sub-tree $\Gamma'$ at its special puncture. Here $\dim(\Gamma)=\dim(\Gamma')$.
\item[$(2_1)$] Join a partial flow tree of type $(2_1)$ to a rigid flow sub-tree $\Gamma'$ at its special puncture. Here $\dim(\Gamma)=\dim(\Gamma')-1$.
\item[$(2_2)$] Join a partial flow tree of type $(2_2)$ to a rigid flow sub-tree $\Gamma'$ at its special puncture. Here $\dim(\Gamma)=\dim(\Gamma')-n+1$.
\item[$(3_0)$] Join two rigid flow sub-trees $\Gamma'$ and $\Gamma''$ at their special puncture to two of the special punctures of a partial tree of type $(3_0)$. Here $\dim(\Gamma)=\dim(\Gamma')+\dim(\Gamma'')-n+1$.
\item[$(3_1)$] Join two rigid flow sub-trees $\Gamma'$ and $\Gamma''$ at their special puncture to two of the special punctures of a partial tree of type $(3_1)$. Here $\dim(\Gamma)=\dim(\Gamma')+\dim(\Gamma'')-n$.
\end{itemize}
\end{lma}

\begin{pf}
Let $\Gamma$ be any rigid flow sub-tree and let $s$ be its special puncture. If the vertex or marked point in the strip region following $s$ is an edge point or a switch then cut $\Gamma$ at the marked point following this point. If there is no such vertex or switch following $s$ then cut $\Gamma$ at the two marked points closest to the $3$-valent vertex of $\Gamma$ adjacent to $s$. The statement on dimensions are immediate from the definition of $\dim(\Gamma)$.
\end{pf}

\subsubsection{Uniform invertibility}\label{ss:unifinv}
Let $\Gamma$ be a rigid flow sub-tree and consider the operator
$$
\bar\pa\colon V_{2,\delta}(\Gamma,\lambda)\oplus V_{\rm sol}(\Gamma,\lambda)\oplus V_{\rm con}(\Gamma,\lambda)\to W_{1,\delta}(\Gamma,\lambda).
$$
Consider the end $[0,\infty)\times[0,1]$ in $\bar\Delta_{p',m'}(\Gamma,\lambda)$ corresponding to a special puncture $s$. Note that any element $v\in V_{2,\delta}(\Gamma,\lambda)$ has a well defined limit $v(s)=\lim_{z\to\infty} v(z)$, $z\in[0,\infty)\times[0,1]$. We call the map ${\rm ev}_s\colon v\mapsto v(s)$ {\em evaluation at $s$}. The restriction of ${\rm ev}_s$ to $\krn(\bar\pa)$ gives a map ${\rm ev}_s\colon \krn(\bar\pa)\to\R^n\approx T_s M$. Let $K_s(\Gamma,\lambda)$ denote the image of this map.

Our transversality assumptions imply that the rigid flow sub-tree $\Gamma$ is a smooth point of a $\dim(\Gamma)$-dimensional family $\Omega$ of partial flow trees. There is an obvious evaluation map $T_\Gamma\Omega\to T_s M$, we call the image of that map the {\em flow tree tangent space of $\Gamma$ at $s$} and denote it $F_s(\Gamma)$, see Remark \ref{r:dimrigidsubtree}. Let
$F_s'(\Gamma)\subset T^\ast_s M$ be the subspace of linear functionals $\la v,\bullet\ra$, $v\in F_s(\Gamma)$ and note that $F_s'(\Gamma)$ can be viewed as a vanishing condition at $s$. We write $V_{F'_s(\Gamma)}(\Gamma,\lambda)\subset V_{2,\delta}(\Gamma,\lambda)$ for the closed subspace of all elements which satisfy this vanishing condition.

\begin{prp}\label{p:rigftpart}
If $\Gamma$ is a rigid flow sub-tree then the index of
$$
\bar\pa\colon V_{2,\delta}(\Gamma,\lambda)\oplus V_{\rm sol}(\Gamma,\lambda)\oplus V_{\rm con}(\Gamma,\lambda)\to W_{1,\delta}(\Gamma,\lambda)
$$
equals $\dim(\Gamma)$. Moreover $K_s(\Gamma,\lambda)$ converges to $F_s(\Gamma)$ as $\lambda\to 0$. In particular, there exists a constant $C_\Gamma$ such that for all sufficiently small $\lambda>0$
\begin{equation}\label{e:estsubtree}
\|v\|\le C_\Gamma\|\bar\pa v\|,\quad v\in
V_{F_s'(\Gamma)}(\Gamma,\lambda).
\end{equation}
\end{prp}

\begin{pf}
It is a consequence of Proposition \ref{p:tfdim=dfdim} that this Fredholm operator has index $\dim(\Gamma)$. (The fact that there are conformal variations supported at each marked point, instead of as in a rigid tree at all but one, is compensated in the dimension formula by the addition of $1$ for the special puncture.)

Note that $\Gamma$ is obtained from other rigid flow sub-trees by the operations described in Lemma \ref{l:build}. Moreover, the cutting procedure in the proof of Lemma \ref{l:build} can be applied repeatedly starting with any rigid flow sub-tree until no further splitting off of partial trees in the list is possible. The only rigid flow sub-trees which cannot be cut any further are the following. Rigid flow sub-trees with one special puncture and one puncture and rigid flow subtrees with one special puncture and one end. We call such pieces {\em minimal} rigid flow sub-trees. Note that Lemmas \ref{l:Iw} and \ref{l:Iw+} imply that Proposition \ref{p:rigftpart} holds for minimal rigid flow sub-trees of the former and latter types respectively. Thus, by induction the proposition follows from the following claim.

\begin{clm}
If Proposition \ref{p:rigftpart} holds for a rigid flow sub-tree $\Gamma'$ (for two rigid flow sub-trees $\Gamma'$ and $\Gamma''$) then it holds also for the rigid flow sub-tree $\Gamma$ obtained from $\Gamma'$ (from $\Gamma'$ and $\Gamma''$) by any  of the operations described in Lemma \ref{l:build}.
\end{clm}

We verify this claim by checking all cases. Consider first the case when $\Gamma=\Gamma'\cup\Gamma_0$ where $\Gamma_0$ has type $(2_0)$. Let $p$ denote the marked point where $\Gamma'$ and $\Gamma_0$ are joined and let $q$ denote the other marked point of $\Gamma_0$. We decompose cut-off solutions and conformal variations of $\Gamma$ into those supported near $p$ and $q$ and those supported far from these points:
\begin{align*}
&V_{\rm sol}(\Gamma,\lambda)\oplus V_{\rm con}(\Gamma,\lambda)=\\
&\bigl(V'_{\rm sol}(\Gamma,\lambda)
\oplus V'_{\rm con}(\Gamma,\lambda)\bigr)
\oplus \bigl(V_{\rm sol}(p,\lambda)\oplus V_{\rm con}(p,\lambda)\bigr)
\oplus V_{\rm sol}(q,\lambda).
\end{align*}
Assume that \eqref{e:estsubtree} does not hold. Then there exists a sequence of functions $v_\lambda\in V_{F_q'(\Gamma)}(\Gamma,\lambda)$, $\lambda\to 0$, such that
\begin{align}\label{e:20contr1}
&\|v_\lambda\|=1,\\\label{e:20contr2}
&\|\bar\pa v_\lambda\|_{W_{1,\delta}(\Gamma,\lambda)}\to 0,
\quad\text{as }\lambda\to 0.
\end{align}
Write
\begin{align*}
v_\lambda&=v_\lambda^0 + v_\lambda^p + v_\lambda^q,\\
v_\lambda^0 &\in V_{2,\delta}(\Gamma,\lambda)
\oplus V'_{\rm sol}(\Gamma,\lambda)\oplus V'_{\rm con}(\Gamma,\lambda),\\
v_\lambda^p &\in V_{\rm sol}(p,\lambda)\oplus V_{\rm con}(p,\lambda),\\
v_\lambda^q &\in V_{\rm sol}(q,\lambda).
\end{align*}

Let $[T_1,T_2]\times[0,1]\subset \bar\Delta_{p',m'}(\Gamma,\lambda)$ be the region which contains $p=(0,0)$.
Let $\beta[-]\colon \bar\Delta_{p',m'}(\Gamma,\lambda)\to \C$ be a cut-off function which is real valued and holomorphic on the boundary and with the following properties: $\beta[-]=1$ in the region of $\bar\Delta_{p',m'}(\Gamma,\lambda)$ to the left of the vertical segment $\{\frac12 T_1\}\times[0,1]$, $\beta[-]=0$ to the right of $\{0\}\times[0,1]$, and $|D^k\beta[-]|=\Ordo(\lambda)$, $k=1,2$. (Recall that $T_j$ grows like $\lambda^{-1}$ as $\lambda\to 0$.) Similarly, let $\beta[+]\colon\bar\Delta_{p',m'}(\Gamma,\lambda)\to\C$ be a cut off function equal to $1$ to the right of $\{\frac12 T_2\}\times[0,1]$ and equal to $0$ to the left of $\{0\}\times[0,1]$ with $|D^k\beta[+]|=\Ordo(\lambda)$.

Let $v_\lambda^p[+]\in V_{\rm sol}(p,\lambda)\subset V_{\rm sol}(\Gamma',\lambda)$ be the cut off constant function which agrees with $v_\lambda^p$ at $p$ and let $v_\lambda^p[-]\in V_{\rm sol}(\Gamma_0,\lambda)$ be the cut off constant function which agrees with $v_\lambda^p$ at $p$. Note that the function $\beta[-]v_\lambda^0+v_\lambda^p[-]+v_\lambda^q$ lies in $V_{F_q'(\Gamma)}(\Gamma,\lambda)\subset V_{2,\delta}(\Gamma_0,\lambda)\oplus V_{\rm sol}(\Gamma_0,\lambda)$ and that the function $\beta[+]v_\lambda^0+v_\lambda^p[+]$ lies in $V_{2,\delta}(\Gamma',\lambda)\oplus V_{\rm sol}(\Gamma',\lambda)\oplus V_{\rm con}(\Gamma',\lambda)$. Moreover,
\begin{align*}
&\bar\pa(\beta[-]v_\lambda^0+v_\lambda^p[-]+v_\lambda^q)
=(\bar\pa\beta[-])v_\lambda^0+\beta[-]\bar\pa v_\lambda,\\
&\bar\pa(\beta[+]v_\lambda^0+v_\lambda^p[+])
=(\bar\pa\beta[+])v_\lambda^0+\beta[+]\bar\pa v_\lambda.
\end{align*}
Thus
\begin{align}\label{e:left20}
&\|\bar\pa(\beta[-]v_\lambda^0+v_\lambda^p[-]+v_\lambda^q)\|_{ W_{1,\delta}(\Gamma_0,\lambda)}=\ordo(1),\\\label{e:right20}
&\|\bar\pa(\beta[+]v_\lambda^0+v_\lambda^p[+])\|_{W_{1,\delta}(\Gamma',\lambda)}
=\ordo(1).
\end{align}

Identify $\Delta_{2,1}(\Gamma_0,\lambda)$ with $\R\times[0,1]$ where $q$ and $p$ corresponds to $-\infty$ and $+\infty$, respectively. It follows from Lemma \ref{l:Is} that the $\bar\pa$-operator on the space of functions $H_A\subset H_{2,\delta}\oplus F_{\rm sol}$, where $F_{\rm sol}=F_{\rm sol}(q)\oplus F_{\rm sol}(p)$, $F(q)=\R^n$, $F(p)=F_p(\Gamma')$, and where the vanishing condition $A$ is given by $F_q'(\Gamma)$ at $q$ is an isomorphism. Thus there exists a function $\eta_\lambda\in H_A$ such that
$$
\bar\pa\eta_\lambda=\bar\pa(\beta[-]v_\lambda^0+v_\lambda^p[-]+v_\lambda^q)
$$
and $\|\eta_\lambda\|_{H_{2,\delta}\oplus F_{\rm sol}}\le C\|\bar\pa(\beta[-]v_\lambda^0+v_\lambda^p[-]+v_\lambda^q)\|_{H_{1,\delta}}$. In particular,
$$
\bar\pa\bigl(\beta[-]v_\lambda^0+v_\lambda^p[-]+v_\lambda^q-\eta_\lambda\bigr)=0,
$$
and hence $\beta[-]v_\lambda+v_\lambda^p[-]+v_\lambda^q-\eta_\lambda$ lies in the kernel of the $\bar\pa$-operator on $H_{2,\delta}\oplus F_{\rm sol}$. Moreover the function satisfies the vanishing condition $F_q'(\Gamma)$ at $q$. Decomposing $\eta_\lambda$ according to the direct sum decomposition of $H_{2,\delta}\oplus F_{\rm sol}(q)\oplus F_{\rm sol}(p)$ we get $\eta_\lambda=\eta_\lambda^0+\eta_\lambda^p+\eta_\lambda^q$. It follows form the characterization of the kernel in Lemma \ref{l:Is} that
$$
\ev_{p}\bigl(\beta[-]v_\lambda^0+v_\lambda^p[-]+v_\lambda^q-\eta_\lambda\bigr)\in F_{q}(\Gamma)^{\perp}.
$$
Thus,
\begin{equation}\label{e:20vpleft}
\ev_p(v_\lambda^0+v_\lambda^p+v_\lambda^q)+\ordo(1)\in F_p(\Gamma')^\perp.
\end{equation}
(There are two contributions to the error term $\ordo(1)$: $\eta_\lambda$ of size $\ordo(1)$ and the difference between $F_q(\Gamma)^\perp$ and $F_p(\Gamma')^\perp$ which arises as the tree tangent space is transported by the flow through the region around the edge point between $p$ and $q$ of size $\Ordo(\lambda)$.)

By assumption the $\bar\pa$-operator on $V_{F'_p(\Gamma')}\subset V_{2,\delta}(\Gamma',\lambda)\oplus V_{\rm sol}(\Gamma',\lambda)\oplus V_{\rm con}(\Gamma',\lambda)$ is an isomorphism. Hence there exists $\xi_\lambda\in V_{F_p'(\Gamma')}$ such that $\bar\pa\xi_\lambda=\bar\pa( \beta[+]v_\lambda^0+v_\lambda^p[+])$. We conclude that $\bar\pa\bigl(\beta[+]v_\lambda^0 + v_\lambda^p[+]-\xi_\lambda\bigr)=0$. Arguing as above and using the inductive assumption about the kernel of the $\bar\pa$-operator on $\Delta_{p',m'}(\Gamma',\lambda)$ we conclude that
\begin{equation}\label{e:20vpright}
\ev_p(v_\lambda^0+v_\lambda^p)+\ordo(1)\in F_p(\Gamma').
\end{equation}
Equations \eqref{e:20vpleft} and \eqref{e:20vpright} imply that $v_\lambda^p=\ordo(1)$. In particular,
\begin{align*}
&\|\bar\pa(\beta[-]v_\lambda^0 + v_\lambda^q[-])\|=\ordo(1),\\
&\|\bar\pa(\beta[+]v_\lambda^0)\|=\ordo(1),
\end{align*}
and we conclude from the estimates for the $\bar\pa$-operator on $H_{2,\delta}\oplus F_{\rm sol}(q)$ and $V_{F'_p(\Gamma')}$, respectively, that
$\|\beta[-]v_\lambda^0\|\to 0$ and $v_\lambda^q\to 0$, and that $\|\beta[+]v_\lambda^0\|\to 0$.

Finally, let $\gamma\colon[T_1,T_2]\times[0,1]\to \C$ be a cut off function which is equal to $1$ on $[\frac{3}{4} T_1, \frac{3}{4} T_2]\times[0,1]$, equal to $0$ outside $[\frac{7}{8} T_1,\frac{7}{8} T_2]\times[0,1]$ and with $|D^k\gamma|=\Ordo(\lambda)$, $k=1,2$. Using Lemma \ref{l:Iw} with negative exponential weights at both ends we see that the $\bar\pa$-operator restricted to the $n$-codimensional subspace of functions vanishing at $0$ is an isomorphism. By definition of $V_{2,\delta}(\Gamma,\lambda)$, $v_\lambda^0(p)=0$ and we conclude from this and $\|\bar\pa(\gamma v_\lambda^0)\|\to 0$ that $\|\gamma v_\lambda^0\|\to 0$. To summarize, \eqref{e:20contr1} and \eqref{e:20contr2} imply that
$$
\|v_\lambda^0\|\le \|\beta[-]v_\lambda^0\|+\|\beta[+]v_\lambda^0\|+\|\gamma v_\lambda^0\|\to 0
$$
and that
$$
\|v_\lambda^p\|\to 0,\quad\quad \|v_\lambda^q\|\to 0.
$$
This however contradicts \eqref{e:20contr1} and we conclude that the estimate holds.

It remains to prove the statement on the kernel of $\Gamma$. This is straightforward after the estimate is established. Pick a function
$v_\lambda'=v_\lambda^0+v^p_\lambda[+]$ in the kernel of $\bar\pa$ on $V_{2,\delta}(\Gamma',\lambda)$. Continue the function $v^p_\lambda[+]$ constantly to the left in $\bar\Delta_{p,m}(\Gamma,\lambda)$. Expressed in terms of the direct sum decomposition used above the continuation is
$$
v_\lambda^p+\sigma_\lambda+v_\lambda^q,
$$
where $\sigma_\lambda$ has finite support near the edge point between $p$ and $q$. Then
$$
\bar\pa(\beta[+]v_\lambda^0 + v_\lambda^p + \sigma_\lambda +v_\lambda^q)=(\bar\pa\beta[+]) v_\lambda^0.
$$
By the estimate just established and the computation of the Fredholm index of $\bar\pa$ on $\bar\Delta_{p,m}(\Gamma,\lambda)$ we find that there exists a function $\rho_\lambda$ satisfying the vanishing condition $F'_q(\Gamma)$ at $q$ such that $\bar\pa \rho_\lambda = \bar\pa\beta[+] v_\lambda^0$ and $\|\rho_\lambda\|=\Ordo(\lambda)$. The kernel is thus spanned by functions of the from $\beta[+] v_\lambda^0 + v_\lambda^p+\sigma_\lambda +v_\lambda^q-\rho_\lambda$. The convergence statement for the kernel follows.

The proofs in the other cases are similar so we just point out the differences from the case just given. In the case $(2_1)$ of Lemma \ref{l:build} we have $\Gamma=\Gamma'\cup\Gamma_0$ where $\Gamma_0$ is of type $(2_1)$. In this case we proceed as above. The main difference from the case considered is that the kernel of the $\bar\pa$-problem on $\Gamma_0$ is smaller. In particular, the analogue of \eqref{e:left20} in this case is
\begin{equation*}
\ev_p(v_\lambda+v_\lambda^p+v_\lambda^q)+\ordo(1)\in (F_p(\Gamma')+W)^\perp,
\end{equation*}
where $W\otimes\C$ is the complex line in which the boundary condition rotates, see Lemma \ref{l:Is-}. The estimate is then established using Lemma \ref{l:Is-} and the inductive hypothesis. Finally, solutions on $\bar\Delta_{p',m'}(\Gamma',\lambda)$ with values in $W^\perp$ can be continued as above.

In the case $(2_2)$ of Lemma \ref{l:build} we have $\Gamma=\Gamma'\cup \Gamma_0$, where $\Gamma_0$ is of type $(2_2)$ with $2$-valent puncture $r$ and special punctures $p$, the special puncture of $\Gamma'$, and $q$. In this case replace the problem on $[T_1,T_2]\times[0,1]$ considered above by the problem on a $3$-punctured disk which is the union of two regions of the form $[T_1,T_2]\times[0,1]$, the region $B(r)$, and $[0,\infty)\times[0,1]$ corresponding to the $2$-valent puncture $r$. Note that the conformal variation associated to $v$ is given by $\pa s\cdot \bar\pa(\beta f)$, where $f=\lambda^{-1}b\pa_\tau$ in $[0,\infty)\times[0,1]$. In particular, the conformal variation behaves close to the limit (when $\theta_\lambda\ll\delta$) just like a cut off constant solution with value $0$ at the $2$-valent puncture and common value at the other two punctures and Lemma \ref{l:Is-} gives the necessary information on the $\bar\pa$-operator in this case. To establish the estimate we again argue by contradiction and write
$$
w_\lambda= v_\lambda + v_\lambda^p + v_\lambda^q + v_\lambda^r.
$$
We conclude from Lemma \ref{l:IIIs2} and the vanishing condition at $q$ that
$$
\ev_p(w_\lambda)=\ordo(1).
$$
This suffices to establish the desired estimate. The kernel function can be glued as above.

In case $(3_0)$ of Lemma \ref{l:build} we have $\Gamma=\Gamma'\cup\Gamma_0\cup\Gamma''$, where $\Gamma_0$ is a tree of type $(3_0)$ with special vertices $p$, $q$, and $r$, where $p$ and $q$ are the special vertices of $\Gamma'$ and $\Gamma''$, respectively. To establish the estimate we again argue by contradiction writing
$$
w_\lambda= v_\lambda + v_\lambda^p + v_\lambda^q + v_\lambda^r.
$$
Lemma \ref{l:IIIs} allows us to conclude from the vanishing condition at $r$ that
$$
\ev_{p}(w_\lambda)+\ordo(1)=\ev_q(w_\lambda)+\ordo(1)\in \bigl(F_p(\Gamma')\cap F_q(\Gamma'')\bigr)^\perp.
$$
On the other hand the inductive hypothesis applied to the $\bar\pa$-operator on $\Delta_{\Gamma'}(\lambda)$ and $\Delta_{\Gamma''}(\lambda)$ shows that
$$
\ev_p(w_\lambda)+\ordo(1)\in F_p(\Gamma')
$$
and
$$
\ev_q(w_\lambda)+\ordo(1)\in F_q(\Gamma'').
$$
These together with the above allows us to conclude $v_\lambda^p\to 0$ and $v_\lambda^q\to 0$. The estimate then follows using the same argument as in the $(2_0)$-case. Also, gluing solutions works much as before.

Finally, the case $(3_1)$ of Lemma \ref{l:build} is obtained from a modification of the $(3_0)$-case which is entirely analogous to the modification of the $(2_0)$-case needed for the $(2_1)$-case. The information on the $\bar\pa$-operator needed follows from Lemma \ref{l:IIIs-}.
\end{pf}

With Proposition \ref{p:rigftpart} established we can now prove the uniform invertibility for rigid flow trees.

\begin{prp}\label{p:rightinv}
Let $\Gamma$ be a rigid flow tree. Then for all sufficiently small $\lambda>0$, the linearization of the $\bar\pa_J$-operator at $u_\lambda^{(0,0)}$
$$
L\bar\pa_J\colon V_{2,\delta}^{(0,0)}(\Gamma,\lambda)\oplus V_{\rm sol}^{(0,0)}(\Gamma,\lambda)\oplus
V_{\rm con}^{(0,0)}(\Gamma,\lambda)\to W_{1,\delta}^{(0,0)}(\Gamma,\lambda)
$$
has a right inverse $Q_\lambda$ and there exists a constant $C$ such that for all
$\lambda>0$ sufficiently small
$$
\|Q_\lambda\|\le C.
$$
\end{prp}

\begin{pf}
It follows by Proposition \ref{p:tfdim=dfdim} that the Fredholm index of $L\bar\pa_J$ equals $0$. Therefore the proposition follows once we establish the estimate
\begin{equation}\label{e:estrinv}
\|v_\lambda\|\le C\|\bar\pa v_\lambda\|_{W_{1,\delta}^{(0,0)}(\Gamma,\lambda)}.
\end{equation}
In order to establish \eqref{e:estrinv} cut the tree at the marked point $p$ where no conformal variation is supported. This produces two rigid flow sub-trees $\Gamma'$ and $\Gamma''$. We now argue by contradiction. Assume that \eqref{e:estrinv} does not hold then there exists a sequence of functions $v_\lambda\in V_{2,\delta}^{(0,0)}(\Gamma,\lambda)\oplus V_{\rm con}^{(0,0)}(\Gamma,\lambda)\oplus V_{\rm sol}^{(0,0)}(\Gamma,\lambda)$ such that
\begin{align}\label{e:rinvcontr1}
&\|v_\lambda\|=1,\\\label{e:rinvcontr2}
&\|\bar\pa v_\lambda\|_{W_{1,\delta}^{(0,0)}(\Gamma,\lambda)}\to 0,
\quad\text{as }\lambda\to 0.
\end{align}
Write
$$
v_\lambda=v_\lambda^0+v^p_\lambda,
$$
where $v_\lambda^0\in W_{2,\delta}^{(0,0)}(\Gamma,\lambda)\oplus V'_{\rm sol}(\Gamma,\lambda)\oplus V'_{\rm con}(\Gamma,\lambda)$ and $v_\lambda^p\in V_{\rm sol}(p,\lambda)$. We chose cut-off functions $\beta[+]$ and $\beta[-]$ as in the proof of Proposition \ref{p:rigftpart}. It is a consequence of Lemma \ref{l:Is} that there exists $\eta_\lambda[+]$ and $\eta_\lambda[-]$ which satisfy the vanishing conditions $F_p'(\Gamma')$ and $F_p'(\Gamma'')$ respectively such that $\bar\pa\eta[\pm]=\bar\pa(\beta[\pm]v_\lambda^0+v^p[\pm])$. We conclude that
$\bar\pa\bigl(\beta[\pm]v_\lambda^0 + v_\lambda^p[\pm] -\eta_\lambda[\pm]\bigr)=0$. It follows from Proposition \ref{p:rigftpart} that
$$
\ev_{p}(\beta[\pm]v_\lambda+v_\lambda[\pm]-\eta_\lambda[\pm])+\ordo(1)\in F_p(\Gamma[\pm]).
$$
Since $F_p(\Gamma[+])\cap F_p(\Gamma[-])=l$ and we have a vanishing condition in direction $l$ we conclude that $v_\lambda^p\to 0$. Proposition \ref{p:rigftpart} then implies that $\beta[\pm]v_\lambda^0\to 0$. Finally, we pick a cut-off function $\gamma$ as in the proof of Proposition \ref{p:rigftpart} and conclude as there that $\|\gamma v_\lambda^0\|\to 0$. This contradicts \eqref{e:rinvcontr1} and we conclude that \eqref{e:estrinv} holds for all sufficiently small $\lambda>0$.
\end{pf}

\subsubsection{Existence and surjectivity}\label{ss:existsurj}
With the uniform invertibility established we are in position to prove Theorem \ref{t:treetodisk}.

\begin{pf}[Proof of Theorem \ref{t:treetodisk}]
After Proposition \ref{p:rightinv} and Remark \ref{r:approxsol}, in order to apply Lemma \ref{l:FloerPicard} to the $\bar\pa_J$-operator, we must only establish the quadratic estimate
\begin{align*}
&\|N(u)-N(v)\|\le C(\|u-v\|)(\|u\|+\|v\|),\\
& u,v\in V_{2,\delta}^{(0,0)}(\Gamma,\lambda)\oplus V_{\rm sol}^{(0,0)}(\Gamma,\lambda)\oplus V_{\rm con}^{(0,0)}(\Gamma,\lambda),
\end{align*}
see \eqref{e:quadest}. The proof of this quadratic estimate is similar to Lemma 7.16 in \cite{EES2} (see also Subsection 4.3 in \cite{EES3}). The new ingredient in the present situation is the space of cut-off solutions and the action of the conformal variations on the function $u_\lambda^{(0,0)}$. The cut-off solutions and the action just mentioned has the following property: altering $u_\lambda^{(0,0)}$ by adding a cut-off solution or by an element in $\B_{p,m}$, we obtain a function $\tilde u_\lambda$ which is $J$-holomorphic except possibly in regions of the form $B(v)$ or $R(p)$, where $v$ is a $3$- or $2$-valent vertex of $\Gamma$ and where $p$ is an edge point. Hence, the non-linear term corresponding to cut off solutions or actions by elements in $\B_{p,m}$ is supported only in the uniformly finite neighborhoods $R(p)$ or $B(v)$ in $\Delta_{p,m}^{(0,0)}(\Gamma,\lambda)$ where the weight function is small. This together with the fact that the diffeomorphisms corresponding to elements in $\B_{p,m}$, used to transport elements in $V_{2,\delta}^{(0,0)}(\Gamma,\lambda)$, satisfies \eqref{e:derssmall} and the fact that the trivialization compensates for the change in weight, see \eqref{e:wigthcomp}, allow us to establish the quadratic estimate by repeating the argument in \cite{EES2}.

This quadratic estimate in combination with Proposition \ref{p:rightinv} and Lemma \ref{l:FloerPicard} imply that for each sufficiently small $\lambda>0$ and each rigid flow tree $\Gamma$, there exists a unique $J$-holomorphic disk in a $\delta$-neighborhood of $u_\lambda^{(0,0)}$ in $\W_{2,\delta}(\Gamma,\lambda)$ for some $\delta>0$. To finish the proof it remains to show that if $w_\lambda$ are any $J$-holomorphic disks with image in some small neighborhood of $\Gamma$ then for $\lambda>0$ small enough $w_\lambda$ is in the image of the exponential map of an element in the $\delta$-neighborhood.

Let $w_\lambda\colon D_m\to T^\ast M$ be a holomorphic disk in a $\delta_0$-neighborhood of $\Gamma$, where $\delta_0>0$ is so small that our metric has the standard form in this neighborhood. It is a consequence of the proof of Theorem \ref{t:disktotree} that for small $\lambda>0$, $w_\lambda$ takes each edge of its domain to an $\Ordo(\lambda)$-neighborhood of a gradient flow line, and that it takes $\log(\lambda^{-1})$-neighborhoods of its boundary minima and tangency points to regions of $\Ordo(\lambda\log(\lambda^{-1}))$-neighborhoods of points. It follows from our transversality assumptions for trees that the image of $w_\lambda$ lies in an $\Ordo(\lambda\log(\lambda^{-1}))$-neighborhood of a flow tree. Moreover, since $w_\lambda$ is a rigid disk, the flow tree close to $w_\lambda$ is rigid as well. Since the only rigid tree in a small finite neighborhood of $\Gamma$ is $\Gamma$ itself it follows that $w_\lambda$ maps to an $\Ordo(\lambda\log(\lambda^{-1}))$ neighborhood of $\Gamma$.

Recall, see Remark \ref{r:domaincontrolI}, that the conformal structure on the domain of $w_\lambda$ is determined up to $\Ordo(\log(\lambda^{-1}))$ by the rigid limit tree $\Gamma$. It follows from this that the conformal structure on $\Delta_{p,m}^{(0,0)}(\Gamma,\lambda)$ is at distance at most $\Ordo(\log(\lambda^{-1}))$ from that of the domain of $w_\lambda$ (where we mark the points where $w_\lambda$ intersects $D_p$). To see this we start at the positive puncture: after suitable translations the two maps agree to within $\Ordo(\lambda\log(\lambda^{-1}))$ at the first marked point. After that we use the fact that both maps stay at $\Ordo(\lambda)$-distance from the double lift of gradient lines which are separated by at most $\Ordo(\lambda\log(\lambda^{-1}))$ (arising from variations near vertices and edge points) to conclude that the marked points in the domain are of distance bounded by $\Ordo(\log(\lambda^{-1}))$.

There thus exists $(\kappa,\nu)$ with $|(\kappa,\nu)|=\Ordo(\lambda\log(\lambda^{-1}))$ such that the domain of $w_\lambda$ can be taken as $\Delta_{p,m}^{(\kappa,\nu)}(\Gamma,\lambda)$. We must estimate the distance between the map $u_\lambda^{(\kappa,\nu)}$ and $w_\lambda$. To do this we use the uniqueness up to translation of the local solutions for estimates in the bounded regions $B(v)$ near $Y_0$- and $Y_1$-vertices, switches, and $2$-valent punctures. In regions complementary to such regions and to edge point regions, where the weight in $V_{2,\delta}^{(\kappa,\nu)}(\Gamma,\lambda)$ is large, we take advantage of the simple boundary conditions and use Fourier analysis to derive the estimates needed. In the remaining edge point regions the weight is small and it is straightforward to estimate the distance using flow tree convergence as in Lemma \ref{l:gradconv}.

Consider first a neighborhood of a boundary minimum corresponding to a $Y_0$-vertex $v$ of $\Gamma$. More precisely we consider such a neighborhood of the form $\Delta_3(T_1,T_2,T_3)$ which is bounded by vertical segments near the nearest edge point regions.  It follows by the argument given above that $w_\lambda\colon\Delta_3(T_1,T_2,T_3)\to T^\ast M$ is a holomorphic map mapping into an $\Ordo(\lambda\log(\lambda^{-1}))$ neighborhood of the image of the local solution at $v$. Moreover, as in Section \ref{S:disktotree}, $|Dw_\lambda|=\Ordo(\lambda)$. It follows by Lemma \ref{l:gradconv} that outside a $\log(\lambda^{-1})$-neighborhood of the boundary minimum, $w_\lambda$ converges at rate $\Ordo(\lambda)$ to a a flow tree.  In particular, the re-scaling $\lambda^{-1} w_\lambda$ converges to a holomorphic map $w_\infty$ on $\Delta_3$ with the same asymptotics as the local solution. Let $c$ be a constant vector in $\R^n$ such that $\lambda^{-1} w_\lambda(p)+c=s(p)$ for some $p\in\Delta_3(T_1,T_2,T_3)$ then also $\lambda^{-1}w_\lambda +c$ converges to a holomorphic map which, by uniqueness of the local solution up to translation in $\R^n$ (see \S\ref{ss:solY0}), must be equal to $s$. Moreover, since the distance from $w_\lambda$ to $\Gamma$ is $\Ordo(\lambda\log(\lambda^{-1}))$ it follows that for sufficiently small $\lambda>0$
$$
w_\lambda +\lambda c= s_\lambda + \eta_\lambda^0,
$$
where $\lambda c=\Ordo(\lambda\log(\lambda^{-1}))$, where $|D^k\eta_\lambda|=\Ordo(\lambda)$, $k=0,1,2$.

Consider next the difference $s_\lambda-(w_\lambda+\lambda c)$ restricted to a strip region $[\tau_1,\tau_2]\times[0,1]\subset\Delta_3(T_1,T_2,T_3)$ where $\{\tau_1\}\times[0,1]$ lies slightly inside $B(v)$ and where $\{\tau_2\}\times[0,1]$ lies near an edge point. The fact that $w_\lambda+\lambda c$ converges to $s_\lambda$ at rate $\Ordo(\lambda)$ and the argument in Lemma \ref{l:gradconv} (see also Remark \ref{r:gradconv}), imply that there are $\Ordo(\lambda\log(\lambda^{-1}))$ $C^2$-bounds for $w_\lambda+\lambda c -s_\lambda$ near both vertical lines bounding the region. For simplicity of notation we assume that $-\tau_1=\tau_2=\tau_0$. Fourier analysis implies that
$$
(w_\lambda+c\lambda)-s_\lambda =\sum_{n\in\Z} c_n e^{n\pi z}.
$$
Denoting this function $(x(\tau+it),y(\tau+it))$, the $\Ordo(\lambda\log(\lambda^{-1}))$ $C^0$-bound at the ends implies that
$$
\int_0^1\Bigl\la y(\pm\tau_0,t)-y(\pm\tau_0,0),
y(\pm\tau_0,t)-y(\pm\tau_0,0)\Bigr\ra\,dt
=\Ordo(\lambda^2\log(\lambda^{-1})^2),
$$
However,
$$
\int_0^1\Bigl\la y(\pm\tau_0,t)-y(\pm\tau_0,0), y(\pm\tau_0,t)-y(\tau_0,0)\Bigr\ra\,dt= c_0^2
\sum_{n\ne 0} (2\pi n)^{-1} c_n^2 (e^{2\pi n\tau_0}-e^{-2\pi n\tau_0}),
$$
and
$$
\int_{-\tau_0}^{\tau_0}\int_0^1 (\sum_{n\ne 0}c_n e^{\pi n z})^2e^{2\delta|\tau_0-\tau|}\,d\tau dt=e^{2\delta\tau_0}\sum_{n\ne 0} (2\pi n-2\delta)^{-1}c_n^2(e^{(2\pi n-2\delta)\tau_0}-e^{(-2\pi n-2\delta)\tau_0}).
$$
Therefore the weighted $L^2$-norm of the non-constant terms in the Fourier expansion is small. Applying a similar argument using the bounds on the higher derivatives at the ends we conclude that $w_\lambda=u_\lambda^{(\kappa,\nu)}+v_\lambda^0+v_\lambda^1$ in $\Delta_3(T_1,T_2,T_3)$ where the norm of $(v_\lambda^0,v_\lambda^1)\in V_{2,\delta}^{(\kappa,\nu)}(\Gamma,\lambda)\oplus V_{\rm sol}^{(\kappa,\nu)}(\Gamma,\lambda)\oplus V_{\rm con}^{(\kappa,\nu)}(\Gamma,\lambda)$ is $\Ordo(\lambda\log(\lambda^{-1}))$. Here the function $v_\lambda^1$ is determined by the constant Fourier components above and by the difference between $s_\lambda$ and $u_\lambda^{(\kappa,\nu)}$. Its norm is therefore $\Ordo(\lambda\log(\lambda^{-1}))$. The function $v_\lambda^0$ is determined by the non-constant Fourier components which was estimated above and the "remainder" of the constant Fourier components supported in $B(v)$ and hence also its norm is $\Ordo(\lambda\log(\lambda^{-1}))$.

The distance between $w_\lambda$ and $u_\lambda^{(\kappa,\nu)}$ in subsets corresponding to $Y_1$-vertices, switches, and $2$-valent punctures bounded by vertical segments near nearby edge points can be bounded in much the same way. The main difference is that in these cases there appear restrictions on the translation vector $c$ above. (For $2$-valent punctures and switches it must be parallel to the cusp edge, for $2$-valent punctures it equals $0$.) We thus have a bound on the distance in such regions. The remaining parts of the disk are a finite number of strip regions. Here the distance can be bounded applying the argument used in the strip region considered above: Flow tree convergence implies a $C^2$-bound near the vertical segments bounding the region. Fourier analysis then gives bounds on the cut-off solutions  corresponding to the constant Fourier modes and also weighted $L^2$-estimates for higher Fourier modes corresponding to elements in $V_{2,\delta}^{(\kappa,\nu)}(\Gamma,\lambda)$. (Note that if the boundary conditions are as in \eqref{e:norform2} then we multiply the the first component of the difference $w_\lambda-s_\lambda$ by $e^{-i\theta_\lambda t}$ and then use the Fourier analysis argument above, where $\theta_\lambda\to 0$ is the angle between the two lines in the boundary condition.)

In regions of the form $[0,\infty)\times[0,1]$ corresponding to  $1$-valent punctures and ends, the local solutions in the limit are unique up to pre-composition with translation in the source and addition of cut-off solutions. By flow tree convergence, the maps $u_\lambda^{(\kappa,\nu)}$ and $w_\lambda$ maps points near the closest edge point at most $\Ordo(\lambda\log(\lambda^{-1}))$ apart, so we conclude  that $w_\lambda=u_\lambda^{(\kappa,\nu)}+v_\lambda^0+v_\lambda^1$ as above.

Finally, in a region $R(p)$ around an edge point $p$ the distance between the maps are at most $\Ordo(\lambda\log(\lambda^{-1}))$ since both are close to the same flow tree, see Lemma \ref{l:gradconv} and Remark \ref{r:gradconv}, and since the size of any change arising above has this size. We conclude that $w_\lambda=\Exp_{(\kappa,\nu)}(v,f)$ where $|(\kappa,\nu)|=\Ordo(\lambda\log(\lambda^{-1}))$ and where $\|v\|+\|f\|=\Ordo(\lambda\log(\lambda^{-1}))$. This finishes the proof.
\end{pf}

\begin{pf}[Proof of Theorem \ref{t:main} (b)]
The theorem follows from Theorems \ref{t:main} (a), \ref{t:disktotree}, and \ref{t:treetodisk}.
\end{pf}

%% file: Sec/7expl.tex
\section{An example}\label{S:expl}
\begin{figure}[htbp]
\begin{center}
\includegraphics[angle=0, width=6cm]{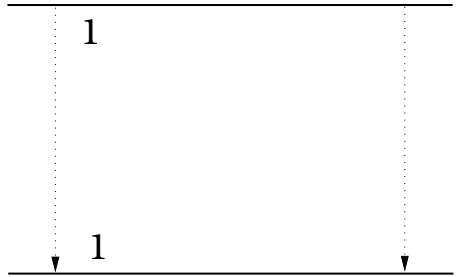}
\end{center}
\caption{The family of trees before isotopy. The flow line $1$ is determined by the sheets $F_1$ and $F_2$.}
\label{f:triv}
\end{figure}

In this section we consider an elementary and local example where we see all vertices of rigid flow trees appearing. Consider two fronts in $J^1(\R^2)$. The first one $F_1$ is simply the $0$-section and the second one $F_2$ is a sheet given by the function $z(x_1,x_2)=K-x_1$, where $K\gg 0$. We consider the parts of these fronts in a square $-R\le |x_j|\le R$, $j=1,2$. There is an obvious $1$-parameter family of holomorphic strips in $T^\ast\R^2=\C^2$ with boundary of the Lagrangian submanifolds corresponding to the two fronts, given by $(x_1+iy_1,x_2+iy_2)=(\tau+it, s)$, where $-R\le s\le R$. The corresponding family of gradient flow lines is as easily described, see Figure \ref{f:triv}.

\begin{figure}[htbp]
\begin{center}
\includegraphics[angle=0, width=6cm]{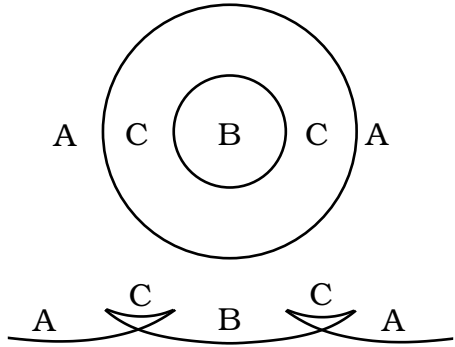}
\end{center}
\caption{A Legendrian isotopy of the $0$-section, top view and cross-section.}
\label{f:isot}
\end{figure}

We next change the $0$-section by a Legendrian isotopy into $F_1'$ as shown in Figure \ref{f:isot}. Label the pieces of the front $F_1'$ by $A$, $B$, and $C$ as shown in Figure \ref{f:isot}, and label the front of $F_2$ by $D$.

\begin{figure}[htbp]
\begin{center}
\includegraphics[angle=0, width=6cm]{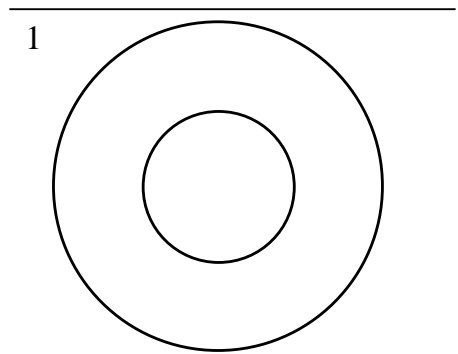}
\end{center}
\caption{Initially members of the family looks as before the isotopy; $1=A|D$.}
\label{f:tree1}
\end{figure}

\begin{figure}[htbp]
\begin{center}
\includegraphics[angle=0, width=6cm]{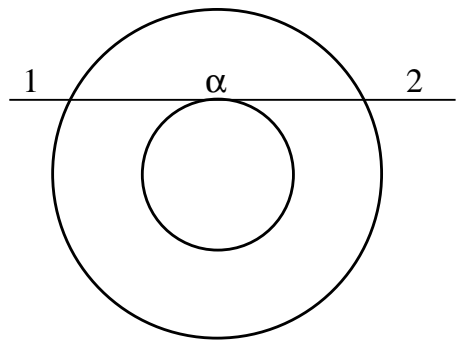}
\end{center}
\caption{The first critical moment; $1=2=A|D$; $\alpha=O$.}
\label{f:tree2}
\end{figure}

\begin{figure}[htbpp]
\begin{center}
\includegraphics[angle=0, width=6cm]{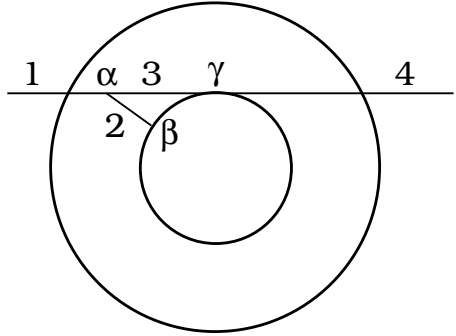}
\end{center}
\caption{After the first critical moment; $1=A|D$, $2=A|C$, $3=C|D$, $4=A|D$; $\alpha=Y_0$, $\beta=E$, $\gamma=S$.}
\label{f:tree3}
\end{figure}

\begin{figure}[htbp]
\begin{center}
\includegraphics[angle=0, width=6cm]{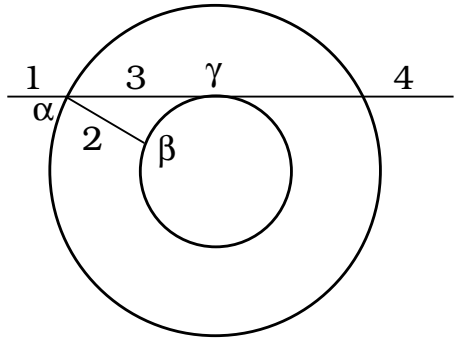}
\end{center}
\caption{The second critical moment; $1=A|D$, $2=A|C$, $3=C|D$, $4=A|D$; $\alpha=O$, $\beta=E$, $\gamma=S$.}
\label{f:tree4}
\end{figure}

\begin{figure}[htbp]
\begin{center}
\includegraphics[angle=0, width=6cm]{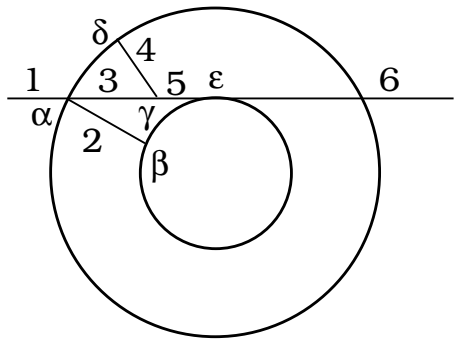}
\end{center}
\caption{After the second critical moment; $1=A|D$, $2=A|C$, $3=B|D$, $4=B|C$, $5=C|D$, $6=A|D$; $\alpha=Y_1$, $\beta=E$, $\gamma=Y_0$, $\delta=E$, $\epsilon=S$.}
\label{f:tree5}
\end{figure}

\begin{figure}[htbp]
\begin{center}
\includegraphics[angle=0, width=6cm]{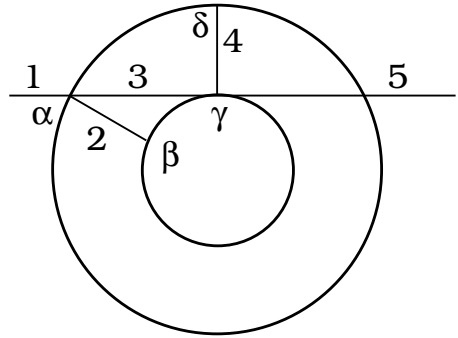}
\end{center}
\caption{The third critical moment; $1=A|D$, $2=A|C$, $3=B|D$, $4=B|C$, $5=A|D$; $\alpha=Y_1$, $\beta=E$, $\gamma=O$, $\delta=E$.}
\label{f:tree6}
\end{figure}

\begin{figure}[htbp]
\begin{center}
\includegraphics[angle=0, width=6cm]{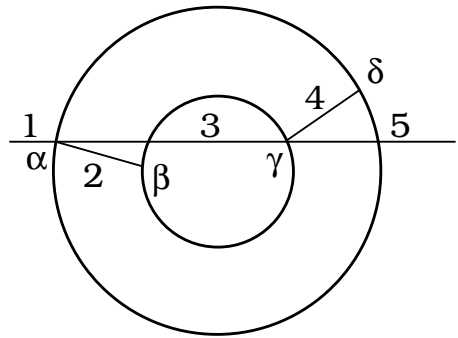}
\end{center}
\caption{After the third critical moment; $1=A|D$, $2=A|C$, $3=B|D$, $4=B|C$, $5=A|D$; $\alpha=\gamma=Y_1$, $\beta=\delta=E$.}
\label{f:tree7}
\end{figure}

We will use the following notation for trees. Vertices are labeled by Greek letters and edges by numbers. If $n$ is a number and $X,Y\in\{A,B,C,D\}$ we write $n=X|Y$ if the edge labeled $n$ is a flow line of the gradient determined by the sheets $X$ and $Y$. If $\alpha$ is the label of a vertex then we write $\alpha=Y_j$, $j=0,1$ if the vertex is a $Y_j$-vertex, $\alpha=S$ if it is a switch, $\alpha=E$ if it is an end, and $\alpha=O$ if the vertex is some other kind of vertex (which does not appear in rigid trees). The first half of the $1$-parameter family of flow trees corresponding to the simple $1$-parameter family is schematically sketched in Figures \ref{f:tree1} - \ref{f:tree7}. (The second half is a copy of the first half reflected in a horizontal line through the center of the circles.) There are three critical instances in the family where we find trees with vertices which does not appear in rigid trees. Also the geometric dimension of these critical trees are one smaller than their formal
dimension.

Finally, imagine that the gradient line labeled $1$ in all figures continues to a minimum. Then the partial tree obtained by cutting, for example, the tree in Figure \ref{f:tree5} in the middle of the edge labeled $5$, cannot be part of any rigid tree since its formal dimension equals $3>2$. Also, this partial tree is an example of a tree for which the evaluation map at the special puncture is not injective, see the proof of Proposition \ref{p:ttv}.

\vspace{1cm}

\vspace{1cm}

\vspace{1cm}

%% file: Sec/ref.tex
\vspace{1cm}

\noindent
{\sc
USC, Department of mathematics\\
3620 S Vermont Ave\\
Los Angeles, CA 90089
}

\vspace{.1cm}
\noindent
{\tt tekholm{\@@}usc.edu}